\documentclass{siamltex1213}

\usepackage{latexsym}
\usepackage{mathrsfs}
\usepackage{amsfonts}
\usepackage{amssymb}
\usepackage{epsfig}
\usepackage{subfigure}
\usepackage{cite}
\usepackage{graphicx}
\usepackage{float}
\usepackage{multirow}
\usepackage{slashbox}
\usepackage{bm}
\usepackage{amsmath}
\usepackage{color}
\usepackage[subnum]{cases}
\usepackage{enumerate}

\def \reff#1{{\rm(\ref{#1})}}
\def \beginproof{\vskip .1cm\par\noindent {\bf Proof.}\ }
\def \endproof{\hskip .5cm $\Box$ \vskip .2cm}

\newtheorem{remark}{Remark}[section]
\newtheorem{assumption}{Assumption}[section]
\newtheorem{example}{Example}[section]

\title{Alternating Direction Method of Multipliers for A Class of Nonconvex and Nonsmooth Problems with Applications to Background/Foreground Extraction}
\author{Lei Yang\footnotemark[1] \and Ting Kei Pong\footnotemark[1] \and Xiaojun Chen\footnotemark[1]}

\begin{document}
\maketitle

\renewcommand{\thefootnote}{\fnsymbol{footnote}}
\footnotetext[1]{Department of Applied Mathematics, The Hong Kong Polytechnic University, Hung Hom, Kowloon, Hong Kong, P.R. China. (\email{lei.yang@connect.polyu.hk}, \email{tk.pong@polyu.edu.hk}, \email{xiaojun.chen@polyu.edu.hk}). The second author's work was supported in part by Hong Kong Research Grants Council PolyU253008/15p. The third author's work was supported in part by Hong Kong Research Grants Council PolyU153001/14p.}
\renewcommand{\thefootnote}{\arabic{footnote}}

\begin{abstract}
In this paper, we study a general optimization model, which covers a large class of existing models for many applications in imaging sciences. To solve the resulting possibly nonconvex, nonsmooth and non-Lipschitz optimization problem, we adapt the alternating direction method of multipliers (ADMM) with a general dual step-size to solve a reformulation that contains three blocks of variables, and analyze its convergence. We show that for any dual step-size less than the golden ratio, there exists a computable threshold such that if the penalty parameter is chosen above such a threshold and the sequence thus generated by our ADMM is bounded, then the cluster point of the sequence gives a stationary point of the nonconvex optimization problem. We achieve this via a potential function specifically constructed for our ADMM. Moreover, we establish the global convergence of the whole sequence if, in addition, this special potential function is a Kurdyka-{\L}ojasiewicz function. Furthermore, we present a simple strategy for initializing the algorithm to guarantee boundedness of the sequence. Finally, we perform numerical experiments comparing our ADMM with the proximal alternating linearized minimization (PALM) proposed in \cite{bst2014} on the background/foreground extraction problem with real data. The numerical results show that our ADMM with a nontrivial dual step-size is efficient.
\end{abstract}


%\begin{AMS}
%90C26, 90C46, 90C90, 65K10
%\end{AMS}

\pagestyle{myheadings}
\thispagestyle{plain}
\markboth{L. YANG, T. K. PONG, X. CHEN}{ADMM FOR NONCONVEX AND NONSMOOTH PROBLEMS}

\graphicspath{{images/}}

\section{Introduction}

In this paper, we consider the following optimization problem:
\begin{eqnarray}\label{generalmodel1}
\min \limits_{L,S}~\Psi(L)+\Phi(S)
+\frac{1}{2}\left\|D-\mathcal{A}\left[\mathcal{B}(L)+\mathcal{C}(S)\right]\right\|_F^2,
\end{eqnarray}
where \vspace{0.5mm}
\begin{itemize}
\item $\Psi, \Phi: \mathbb{R}^{m\times n} \rightarrow \mathbb{R}_{+}\cup\{\infty\}$ are proper closed nonnegative functions, and $\Psi$ is \textit{convex}, while $\Phi$ is possibly \textit{nonconvex}, \textit{nonsmooth} and \textit{non-Lipschitz}; \vspace{0.5mm}
\item $\mathcal{A}, \mathcal{B}, \mathcal{C}: \mathbb{R}^{m\times n} \rightarrow \mathbb{R}^{m\times n}$ are linear maps and $\mathcal{B}$, $\mathcal{C}$ are injective.
\end{itemize}
In particular, $\Psi(L)$ and $\Phi(S)$ in \eqref{generalmodel1} can be regularizers used for inducing the desired structures. For instance, $\Psi(L)$ can be used for inducing low rank in $L$. One possible choice is $\Psi(L)=\|L\|_*$ (see next section for notation and definitions). Alternatively, one may consider $\Psi(L)=\delta_{\Omega}(L)$, where $\Omega$ is a compact convex set such as $\Omega = \{L \in {\mathbb{R}}^{m\times n}~|~ \|L\|_\infty \le l,\ L_{:1} = L_{:2} = \cdots = L_{:n}\}$ with $l>0$, or $\Omega=\{L \in {\mathbb{R}}^{m\times n}~|~ \|L\|_* \le r\}$ with $r>0$; the former choice restricts $L$ to have rank at most $1$ and makes \eqref{generalmodel1} nuclear-norm-free (see \cite{lny2013,lwhc2014}). On the other hand, $\Phi(S)$ can be used for inducing sparsity. In the literature, $\Phi(S)$ is typically separable, i.e., taking the form
\begin{equation}\label{sepphi}
\Phi(S)=\mu\sum^m_{i=1}\sum^n_{j=1}\phi(s_{ij}),
\end{equation}
where $\phi$ is a nonnegative continuous function with $\phi(0)=0$ and $\mu>0$ is a regularization parameter. Some concrete examples of $\phi$ include:\vspace{0.5mm}
\begin{itemize}
\item[1.] bridge penalty \cite{hhjm2008,kf2000}: $\phi(t)=|t|^p$ for $0 < p \leq 1$; \vspace{0.5mm}
\item[2.] fraction penalty \cite{gr1992}: $\phi(t)=\alpha |t|/(1+\alpha|t|)$ for $\alpha>0$; \vspace{0.5mm}
\item[3.] logistic penalty \cite{nnzc2008}: $\phi(t)=\log(1+\alpha|t|)$ for $\alpha>0$; \vspace{0.5mm}
\item[4.] smoothly clipped absolute deviation \cite{f1997}: $\phi(t)=\int^{|t|}_{0} \min(1, (\alpha - s/\mu)_+/(\alpha-1))\,\mathrm{d}s$ for $\alpha>2$; \vspace{0.5mm}
\item[5.] minimax concave penalty \cite{z2010}: $\phi(t) = \int^{|t|}_0 (1-s/(\alpha\mu))_+\,\mathrm{d}s$ for $\alpha>0$; \vspace{0.5mm}
\item[6.] hard thresholding penalty function \cite{fl2001}: $\phi(t)=\mu - (\mu - |t|)^2_+ / \mu$. \vspace{0.5mm}
\end{itemize}
\vspace{0.5mm}
The bridge penalty and the logistic penalty have also been considered in \cite{DDLZH13}. Finally, the linear map $\mathcal{A}$ can be suitably chosen to model different scenarios. For example, $\mathcal{A}$ can be chosen to be the identity map for extracting $L$ and $S$ from a noisy data $D$, and the blurring map for a blurred data $D$. The linear map $\mathcal{B}$ can be the identity map or some ``dictionary" that spans the data space (see, for example, \cite{llysym2013}), and $\mathcal{C}$ can be chosen to be the identity map or the inverse of certain sparsifying transform (see, for example, \cite{ocs2015}). More examples of \eqref{generalmodel1} can be found in \cite{bz2014,bsjjz2015,CLMW11,XCS2012,DDLZH13,pgwxm2012}.

One representative application that is frequently modeled by \eqref{generalmodel1} via a suitable choice of $\Phi$, $\Psi$, $\cal A$, $\cal B$ and $\cal C$ is the background/foreground extraction problem, which is an important problem in video processing; see \cite{Bou11,Bou14} for recent surveys. In this problem, one attempts to separate the relatively static information called ``background" and the moving objects called ``foreground" in a video. The problem can be modeled by \eqref{generalmodel1}, and such models are typically referred to as RPCA-based models. In these models, each image is stacked as a column of a data matrix $D$, the relatively static background is then modeled as a low rank matrix, while the moving foreground is modeled as sparse outliers. The data matrix $D$ is then decomposed (approximately) as the sum of a low rank matrix $L\in\mathbb{R}^{m\times n}$ modeling the background and a sparse matrix $S\in\mathbb{R}^{m\times n}$ modeling the foreground. Various approximations are then used to induce low rank and sparsity, resulting in different RPCA-based models, most of which take the form of \eqref{generalmodel1}. One example is to set $\Psi$ to be the nuclear norm of $L$, i.e., the sum of singular values of $L$, to promote low rank in $L$ and $\Phi$ to be the $\ell_1$ norm of $S$ to promote sparsity in $S$, as in \cite{CLMW11}. Besides convex regularizers, nonconvex models have also been widely studied recently and their performances are promising; see \cite{DDLZH13,wcx2015} for background/foreground extraction and \cite{bc2015,cz2010,gzzf2014,n2014,nnzc2008,zrglz2015} for other problems in image processing. There are also nuclear-norm-free models that do not require matrix decomposition of the matrix variable $L$ when solving them, making the model more practical especially when the size of matrix is large. For instance, in \cite{lny2013}, the authors set $\Phi$ to be the $\ell_1$ norm of $S$ and $\Psi$ to be the indicator function of $\Omega = \{L \in {\mathbb{R}}^{m\times n}~|~ \ L_{:1} = L_{:2} = \cdots = L_{:n}\}$. A similar approach was also adopted in \cite{lwhc2014} with promising performances. Clearly, for nuclear-norm-free models, one can also take $\Phi$ to be some nonconvex sparsity inducing regularizers, resulting in a special case of \eqref{generalmodel1} that has not been explicitly considered in the literature before; we will consider these models in our numerical experiments in Section~\ref{sec5}. The above discussion shows that problem \eqref{generalmodel1} is flexible enough to cover a wide range of RPCA-based models for background/foreground extraction.

Problem \eqref{generalmodel1}, though nonconvex in general, as we will show later in Section~\ref{sec3}, can be reformulated into an optimization problem with three blocks of variables. This kind of problems containing several blocks of variables has been widely studied in the literature; see, for example, \cite{lny2013,mlfwl2014,pgwxm2012}. Hence, it is natural to adapt the algorithm used there, namely, the alternating direction method of multipliers (ADMM), for solving \eqref{generalmodel1}. Classically, the ADMM can be applied to solving problems of the following form that contains 2 blocks of variables:
\begin{eqnarray}\label{admmmodel}
\min\limits_{x_1,x_2} \left\{f_1(x_1) + f_2(x_2)~|~\mathcal{A}_1(x_1) + \mathcal{A}_2(x_2) = b \right\},
\end{eqnarray}
where $f_1$ and $f_2$ are proper closed convex functions, $\mathcal{A}_1$ and $\mathcal{A}_2$ are linear operators. The iterative scheme of ADMM is
\begin{eqnarray*}
\left\{\begin{aligned}
&x_1^{k+1} \in \mathop{\mathrm{Argmin}}
\limits_{x_1}\left\{\mathcal{L}_{\beta}(x_1,x_2^k,z^k)\right\},\\
&x_2^{k+1} \in \mathop{\mathrm{Argmin}}
\limits_{x_2}\left\{\mathcal{L}_{\beta}(x_1^{k+1},x_2,z^k)\right\},\\
&z^{k+1} = z^k - \tau\beta(\mathcal{A}_1(x_1^{k+1})+\mathcal{A}_2(x_2^{k+1})-b),
\end{aligned}\right.
\end{eqnarray*}
where $\tau \in (0,\frac{\sqrt{5}+1}2)$ is the dual step-size and $\mathcal{L}_{\beta}$ is the augmented Lagrangian function for \eqref{admmmodel} defined as
\begin{eqnarray*}
\begin{aligned}
\mathcal{L}_{\beta}(x_1,x_2,z)~:=&~f_1(x_1) + f_2(x_2)-\langle z, \mathcal{A}_1(x_1) + \mathcal{A}_2(x_2) - b \rangle \\
&~+\frac{\beta}{2}\|\mathcal{A}_1(x_1) + \mathcal{A}_2(x_2) - b\|^2
\end{aligned}
\end{eqnarray*}
with $\beta>0$ being the penalty parameter. Under some mild conditions, the sequence $\{(x_1^k, x_2^k)\}$ generated by the above ADMM can be shown to converge to an optimal solution of \eqref{admmmodel}; see for example, \cite{bt1989,eb1992,GaM76,GlM75}. However, the ADMM used in \cite{lny2013,mlfwl2014,pgwxm2012} does not have a convergence guarantee; indeed, it is shown recently in \cite{CHYY14} that the ADMM, when applied to a convex optimization problem with $3$ blocks of variables, can be divergent in general. This motivates the study of many provably convergent variants of the ADMM for convex problems with more than 2 blocks of variables; see, for example, \cite{hty2012,hy2013,lst2015,lll2015}. Recently, Hong et al. \cite{hlr2014} established the convergence of the multi-block ADMM for certain types of nonconvex problems whose objective is a sum of a possibly {\em nonconvex} Lipschitz differentiable function and a bunch of convex nonsmooth functions when the penalty parameter is chosen above a computable threshold. The problem they considered covers \eqref{generalmodel1} when $\Phi$ is convex, or smooth and possibly nonconvex. Later, Wang et al. \cite{wcx2015} considered a more general type of nonconvex problems that contains \eqref{generalmodel1} as a special case and allows some nonconvex nonsmooth functions in the objective. To solve this type of problems, they considered a variant of the ADMM whose subproblems are simplified by adding a Bregman proximal term. However, their results cannot be applied to the direct adaptation of the ADMM for solving \eqref{generalmodel1}.

In this paper, following the studies in \cite{hlr2014,wcx2015} on convergence of nonconvex ADMM and its variant, and the recent studies in \cite{ah2014,lp2014,wxx2014}, we manage to analyze the convergence of the ADMM applied to solving the possibly nonconvex problem \eqref{generalmodel1}. In addition, we would like to point out that all the aforementioned nonconvex ADMM have a dual step-size of $\tau = 1$. While it is known that the classical ADMM converges for any $\tau\in (0,\frac{\sqrt{5}+1}2)$ for convex problems, and that empirically $\tau\approx \frac{\sqrt{5}+1}2$ works best (see, for example, \cite{FPST12,GaM76,GlM75,lst2015}), to our knowledge, the algorithm with a dual step-size $\tau \neq 1$ has never been studied in the nonconvex scenarios. Thus, we also study the ADMM with a general dual step-size, which will allow more flexibility in the design of algorithms.

The contributions of this paper are as follows:
\begin{enumerate}
  \item We show that for any positive dual step-size $\tau$ less than the golden ratio, the cluster point of the sequence generated by our ADMM gives a stationary point of \eqref{generalmodel1} if the penalty parameter is chosen above a computable threshold depending on $\tau$, whenever the sequence is bounded. We achieve this via a potential function specifically constructed for our ADMM. To the best of our knowledge, this is the first convergence result for the ADMM in the nonconvex scenario with a possibly nontrivial dual step-size ($\tau\neq 1$). This result is also new for the convex scenario for the multi-block ADMM.
  \item We establish global convergence of the whole sequence generated by the ADMM under the additional assumption that the special potential function is a Kurdyka-{\L}ojasiewicz function. Following the discussions in \cite[Section~4]{AtBoReSo10}, one can check that this condition is satisfied for all the aforementioned $\phi$.
  \item Furthermore, we discuss an initialization strategy to guarantee the boundedness of the sequence generated by the ADMM.
\end{enumerate}

We also conduct numerical experiments to evaluate the performance of our ADMM by using different nonconvex regularizers and real data. Our computational results illustrate the efficiency of our ADMM with a nontrivial dual step-size.

The rest of this paper is organized as follows. We present notation and preliminaries in Section~\ref{sec2}. The ADMM for \eqref{generalmodel1} is described in Section~\ref{sec3}. We analyze the convergence of the method in Section~\ref{sec4}. Numerical results are presented in Section~\ref{sec5}, with some concluding remarks given in Section~\ref{sec6}.

%------------------------------------------------------------------------------------------- Section 2
\section{Notation and preliminaries}\label{sec2}

In this paper, we use $\mathbb{R}^{m \times n}$ to denote the set of all $m\times n$ matrices. For a matrix $X \in \mathbb{R}^{m \times n}$, we let $x_{ij}$ denote its $(i,j)$th entry and $X_{:j}$ denote its $j$th column. The number of nonzero entries in $X$ is denoted by $\|X\|_0$ and the largest entry in magnitude is denoted by $\|X\|_{\infty}$. Moreover, the Fr\"{o}benius norm is denoted by $\|X\|_F$, the nuclear norm is denoted by $\|X\|_*$, which is the sum of singular values of $X$; and $\ell_1$-norm and $\ell_p$-quasi-norm ($0 < p < 1$) are given by $\|X\|_1:=\sum_{i=1}^{m}\sum_{j=1}^{n} |x_{ij}|$ and $\|X\|_p:=\left( \sum_{i=1}^{m}\sum_{j=1}^{n} |x_{ij}|^p \right)^{\frac{1}{p}}$, respectively. Furthermore, for two matrices $X$ and $Y$ of the same size, we denote their trace inner product by $\langle X,Y\rangle := \sum_{i=1}^m\sum_{j=1}^nx_{ij}y_{ij}$. Finally, for the linear map $\mathcal{A}:\mathbb{R}^{m \times n}\to \mathbb{R}^{m \times n}$ in \eqref{generalmodel1}, its adjoint is denoted by $\mathcal{A}^*$, while the largest (resp., smallest) eigenvalue of the linear map $\mathcal{A}^*\mathcal{A}$ is denoted by $\lambda_{\max}$ (resp., $\lambda_{\min}$). The identity map is denoted by $\mathcal{I}$.

For an extended-real-valued function $f: \mathbb{R}^{m \times n} \rightarrow [-\infty,\infty]$, we say that it is \textit{proper} if $f(X) > -\infty$ for all $X \in \mathbb{R}^{m \times n}$ and its domain ${\rm dom}f:=\{X\in \mathbb{R}^{m \times n} ~|~f(X) < \infty\}$ is nonempty. %Such a function is \textit{lower semicontinuous} at a point $X\in {\rm dom}f$ if $f(X) \leq \liminf_{Y \rightarrow X} f(Y)$.
For a proper function $f$, we use the notation $Y \xrightarrow{f} X$ to denote $Y \rightarrow X$ and $f(Y) \rightarrow f(X)$. Our basic \textit{(limiting-)subdifferential} \cite[Definition~8.3]{rw1998} of $f$ at $X\in \mathrm{dom}f$ used in this paper, denoted by $\partial f(X)$, is defined as
\begin{eqnarray*}%\label{def_sub}
\partial f(X):=\left\{ D \in \mathbb{R}^{m\times n}: \exists~ X^k \xrightarrow{f} X~\mathrm{and}~D^k \rightarrow D ~\mathrm{with}~D^k \in \widehat{\partial} f(X^k) ~\mathrm{for~all}~k\right\},
\end{eqnarray*}
where $\widehat{\partial} f(U)$ denotes the Fr\'{e}chet subdifferential of $f$ at $U\in \mathrm{dom}f$, which is the set of all $D \in \mathbb{R}^{m\times n}$ satisfying
\begin{eqnarray*}%\label{def_sub}
\liminf\limits_{Y\neq U, Y \rightarrow U} \frac{f(Y)-f(U)-\langle D, Y-U\rangle}{\|Y-U\|_F} \geq 0.
\end{eqnarray*}
From the above definition, we can easily observe that
\begin{eqnarray}\label{robust}
\left\{ D\in\mathbb{R}^{m\times n}: \exists\, X^k \xrightarrow{f} X,~ D^k \rightarrow D,~ D^k \in \partial f(X^k) \right\} \subseteq \partial f(X).
\end{eqnarray}
We also recall that when $f$ is continuously differentiable or convex, the above subdifferential coincides with the classical concept of derivative or convex subdifferential of $f$; see, for example, \cite[Exercise~8.8]{rw1998} and \cite[Proposition~8.12]{rw1998}. Moreover, from the generalized Fermat's rule \cite[Theorem~10.1]{rw1998}, we know that if $X\in\mathbb{R}^{m\times n}$ is a local minimizer of $f$, then $0 \in \partial f(X)$. Additionally, for a function $f$ with several groups of variables, we write $\partial_X f$ (resp., $\nabla_X f$) for the subdifferential (resp., derivative) of $f$ with respect to the group of variables $X$.

For a compact convex set $\Omega\subseteq\mathbb{R}^{m \times n}$, its indicator function $\delta_{\Omega}$ is defined by
\begin{eqnarray*}
\delta_{\Omega}(X) = \left\{
\begin{array}{ll}
0 &\quad\mathrm{if}~X\in\Omega, \\
+\infty &\quad\mathrm{otherwise}.
\end{array}\right.
\end{eqnarray*}
The normal cone of $\Omega$ at the point $X \in \Omega$ is given by $\mathcal{N}_{\Omega}(X)=\partial\delta_{\Omega}(X)$. We also use $\mathrm{dist}(X, \Omega)$ to denote the distance from $X$ to $\Omega$, i.e., $\mathrm{dist}(X, \Omega) := \inf_{Y\in\Omega}\|X-Y\|_F$, and ${\cal P}_\Omega(X)$ to denote the unique closest point to $X$ in $\Omega$.

Next, we recall the Kurdyka-{\L}ojasiewicz (KL) property, which plays an important role in our global convergence analysis. For notational simplicity, we use $\Xi_{\eta}$ ($\eta>0$) to denote the class of concave functions $\varphi:[0,\eta) \rightarrow \mathbb{R}_{+}$ satisfying: (1) $\varphi(0)=0$; (2) $\varphi$ is continuously differentiable on $(0,\eta)$ and continuous at $0$; (3) $\varphi'(x)>0$ for all $x\in(0,\eta)$. Then the KL property can be described as follows.
\begin{definition}[\textbf{KL property and KL function}]
Let $f$ be a proper lower semicontinuous function.
\begin{itemize}
\item[(i)] For $\tilde{X}\in{\rm dom}\,\partial f:=\{X\in \mathbb{R}^{m\times n}: \partial f(X) \neq \emptyset\}$, if there exist an $\eta\in(0, +\infty]$, a neighborhood $V$ of $\tilde{X}$ and a function $\varphi \in \Xi_{\eta}$ such that for all $X \in V \cap \{X\in \mathbb{R}^{m\times n}:f(\tilde{X})<f(X)<f(\tilde{X})+\eta\}$, it holds that
    \begin{eqnarray*}
    \varphi'(f(X)-f(\tilde{X}))\mathrm{dist}(0, \partial f(X)) \geq 1,
    \end{eqnarray*}
    then $f$ is said to have the \textbf{Kurdyka-{\L}ojasiewicz (KL)} property at $\tilde{X}$.

\item[(ii)] If $f$ satisfies the KL property at each point of ${\rm dom}\,\partial f$, then $f$ is called a KL function.
\end{itemize}
\end{definition}

We refer the interested readers to \cite{AtBoReSo10} and references therein for examples of KL functions. We also recall the following uniformized KL property, which was established in \cite[Lemma 6]{bst2014}.

\begin{proposition}[\textbf{Uniformized KL property}]\label{uniKL}
Suppose that $f$ is a proper lower semicontinuous function and $\Gamma$ is a compact set. If $f \equiv f^*$ on $\Gamma$ for some constant $f^*$ and satisfies the KL property at each point of $\Gamma$, then there exist $\varepsilon>0$, $\eta>0$ and $\varphi \in \Xi_{\eta}$ such that
\begin{eqnarray*}
\varphi'(f(X)-f^*)\mathrm{dist}(0, \partial f(X)) \geq 1
\end{eqnarray*}
for all $X \in \{X\in\mathbb{R}^{m\times n}: \mathrm{dist}(X,\Gamma)<\varepsilon\} \cap \{X\in \mathbb{R}^{m\times n}:f^*<f(X)<f^*+\eta\}$.
\end{proposition}

Before ending this section, we discuss first-order necessary conditions for \eqref{generalmodel1}. First, recall that \eqref{generalmodel1} is the same as
\begin{equation*}
    \min_{L,S} \ \mathcal{F}(L, S): = \Psi(L)+\Phi(S)
+\frac{1}{2}\left\|D-\mathcal{A}\left[\mathcal{B}(L)+\mathcal{C}(S)\right]\right\|_F^2.
\end{equation*}
Hence, from \cite[Theorem~10.1]{rw1998}, we have $0\in \partial \mathcal{F}(\bar L,\bar S)$ at any local minimizer $(\bar L,\bar S)$ of \eqref{generalmodel1}. On the other hand, from \cite[Exercise~8.8]{rw1998} and \cite[Proposition~10.5]{rw1998}, we see that
\begin{eqnarray*}
\partial \mathcal{F}(L,\,S) = \left(
\begin{array}{c}
\partial \Psi(L)+\mathcal{B}^*\mathcal{A}^*\left(\mathcal{A}(\mathcal{B}(L)+\mathcal{C}(S))-D\right) \\
\partial \Phi(S)+\mathcal{C}^*\mathcal{A}^*\left(\mathcal{A}(\mathcal{B}(L)+\mathcal{C}(S))-D\right)
\end{array}\right).
\end{eqnarray*}
Consequently, the first-order necessary conditions of \eqref{generalmodel1} at the local minimizer $(\bar L,\bar S)$ is given by:
\begin{eqnarray}\label{optcond}
\left\{\begin{aligned}
&0\in \partial \Psi(\bar{L})+\mathcal{B}^*\mathcal{A}^*\left(\mathcal{A}(\mathcal{B}(\bar{L})+\mathcal{C}(\bar{S}))-D\right),\\
&0\in \partial \Phi(\bar{S})+\mathcal{C}^*\mathcal{A}^*\left(\mathcal{A}(\mathcal{B}(\bar{L})+\mathcal{C}(\bar{S}))-D\right).
\end{aligned}\right.
\end{eqnarray}
In this paper, we say that $(L^*, S^*)$ is a stationary point of \eqref{generalmodel1} if $(L^*, S^*)$ satisfies \eqref{optcond} in place of $(\bar{L}, \bar{S})$.

%--------------------------------------------------------------------------------------------------------section 3
\section{Alternating direction method of multipliers}\label{sec3}

In this section, we present an ADMM for solving \eqref{generalmodel1}, which can be equivalently written as
\begin{eqnarray}\label{generalmodel2}
\begin{aligned}
&\min \limits_{L,S,Z}~~\Psi(L)+\Phi(S)+\frac{1}{2}\left\|D-\mathcal{A}(Z)\right\|_F^2 \\
&~~\mathrm{s.t.}~~~~\mathcal{B}(L)+\mathcal{C}(S) = Z.
\end{aligned}
\end{eqnarray}
To describe the iterates of the ADMM, we first introduce the augmented Lagrangian function of the above optimization problem:
\begin{eqnarray*}
\begin{aligned}
&\mathcal{L}_{\beta} (L,S,Z,\Lambda)=\Psi(L)+\Phi(S)+\frac{1}{2}\|D-\mathcal{A}(Z)\|_F^2 \\
&\hspace{10mm}-\langle\Lambda, ~\mathcal{B}(L)+\mathcal{C}(S)-Z\rangle + \frac{\beta}{2}\|\mathcal{B}(L)+\mathcal{C}(S)-Z\|_F^2,
\end{aligned}
\end{eqnarray*}
where $\Lambda \in \mathbb{R}^{m \times n}$ is the Lagrangian multiplier and $\beta > 0$ is the penalty parameter. The ADMM for solving \eqref{generalmodel2} (equivalently \eqref{generalmodel1}) is then presented as follows:
\vspace{2mm}
%--------------------------------------------------------------------------------------------- Algorithm 1
\begin{center}
\fbox{\parbox{5in}{\vspace{1mm}
~\textbf{Algorithm 1} ~ADMM for solving \eqref{generalmodel2}
\begin{description}
\item{\textbf{Input}:} Initial point $(S^0, Z^0, \Lambda^0)$, dual step-size parameter $\tau>0$, penalty parameter $\beta>0$, $k=0$

\item \hspace{10mm} \textbf{while} a termination criterion is not met, \textbf{do} \vspace{1mm}
\begin{itemize}
\item[] \textbf{Step 1}. Set
\begin{numcases}{}
L^{k+1} \in \mathop{\mathrm{Argmin}} \limits_{L}~\mathcal{L}_{\beta} (L, S^k, Z^k, \Lambda^k) \label{Lsubpro}\\
S^{k+1} \in \mathop{\mathrm{Argmin}} \limits_{S}~ \mathcal{L}_{\beta} (L^{k+1}, S, Z^k, \Lambda^k) \label{Ssubpro}\\
Z^{k+1} = \mathop{\mathrm{argmin}} \limits_{Z}~\mathcal{L}_{\beta} (L^{k+1}, S^{k+1}, Z, \Lambda^k) \label{Zsubpro}\\
\Lambda^{k+1} = \Lambda^k - \tau\beta (\mathcal{B}(L^{k+1})+\mathcal{C}(S^{k+1})-Z^{k+1})  \nonumber
\end{numcases}
\item[] \textbf{Step 2}. Set $k := k + 1$
\end{itemize}
\item \vspace{1mm}\hspace{10mm} \textbf{end while}
\item{\textbf{Output}:} $(L^{k}, S^{k})$
\end{description}}}
\end{center}
\vspace{2mm}

Comparing with the ADMM considered in \cite{hlr2014}, the above algorithm has an extra dual step-size parameter $\tau > 0$ in the $\Lambda$-update. Such a dual step-size was introduced in \cite{GaM76,GlM75} for the classical ADMM (i.e., for convex problems with two separate blocks of variables), and was further studied in \cite{xw2011,FPST12,sty2015,lst2015} for other variants of the ADMM. Numerically, it was also demonstrated in \cite{sty2015} that a larger dual step-size ($\tau \approx \frac{\sqrt{5}+1}{2}$) results in faster convergence for the convex problems they consider. Thus, we adapt this dual step-size $\tau$ in our algorithm above. Surprisingly,
in our numerical experiments, a parameter choice of $\tau \approx \frac{\sqrt{5}+1}{2}$ leads to the worst performance for our nonconvex problems.

When $\tau = 1$, the above algorithm is a special case of the general algorithm studied in \cite{hlr2014} when $\Psi$ and $\Phi$ are smooth functions, or convex nonsmooth functions. The algorithm is shown to converge when $\beta$ is chosen above a computable threshold. However, their convergence result cannot be directly applied when $\tau \neq 1$ or when $\Phi$ is nonsmooth and nonconvex. Nevertheless, following their analysis and the related studies \cite{lp2014,wxx2014,wcx2015}, the above algorithm can be shown to be convergent under suitable assumptions. We will present the convergence analysis in Section~\ref{sec4}.

Before ending this section, we further discuss the three subproblems in Algorithm 1. First, notice that the $L$-update and $S$-update are given by
\begin{eqnarray*}
\left\{\begin{aligned}
&L^{k+1} \in \mathop{\mathrm{Argmin}}\limits_{L}~\left\{ \Psi(L)+\frac{\beta}{2}\|\mathcal{B}(L)+\mathcal{C}(S^k)-Z^k-\frac{1}{\beta}\Lambda^k\|_F^2\right\}, \\
&S^{k+1} \in \mathop{\mathrm{Argmin}}\limits_{S}~\left\{ \Phi(S)+\frac{\beta}{2}\|\mathcal{B}(L^{k+1})+\mathcal{C}(S)-Z^k-\frac{1}{\beta}\Lambda^k\|_F^2\right\}.
\end{aligned}\right.
\end{eqnarray*}
In general, these two subproblems are not easy to solve. However, when $\Psi$ and $\Phi$ are chosen to be some common regularizers used in the literature, for example, $\Psi(L) = \|L\|_*$ and $\Phi(S) = \|S\|_1$, then these subproblems can be solved efficiently via the proximal gradient method. Additionally, when $\Psi(L)=\delta_{\Omega}(L)$ with $\Omega$ being a closed convex set and $\mathcal{B}=\mathcal{I}$, the $L$-update can be given explicitly by
\[
L^{k+1} = {\cal P}_\Omega\left(-\mathcal{C}(S^k)+Z^k+\frac{1}{\beta}\Lambda^k\right),
\]
which can be computed efficiently if $\Omega$ is simple, for example, when $\Omega = \{L \in {\mathbb{R}}^{m\times n}~|~ \|L\|_\infty$ $\le l,\ L_{:1} = L_{:2} = \cdots = L_{:n}\}$ for some $l > 0$. For the $S$-update, when $\Phi$ is given by \eqref{sepphi} with $\phi$ being one of the penalty functions presented in the introduction and $\mathcal{C}=\mathcal{I}$, it can be solved efficiently via a simple root-finding procedure. Finally, from the optimality conditions of \eqref{Zsubpro}, the $Z^{k+1}$ can be obtained by solving the following linear system
\[
\mathcal{A}^{*}\!\mathcal{A}(Z)+\beta Z=\mathcal{A}^{*}(D)-\Lambda^k+\beta\left(\mathcal{B}(L^{k+1})+\mathcal{C}(S^{k+1})\right),
\]
whose complexity would depend on the choice of ${\mathcal{A}}$ in our model \eqref{generalmodel1}. For example, when ${\mathcal{A}}$ is just the identity map, the $Z^{k+1}$ is given explicitly by
\[
Z^{k+1}=\frac{1}{1+\beta}\left[D-\Lambda^k+\beta\left(\mathcal{B}(L^{k+1})+\mathcal{C}(S^{k+1})\right)\right].
\]

\section{Convergence analysis}\label{sec4}

In this section, we discuss the convergence of Algorithm 1 for $0<\tau<\frac{1+\sqrt{5}}{2}$. We first present the first-order optimality conditions for the subproblems in Algorithm 1 as follows, which will be used repeatedly in our convergence analysis below.
\begin{numcases}{}
0 \in \partial \Psi(L^{k+1})-\mathcal{B}^{*}(\Lambda^k)+\beta\mathcal{B}^{*}\left(\mathcal{B}(L^{k+1})
+\mathcal{C}(S^k)-Z^k\right), \label{lopt}\\
0 \in \partial \Phi(S^{k+1})-\mathcal{C}^{*}(\Lambda^k)+\beta\mathcal{C}^{*}\left(\mathcal{B}(L^{k+1})
+\mathcal{C}(S^{k+1})-Z^k\right), \label{sopt} \\
0 = \mathcal{A}^{*}(\mathcal{A}(Z^{k+1})-D) + \Lambda^k - \beta(\mathcal{B}(L^{k+1})+\mathcal{C}(S^{k+1})-Z^{k+1}), \label{zopt}\\
\Lambda^{k+1}-\Lambda^k = -\tau\beta\left(\mathcal{B}(L^{k+1})+\mathcal{C}(S^{k+1})-Z^{k+1}\right). \label{lamopt}
\end{numcases}
Our convergence analysis is largely based on the following potential function:
\begin{eqnarray*}
\Theta_{\tau,\beta}(L, S, Z, \Lambda) = \mathcal{L}_{\beta} (L, S, Z, \Lambda) + \theta(\tau)\beta\|\mathcal{B}(L)+\mathcal{C}(S)-Z\|_F^2,
\end{eqnarray*}
where
\begin{eqnarray}\label{deftheta}
\theta(\tau):=\max\left\{1-\tau, ~\frac{(\tau-1)\tau^2}{1+\tau-\tau^2}\right\}, \quad \mathrm{for}~~0<\tau<\frac{1+\sqrt{5}}{2}.
\end{eqnarray}
Note that $\theta(\cdot)$ is a convex and nonnegative function on $\left(0,~ \frac{1+\sqrt{5}}{2}\right)$. Thus,  for any $(L, S, Z, \Lambda)$, we have $\Theta_{\tau,\beta}(L, S, Z, \Lambda) \geq \mathcal{L}_{\beta} (L, S, Z, \Lambda)$ for $0<\tau<\frac{1+\sqrt{5}}{2}$, and the equality holds when $\tau = 1$ (so that $\theta(\tau) = 0$).

Our convergence analysis also relies on the following assumption.
\begin{assumption}\label{assum}
$\Psi$, $\Phi$, $\mathcal{B}$, $\mathcal{C}$, $\beta$ and $\tau$ satisfy
\begin{itemize}
\item[(a1)] $\mathcal{B}^*\mathcal{B}\succeq\sigma\mathcal{I}$ for some $\sigma>0$ and $\mathcal{C}^*\mathcal{C}\succeq\sigma'\mathcal{I}$ for some $\sigma'>0$;
\item[(a2)] $\Psi$ is continuous in its domain;
\item[(a3)] the first iterate $(L^1, S^1, Z^1, \Lambda^1)$ satisfies
\begin{eqnarray*}
\Theta_{\tau,\beta}(L^1, S^1, Z^1, \Lambda^1) < h_0:=\liminf\limits_{\|L\|_F+\|S\|_F\rightarrow\infty} \Psi(L)+\Phi(S).
\end{eqnarray*}
\end{itemize}
\end{assumption}

\begin{remark}[\textbf{Note on Assumption~\ref{assum}}]
(i) Since $\mathcal{B}$ and $\mathcal{C}$ in \eqref{generalmodel1} are injective, (a1) holds trivially; (ii) (a2) holds for many common regularizers (for example, the nuclear norm) or the indicator function of a set; (iii) (a3) places conditions on the {\em first} iterate of the algorithm. It is not hard to observe that this assumption holds trivially if both $\Psi$ and $\Phi$ are coercive, i.e., if
$\liminf\limits_{\|L\|_F+\|S\|_F\rightarrow\infty} \Psi(L)+\Phi(S) = \infty$. We will discuss more sufficient conditions for this assumption after our convergence results, i.e., after Theorem~\ref{thm42}.
\end{remark}

We now start our convergence analysis by proving the following preparatory lemma, which states that the potential function is decreasing along the sequence generated from Algorithm 1 if the penalty parameter $\beta$ is chosen above a computable threshold.

\begin{lemma}\label{Phidelemma}
Suppose that $0 < \tau < \frac{1+\sqrt{5}}{2}$ and $\{(L^k, S^k, Z^k, \Lambda^k)\}$ is a sequence generated by Algorithm 1. If (a1) in Assumption \ref{assum} holds, then for $k \geq 1$, we have
\begin{eqnarray}\label{difofp}
\hspace{-8mm}\begin{aligned}
&\quad\Theta_{\tau,\beta} (L^{k+1}, S^{k+1}, Z^{k+1}, \Lambda^{k+1}) - \Theta_{\tau,\beta} (L^k, S^k, Z^k, \Lambda^k)\\
&\leq \textstyle{\left(\max\left\{\frac{1}{\tau}, \frac{\tau^2}{1+\tau-\tau^2}\right\}\cdot\frac{\lambda_{\max}^2}{\beta} - \frac{\lambda_{\min}+\beta}{2}\right)}\|Z^{k+1} - Z^k\|_F^2 -\frac{\sigma\beta}{2}\|L^{k+1} - L^k\|_F^2.
\end{aligned}
\end{eqnarray}
Moreover, if $\beta \geq -\frac{\lambda_{\min}}{2}+\frac{1}{2}\sqrt{\lambda_{\min}^2+\max\left\{\frac{1}{\tau}, \frac{\tau^2}{1+\tau-\tau^2}\right\} \cdot 8\lambda_{\max}^2}$, then the sequence $\{\Theta_{\tau,\beta}(L^k, S^k, Z^k$, $\Lambda^k)\}_{k=1}^{\infty}$ is decreasing.
\end{lemma}
\beginproof
We start our proof by noticing that
\begin{eqnarray}\label{first_term_tmp}
&&\Theta_{\tau,\beta} (L^{k+1}, S^{k+1}, Z^{k+1}, \Lambda^{k+1}) - \Theta_{\tau,\beta} (L^{k+1}, S^{k+1}, Z^{k+1}, \Lambda^k) \nonumber\\
&=& -\langle \Lambda^{k+1} - \Lambda^k, L^{k+1} + S^{k+1} - Z^{k+1} \rangle = \frac{1}{\tau\beta}\|\Lambda^{k+1} - \Lambda^k\|_F^2,
\end{eqnarray}
where the last equality follows from \eqref{lamopt}. We next derive an upper bound of $\|\Lambda^{k+1} - \Lambda^k\|_F^2$.
To proceed, we first note from \eqref{zopt} that
\begin{eqnarray*}
0&=& \mathcal{A}^{*}(\mathcal{A}(Z^{k+1})-D) + \Lambda^k - \beta(\mathcal{B}(L^{k+1}) + \mathcal{C}(S^{k+1}) - Z^{k+1})  \nonumber\\
 &=& \mathcal{A}^{*}(\mathcal{A}(Z^{k+1})-D) + \Lambda^k + \frac{1}{\tau}(\Lambda^{k+1}-\Lambda^k) \nonumber\\
\Longrightarrow && \Lambda^{k+1} = \tau \mathcal{A}^*(D - \mathcal{A}(Z^{k+1})) + (1-\tau)\Lambda^k,
\end{eqnarray*}
where the second equality follows from \eqref{lamopt}. Hence, for $k\ge 1$,
\begin{eqnarray}
&&\Lambda^{k+1} - \Lambda^{k}  \nonumber \\
&=& [\tau \mathcal{A}^*(D - \mathcal{A}(Z^{k+1})) + (1-\tau)\Lambda^k] - [\tau \mathcal{A}^*(D - \mathcal{A}(Z^{k})) + (1-\tau)\Lambda^{k-1}] \nonumber\\
&=& \tau \mathcal{A}^*\mathcal{A}(Z^{k} - Z^{k+1}) + (1-\tau)(\Lambda^k - \Lambda^{k-1}). \label{lamdiff}
\end{eqnarray}
We now consider two separate cases: $0<\tau\leq 1$ and $1<\tau<\frac{1+\sqrt{5}}{2}$. \vspace{1mm}
\begin{itemize}
\item For $0<\tau\leq 1$, it follows from the convexity of $\|\cdot\|_F^2$ that
      \begin{eqnarray*}
      \|\Lambda^{k+1} - \Lambda^k\|_F^2 &=& \left\|\tau \mathcal{A}^*\mathcal{A}(Z^{k} - Z^{k+1}) + (1-\tau)(\Lambda^k - \Lambda^{k-1})\right\|_F^2 \\
      &\leq& \tau\lambda^2_{\max}\|Z^{k+1} - Z^k\|_F^2 + (1-\tau)\|\Lambda^k - \Lambda^{k-1}\|_F^2.
      \end{eqnarray*}
      We further add $-(1-\tau)\left\|\Lambda^{k+1} - \Lambda^{k}\right\|_F^2$ to both sides of the above inequality and simplify the resulting inequality to get
      \begin{eqnarray}\label{lamboundbyZ1}
      \hspace{-4mm}\begin{aligned}
      &~~\|\Lambda^{k+1} - \Lambda^k\|_F^2  \\
      &\leq \lambda^2_{\max}\|Z^{k+1} - Z^k\|_F^2 + {\textstyle\frac{1-\tau}{\tau}}\left(\|\Lambda^k - \Lambda^{k-1}\|_F^2 - \|\Lambda^{k+1} - \Lambda^k\|_F^2 \right)  \\
      &= (1-\tau)\tau\beta^2\left(\|\mathcal{B}(L^k)+\mathcal{C}(S^k)-Z^k\|_F^2 - \|\mathcal{B}(L^{k+1})+\mathcal{C}(S^{k+1})-Z^{k+1}\|_F^2\right)  \\
      &~~+\lambda^2_{\max}\|Z^{k+1} - Z^k\|_F^2.
      \end{aligned}
      \end{eqnarray}
      where the last equality follows from \eqref{lamopt}.

\item For $1 < \tau < \frac{1+\sqrt{5}}{2}$, dividing $\tau$ from both sides of \eqref{lamdiff}, we have
      \begin{eqnarray*}
      \frac{1}{\tau}\left(\Lambda^{k+1} - \Lambda^{k}\right) &=& \mathcal{A}^*\mathcal{A}\left(Z^{k} - Z^{k+1}\right)
      + \left(\frac{1}{\tau}-1\right)(\Lambda^k - \Lambda^{k-1}) \\
      &=& \frac{1}{\tau}~ \tau \mathcal{A}^*\mathcal{A}\left(Z^{k} - Z^{k+1}\right) + \left(1-\frac{1}{\tau}\right)(\Lambda^{k-1} - \Lambda^{k}).
      \end{eqnarray*}
      This together with $0<\frac{1}{\tau}<1$ and the convexity of $\|\cdot\|_F^2$, implies that
      \begin{eqnarray*}
      \left\|{\textstyle\frac{1}{\tau}}\left(\Lambda^{k+1} - \Lambda^{k}\right)\right\|_F^2 &\leq& {\textstyle\frac{1}{\tau}}\|\tau \mathcal{A}^*\mathcal{A}\left(Z^{k} - Z^{k+1}\right)\|_F^2+{\textstyle\left(1-\frac{1}{\tau}\right)}\|\Lambda^{k-1} - \Lambda^{k}\|_F^2  \nonumber\\
      &\le&\tau\lambda_{\max}^2\|Z^{k+1} - Z^k\|_F^2 + {\textstyle\left(1-\frac{1}{\tau}\right)}\|\Lambda^{k} - \Lambda^{k-1}\|_F^2 \nonumber\vspace{2mm}\\
      \Longrightarrow \quad \|\Lambda^{k+1} - \Lambda^k\|_F^2 &\leq& \tau^3\lambda_{\max}^2\|Z^{k+1} - Z^k\|_F^2 + \left(\tau^2-\tau\right)\|\Lambda^{k} - \Lambda^{k-1}\|_F^2.
      \end{eqnarray*}
      Then, adding $-\left(\tau^2-\tau\right)\left\|\Lambda^{k+1}-\Lambda^{k}\right\|_F^2$ to both sides of the above inequality, simplifying the resulting inequality and using the fact that $1+\tau-\tau^2>0$ for $1<\tau<\frac{1+\sqrt{5}}{2}$, we see that
      \begin{eqnarray}\label{lamboundbyZ2}
      \hspace{-4mm}\begin{aligned}
      &~~\|\Lambda^{k+1} - \Lambda^k\|_F^2 \\
      &\leq {\textstyle\frac{\tau^3\lambda_{\max}^2}{1+\tau-\tau^2}}\|Z^{k+1} - Z^k\|_F^2 + {\textstyle\frac{\tau^2-\tau}{1+\tau-\tau^2}}\left(\|\Lambda^k - \Lambda^{k-1}\|_F^2-\|\Lambda^{k+1} - \Lambda^k\|_F^2\right)  \\
      &={\textstyle\frac{\tau^3\lambda_{\max}^2}{1+\tau-\tau^2}}\|Z^{k+1} - Z^k\|_F^2+{\textstyle\frac{(\tau-1)\tau^3\beta^2}{1+\tau-\tau^2}}\left(\|\mathcal{B}(L^k)+\mathcal{C}(S^k)-Z^k\|_F^2\right.\\
      &\quad-\left.\|\mathcal{B}(L^{k+1})+\mathcal{C}(S^{k+1})-Z^{k+1}\|_F^2\right),
      \end{aligned}
      \end{eqnarray}
      where the equality follows from \eqref{lamopt}.
\end{itemize}
\vspace{2mm}
Thus, for $0<\tau<\frac{1+\sqrt{5}}{2}$, combining \eqref{lamboundbyZ1}, \eqref{lamboundbyZ2} and recalling the definition of $\theta(\tau)$ in \eqref{deftheta}, we have
\begin{eqnarray}\label{lamboundbyZ}
\begin{aligned}
&\frac{1}{\tau\beta}\|\Lambda^{k+1} - \Lambda^k\|_F^2
\leq {\textstyle\max\left\{\frac{1}{\tau}, \frac{\tau^2}{1+\tau-\tau^2}\right\}\cdot\frac{\lambda_{\max}^2}{\beta}}\|Z^{k+1} - Z^k\|_F^2 \\
&~~+\theta(\tau)\beta\left(\|\mathcal{B}(L^k)+\mathcal{C}(S^k)-Z^k\|_F^2 - \|\mathcal{B}(L^{k+1})+\mathcal{C}(S^{k+1})-Z^{k+1}\|_F^2\right).
\end{aligned}
\end{eqnarray}
%Substituting \eqref{lamboundbyZ} into \eqref{first_term_tmp}, we obtain further that
%\begin{eqnarray}\label{first_term}
%\begin{aligned}
%&\quad\Theta_{\tau,\beta}(L^{k+1}, S^{k+1}, Z^{k+1}, \Lambda^{k+1}) - \Theta_{\tau,\beta} (L^{k+1}, S^{k+1}, Z^{k+1}, \Lambda^k) \\
%&\leq {\textstyle\max\left\{\frac{1}{\tau}, \frac{\tau^2}{1+\tau-\tau^2}\right\}\cdot\frac{\lambda_{\max}^2}{\beta}}\|Z^{k+1} - Z^k\|_F^2\\
%&\quad+\theta(\tau)\beta\left(\|\mathcal{B}(L^k)+\mathcal{C}(S^k)-Z^k\|_F^2 - \|\mathcal{B}(L^{k+1})+\mathcal{C}(S^{k+1})-Z^{k+1}\|_F^2\right).
%\end{aligned}
%\end{eqnarray}

Next, note that the function $Z \mapsto \mathcal{L}_{\beta} (L^{k+1}, S^{k+1}, Z, \Lambda^{k})$ is strongly convex with modulus at least $\lambda_{\min}+\beta$. Using this fact and the definition of $Z^{k+1}$ as a minimizer in \eqref{Zsubpro}, we see that
\begin{eqnarray}\label{second_term}
\hspace{-6mm}\begin{aligned}
&\quad\Theta_{\tau,\beta} (L^{k+1}, S^{k+1}, Z^{k+1}, \Lambda^{k}) - \Theta_{\tau,\beta}(L^{k+1}, S^{k+1}, Z^k, \Lambda^k) \\
&=\mathcal{L}_{\beta} (L^{k+1}, S^{k+1}, Z^{k+1}, \Lambda^{k}) - \mathcal{L}_{\beta} (L^{k+1}, S^{k+1}, Z^k, \Lambda^k)  \\
&~~+\theta(\tau)\beta\left(\|\mathcal{B}(L^{k+1})+\mathcal{C}(S^{k+1})-Z^{k+1}\|_F^2
-\|\mathcal{B}(L^{k+1})+\mathcal{C}(S^{k+1})-Z^{k}\|_F^2\right)   \\
&\leq\theta(\tau)\beta\left(\|\mathcal{B}(L^{k+1})+\mathcal{C}(S^{k+1})-Z^{k+1}\|_F^2
-\|\mathcal{B}(L^{k+1})+\mathcal{C}(S^{k+1})-Z^{k}\|_F^2\right) \\
&~~-{\textstyle\frac{\lambda_{\min}+\beta}{2}}\|Z^{k+1} - Z^k\|_F^2.
\end{aligned}
\end{eqnarray}

Moreover, using the fact that $S^{k+1}$ is a minimizer in \eqref{Ssubpro}, we have
\begin{eqnarray}\label{third_term}
\begin{aligned}
&\quad\Theta_{\tau,\beta} (L^{k+1}, S^{k+1}, Z^{k}, \Lambda^{k}) - \Theta_{\tau,\beta}(L^{k+1}, S^k, Z^k, \Lambda^k) \\
&= \mathcal{L}_{\beta} (L^{k+1}, S^{k+1}, Z^{k}, \Lambda^{k}) - \mathcal{L}_{\beta} (L^{k+1}, S^k, Z^k, \Lambda^k)  \\
&\quad+\theta(\tau)\beta\left(\|\mathcal{B}(L^{k+1})+\mathcal{C}(S^{k+1})-Z^{k}\|_F^2
-\|\mathcal{B}(L^{k+1})+\mathcal{C}(S^{k})-Z^{k}\|_F^2\right)   \\
&\leq \theta(\tau)\beta\left(\|\mathcal{B}(L^{k+1})+\mathcal{C}(S^{k+1})-Z^{k}\|_F^2
-\|\mathcal{B}(L^{k+1})+\mathcal{C}(S^{k})-Z^{k}\|_F^2\right).
\end{aligned}
\end{eqnarray}

Finally, note that $L \mapsto \mathcal{L}_{\beta}(L, S^{k}, Z^{k}, \Lambda^{k})$ is strongly convex with modulus at least $\sigma\beta$ from (a1) in Assumption \ref{assum}. From this, we can similarly obtain
\begin{eqnarray}\label{fourth_term}
\begin{aligned}
&\quad\Theta_{\tau,\beta} (L^{k+1}, S^{k}, Z^{k}, \Lambda^{k}) - \Theta_{\tau,\beta}(L^{k}, S^{k}, Z^k, \Lambda^k) \\
&\leq\theta(\tau)\beta\left(\|\mathcal{B}(L^{k+1})+\mathcal{C}(S^{k})-Z^{k}\|_F^2-\|\mathcal{B}(L^{k})
+\mathcal{C}(S^{k})-Z^{k}\|_F^2\right) \\
&~~-\frac{\sigma\beta}{2}\|L^{k+1}-L^k\|_F^2.
\end{aligned}
\end{eqnarray}
Thus, summing \eqref{first_term_tmp}, \eqref{lamboundbyZ}, \eqref{second_term}, \eqref{third_term} and \eqref{fourth_term}, we obtain \eqref{difofp}.

Now, suppose in addition that $\beta \geq -\frac{\lambda_{\min}}{2}+\frac{1}{2}\sqrt{\lambda_{\min}^2+\max\left\{\frac{1}{\tau}, \frac{\tau^2}{1+\tau-\tau^2}\right\} \cdot 8\lambda_{\max}^2}$. Then it is easy to check that
\begin{eqnarray*}
\max\left\{\frac{1}{\tau}, \frac{\tau^2}{1+\tau-\tau^2}\right\}\cdot\frac{\lambda_{\max}^2}{\beta} - \frac{\lambda_{\min}+\beta}{2}\leq0.
\end{eqnarray*}
Hence we see from \eqref{difofp} that
\begin{eqnarray*}
\Theta_{\tau,\beta} (L^{k+1}, S^{k+1}, Z^{k+1}, \Lambda^{k+1}) - \Theta_{\tau,\beta} (L^k, S^k, Z^k, \Lambda^k) \leq 0,
\end{eqnarray*}
which means that $\{\Theta_{\tau,\beta} (L^k, S^k, Z^k, \Lambda^k)\}_{k=1}^\infty$ is decreasing. This completes the proof.
\endproof

We next show that the sequence generated by Algorithm 1 is bounded if $\beta$ is chosen above a computable threshold, under (a1) and (a3) in Assumption~\ref{assum}. For notational simplicity, from now on, we let
\begin{eqnarray}\label{betalowbound}
\hspace{-1cm}{\textstyle\bar{\beta}:= \max\left\{\max\{1/\tau,\tau\}\cdot\lambda_{\max},
-\frac{\lambda_{\min}}{2}\!+\frac{1}{2}\sqrt{\lambda_{\min}^2\!+\!\max\left\{\frac{1}{\tau}, \frac{\tau^2}{1+\tau-\tau^2}\right\} \cdot 8\lambda_{\max}^2}\right\}.}
\end{eqnarray}

\begin{proposition}[\textbf{Boundedness of sequence generated by ADMM}]\label{boundnesslemma}
Suppose that $0 < \tau < \frac{1+\sqrt{5}}{2}$ and $\beta>\bar{\beta}$. If (a1) and (a3) in Assumption \ref{assum} hold, then a sequence $\{(L^k, S^k, Z^k, \Lambda^k)\}_{k=1}^{\infty}$ generated by Algorithm 1 is bounded.
\end{proposition}
\beginproof
With our choice of $\beta$ and (a1) in Assumption \ref{assum}, we see immediately from Lemma \ref{Phidelemma} that the sequence $\{\Theta_{\tau,\beta}(L^k, S^k, Z^k, \Lambda^k)\}_{k=1}^{\infty}$ is decreasing. This together with (a3) in Assumption \ref{assum} shows that, for $k \geq 1$,
\begin{eqnarray}\label{Phivalbound1}
\begin{aligned}
&\quad h_0~>~\Theta_{\tau,\beta} (L^1, S^1, Z^1, \Lambda^1) ~\geq ~\Theta_{\tau,\beta} (L^k, S^k, Z^k, \Lambda^k) \\
&= \Psi(L^k)+\Phi(S^k)+\frac{1}{2}\|D-\mathcal{A}(Z^k)\|_F^2-\langle \Lambda^k, \mathcal{B}(L^k)+\mathcal{C}(S^k)-Z^k \rangle \\
&~~~+\left(1+2\theta(\tau)\right)\frac{\beta}{2}\|\mathcal{B}(L^k)+\mathcal{C}(S^k)-Z^k\|_F^2 \\
&= \Psi(L^k)+\Phi(S^k)+\frac{1}{2}\|D-\mathcal{A}(Z^k)\|_F^2 + \frac{\beta}{2}\|\mathcal{B}(L^k)+\mathcal{C}(S^k)-Z^k-\frac{1}{\beta}\Lambda^k\|_F^2 \\
&~~~-\frac{1}{2\beta}\|\Lambda^k\|_F^2+\theta(\tau)\beta\|\mathcal{B}(L^k)+\mathcal{C}(S^k)-Z^k\|_F^2,
\end{aligned}
\end{eqnarray}
where the last equality is obtained by completing the square. We next derive an upper bound for $\left\|\Lambda^k\right\|_F^2$. We start by substituting \eqref{lamopt} into \eqref{zopt} and rearranging terms to obtain
\begin{eqnarray}\label{z_optcon}
&0 = \mathcal{A}^*(\mathcal{A}(Z^k)-D) + \Lambda^{k-1} + \frac{1}{\tau}(\Lambda^k - \Lambda^{k-1})\nonumber\\
\Longrightarrow&-\tau\Lambda^k = \tau \mathcal{A}^*(\mathcal{A}(Z^k)-D) + \left(1 - \tau\right)(\Lambda^k - \Lambda^{k-1}).
\end{eqnarray}
We now consider two different cases:
\begin{itemize}
\item For $0<\tau\leq 1$, it follows from the convexity of $\|\cdot\|_F^2$ and \eqref{z_optcon} that
      \begin{eqnarray*}
      \|-\tau\Lambda^k\|_F^2 &\leq& \tau\|\mathcal{A}^{*}(\mathcal{A}(Z^k)-D)\|_F^2 + (1-\tau)\|\Lambda^k - \Lambda^{k-1}\|_F^2 \\
      &\leq& \tau\lambda_{\max}\|\mathcal{A}(Z^k)-D\|_F^2 + (1-\tau)\|\Lambda^k - \Lambda^{k-1}\|_F^2 \\
      &=& \tau\lambda_{\max}\|\mathcal{A}(Z^k)-D\|_F^2 + (1-\tau)\tau^2\beta^2\|\mathcal{B}(L^k) + \mathcal{C}(S^k) - Z^k\|_F^2,
      \end{eqnarray*}
      where the equality follows from \eqref{lamopt}. Then, we have
      \begin{eqnarray}\label{lambound1}
      \hspace{-1cm}\|\Lambda^k\|_F^2 \leq \frac{\lambda_{\max}}{\tau}\|\mathcal{A}(Z^k)-D\|_F^2 + (1-\tau)\beta^2\|\mathcal{B}(L^k) + \mathcal{C}(S^k) - Z^k\|_F^2.
      \end{eqnarray}

\item For $1 < \tau < \frac{1+\sqrt{5}}{2}$, by dividing $-\tau$ from both sides of \eqref{zopt}, we obtain
      \begin{eqnarray*}
      \Lambda^k = \frac{1}{\tau}~\tau \mathcal{A}^*(D-\mathcal{A}(Z^k)) + \left(1-\frac{1}{\tau}\right)(\Lambda^{k}-\Lambda^{k-1}).
      \end{eqnarray*}
      Then, since $0<\frac{1}{\tau}<1$, using the convexity of $\|\cdot\|_F^2$ and \eqref{lamopt}, we have
      \begin{eqnarray}\label{lambound2}
      \begin{aligned}
      \|\Lambda^k\|_F^2 &\leq~ \frac{1}{\tau}\|\tau \mathcal{A}^*(D-\mathcal{A}(Z^k))\|_F^2 + \left(1-\frac{1}{\tau}\right)\|\Lambda^k - \Lambda^{k-1}\|_F^2 \\
      &\leq~\tau\lambda_{\max}\|D - \mathcal{A}(Z^k)\|_F^2 + (\tau-1)\tau\beta^2\|\mathcal{B}(L^k)+\mathcal{C}(S^k)-Z^k\|_F^2.
      \end{aligned}
      \end{eqnarray}
\end{itemize}
Thus, combining \eqref{lambound1} and \eqref{lambound2}, we have
\begin{eqnarray}\label{lambound}
\|\Lambda^k\|_F^2 &\leq& \max\{1/\tau,\tau\}\cdot\lambda_{\max}\|D - \mathcal{A}(Z^k)\|_F^2 \nonumber\\
&&+ \max\{1-\tau, (\tau-1)\tau\}\beta^2\|\mathcal{B}(L^k)+\mathcal{C}(S^k)-Z^k\|_F^2 \nonumber\\
\Longrightarrow ~~ -\frac{1}{2\beta}\|\Lambda^k\|_F^2 &\geq& {\textstyle-\frac{\max\{1/\tau,\tau\}\lambda_{\max}}{2\beta}}\|D - \mathcal{A}(Z^k)\|_F^2 \\
&&\quad-{\textstyle\frac{\max\{1-\tau, (\tau-1)\tau\}\beta}{2}}\|\mathcal{B}(L^k)+\mathcal{C}(S^k)-Z^k\|_F^2. \nonumber
\end{eqnarray}
Substituting \eqref{lambound} into \eqref{Phivalbound1}, we have
\begin{eqnarray}\label{Phivalbound}
\begin{aligned}
&h_0 > \Theta_{\tau,\beta} (L^k, S^k, Z^k, \Lambda^k) \ge \Psi(L^k)+\Phi(S^k)\\
&\quad+\frac{1}{2}\left(1-\max\{1/\tau,\tau\}\cdot\frac{\lambda_{\max}}{\beta}\right)\|D-\mathcal{A}(Z^k)\|_F^2\\
&\quad+\frac{\beta}{2}\|\mathcal{B}(L^k)+\mathcal{C}(S^k)-Z^k-\frac{1}{\beta}\Lambda^k\|_F^2\\
&\quad+\left[2\theta(\tau)-\max\{1-\tau, (\tau-1)\tau\}\right]\cdot\frac{\beta}{2}\|\mathcal{B}(L^k)+\mathcal{C}(S^k)-Z^k\|_F^2.
\end{aligned}
\end{eqnarray}

With \eqref{Phivalbound} established, we are now ready to prove the boundedness of the sequence. We start with the observation that for $0 < \tau < \frac{1+\sqrt{5}}{2}$ and $\beta>\bar{\beta}$, we always have
\begin{eqnarray}\label{coef1}
1-\max\{1/\tau,\tau\}\cdot\frac{\lambda_{\max}}{\beta} > 0
\end{eqnarray}
and
\begin{eqnarray}\label{coef2}
{\textstyle2\theta(\tau)-\max\{1-\tau, (\tau-1)\tau\}=\left\{
\begin{aligned}
&1-\tau > 0,  &&\mathrm{for}~0<\tau<1, \\
&0,  &&\mathrm{for}~\tau=1, \\
&{\textstyle\frac{\tau(\tau-1)(\tau^2+\tau-1)}{1+\tau-\tau^2}} > 0,  &&\mathrm{for}~1 < \tau < {\textstyle\frac{1+\sqrt{5}}{2}},
\end{aligned}\right.}
\end{eqnarray}
where $\theta(\tau)$ is defined in \eqref{deftheta}. Then we consider two cases:
\begin{itemize}
\item For $\tau \in \left(0, 1\right)\cup\left(1,\frac{1+\sqrt{5}}{2}\right)$, it follows from \eqref{Phivalbound}, \eqref{coef1}, \eqref{coef2}, and the nonnegativity of $\Psi$ and $\Phi$ that
    $\{\|D-\mathcal{A}(Z^k)\|_F\}$, $\{\|\mathcal{B}(L^k)+\mathcal{C}(S^k)-Z^k-\frac{1}{\beta}\Lambda^k\|_F\}$ and $\{\|\mathcal{B}(L^k)+\mathcal{C}(S^k)-Z^k\|_F\}$ are bounded; and moreover,
    \begin{eqnarray*}%\label{Sbound}
    \Psi(L^k)+\Phi(S^k)<h_0.
    \end{eqnarray*}
    The boundedness of $\{L^k\}$ and $\{S^k\}$ follows immediately from this last relation. Furthermore, $\{\Lambda^k\}$ is bounded since
    \begin{eqnarray*}%\label{lambdabound}
    \|\Lambda^k\|_F \leq \beta\|\mathcal{B}(L^k)+\mathcal{C}(S^k)-Z^k-\frac{1}{\beta}\Lambda^k\|_F + \beta\|\mathcal{B}(L^k)+\mathcal{C}(S^k)-Z^k\|_F.
    \end{eqnarray*}
    Finally, we obtain the boundedness of $\{Z^k\}$ from
    \begin{eqnarray}\label{Zbound}
    \begin{aligned}
    \|Z^k\|_F &\leq \|\mathcal{B}(L^k)+\mathcal{C}(S^k)-Z^k-\frac{1}{\beta}\Lambda^k\|_F + \|\mathcal{B}(L^k)\|_F\\
    &\quad+\|\mathcal{C}(S^k)\|_F + \frac{1}{\beta}\|\Lambda^k\|_F.
    \end{aligned}
    \end{eqnarray}

\item For $\tau=1$, it follows from \eqref{Phivalbound}, \eqref{coef1}, \eqref{coef2}, and the nonnegativity of $\Psi$ and $\Phi$ that $\{\|D-\mathcal{A}(Z^k)\|_F\}$ and $\{\|\mathcal{B}(L^k)+\mathcal{C}(S^k)-Z^k-\frac{1}{\beta}\Lambda^k\|_F\}$ are bounded; and moreover $\Psi(L^k)+\Phi(S^k)<h_0$, from which we see immediately that $\{L^k\}$ and $\{S^k\}$ are bounded. The boundedness of $\{\Lambda^k\}$ now follows from \eqref{z_optcon} with $\tau=1$, i.e., $\Lambda^k = \mathcal{A}^*(D-\mathcal{A}(Z^k))$. The boundedness of $\{Z^k\}$ again follows from \eqref{Zbound}.

\end{itemize}

This completes the proof.
\endproof

We are now ready to prove our first global convergence result for Algorithm 1, which also characterizes the cluster point of the sequence generated.

\begin{theorem}[\textbf{Global subsequential convergence}]\label{convergencethe}
Suppose that $0 < \tau < \frac{1+\sqrt{5}}{2}$ and $\beta > \bar{\beta}$. If Assumption \ref{assum} holds, then
\begin{enumerate}[{\rm (i)}]
\item $\lim\limits_{k\rightarrow \infty} \|L^{k+1} - L^k\|_F + \|S^{k+1} - S^k\|_F + \|Z^{k+1} - Z^k\|_F+\|\Lambda^{k+1} - \Lambda^k\|_F  = 0$;

\item Any cluster point $(L^*, S^*, Z^*, \Lambda^*)$ of a sequence $\{(L^k, S^k, Z^k, \Lambda^k)\}$ generated by Algorithm 1 is a stationary point of \eqref{generalmodel1}.
\end{enumerate}
\end{theorem}
\beginproof
The boundedness of the sequence $\{(L^k, S^k, Z^k, \Lambda^k)\}$ follows immediately from Proposition \ref{boundnesslemma} and thus a cluster point exists. We now prove statement (i).

Suppose that $(L^*, S^*, Z^*, \Lambda^*)$ is a cluster point of the sequence $\{(L^k, S^k, Z^k, \Lambda^k)\}$ and let $\{(L^{k_i}, S^{k_i}$, $Z^{k_i}, \Lambda^{k_i})\}$ be a convergent subsequence such that
\begin{eqnarray*}
\lim_{i\rightarrow\infty} (L^{k_i}, S^{k_i}, Z^{k_i}, \Lambda^{k_i}) = (L^*, S^*, Z^*, \Lambda^*).
\end{eqnarray*}
By summing \eqref{difofp} from $k=1$ to $k = k_i-1$, we have
\begin{eqnarray}\label{Phivalde}
\begin{aligned}
&\quad \Theta_{\tau,\beta} (L^{k_i}, S^{k_i}, Z^{k_i}, \Lambda^{k_i}) - \Theta_{\tau,\beta} (L^1, S^1, Z^1, \Lambda^1) \\
&\leq -C\sum\limits^{k_i-1}_{k=1}\|Z^{k+1} - Z^k\|_F^2-\frac{\sigma\beta}{2}\sum\limits^{k_i-1}_{k=1}\|L^{k+1} - L^k\|_F^2,
\end{aligned}
\end{eqnarray}
where $C:=\frac{\lambda_{\min}+\beta}{2}-\max\left\{\frac{1}{\tau}, \frac{\tau^2}{1+\tau-\tau^2}\right\}\cdot\frac{\lambda_{\max}^2}{\beta} > 0$ (since $\beta>\bar{\beta}$). Passing to the limit in \eqref{Phivalde} and rearranging terms in the resulting relation, we obtain
\begin{eqnarray*}
\begin{aligned}
&C\sum\limits^{\infty}_{k=1}\|Z^{k+1} - Z^k\|_F^2 + \frac{\sigma\beta}{2}\sum\limits^{\infty}_{k=1}\|L^{k+1} - L^k\|_F^2 \\
&\leq \Theta_{\tau,\beta} (L^1, S^1, Z^1, \Lambda^1) - \Theta_{\tau,\beta} (L^*, S^*, Z^*, \Lambda^*) < \infty,
\end{aligned}
\end{eqnarray*}
where the last inequality follows from the properness of $\Psi$ and $\Phi$. This together with $C > 0$ and $\sigma>0$ implies that
\begin{eqnarray*}
\sum\limits^{\infty}_{k=1}\|Z^{k+1} - Z^k\|_F^2 < \infty \quad \mathrm{and} \quad \sum\limits^{\infty}_{k=1}\|L^{k+1} - L^k\|_F^2 <\infty.
\end{eqnarray*}
Hence, we have
\begin{eqnarray}\label{limit1}
Z^{k+1} - Z^k  \rightarrow 0, ~L^{k+1} - L^k \rightarrow 0.
\end{eqnarray}
Next, by summing both sides of \eqref{lamboundbyZ} from $k=1$ to $k=k_i$ and passing to the limit, we have
\begin{eqnarray*}
\begin{aligned}
&\sum\limits^{\infty}_{k=1}\|\Lambda^{k+1} - \Lambda^k\|_F^2
\leq{\textstyle\max\left\{\frac{1}{\tau}, \frac{\tau^2}{1+\tau-\tau^2}\right\}\cdot\tau\lambda_{\max}^2}\sum\limits^{\infty}_{k=1}\|Z^{k+1} - Z^k\|_F^2 \\
&\quad +\theta(\tau)\tau\beta^2\left(\|\mathcal{B}(L^1)+\mathcal{C}(S^1)-Z^1\|_F^2 - \liminf\limits_{k\rightarrow\infty}\|\mathcal{B}(L^{k+1})+\mathcal{C}(S^{k+1})-Z^{k+1}\|_F^2\right),
\end{aligned}
\end{eqnarray*}
from which we conclude that
\begin{eqnarray}\label{limit2}
\Lambda^{k+1} - \Lambda^k  \rightarrow 0.
\end{eqnarray}
Finally, we have $S^{k+1} - S^k  \rightarrow 0$ from \eqref{limit1}, \eqref{limit2}, \eqref{lamopt} and (a1) in Assumption \ref{assum}. This proves statement (i).

We next prove statement (ii). From the lower semicontinuity of $\Theta_{\tau,\beta}$ (since $\Psi$ and $\Phi$ are lower semicontinuous), we have
\begin{eqnarray}\label{semicont1}
\begin{aligned}
&\liminf\limits_{i \rightarrow \infty}~\Theta_{\tau,\beta} (L^{k_i+1}, S^{k_i+1}, Z^{k_i}, \Lambda^{k_i})
\geq~\Psi(L^*)+\Phi(S^*)+\frac{1}{2}\|D-\mathcal{A}(Z^*)\|_F^2  \\
&\quad-\langle \Lambda^*,~\mathcal{B}(L^*)+\mathcal{C}(S^*)-Z^*\rangle+ \left(1+2\theta(\tau)\right)\frac{\beta}{2}\|\mathcal{B}(L^*)+\mathcal{C}(S^*)-Z^*\|_F^2.
\end{aligned}
\end{eqnarray}
On the other hand, from the definition of $S^{k_i+1}$ as a minimizer in \eqref{Ssubpro}, we have
\begin{eqnarray*}
\begin{aligned}
&\Theta_{\tau,\beta} (L^{k_i+1}, S^{k_i+1}, Z^{k_i}, \Lambda^{k_i})\leq \Theta_{\tau,\beta} (L^{k_i+1}, S^{*}, Z^{k_i}, \Lambda^{k_i}) \\
&\quad+\theta(\tau)\beta\left(\|\mathcal{B}(L^{k_i+1})+\mathcal{C}(S^{k_i+1})-Z^{k_i}\|_F^2
-\|\mathcal{B}(L^{k_i+1})+\mathcal{C}(S^{*})-Z^{k_i}\|_F^2\right).
\end{aligned}
\end{eqnarray*}
Taking limit in above equality, and invoking statement (i) and (a2) in Assumption \ref{assum}, we see that
\begin{eqnarray}\label{semicont2}
\begin{aligned}
&\limsup\limits_{i \rightarrow \infty}~\Theta_{\tau,\beta} (L^{k_i+1}, S^{k_i+1}, Z^{k_i}, \Lambda^{k_i})\leq \Psi(L^*)+\Phi(S^*) + \frac{1}{2}\|D-\mathcal{A}(Z^*)\|_F^2 \\
&\quad -\langle \Lambda^*,~\mathcal{B}(L^*)+\mathcal{C}(S^*)-Z^* \rangle + \left(1+2\theta(\tau)\right)\frac{\beta}{2}\|\mathcal{B}(L^*)+\mathcal{C}(S^*)-Z^*\|_F^2.
\end{aligned}
\end{eqnarray}
Then, combining \eqref{semicont1} and \eqref{semicont2}, we see that
\begin{eqnarray*}
\begin{aligned}
&\lim\limits_{i \rightarrow \infty}~\Theta_{\tau,\beta} (L^{k_i+1}, S^{k_i+1}, Z^{k_i}, \Lambda^{k_i})= \Psi(L^*)+\Phi(S^*) + \frac{1}{2}\|D-\mathcal{A}(Z^*)\|_F^2 \\
&\quad -\langle \Lambda^*,~\mathcal{B}(L^*)+\mathcal{C}(S^*)-Z^* \rangle + \left(1+2\theta(\tau)\right)\frac{\beta}{2}\|\mathcal{B}(L^*)+\mathcal{C}(S^*)-Z^*\|_F^2,
\end{aligned}
\end{eqnarray*}
which, together with (a2) in Assumption \ref{assum}, $L^{k+1}-L^k\rightarrow0$, $S^{k+1}-S^k\rightarrow0$ and the definition of $\Theta_{\tau,\beta}$, implies that
\begin{eqnarray}\label{limitS}
\lim\limits_{i \rightarrow \infty} \Phi(S^{k_i+1}) = \Phi(S^{*}).
\end{eqnarray}

Thus, passing to the limit in \eqref{lopt}-\eqref{lamopt} along $\{(L^{k_i}, S^{k_i}, Z^{k_i}, \Lambda^{k_i})\}$ and invoking statement (i), \eqref{limitS} and \eqref{robust}, we see that
\begin{eqnarray}\label{optlimit}
\left\{\begin{aligned}
&0 \in \partial \Psi(L^*)-\mathcal{B}^*(\Lambda^*)+\beta\mathcal{B}^*\left(\mathcal{B}(L^{*})+\mathcal{C}(S^*)-Z^*\right), \\
&0 \in \partial \Phi(S^{*})-\mathcal{C}^*(\Lambda^*)+\beta\mathcal{C}^*\left(\mathcal{B}(L^{*})+\mathcal{C}(S^*)-Z^*\right), \\
&0 = \mathcal{A}^{*}(\mathcal{A}(Z^{*})-D) + \Lambda^* - \beta(\mathcal{B}(L^{*})+\mathcal{C}(S^{*})-Z^{*}),\\
&\mathcal{B}(L^*)+\mathcal{C}(S^*) = Z^*.
\end{aligned}\right.
\end{eqnarray}
Rearranging terms in \eqref{optlimit}, it is not hard to obtain
\begin{eqnarray*}
\left\{\begin{aligned}
&0\in \partial \Psi(L^*)+\mathcal{B}^*\mathcal{A}^*\left(\mathcal{A}(\mathcal{B}(L^*)+\mathcal{C}(S^*))-D\right),\\
&0\in \partial \Phi(S^*)+\mathcal{C}^*\mathcal{A}^*\left(\mathcal{A}(\mathcal{B}(L^*)+\mathcal{C}(S^*))-D\right).
\end{aligned}\right.
\end{eqnarray*}
This shows that $(L^*, S^*, Z^*, \Lambda^*)$ is a stationary point of \eqref{generalmodel1}. This completes the proof.
\endproof

\begin{remark}[\textbf{Comments on the computable threshold}]
From the above discussions, we establish under Assumption~\ref{assum} the convergence of the ADMM with $0<\tau<\frac{1+\sqrt{5}}{2}$ when the penalty parameter $\beta$ is chosen above a computable threshold $\bar{\beta}$ which depends on $\tau$. The existence of this kind of threshold is also obtained in the recent studies \cite{ah2014,hlr2014,lp2014,wxx2014,wcx2015} on the nonconvex ADMM and its variants with $\tau=1$. In Fig.\,\ref{beta0}, we plot $\bar{\beta}$ against $\tau$ with $\mathcal{A}$ being the identity map (hence, $\lambda_{\max}=\lambda_{\min}=1$). It is not hard to see from Fig.\,\ref{beta0} that for a given penalty parameter $\beta > 1$, we can always choose a dual step-size $\tau$ from an interval containing 1 so that the corresponding ADMM is convergent.
\begin{figure}[ht]
\centering
\includegraphics[height=7cm]{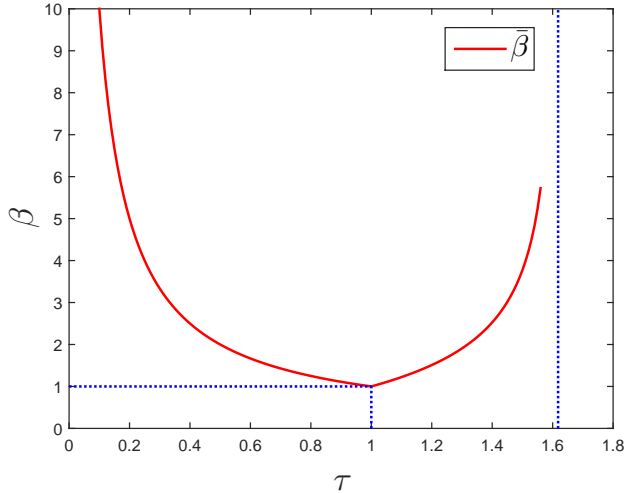}
\caption{The computable threshold $\bar\beta$ for $0<\tau<\frac{1+\sqrt{5}}{2}$.}\label{beta0}
\end{figure}
\end{remark}

\begin{remark}[\textbf{Practical computation consideration on penalty parameter}]\label{rem4penalty}
In computation, for a $0<\tau<\frac{1+\sqrt{5}}{2}$, the $\bar{\beta}$ in \eqref{betalowbound} may be too large and hence fixing a $\beta$ close to it can lead to slow convergence. As in \cite{sty2015,lp2014Douglas}, one could possibly accelerate the algorithm by initializing the algorithm with a small $\beta$ (less than $\bar{\beta}$) and then increasing the $\beta$ by a constant ratio until $\beta>\bar{\beta}$ if the sequence generated becomes unbounded or the successive change does not vanish sufficiently fast. Clearly, after at most finitely many increases, the penalty parameter $\beta$ gets above the threshold $\bar\beta$ and the convergence of the resulting algorithm is guaranteed by Theorem~\ref{convergencethe} under Assumption~\ref{assum}. On the other hand, if $\beta$ is never increased, this means that the successive change goes to zero and the sequence is bounded. Then it is routine to show that any cluster point is a stationary point if $\Phi$ is continuous in its domain.
\end{remark}

Under the additional assumption that the potential function $\Theta_{\tau,\beta}$ is a KL function, we show in the next theorem that the whole sequence generated by Algorithm 1 is convergent if $\beta$ is greater than a computable threshold, again under Assumption~\ref{assum}. Our proof makes use of the uniformized KL property; see Proposition~\ref{uniKL}. This technique was previously used in \cite{bst2014} to prove the convergence of the proximal alternating linearized minimization algorithm for nonconvex and nonsmooth problems, and later in \cite{wcx2015,wxx2014} to prove the global convergence of the Bregman ADMM with $\tau = 1$. Our analysis, though follows a similar line of arguments as in \cite{wcx2015,wxx2014}, is much more intricate. This is because when $\tau \neq 1$, the successive change in the dual variable cannot be controlled solely by the successive changes in the primal variables.

\begin{theorem}[\textbf{Global convergence of the whole sequence}]\label{thm42}
Let $0<\tau<\frac{1+\sqrt{5}}{2}$ and $\beta > \bar \beta$.
%\begin{eqnarray*}
%\beta > \max\left\{\max\{1/\tau,\tau\}\cdot\lambda_{\max},~-\frac{\lambda_{\min}}{2}+\frac{1}{2}\sqrt{\lambda_{\min}^2+\max\left\{\frac{1}{\tau}, \frac{\tau^2}{1+\tau-\tau^2}\right\} \cdot 8\lambda_{\max}^2}\right\}.
%\end{eqnarray*}
Suppose in addition that Assumption~\ref{assum} holds and the potential function $\Theta_{\tau,\beta}(\cdot)$ is a KL function. Then, the sequence $\{(L^k$, $S^k$, $Z^k$, $\Lambda^k)\}_{k=1}^{\infty}$ generated by Algorithm 1 converges to a stationary point of \eqref{generalmodel1}.
\end{theorem}
\beginproof
In view of Theorem~\ref{convergencethe}, we only need to show that the sequence is convergent. We start by noting from \eqref{Phivalbound}, \eqref{coef1} and \eqref{coef2} that $\{\Theta_{\tau,\beta} (L^k$, $S^k$, $Z^k$, $\Lambda^k)\}_{k=1}^{\infty}$ is bounded below. Since this sequence is also decreasing from Theorem~\ref{Phidelemma}, we conclude that $\lim_{k\rightarrow\infty}\Theta_{\tau,\beta}(L^k, S^k, Z^k, \Lambda^k) =: \theta^*$ exists. In the following, we will consider two cases.

\textbf{Case 1)} Suppose first that $\Theta_{\tau,\beta}(L^{N}, S^{N}, Z^{N}, \Lambda^{N})=\theta^*$ for some $N \geq 1$. Since $\{\Theta_{\tau,\beta} (L^k$, $S^k$, $Z^k$, $\Lambda^k)\}_{k=1}^{\infty}$ is decreasing, we must have $\Theta_{\tau,\beta}(L^k, S^k, Z^k, \Lambda^k)=\theta^*$ for all $k \geq N$. Then, it follows from \eqref{difofp} that $L^{N+t}=L^{N}$ and $Z^{N+t}=Z^{N}$ for all $t \geq 0$. Hence, $\{L^k\}$ and $\{Z^k\}$ converge finitely. Moreover, from \eqref{lamdiff}, we have
\begin{eqnarray*}
\|\Lambda^{k+1} - \Lambda^{k}\|_F = |1-\tau|\cdot\|\Lambda^k - \Lambda^{k-1}\|_F = \cdots = |1-\tau|^{k+1-N}\cdot\|\Lambda^N - \Lambda^{N-1}\|_F
\end{eqnarray*}
for all $k\geq N$. Since $0<\tau<\frac{1+\sqrt{5}}{2}$, we have $0<1-|1-\tau|\leq1$ and hence we see further that
\begin{eqnarray}\label{sumoflam}
\sum\limits_{k=N}^{\infty} \|\Lambda^{k+1}-\Lambda^{k}\|_F \leq \frac{1}{1-|1-\tau|}\|\Lambda^N-\Lambda^{N-1}\|_F < \infty,
\end{eqnarray}
which implies the convergence of $\{\Lambda^k\}$. Additionally, for all $k\geq N$, we have
\begin{eqnarray*}%\label{SboundbyLam}
\|S^{k+1}-S^k\|_F
&\leq&\frac{1}{\sqrt{\sigma'}}\|\mathcal{C}(S^{k+1})-\mathcal{C}(S^k)\|_F \nonumber \\
&=&\frac{1}{\sqrt{\sigma'}}\left\|\frac{1}{\tau\beta}(\Lambda^{k}-\Lambda^{k+1}) - \frac{1}{\tau\beta}(\Lambda^{k-1}-\Lambda^k)\right\|_F\nonumber\\
&\leq&\frac{1}{\tau\beta\sqrt{\sigma'}}\|\Lambda^{k+1}-\Lambda^k\|_F
+\frac{1}{\tau\beta\sqrt{\sigma'}}\|\Lambda^{k}-\Lambda^{k-1}\|_F,
\end{eqnarray*}
where the first inequality follows from (a1) in Assumption \ref{assum} and the equality follows from \eqref{lamopt}. This together with \eqref{sumoflam}, implies that $\sum_{k=N}^{\infty} \|S^{k+1}-S^{k}\|_F<\infty$. Thus, $\{S^k\}$ is also convergent. Consequently, we see that $\{(L^k$, $S^k$, $Z^k$, $\Lambda^k)\}_{k=1}^{\infty}$ is a convergent sequence in this case.

\textbf{Case 2)} From now on, we consider the case where $\Theta_{\tau,\beta}(L^k$, $S^k$, $Z^k$, $\Lambda^k)>\theta^*$ for all $k\geq1$. In this case, we will divide the proof into three steps: \textbf{1.} we first prove that $\Theta_{\tau,\beta}$ is constant on the set of cluster points of the sequence $\{(L^k$, $S^k$, $Z^k$, $\Lambda^k)\}_{k=1}^{\infty}$ and then apply the uniformized KL property; \textbf{2.} we bound the distance from 0 to $\partial\Theta_{\tau,\beta}(L^{k}, S^{k}, Z^{k}, \Lambda^{k})$; \textbf{3.} we show that the sequence $\{(L^k$, $S^k$, $Z^k$, $\Lambda^k)\}_{k=1}^{\infty}$ is a Cauchy sequence and hence is convergent. The complete proof is presented as follows.

\textit{Step \textbf{1.}} We recall from Proposition~\ref{boundnesslemma} that the sequence $\{(L^k$, $S^k$, $Z^k$, $\Lambda^k)\}_{k=1}^{\infty}$ generated by Algorithm 1 is bounded and hence must have at least one cluster point. Let $\Gamma$ denote the set of cluster points of $\{(L^k$, $S^k$, $Z^k$, $\Lambda^k)\}_{k=1}^{\infty}$. We will show that $\Theta_{\tau,\beta}$ is constant on $\Gamma$.

To this end, take any $(L^*, S^*, Z^*, \Lambda^*)\in \Gamma$ and consider a convergent subsequence $\{(L^{k_i}, S^{k_i}, Z^{k_i}, \Lambda^{k_i})\}$ with $\lim_{i\rightarrow\infty}(L^{k_i}, S^{k_i}, Z^{k_i}, \Lambda^{k_i}) = (L^*, S^*, Z^*, \Lambda^*)$. Then from the lower semicontinuity of $\Theta_{\tau,\beta}$ (since $\Psi$ and $\Phi$ are lower semicontinuous) and the definition of $\theta^*$, we have
\begin{eqnarray}\label{liminf}
\theta^* = \lim_{i \rightarrow \infty}~\Theta_{\tau,\beta}(L^{k_i}, S^{k_i}, Z^{k_i}, \Lambda^{k_i}) \geq \Theta_{\tau,\beta}(L^*, S^*, Z^*, \Lambda^*).
\end{eqnarray}
On the other hand, notice from the definition of $S^{k+1}$ as a minimizer in \eqref{Ssubpro} that
\begin{eqnarray*}
\begin{aligned}
&\quad\Theta_{\tau,\beta} (L^{k_i}, S^{k_i}, Z^{k_i-1}, \Lambda^{k_i-1}) - \Theta_{\tau,\beta}(L^{k_i}, S^*, Z^{k_i-1}, \Lambda^{k_i-1}) \\
&= \mathcal{L}_{\beta} (L^{k_i}, S^{k_i}, Z^{k_i-1}, \Lambda^{k_i-1}) - \mathcal{L}_{\beta} (L^{k_i}, S^*, Z^{k_i-1}, \Lambda^{k_i-1})  \\
&\qquad+\theta(\tau)\beta\left(\|\mathcal{B}(L^{k_i})+\mathcal{C}(S^{k_i})-Z^{k_i-1}\|_F^2
-\|\mathcal{B}(L^{k_i})+\mathcal{C}(S^{*})-Z^{k_i-1}\|_F^2\right)   \\
&\leq \theta(\tau)\beta\left(\|\mathcal{B}(L^{k_i})+\mathcal{C}(S^{k_i})-Z^{k_i-1}\|_F^2
-\|\mathcal{B}(L^{k_i})+\mathcal{C}(S^{*})-Z^{k_i-1}\|_F^2\right).
\end{aligned}
\end{eqnarray*}
This together with Theorem \ref{convergencethe}(i), the continuity of $\Theta_{\tau,\beta}$ with respect to $L$ (from (a2) in Assumption \ref{assum}), $Z$ and $\Lambda$; and the definition of $\theta^*$ implies that
\begin{eqnarray}\label{limsup}
\theta^* = \lim\limits_{i\rightarrow\infty}~\Theta_{\tau,\beta}(L^{k_i}, S^{k_i}, Z^{k_i}, \Lambda^{k_i}) \leq\Theta_{\tau,\beta}(L^*, S^*, Z^*, \Lambda^*).
\end{eqnarray}
Combining \eqref{liminf} and \eqref{limsup}, we conclude that $\Theta_{\tau,\beta}(L^*, S^*, Z^*, \Lambda^*) = \theta^*$. Since $(L^*, S^*, Z^*, \Lambda^*)\in \Gamma$ is arbitrary, we conclude further that the potential function $\Theta_{\tau,\beta}$ is constant on $\Gamma$.

The fact that $\Theta_{\tau,\beta}\equiv \theta^*$ on $\Gamma$ together with our assumption that $\Theta_{\tau,\beta}(\cdot)$ is a KL function and Proposition \ref{uniKL} implies that there exist $\varepsilon>0$, $\eta>0$ and $\varphi\in\Xi_{\eta}$, such that
\begin{eqnarray*}
\varphi'\left(\Theta_{\tau,\beta}(L, S, Z, \Lambda)-\theta^*\right)\mathrm{dist}\left(0,~\partial\Theta_{\tau,\beta}(L, S, Z, \Lambda)\right) \geq 1
\end{eqnarray*}
for all $(L, S, Z, \Lambda)$ satisfying $\mathrm{dist}((L, S, Z, \Lambda),\Gamma) < \varepsilon$ and $\theta^*<\Theta_{\tau,\beta}(L, S, Z, \Lambda)<\theta^*+\eta$. On the other hand, since $\lim_{k\rightarrow \infty}\mathrm{dist}((L^{k}, S^{k}, Z^{k}, \Lambda^{k}),\Gamma)=0$ by the definition of $\Gamma$, and $\Theta_{\tau,\beta}(L^{k}, S^{k}$, $Z^{k}$, $\Lambda^{k})\rightarrow\theta^*$, then for such $\varepsilon$ and $\eta$, there exists $k_1\ge 3$ such that $\mathrm{dist}((L^{k}, S^{k}, Z^{k}, \Lambda^{k}),\Gamma) < \varepsilon$ and $\theta^*<\Theta_{\tau,\beta}(L^{k}, S^{k}, Z^{k}, \Lambda^{k})<\theta^*+\eta$ for all $k\geq k_1$. Thus, for $k\geq k_1$, we have
\begin{eqnarray}\label{klpotential}
\varphi'\left(\Theta_{\tau,\beta}(L^{k}, S^{k}, Z^{k}, \Lambda^{k})-\theta^*\right)\mathrm{dist}\left(0,~\partial\Theta_{\tau,\beta}(L^{k}, S^{k}, Z^{k}, \Lambda^{k})\right) \geq 1.
\end{eqnarray}

\textit{Step \textbf{2.}} We next consider the subdifferential $\partial\Theta_{\tau,\beta}(L^{k}, S^{k}, Z^{k}, \Lambda^{k})$. Looking at the partial subdifferential with respect to $L$, we have
\begin{eqnarray*}
\begin{aligned}
&\quad \partial_{L}\Theta_{\tau,\beta}(L^{k}, S^{k}, Z^{k}, \Lambda^{k}) \\
&=\partial \Psi(L^k)-\mathcal{B}^*(\Lambda^{k})+(1+2\theta(\tau))\beta\mathcal{B}^*(\mathcal{B}(L^k)+\mathcal{C}(S^k)-Z^k) \\
&=\partial \Psi(L^k)-\mathcal{B}^*(\Lambda^{k-1})+\beta\mathcal{B}^*(\mathcal{B}(L^k)+\mathcal{C}(S^{k-1})-Z^{k-1})
+2\theta(\tau)\beta\mathcal{B}^*(\mathcal{B}(L^k)+\mathcal{C}(S^k)-Z^k)\\
&\quad~~-\mathcal{B}^*(\Lambda^{k}-\Lambda^{k-1})+\beta\mathcal{B}^*(\mathcal{C}(S^{k})-Z^{k}-\mathcal{C}(S^{k-1})+Z^{k-1}) \\
&\ni 2\theta(\tau)\beta\mathcal{B}^*(\mathcal{B}(L^k)+\mathcal{C}(S^k)-Z^k)
-\mathcal{B}^*(\Lambda^{k}-\Lambda^{k-1})+\beta\mathcal{B}^*(\mathcal{C}(S^{k})-Z^{k}-\mathcal{C}(S^{k-1})+Z^{k-1})  \\
&\stackrel{(\mathrm{i})}{=} -{\textstyle\left(1+\frac{2\theta(\tau)}{\tau}\right)}\mathcal{B}^*(\Lambda^{k}-\Lambda^{k-1})
+\beta\mathcal{B}^*[(\mathcal{C}(S^{k})-Z^{k})-(\mathcal{C}(S^{k-1})-Z^{k-1})] \\
&\stackrel{(\mathrm{ii})}{=} -{\textstyle\left(1+\frac{2\theta(\tau)}{\tau} \right)\mathcal{B}^*(\Lambda^{k}-\Lambda^{k-1})
+\beta\mathcal{B}^*\left[\left(-\mathcal{B}(L^k)-\frac{\Lambda^k-\Lambda^{k-1}}{\tau\beta}\right)-\left(-\mathcal{B}(L^{k-1}) -\frac{\Lambda^{k-1}-\Lambda^{k-2}}{\tau\beta}\right)\right]}\\
&= -{\textstyle\left(1+\frac{2\theta(\tau)+1}{\tau}\right)}\mathcal{B}^*(\Lambda^{k}-\Lambda^{k-1})
+{\textstyle\frac{1}{\tau}}\mathcal{B}^*(\Lambda^{k-1}-\Lambda^{k-2})-\beta\mathcal{B}^*\mathcal{B}(L^{k}-L^{k-1}),
\end{aligned}
\end{eqnarray*}
where the inclusion follows from \eqref{lopt}, and the equalities (i) and (ii) follow from \eqref{lamopt}. Similarly,
\begin{eqnarray*}
\begin{aligned}
&\quad \partial_{S}\Theta_{\tau,\beta}(L^{k}, S^{k}, Z^{k}, \Lambda^{k}) \\
&=\partial \Phi(S^k)-\mathcal{C}^*(\Lambda^{k})+(1+2\theta(\tau))\beta\mathcal{C}^*(\mathcal{B}(L^k)+\mathcal{C}(S^k)-Z^k) \\
&=\partial \Phi(S^k)-\mathcal{C}^*(\Lambda^{k-1})+\beta\mathcal{C}^*(\mathcal{B}(L^{k})+\mathcal{C}(S^k)-Z^{k-1})\\
&\quad~~+2\theta(\tau)\beta\mathcal{C}^*(\mathcal{B}(L^k)+\mathcal{C}(S^k)-Z^k)-\mathcal{C}^*(\Lambda^k-\Lambda^{k-1})
-\beta\mathcal{C}^*(Z^k-Z^{k-1}) \\
&\ni 2\theta(\tau)\beta\mathcal{C}^*(\mathcal{B}(L^k)+\mathcal{C}(S^k)-Z^k)
-\mathcal{C}^*(\Lambda^k-\Lambda^{k-1})-\beta\mathcal{C}^*(Z^k-Z^{k-1}) \\
&=-\left(1+\frac{2\theta(\tau)}{\tau}\right)\mathcal{C}^*(\Lambda^{k}-\Lambda^{k-1})-\beta\mathcal{C}^*(Z^k-Z^{k-1}),
\end{aligned}
\end{eqnarray*}
where the inclusion follows from \eqref{sopt} and the last equality follows from \eqref{lamopt}. Moreover,
\begin{eqnarray*}
\begin{aligned}
&\quad \nabla_{Z}\Theta_{\tau,\beta}(L^{k}, S^{k}, Z^{k}, \Lambda^{k}) \\
&=\mathcal{A}^*(\mathcal{A}(Z^k)-D)+\Lambda^k-\beta(\mathcal{B}(L^k)+\mathcal{C}(S^k)-Z^k) \\
&\quad~~-2\theta(\tau)\beta(\mathcal{B}(L^k)+\mathcal{C}(S^k)-Z^k) \\
&=\mathcal{A}^*(\mathcal{A}(Z^k)-D)+\Lambda^{k-1}-\beta(\mathcal{B}(L^k)+\mathcal{C}(S^k)-Z^k) \\
&\quad~~-2\theta(\tau)\beta(\mathcal{B}(L^k)+\mathcal{C}(S^k)-Z^k) +(\Lambda^k - \Lambda^{k-1}) \\
&=-2\theta(\tau)\beta(\mathcal{B}(L^k)+\mathcal{C}(S^k)-Z^k) +(\Lambda^k - \Lambda^{k-1}) \\
&=\left(1 +\frac{2\theta(\tau)}{\tau}\right)(\Lambda^{k}-\Lambda^{k-1}),
\end{aligned}
\end{eqnarray*}
where the third equality follows from \eqref{zopt} and the last equality follows from \eqref{lamopt}. Finally,
\begin{eqnarray*}
\nabla_{\lambda}\Theta_{\tau,\beta}(L^{k}, S^{k}, Z^{k}, \Lambda^{k}) = -(\mathcal{B}(L^k)+\mathcal{C}(S^k)-Z^k)=\frac{1}{\tau\beta}(\Lambda^k-\Lambda^{k-1}),
\end{eqnarray*}
where the last equality follows from \eqref{lamopt}. Thus, from the above relations, there exists $a>0$ so that
\begin{eqnarray}\label{distbound}
\hspace{-4mm}\begin{aligned}
&\mathrm{dist}\left(0,~\partial\Theta_{\tau,\beta}\left(L^{k}, S^{k}, Z^{k}, \Lambda^{k}\right)\right) \\
&\leq a\left(\|L^{k}-L^{k-1}\|_F+\|Z^{k}-Z^{k-1}\|_F+\|\Lambda^{k}-\Lambda^{k-1}\|_F+\|\Lambda^{k-1}-\Lambda^{k-2}\|_F\right).
\end{aligned}
\end{eqnarray}

\textit{Step \textbf{3.}} We now prove the convergence of the sequence by combining \eqref{distbound} with \eqref{klpotential}. For notational simplicity, define
\[
\Delta^k:=\varphi\left(\Theta_{\tau,\beta}(L^{k}, S^{k}, Z^{k}, \Lambda^{k})-\theta^*\right)-\varphi\left(\Theta_{\tau,\beta}(L^{k+1}, S^{k+1}, Z^{k+1}, \Lambda^{k+1})-\theta^*\right).
\]
Since $\Theta_{\tau,\beta}$ is decreasing and $\varphi$ is monotonic, it is easy to see $\Delta^k\geq0$ for $k\geq1$. Then we have for all $k\geq k_1$ that
\begin{eqnarray}\label{haha}
\begin{aligned}
&\quad a\left(\|L^{k}-L^{k-1}\|_F+\|Z^{k}-Z^{k-1}\|_F+\|\Lambda^{k}-\Lambda^{k-1}\|_F+\|\Lambda^{k-1}-\Lambda^{k-2}\|_F\right)\cdot\Delta^k \\
&\geq \mathrm{dist}(0,~\partial\Theta_{\tau,\beta}(L^{k}, S^{k}, Z^{k}, \Lambda^{k}))\cdot\Delta^k \\
&\geq \mathrm{dist}(0,~\partial\Theta_{\tau,\beta}(L^{k}, S^{k}, Z^{k}, \Lambda^{k}))\cdot\varphi'\left(\Theta_{\tau,\beta}(L^{k}, S^{k}, Z^{k}, \Lambda^{k})-\theta^*\right)\\
&\quad \cdot\left[\Theta_{\tau,\beta} (L^k, S^k, Z^k, \Lambda^k)-\Theta_{\tau,\beta} (L^{k+1}, S^{k+1}, Z^{k+1}, \Lambda^{k+1})\right] \\
&\geq \Theta_{\tau,\beta} (L^k, S^k, Z^k, \Lambda^k)-\Theta_{\tau,\beta} (L^{k+1}, S^{k+1}, Z^{k+1}, \Lambda^{k+1}) \\
&\geq b_1\|L^{k+1} - L^k\|_F^2 + b_2\|Z^{k+1} - Z^k\|_F^2   \\
&\geq \frac{1}{2}\min\{b_1, b_2\}\cdot\left[\|L^{k+1}-L^k\|_F+\|Z^{k+1}-Z^k\|_F\right]^2,
\end{aligned}
\end{eqnarray}
where the first inequality follows from \eqref{distbound}, the second inequality follows from the concavity of $\varphi$, the third inequality follows from \eqref{klpotential}, the fourth inequality follows from \eqref{difofp} with $b_1:=\frac{\sigma\beta}{2}$ and $b_2:=\frac{\lambda_{\min}+\beta}{2}-\max\left\{\frac{1}{\tau}, \frac{\tau^2}{1+\tau-\tau^2}\right\}\cdot\frac{\lambda_{\max}^2}{\beta}$.

Dividing both sides of \eqref{haha} by $c:=\frac{1}{2}\min\{b_1, b_2\}$, taking the square root and using the inequality $\sqrt{uv}\leq \frac{u + v}{2}$ for $u, v \ge 0$ to further upper bound the left hand side of the resulting inequality, we obtain that
\begin{eqnarray}\label{klresult1}
\hspace{-10mm}
\begin{aligned}
&\textstyle{\frac{1}{2\gamma}}\left(\|L^{k}-L^{k-1}\|_F+\|Z^{k}-Z^{k-1}\|_F+\|\Lambda^{k}-\Lambda^{k-1}\|_F+\|\Lambda^{k-1}-\Lambda^{k-2}\|_F\right)
+\textstyle{\frac{\gamma a}{2c}}\Delta^k\\
& \ge \|L^{k+1}-L^k\|_F+\|Z^{k+1}-Z^k\|_F ,
\end{aligned}
\end{eqnarray}
where $\gamma$ is an arbitrary positive constant. On the other hand, it follows from \eqref{lamdiff} that
\begin{eqnarray*}
\|\Lambda^{k}-\Lambda^{k-1}\|_F &=& \|\tau \mathcal{A}^*\mathcal{A}(Z^{k-1}-Z^{k}) + (1-\tau)(\Lambda^{k-1} - \Lambda^{k-2})\|_F \\
&\leq& \tau\lambda_{\max}\|Z^{k}-Z^{k-1}\|_F+|1-\tau|\cdot\|\Lambda^{k-1} - \Lambda^{k-2}\|_F.
\end{eqnarray*}
Adding $-|1-\tau|\cdot\|\Lambda^{k} - \Lambda^{k-1}\|_F$ to both sides of the above inequality and simplifying the resulting inequality, we obtain that
\begin{eqnarray}\label{lambo1}
\begin{aligned}
&\quad\|\Lambda^{k}-\Lambda^{k-1}\|_F \\
&\leq \textstyle{\frac{\tau\lambda_{\max}}{1-|1-\tau|}}\|Z^{k}-Z^{k-1}\|_F + \textstyle{\frac{|1-\tau|}{1-|1-\tau|}}\left(\|\Lambda^{k-1} - \Lambda^{k-2}\|_F-\|\Lambda^{k}-\Lambda^{k-1}\|_F\right)\\
&= d_1\|Z^{k}-Z^{k-1}\|_F + d_2\left(\|\Lambda^{k-1} - \Lambda^{k-2}\|_F-\|\Lambda^{k}-\Lambda^{k-1}\|_F\right),
\end{aligned}
\end{eqnarray}
where we write $d_1:=\frac{\tau\lambda_{\max}}{1-|1-\tau|}$ and $d_2:=\frac{|1-\tau|}{1-|1-\tau|}$ for notational simplicity. Similarly,
\begin{eqnarray}\label{lambo2}
\begin{aligned}
&\|\Lambda^{k-1}-\Lambda^{k-2}\|_F \leq  d_1\|Z^{k-1}-Z^{k-2}\|_F \\
&\qquad+d_2\left(\|\Lambda^{k-2} - \Lambda^{k-3}\|_F-\|\Lambda^{k-1}-\Lambda^{k-2}\|_F\right).
\end{aligned}
\end{eqnarray}
Then substituting \eqref{lambo1} and \eqref{lambo2} into \eqref{klresult1} and rearranging terms, we have
\begin{eqnarray}\label{klresult2}
\begin{aligned}
&{\textstyle\left(1-\frac{1}{2\gamma}\right)}\|L^{k+1}-L^k\|_F + {\textstyle\left(1-\frac{1}{2\gamma}-\frac{d_1}{\gamma}\right)}\|Z^{k+1}-Z^k\|_F \\
&\leq{\textstyle\frac{1}{2\gamma}}\left(\|L^{k}-L^{k-1}\|_F-\|L^{k+1}-L^k\|_F\right) \\
&~~+{\textstyle\left(\frac{1}{2\gamma}+\frac{d_1}{\gamma}\right)}\left(\|Z^{k}-Z^{k-1}\|_F-
\|Z^{k+1}-Z^{k}\|_F\right) \\
&~~+{\textstyle\frac{d_1}{2\gamma}}\left(\|Z^{k-1}-Z^{k-2}\|_F-\|Z^{k}-Z^{k-1}\|_F\right) \\
&~~+{\textstyle\frac{d_2}{2\gamma}}\left(\|\Lambda^{k-1}-\Lambda^{k-2}\|_F-\|\Lambda^{k}
-\Lambda^{k-1}\|_F\right) \\
&~~+{\textstyle\frac{d_2}{2\gamma}}\left(\|\Lambda^{k-2}-\Lambda^{k-3}\|_F-\|\Lambda^{k-1}-\Lambda^{k-2}\|_F\right)
+{\textstyle\frac{\gamma a}{2c}}\Delta^k.
\end{aligned}
\end{eqnarray}
Thus, summing \eqref{klresult2} from $k=k_1$ to $\infty$, we have
\begin{eqnarray*}
\begin{aligned}
&\quad{\textstyle\left(1-\frac{1}{2\gamma}\right)\sum_{k=k_1}^{\infty}}\|L^{k+1}-L^k\|_F + {\textstyle\left(1-\frac{1}{2\gamma}-\frac{d_1}{\gamma}\right)
\sum_{k=k_1}^{\infty}}\|Z^{k+1}-Z^k\|_F \\
&\leq{\textstyle\frac{1}{2\gamma}}\|L^{k_1}-L^{k_1-1}\|_F
+{\textstyle\left(\frac{1}{2\gamma}+\frac{d_1}{\gamma}\right)}\|Z^{k_1}-Z^{k_1-1}\|_F
+{\textstyle\frac{d_1}{2\gamma}}\|Z^{k_1-1}-Z^{k_1-2}\|_F \\
&~~+{\textstyle\frac{d_2}{2\gamma}}\|\Lambda^{k_1-1}-\Lambda^{k_1-2}\|_F+{\textstyle\frac{d_2}{2\gamma}}\|\Lambda^{k_1-2} - \Lambda^{k_1-3}\|_F+{\textstyle\frac{a\gamma}{2c}}\,\varphi\left(\Theta_{\tau,\beta}(L^{k_1}, S^{k_1}, Z^{k_1}, \Lambda^{k_1})-\theta^*\right) \\
&< \infty.
\end{aligned}
\end{eqnarray*}
Recall that $\gamma$ introduced in \eqref{klresult1} is an arbitrary positive constant. Taking $\gamma>\frac{1+2d_1}{2}$ and hence $1-\frac{1}{2\gamma}>1-\frac{1}{2\gamma}-\frac{d_1}{\gamma}>0$, we have from the above inequality that
\begin{eqnarray*}
\sum\limits_{k=k_1}^{\infty}\|L^{k+1}-L^k\|_F < \infty \quad \mathrm{and} \quad \sum\limits_{k=k_1}^{\infty}\|Z^{k+1}-Z^k\|_F<\infty.
\end{eqnarray*}
Hence $\{L^k\}$ and $\{Z^k\}$ are convergent. Additionally, summing \eqref{lambo1} from $k=k_1$ to $\infty$, we have
\begin{eqnarray*}%\label{sumoflambo}
\begin{aligned}
\sum\limits_{k=k_1}^{\infty}\|\Lambda^{k}-\Lambda^{k-1}\|_F&\leq d_1\sum\limits_{k=k_1}^{\infty}\|Z^{k}-Z^{k-1}\|_F + d_2\|\Lambda^{k_1-1} - \Lambda^{k_1-2}\|_F<\infty,\\
\end{aligned}
\end{eqnarray*}
which implies that $\{\Lambda^k\}$ is convergent. Finally, from \eqref{lamopt} and (a1) in Assumption \ref{assum}, we see that $\{S^k\}$ is also convergent. Consequently, we conclude that $\{(L^k$, $S^k$, $Z^k$, $\Lambda^k)\}_{k=1}^{\infty}$ is a convergent sequence. This completes the proof.
\endproof

Our convergence analysis relies on Assumption \ref{assum}. While (a3) in Assumption~\ref{assum} appears restrictive since it makes assumptions on the first iterate of Algorithm 1, we show below that this assumption would hold upon a suitable choice of initialization. Specifically, if we initialize at $(L^0, S^0, Z^0, \Lambda^0)$ satisfying
\begin{numcases}{}
\Theta_{\tau,\beta}(L^1, S^1, Z^1, \Lambda^1) \leq \Theta_{\tau,\beta}(L^0, S^0, Z^0, \Lambda^0), \label{initializationB}\\
\Theta_{\tau,\beta}(L^0, S^0, Z^0, \Lambda^0) < h_0, \label{initializationA}
\end{numcases}
%where $Z^0:=\mathcal{B}(L^0)-\mathcal{C}(S^0)$,
then it is easy to check that (a3) in Assumption~\ref{assum} holds. In the next proposition, we demonstrate that \eqref{initializationB} can always be satisfied with a suitable initialization. After this, we will propose a specific way to initialize Algorithm 1 for a wide range of problems so that both \eqref{initializationB} and \eqref{initializationA} are satisfied.

\begin{proposition}\label{initial_pro}
Suppose that $0<\tau<\frac{1+\sqrt{5}}{2}$ and $\beta > \bar\beta$. If the initialization $(L^0, S^0, Z^0, \Lambda^0)$ is chosen as $(L^0, S^0)\in\mathrm{dom}\,\Psi \times \mathrm{dom}\,\Phi$ and
\begin{eqnarray}\label{initial_con1}
\Lambda^0 = \mathcal{A}^*(D - \mathcal{A}(Z^0)),
\end{eqnarray}
then we have
\begin{eqnarray*}
\Theta_{\tau,\beta}(L^1, S^1, Z^1, \Lambda^1) \leq \Theta_{\tau,\beta}(L^0, S^0, Z^0, \Lambda^0).
\end{eqnarray*}
\end{proposition}
\beginproof
First, from \eqref{zopt}, we have
\begin{eqnarray}\label{initial_optcon}
\hspace{-8mm}\begin{aligned}
0&= \mathcal{A}^{*}(\mathcal{A}(Z^{1})-D) + \Lambda^0 - \beta(\mathcal{B}(L^{1}) + \mathcal{C}(S^{1}) - Z^{1}) \\
\Longrightarrow~&
\mathcal{B}(L^{1})+\mathcal{C}(S^{1})-Z^{1} = \frac{1}{\beta}\Lambda^0 + \frac{1}{\beta}\mathcal{A}^{*}(\mathcal{A}(Z^{1})-D)=\frac{1}{\beta} \mathcal{A}^*\mathcal{A}(Z^1 - Z^0),
\end{aligned}
\end{eqnarray}
where the last equality follows from \eqref{initial_con1}. Then,
\begin{eqnarray}\label{first_term_initial1}
&&\Theta_{\tau,\beta} (L^{1}, S^{1}, Z^{1}, \Lambda^{1}) - \Theta_{\tau,\beta} (L^{1}, S^{1}, Z^{1}, \Lambda^0) \nonumber\\
&=& -\langle\Lambda^{1}-\Lambda^0, \mathcal{B}(L^{1})+\mathcal{C}(S^{1})-Z^{1}\rangle = \tau\beta\|\mathcal{B}(L^{1})+\mathcal{C}(S^{1})-Z^{1}\|_F^2   \nonumber\\
&=& \left(\tau+\theta(\tau)\right)\beta\|\mathcal{B}(L^{1})+\mathcal{C}(S^{1})-Z^{1}\|_F^2 - \theta(\tau)\beta\|\mathcal{B}(L^{1})+\mathcal{C}(S^{1})-Z^{1}\|_F^2  \nonumber\\
&=& \left(\tau+\theta(\tau)\right)\beta\left\|{\textstyle\frac{1}{\beta} \mathcal{A}^*\mathcal{A}(Z^1 - Z^0)}\right\|_F^2
-\theta(\tau)\beta\|\mathcal{B}(L^{1})+\mathcal{C}(S^{1})-Z^{1}\|_F^2 \nonumber\\
&\leq& \left(\tau+\theta(\tau)\right)\frac{\lambda^2_{\max}}{\beta}\|Z^{1} - Z^0\|_F^2 - \theta(\tau)\beta\|\mathcal{B}(L^{1})+\mathcal{C}(S^{1})-Z^{1}\|_F^2,
\end{eqnarray}
where the second equality follows from \eqref{lamopt} and the fourth equality follows from \eqref{initial_optcon}. Additionally, using the same arguments as in the proof of Lemma \ref{Phidelemma} leading to \eqref{second_term}, \eqref{third_term} and \eqref{fourth_term}, it is easy to see that
\begin{eqnarray}
&&\Theta_{\tau,\beta} (L^{1}, S^{1}, Z^{1}, \Lambda^{0})-\Theta_{\tau,\beta}(L^{1}, S^{1}, Z^0, \Lambda^0)\leq \textstyle{-\frac{\lambda_{\min}+\beta}{2}}\|Z^{1} - Z^0\|_F^2\nonumber\\
&&\qquad+~\theta(\tau)\beta\left(\|\mathcal{B}(L^{1})+\mathcal{C}(S^{1})-Z^{1}\|_F^2
-\|\mathcal{B}(L^{1})+\mathcal{C}(S^{1})-Z^{0}\|_F^2\right), \label{second_term_initial1} \vspace{8mm}\\
&&\Theta_{\tau,\beta} (L^{1}, S^{1}, Z^{0}, \Lambda^{0}) - \Theta_{\tau,\beta}(L^{1}, S^0, Z^0, \Lambda^0)  \nonumber\\
&&\qquad\leq \theta(\tau)\beta\left(\|\mathcal{B}(L^{1})+\mathcal{C}(S^{1})-Z^{0}\|_F^2
-\|\mathcal{B}(L^{1})+\mathcal{C}(S^{0})-Z^{0}\|_F^2\right), \label{third_term_initial1} \vspace{8mm}\\
&&\Theta_{\tau,\beta} (L^{1}, S^{0}, Z^{0}, \Lambda^{0})-\Theta_{\tau,\beta}(L^{0}, S^{0}, Z^0, \Lambda^0) \nonumber\\
&&\qquad\leq\theta(\tau)\beta\left(\|\mathcal{B}(L^{1})+\mathcal{C}(S^{0})-Z^{0}\|_F^2
-\|\mathcal{B}(L^{0})+\mathcal{C}(S^{0})-Z^{0}\|_F^2\right).\label{fourth_term_initial1}
\end{eqnarray}
Summing \eqref{first_term_initial1}, \eqref{second_term_initial1}, \eqref{third_term_initial1} and \eqref{fourth_term_initial1}, we obtain
\begin{eqnarray}\label{initial_diff}
\hspace{-6mm}\begin{aligned}
&\quad\Theta_{\tau,\beta} (L^{1}, S^{1}, Z^{1}, \Lambda^{1}) - \Theta_{\tau,\beta} (L^0, S^0, Z^0, \Lambda^0) \\
&\leq {\textstyle\left(\left(\tau+\theta(\tau)\right)\frac{\lambda^2_{\max}}{\beta} - \frac{\lambda_{\min}+\beta}{2}\right)}\|Z^{1}-Z^0\|_F^2-\theta(\tau)\beta\|\mathcal{B}(L^{0})+\mathcal{C}(S^{0})-Z^{0}\|_F^2.
\end{aligned}
\end{eqnarray}
We now consider two cases:
\begin{itemize}
\item For $0< \tau \leq 1$, it is easy to see $\theta(\tau)=1-\tau$ and
\[
\beta > \max\left\{\frac{\lambda_{\max}}{\tau},
-\frac{\lambda_{\min}}{2}+\frac{1}{2}\sqrt{\lambda_{\min}^2+\frac{8}{\tau}\lambda_{\max}^2}\right\}.
\]
Then, we have
\begin{eqnarray*}
(\tau+\theta(\tau))\frac{\lambda^2_{\max}}{\beta} - \frac{\lambda_{\min}+\beta}{2} = \frac{\lambda^2_{\max}}{\beta} - \frac{\lambda_{\min}+\beta}{2} \leq \frac{\lambda^2_{\max}}{\tau\beta} - \frac{\lambda_{\min}+\beta}{2} < 0.
\end{eqnarray*}

\item For $1<\tau<\frac{1+\sqrt{5}}{2}$, it is easy to see $\theta(\tau)=\frac{(\tau-1)\tau^2}{1+\tau-\tau^2}$ and
\[\beta > \max\left\{\tau\lambda_{\max}, -\frac{\lambda_{\min}}{2}+\frac{1}{2}\sqrt{\lambda^2_{\min}+\frac{8\tau^2}{1+\tau-\tau^2}\lambda_{\max}^2}\right\}.
\]
Then, we have
    \begin{eqnarray*}
    \begin{aligned}
    &(\tau+\theta(\tau))\frac{\lambda^2_{\max}}{\beta} - \frac{\lambda_{\min}+\beta}{2} = \frac{\tau\lambda^2_{\max}}{(1+\tau-\tau^2)\beta} - \frac{\lambda_{\min}+\beta}{2} \\
    &\quad<\frac{\tau^2\lambda^2_{\max}}{(1+\tau-\tau^2)\beta} - \frac{\lambda_{\min}+\beta}{2} < 0.
    \end{aligned}
    \end{eqnarray*}
\end{itemize}
Thus, combining the above with \eqref{initial_diff} and $\theta(\tau)\geq0$, we conclude that
\begin{eqnarray*}
\Theta_{\tau,\beta}(L^1, S^1, Z^1, \Lambda^1) \leq \Theta_{\tau,\beta}(L^0, S^0, Z^0, \Lambda^0).
\end{eqnarray*}
This completes the proof.
\endproof

From Proposition \ref{initial_pro}, we see that if the initialization $(L^0, S^0, Z^0, \Lambda^0)$ is chosen to satisfy the conditions in Proposition \ref{initial_pro}, then \eqref{initializationB} holds. Based on this, we can now present one specific way to initialize Algorithm 1 so that both \eqref{initializationB} and \eqref{initializationA} are satisfied for a class of problems, whose objective functions $\Psi(L)$ and $\Phi(S)$ take forms $\delta_{\Omega}(L)$ and \eqref{sepphi}, respectively; here, $\Omega$ is a compact convex set.

The initialization we consider is:
\begin{eqnarray}\label{initialpoint}
L^0=\mathcal{P}_{\Omega}(\kappa D), ~S^0 = 0, ~Z^0 = \mathcal{B}(L^0), ~
\Lambda^0 = \mathcal{A}^*\left(D - \mathcal{A}(Z^0)\right),
\end{eqnarray}
where $\kappa$ is a scaling parameter. One can easily check that this initialization satisfies \eqref{initial_con1}. Moreover,
\begin{eqnarray*}
\Theta_{\tau,\beta}(L^0, S^0, Z^0, \Lambda^0)=\frac{1}{2}\left\|D-\mathcal{A}\left(Z^0\right)\right\|_F^2
=\frac{1}{2}\left\|D-\mathcal{A}\left(\mathcal{B}(\mathcal{P}_{\Omega}(\kappa D))\right)\right\|_F^2.
\end{eqnarray*}
Thus, the condition \eqref{initializationA} is equivalent to
\begin{equation}\label{ineq}
\frac{1}{2}\left\|D-\mathcal{A}\left(\mathcal{B}(\mathcal{P}_{\Omega}(\kappa D))\right)\right\|_F^2<\liminf\limits_{\|L\|_F+\|S\|_F\rightarrow\infty} \Psi(L)+\Phi(S)=\liminf\limits_{\|S\|_F\rightarrow\infty} \Phi(S).
\end{equation}
We further discuss this inequality for some concrete examples of $\Phi$ presented in the introduction.

\begin{example}
Suppose that $\Phi$ is coercive. Then $\liminf\limits_{\|S\|_F\rightarrow \infty} \Phi(S) = \infty$ and hence \eqref{ineq} holds trivially for any choice of $\kappa$.
\end{example}

\begin{example}
Suppose that $\Phi(S)=\mu\sum^m_{i=1}\sum^n_{j=1} \frac{\alpha |s_{ij}|}{1+\alpha|s_{ij}|}$ for $\alpha>0$. Then $\liminf\limits_{\|S\|_F\rightarrow \infty} \Phi(S) = \mu$.
Hence \eqref{ineq} holds if the parameter $\kappa$ can be chosen so that $\frac{1}{2}\left\|D-\mathcal{A}\left(\mathcal{B}(\mathcal{P}_{\Omega}(\kappa D))\right)\right\|_F^2<\mu$.
\end{example}

\begin{example}
Suppose that $\Phi(S)=\mu\sum^m_{i=1}\sum^n_{j=1}\int^{|s_{ij}|}_{0} \min(1, (\alpha - t/\mu)_+/(\alpha-1))~\mathrm{d}t$ for $\alpha>2$. Then $\liminf\limits_{\|S\|_F\rightarrow \infty} \Phi(S) = \frac{1}{2}(\alpha+1)\mu^2$. Hence \eqref{ineq} holds if $\kappa$ can be chosen so that $\frac{1}{2}\left\|D-\mathcal{A}\left(\mathcal{B}(\mathcal{P}_{\Omega}(\kappa D))\right)\right\|_F^2<\frac{1}{2}(\alpha+1)\mu^2$.
\end{example}

\begin{example}
Suppose that $\Phi(S)=\mu\sum^m_{i=1}\sum^n_{j=1} \int^{|s_{ij}|}_0 (1-t/(\alpha\mu))_+~\mathrm{d}t$ for $\alpha>0$. Then, $\liminf\limits_{\|S\|_F\rightarrow \infty} \Phi(S) = \frac{1}{2}\alpha\mu^2$.
Hence \eqref{ineq} holds if $\kappa$ can be chosen so that $\frac{1}{2}\left\|D-\mathcal{A}\left(\mathcal{B}(\mathcal{P}_{\Omega}(\kappa D))\right)\right\|_F^2<\frac{1}{2}\alpha\mu^2$.
\end{example}

\begin{example}
Suppose that $\Phi(S)=\mu\sum^m_{i=1}\sum^n_{j=1} \mu - (\mu - |s_{ij}|)^2_+ / \mu$. Then it is not hard to show that $\liminf\limits_{\|S\|_F\rightarrow \infty} \Phi(S) = \mu^2$.
Hence \eqref{ineq} holds if $\kappa$ can be chosen so that $\frac{1}{2}\left\|D-\mathcal{A}\left(\mathcal{B}(\mathcal{P}_{\Omega}(\kappa D))\right)\right\|_F^2<\mu^2$.
\end{example}

\section{Numerical experiments}\label{sec5}

In this section, we conduct numerical experiments to show the performances of our algorithm. All experiments are run in Matlab R2014b on a 64-bit PC with an Intel(R) Core(TM) i7-4790 CPU (3.60GHz) and 16GB of RAM equipped with Windows 8.1 OS.

\subsection{Implementation details}

\paragraph{Testing model} We consider the problem of extracting background/foreground from a given video under different scenarios. Specifically, we consider:
\begin{eqnarray}\label{testing_model}
\begin{array}{ll}
\qquad \min \limits_{L,S}~\Phi(S) + \frac{1}{2}\|D-\mathcal{A}(L+S)\|_F^2  \\
\qquad ~\mathrm{s.t.} ~~L \in \Omega,
\end{array}
\end{eqnarray}
where $\Omega = \left\{L \in {\mathbb{R}}^{m\times n}~|~ \|L\|_\infty \le 1,\ L_{:1} = L_{:2} = \cdots = L_{:n}\right\}$ and $\mathcal{A}$ is a linear map. This model corresponds to \eqref{generalmodel1} with $\Psi(L)=\delta_{\Omega}(L)$ and $\mathcal{B}=\mathcal{C}=\mathcal{I}$. We compare the performances of the ADMM with different choices of $\tau$, as well as the proximal alternating linearized minimization (PALM) proposed in \cite{bst2014}, on solving \eqref{testing_model}. For ease of future reference, we recall that the PALM for solving \eqref{testing_model} is given by
\begin{eqnarray*}%\label{PALM}
\left\{\begin{aligned}
&L^{k+1}=\mathcal{P}_{\Omega}\left(L^k-\frac{1}{c_k}\mathcal{A}^*(\mathcal{A}(L^k+S^k)-D)\right), \\
&S^{k+1}\in \mathop{\mathrm{Argmin}}\limits_{S}\left\{\Phi(S) + \frac{d_k}{2}\left\|S-S^k+\frac{1}{d_k}\mathcal{A}^*(\mathcal{A}(L^{k+1}+S^k)-D)\right\|_F^2\right\},
\end{aligned}\right.
\end{eqnarray*}
where $c_k$ and $d_k$ are positive numbers.

In our experiments, we consider the following three choices of sparse regularizers $\Phi(S)$: \vspace{2mm}
\begin{itemize}
\item bridge regularizer: $\Phi(S)=\mu\|S\|_p^p$ for $0<p\leq1$;
\item fraction regularizer: $\Phi(S)=\mu\sum^{m}_{i=1}\sum^n_{j=1}\frac{\alpha |s_{ij}|}{1+\alpha|s_{ij}|}$ for $\alpha>0$;
\item logistic regularizer: $\Phi(S)=\mu\sum^{m}_{i=1}\sum^n_{j=1}\log(1+\alpha|s_{ij}|)$ for $\alpha>0$;
\end{itemize}
\vspace{2mm}
and two choices of linear map $\mathcal{A}$: \vspace{2mm}
\begin{itemize}
\item $\mathcal{A}(L+S):=L+S$: in this case, model \eqref{testing_model} can be applied to extracting background/foreground from a surveillance video with noise. \vspace{2mm}

\item $\mathcal{A}(L+S):=H(L+S)$ with $H\in\mathbb{R}^{m\times m}$ being the matrix representation of a regular blurring operator (the blurring is assumed to occur frame-wise): in this case, model \eqref{testing_model} can be applied to extracting background/foreground from a blurred and noisy surveillance video. \vspace{2mm}
\end{itemize}

\begin{figure}[ht]
\centering
\hspace{1.7cm}\small{Hall}\hspace{1.1cm}\small{Bootstrap}\hspace{0.7cm}\small{Fountain}\hspace{6mm}\small{ShoppingMall} \vspace{1mm}\\
\hspace{-0.7cm}\begin{minipage}[t]{0.12\textwidth}\begin{center}\vspace{-12mm}noisy\end{center}\end{minipage}\includegraphics[height=2cm,width=2cm]{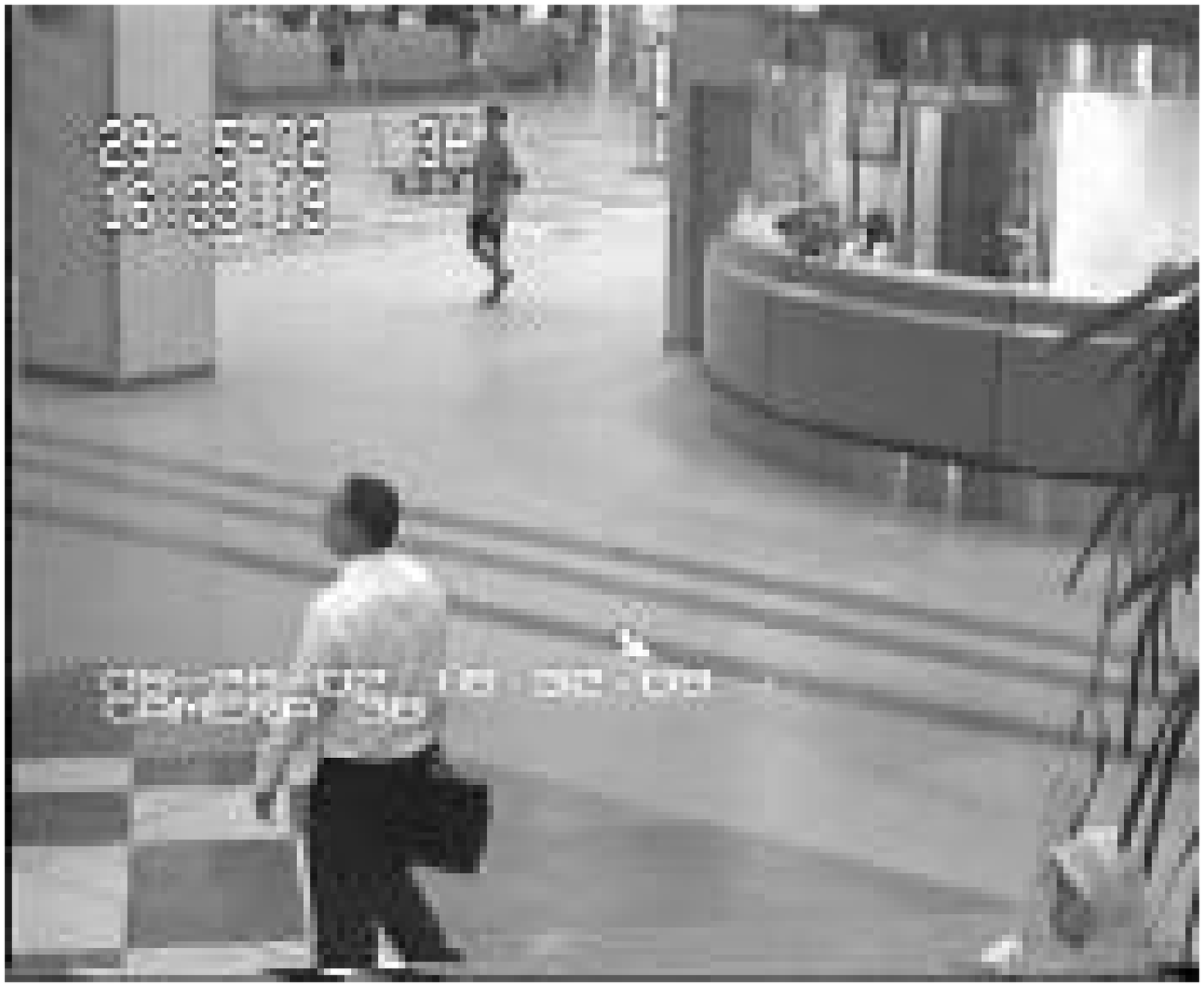}~\includegraphics[height=2cm,width=2cm]{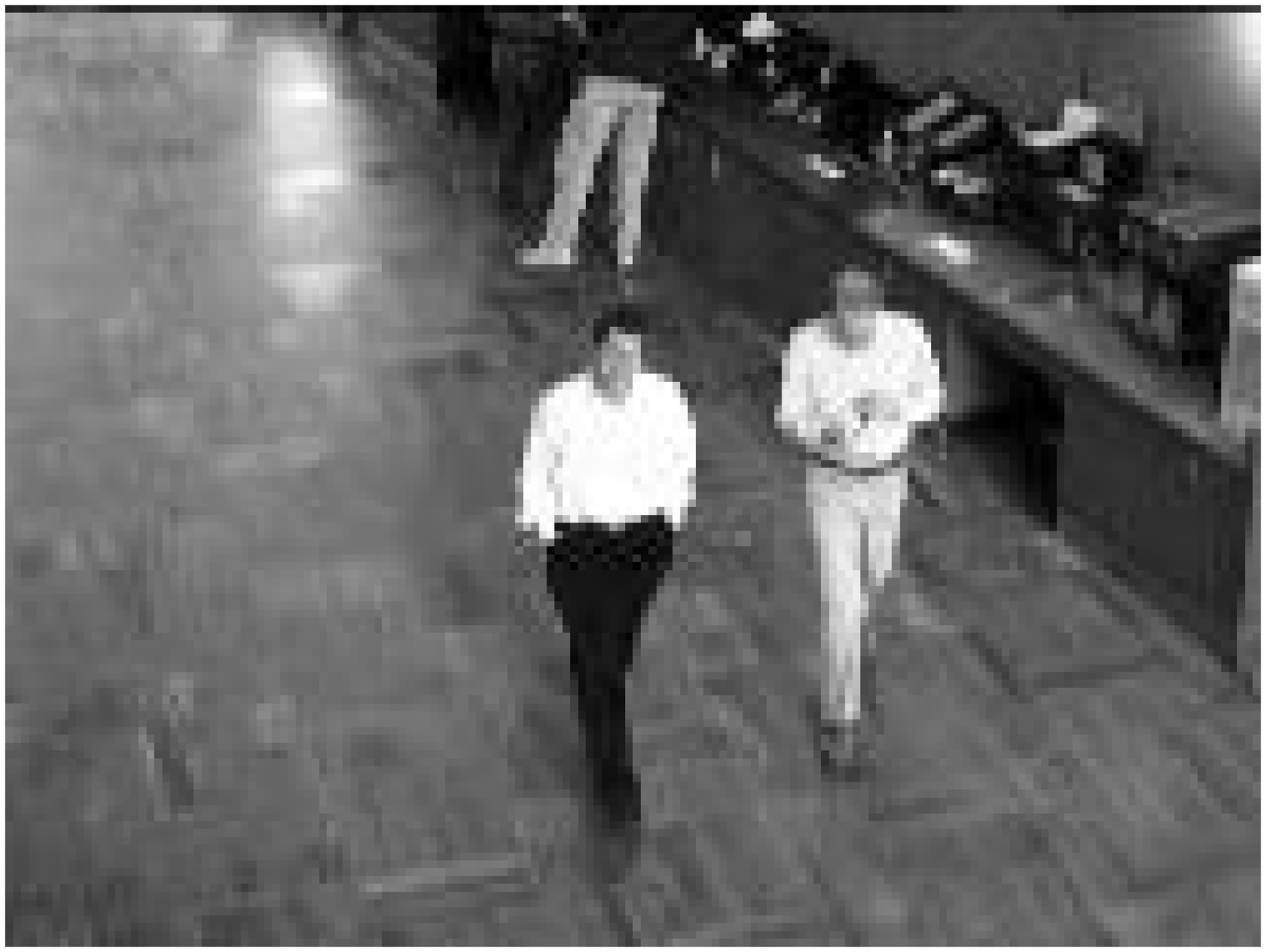}~\includegraphics[height=2cm,width=2cm]{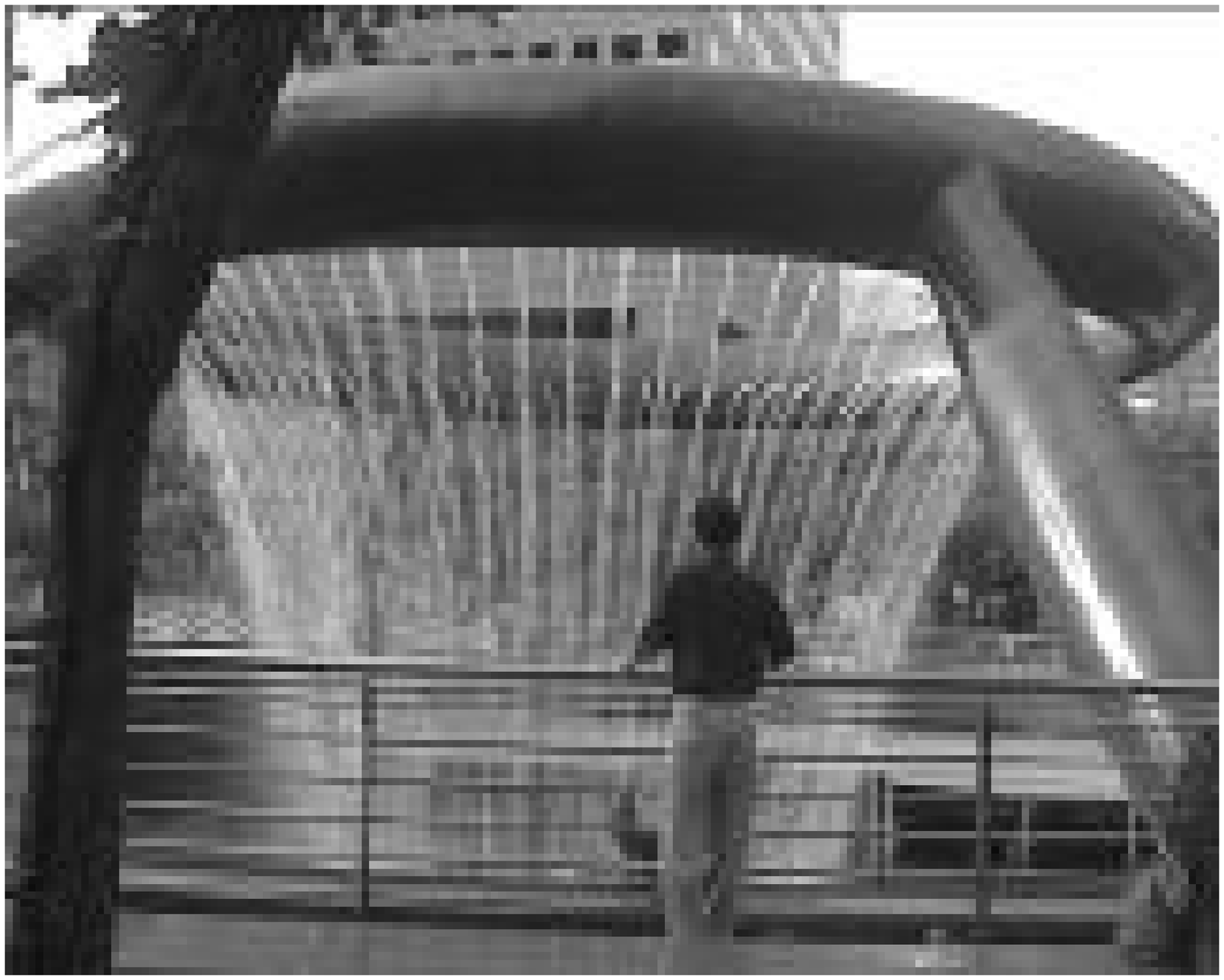}~\includegraphics[height=2cm,width=2cm]{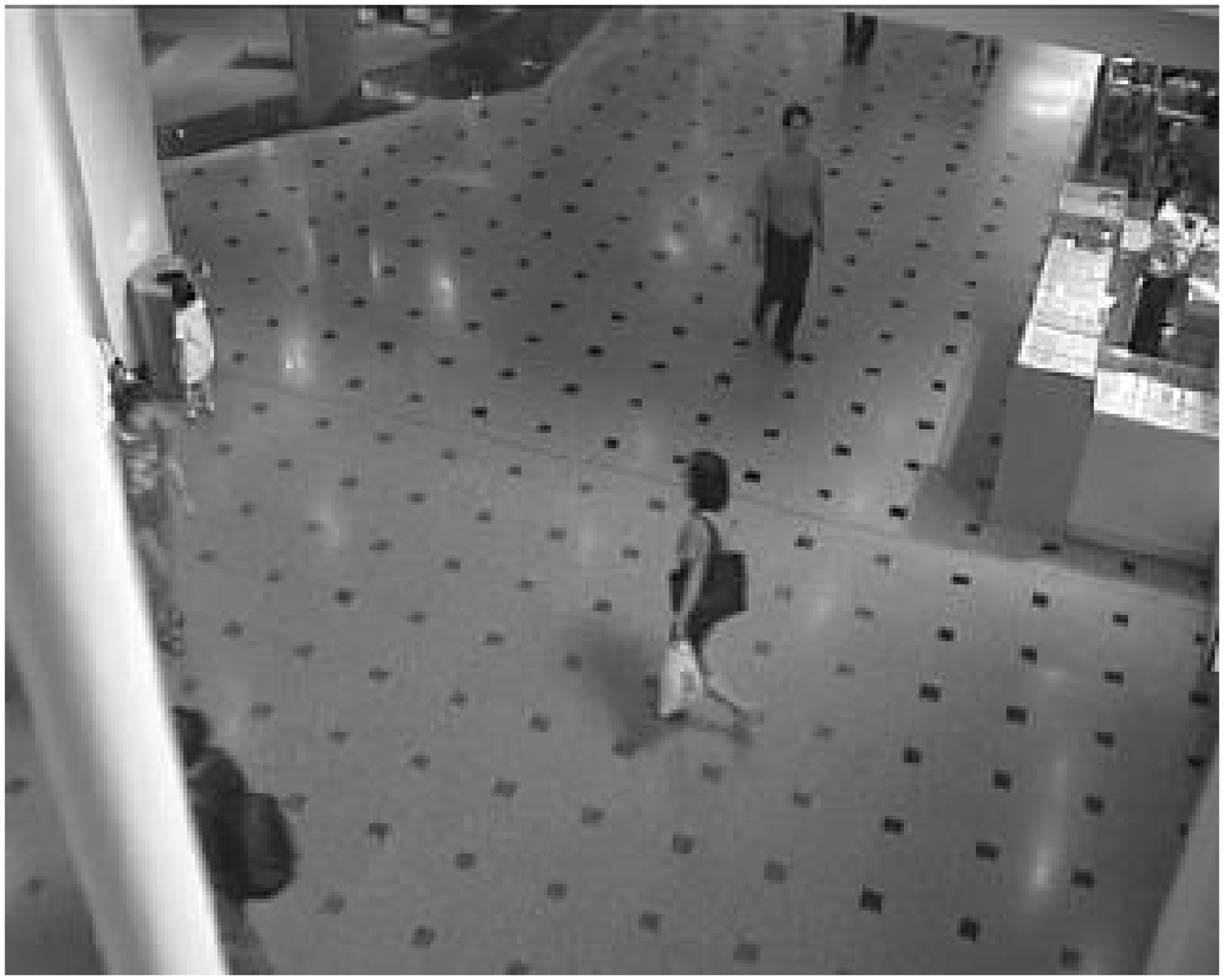}\vspace{1mm}\\
\hspace{-0.7cm}\begin{minipage}[t]{0.12\textwidth}\begin{center}\vspace{-14mm}noisy\\blurred\end{center}\end{minipage}\includegraphics[height=2cm,width=2cm]{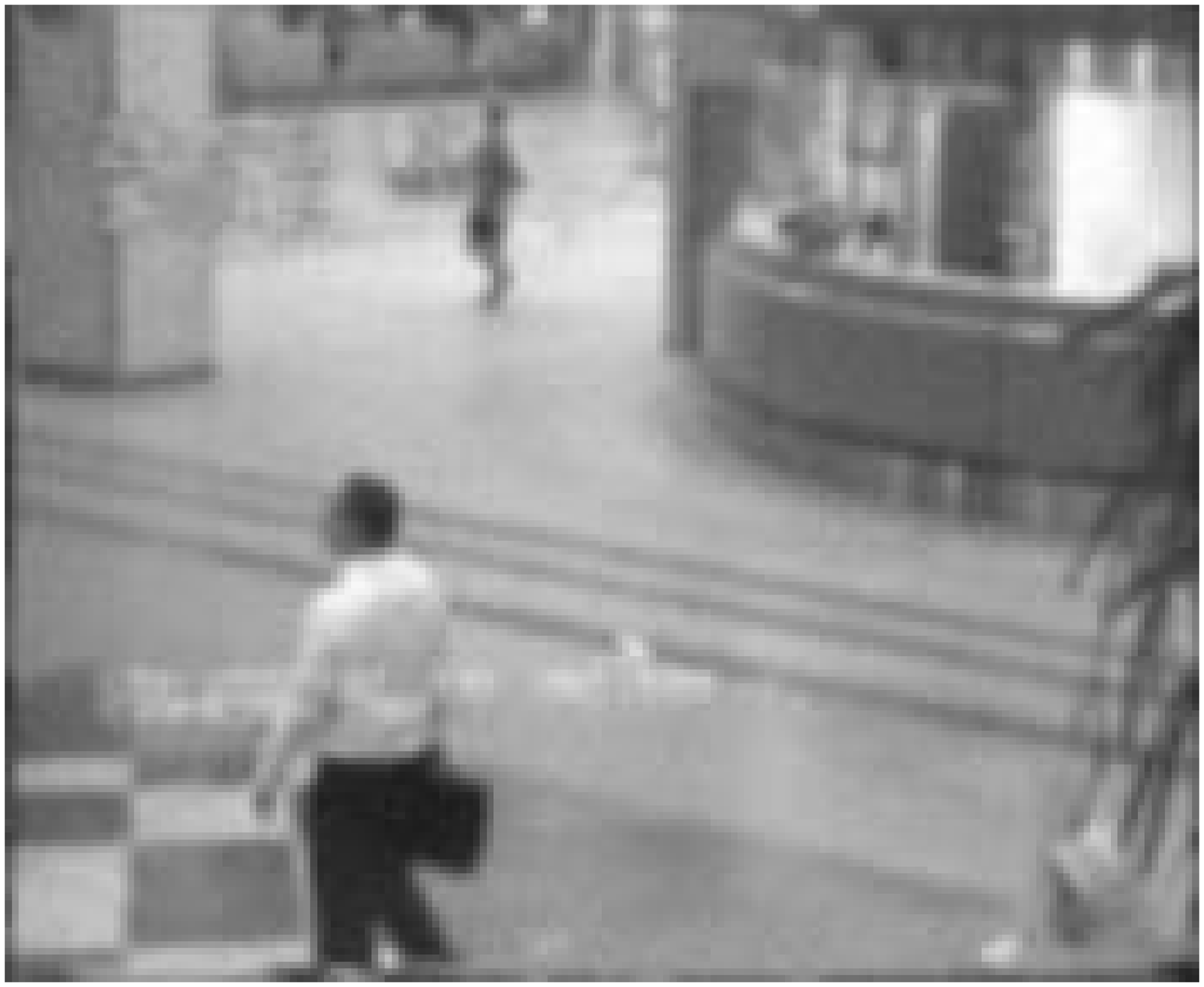}~\includegraphics[height=2cm,width=2cm]{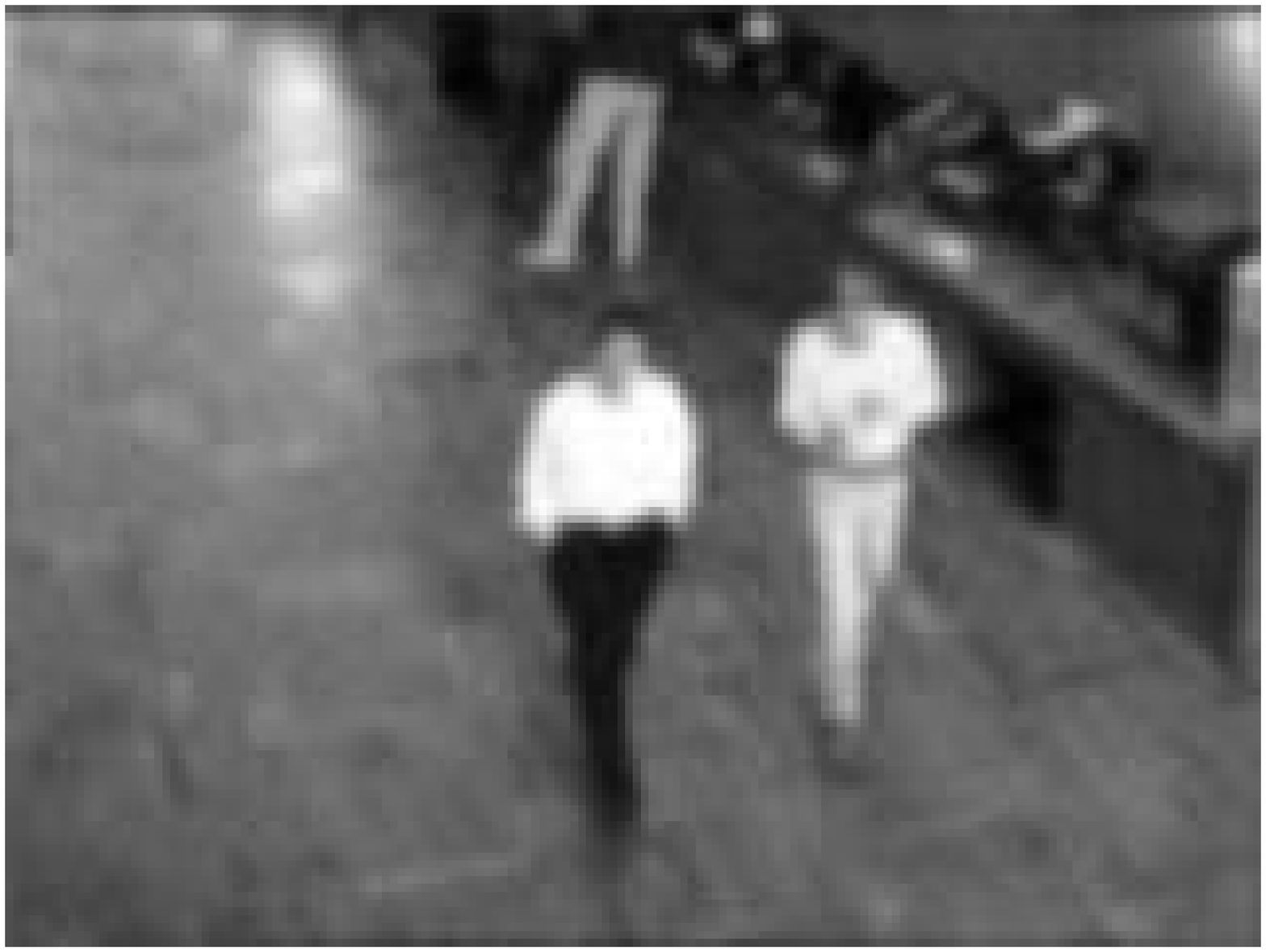}~\includegraphics[height=2cm,width=2cm]{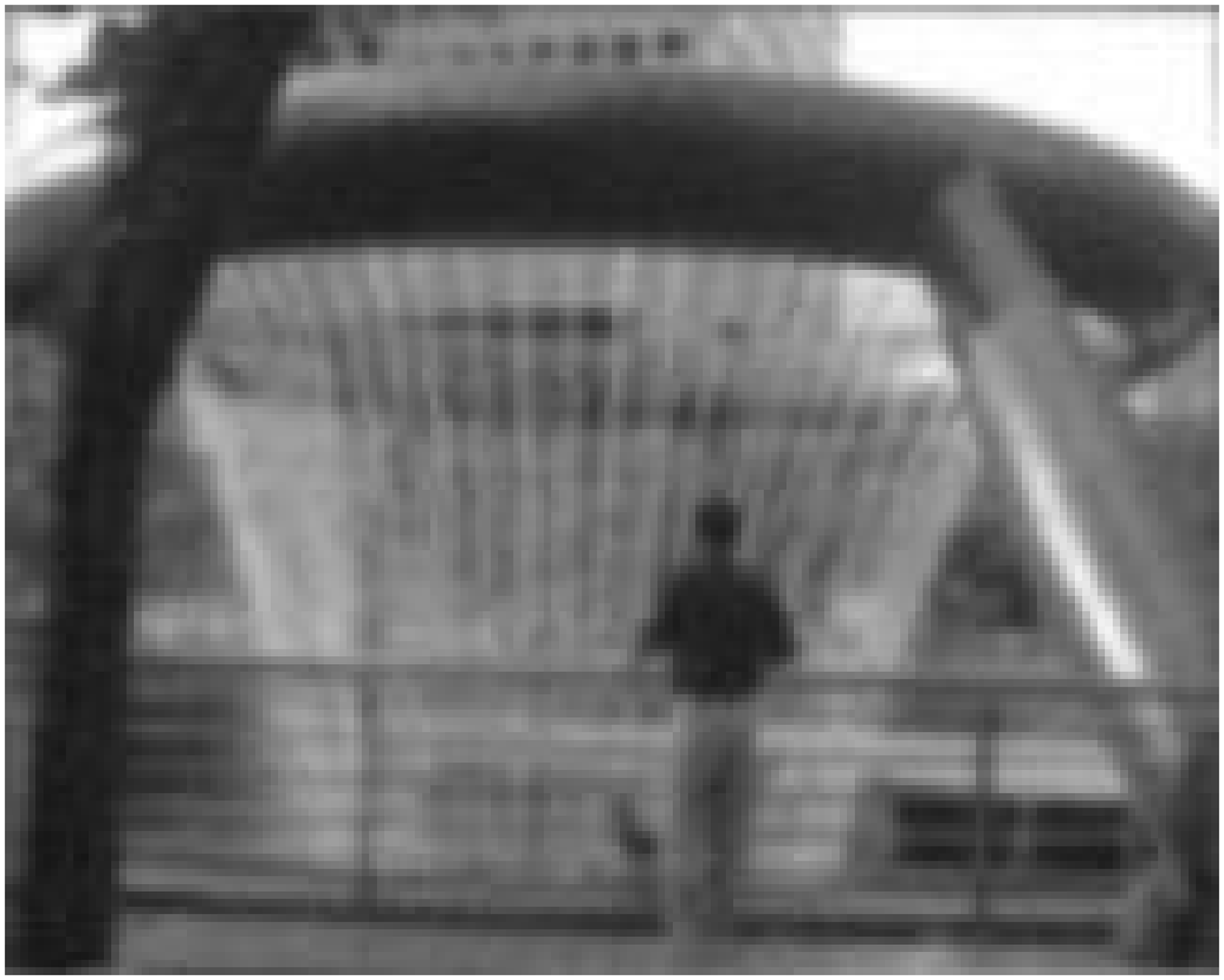}~\includegraphics[height=2cm,width=2cm]{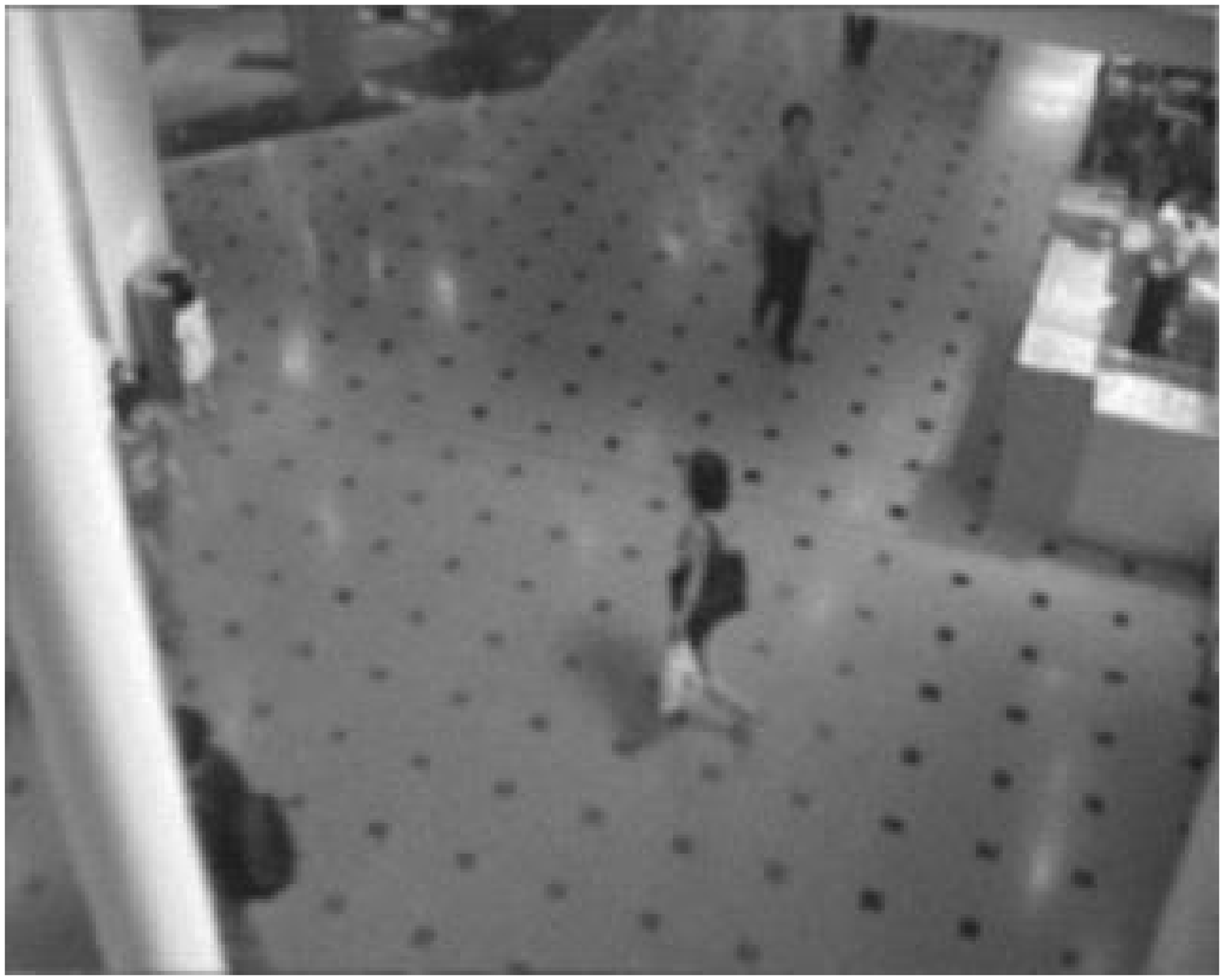}\vspace{1mm}\\
\hspace{-0.7cm}\begin{minipage}[t]{0.12\textwidth}\begin{center}\vspace{-14mm}ground\\truth\end{center}\end{minipage}\includegraphics[height=2cm,width=2cm]{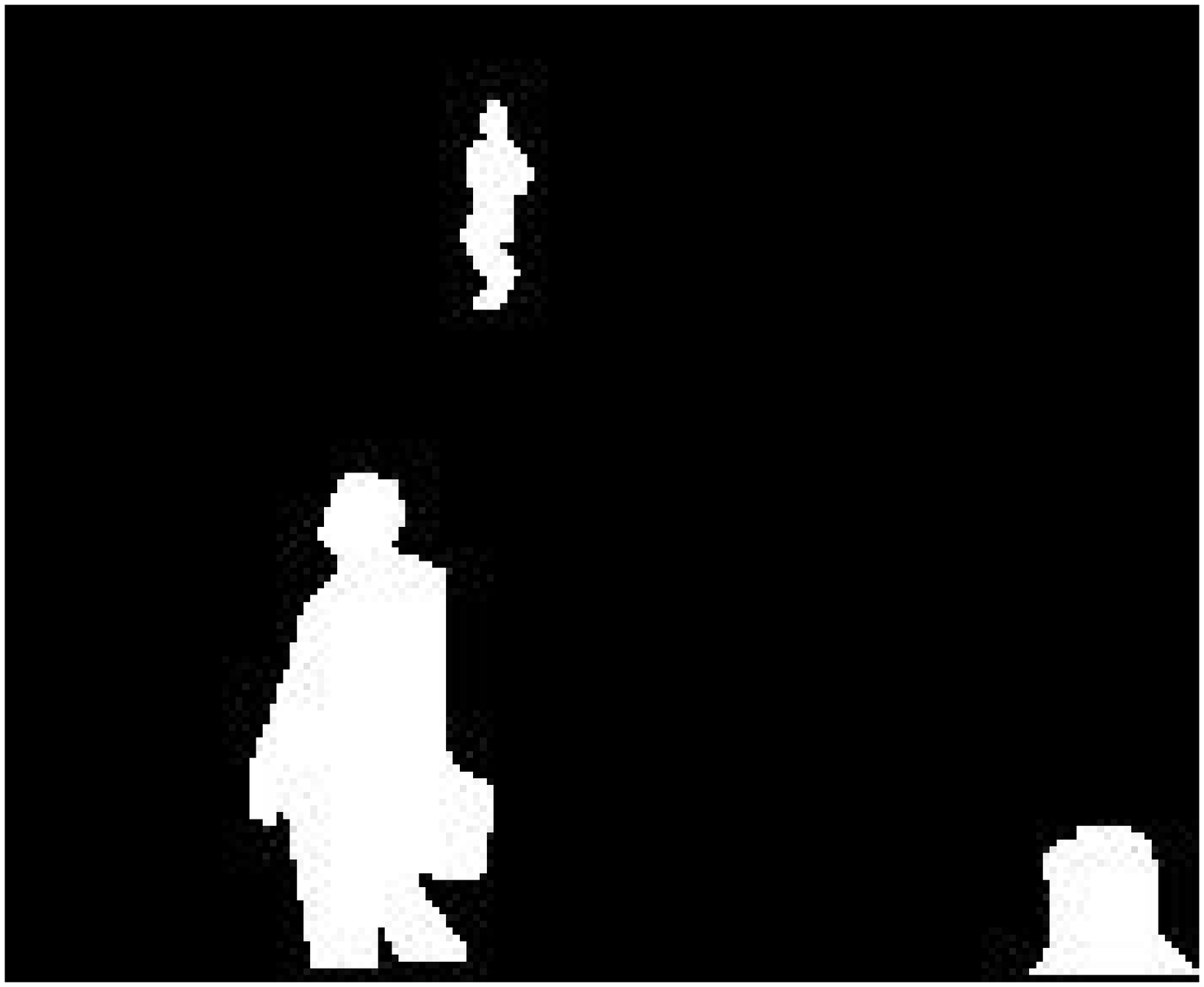}~\includegraphics[height=2cm,width=2cm]{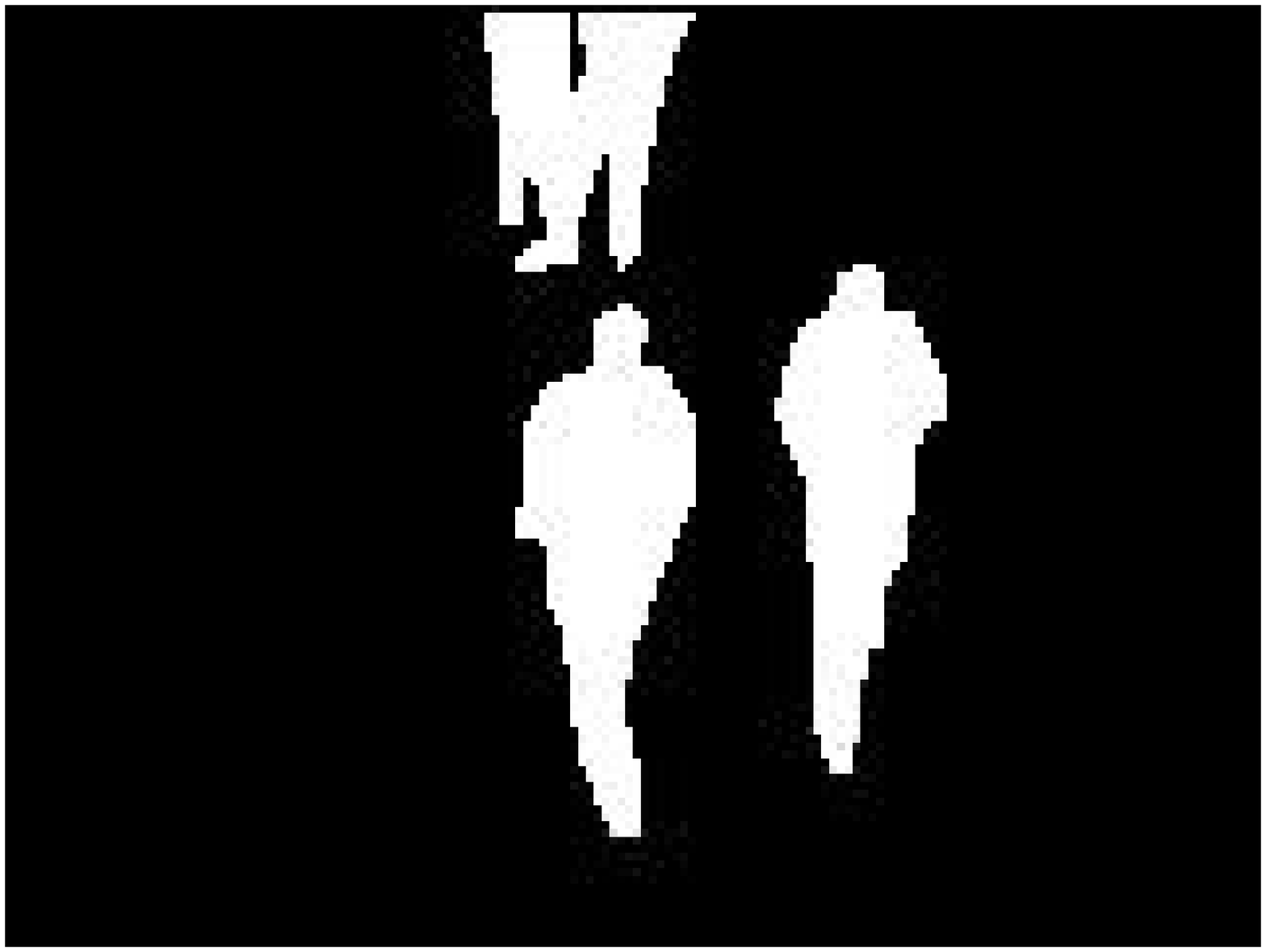}~\includegraphics[height=2cm,width=2cm]{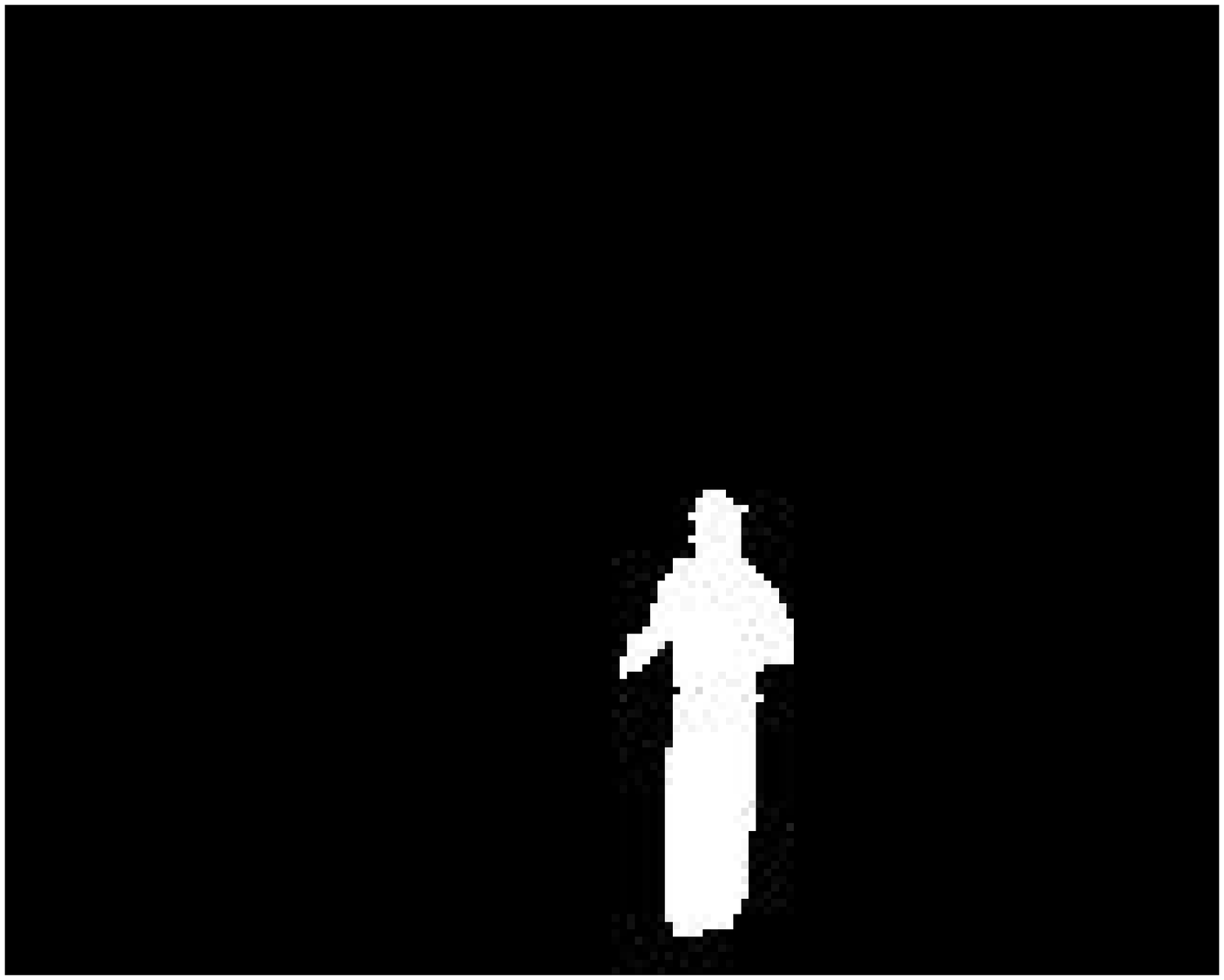}~\includegraphics[height=2cm,width=2cm]{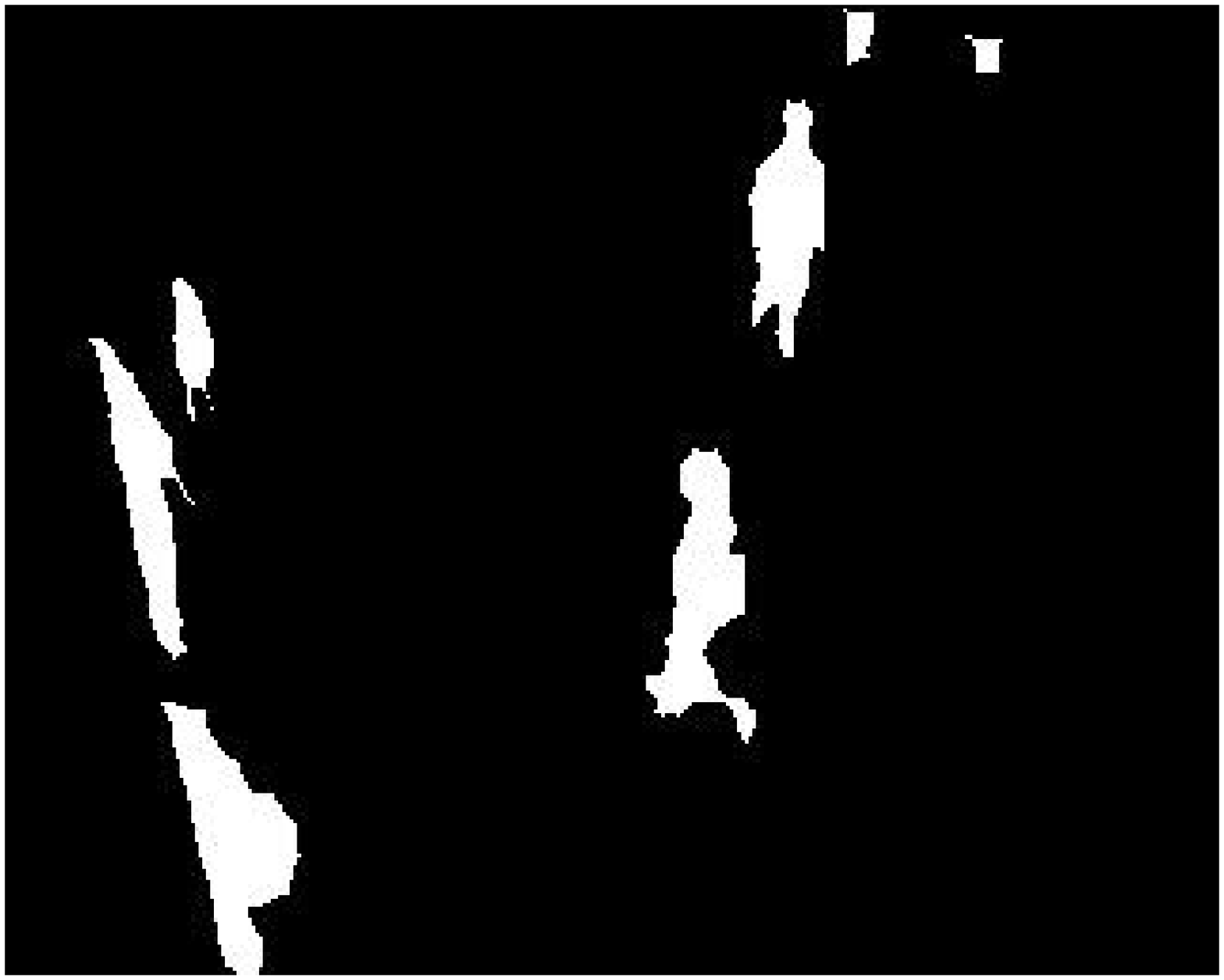}
\caption{One frame (from left to right: airport2180, b01842, Fountain1440 and ShoppingMall1535) of each testing video under different scenarios (the first two rows) and the ground-truth image of foreground of each testing video (the last row).}\label{fig_org}
\end{figure}

\paragraph{Testing videos} We choose four real videos, ``Hall", ``Bootstrap", ``Fountain" and ``ShoppingMall", from the dataset I2R\footnote{This dataset is available in http://perception.i2r.a-star.edu.sg/bk\_model/bk\_index.html. The authors also provide 20 ground-truth images of foregrounds for each video in this dataset.} provided by Li et al. \cite{lhgt2004}. The details of these videos are as follows: \vspace{1mm}
\begin{itemize}
\item \textbf{Hall} video contains 200 $144 \times 176$ frames (from airport2001 to airport2200); \vspace{1mm}
\item \textbf{Bootstrap} video contains 200 $120 \times 160$ frames (from b01801 to b02000); \vspace{1mm}
\item \textbf{Fountain} video contains 200 $128 \times 160$ frames (from Fountain1301 to Fountain1500); \vspace{1mm}
\item \textbf{ShoppingMall} video contains 200 $256 \times 320$ frames (from ShoppingMall1501 to ShoppingMall1700). \vspace{1mm}
\end{itemize}
We show one frame of each testing video under two different scenarios (noisy and noisy blurred), and their ground-truth images of foregrounds in Fig.\,\ref{fig_org}. Additionally, all pixel values of the testing videos are re-scaled into $[0,1]$ in our numerical experiments.

\paragraph{Parameters setting} For the ADMM, we use the following heuristics\footnote{Note from Theorem \ref{convergencethe}(i) that the successive change of each variable goes to zero as $k\rightarrow \infty$. Thus, intuitively, it is more favorable to see a decrease in the successive change as $k$ increases. This heuristic is designed based on this intuition.} to update $\beta$: we initialize $n_s=0$ and $\beta=0.6\bar{\beta}$, where $\bar{\beta}$ is given in \reff{betalowbound}. In the $k$-th iteration, we compute
\begin{eqnarray*}
fnorm^k &=&\|L^k\|_F+\|Z^k\|_F,\\
succ\_chg^{k}&=&\|L^{k}-L^{k-1}\|_F+\|Z^{k}-Z^{k-1}\|_F.
\end{eqnarray*}
Then, we increase $n_s$ by 1 if $succ\_chg^{k} > 0.99\cdot succ\_chg^{k-1}$. Obviously, $n_s$ is non-decreasing in this procedure. We then update $\beta$ as $1.1\beta$ whenever $\beta\leq1.01\bar{\beta}$ and the sequence satisfies either $n_s \geq 0.3k$ or $fnorm^k > 10^{10}$. On the other hand, for PALM, we set $c_k=d_k=\frac{\lambda_{\max}}{0.99}$.

We initialize our algorithm and the PALM at the point specified in \reff{initialpoint} with $\kappa=1$. Moreover, we terminate our ADMM by the following two-stage criterion\footnote{We use this two-stage criterion rather than computing the relative errors of all four variables ($L$, $S$, $Z$, $\Lambda$) in each iteration of our algorithm because computing matrix Frobenius norms can be expensive, especially for large scale problems. This strategy will help reduce the cost per iteration. We examine $\|L^k-L^{k-1}\|_F$ and $\|Z^k-Z^{k-1}\|_F$ in the first stage because these quantities being small intuitively implies that $\|S^k-S^{k-1}\|_F$ and $\|\Lambda^k - \Lambda^{k-1}\|_F$ are small; see the proof of Theorem~\ref{convergencethe}, particularly \eqref{limit1}, \eqref{limit2} and the discussions that follow.}: in each iteration, we check if
\begin{eqnarray*}
\frac{\|L^k-L^{k-1}\|_F+\|Z^k-Z^{k-1}\|_F}{\|L^k\|_F+\|Z^k\|_F+1} < \mathrm{Tol}_{A,1}
\end{eqnarray*}
for some $\mathrm{Tol}_{A,1} > 0$; if it holds, then we further check if
\begin{eqnarray*}
\frac{\|S^k-S^{k-1}\|_F+\|\Lambda^k-\Lambda^{k-1}\|_F}{\|S^k\|_F+\|\Lambda^k\|_F+1} < \mathrm{Tol}_{A,2}
\end{eqnarray*}
for some $\mathrm{Tol}_{A,2} > 0$. We terminate the algorithm if this latter condition is also satisfied. For the PALM, we terminate it when
\begin{eqnarray*}
\frac{\|L^k-L^{k-1}\|_F+\|S^k-S^{k-1}\|_F}{\|L^k\|_F+\|S^k\|_F+1} < \mathrm{Tol}_{P}
\end{eqnarray*}
for some $\mathrm{Tol}_{P} > 0$. The specific values of $\mathrm{Tol}_{A,1}$, $\mathrm{Tol}_{A,2}$ and $\mathrm{Tol}_{P}$ are given in the following experiments.

%{\color{red} TK: 1. Give the values of tol in appropriate places below. 2. The initialization in the code for PALM does not match the description above.}

\subsection{Comparisons between ADMM with different $\tau$ and PALM}
\label{comparisons}
In this subsection, we use the performance profile to evaluate the performances of the ADMM with different $\tau$ and the PALM for extraction under different scenarios. The performance profile is proposed by Dolan and Mor\'{e} \cite{dm2012} as a tool for evaluating and comparing the performance of a collection of solvers $\mathcal{K}$ on a set of test problems $\mathcal{J}$.

To describe this method, we assume that we have $K$ solvers and $J$ problems, and we use the iteration number as a performance measure. Then, for each problem $j$ and solver $k$, we set
\begin{eqnarray*}
\mathrm{iter}_{j,k} = \mathrm{the~iteration~number~required~to~solve~problem}~j~\mathrm{by~solver}~k.
\end{eqnarray*}
and compute the performance ratio
\begin{eqnarray}\label{perratio}
r_{j,k} = \frac{\mathrm{iter}_{j,k}}{\min\{\mathrm{iter}_{j,k}: k\in\mathcal{K}\}}.
\end{eqnarray}
The performance profile of iteration numbers is then defined as the distribution function for the performance ratio, i.e.,
\begin{eqnarray*}
\rho_{k}(\nu)=\frac{1}{J\,}\,\sharp\{j\in\mathcal{J}: r_{j,k}\leq\nu\}
\end{eqnarray*}
for $\nu\geq1$. Similarly, the performance profile of function values is obtained by using $\mathrm{fval}_{j,k}$ in place of $\mathrm{iter}_{j,k}$ in \eqref{perratio}, where $\mathrm{fval}_{j,k}$ denotes the function value at the solution given by solver $k$ for solving problem $j$. Generally speaking, for solver $k\in\mathcal{K}$, the higher $\rho_{k}(\nu)$ indicates a better performance within the factor $\nu$.

In our experiments, we evaluate the following solvers: the ADMM with $\tau=0.8$, the ADMM with $\tau=1$, the ADMM with $\tau=1.6$ and the PALM.

\begin{figure}[ht]
\centering
\subfigure[bridge regularizer]{\includegraphics[width=6.2cm,height=5cm]{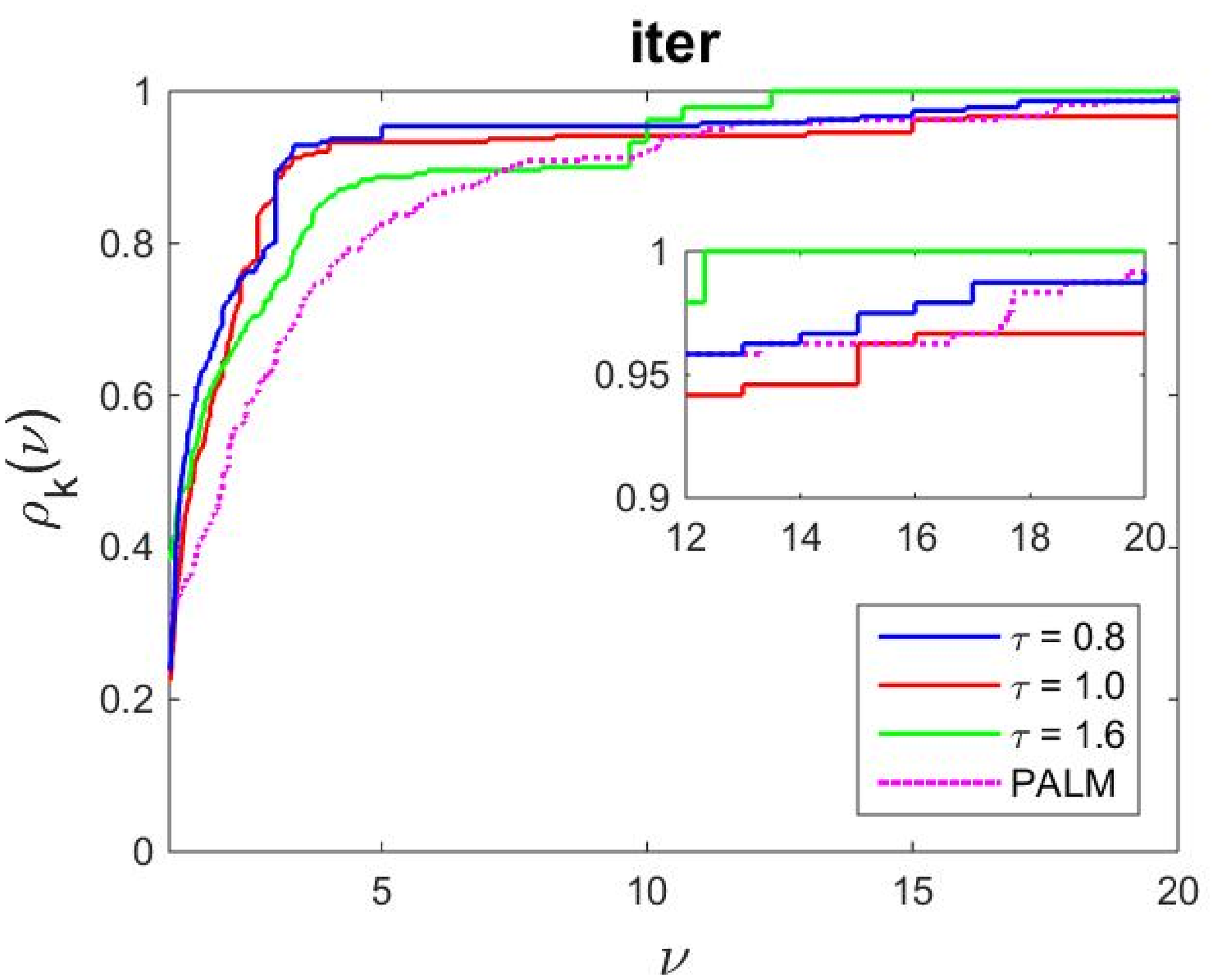}\includegraphics[width=6.2cm,height=5cm]{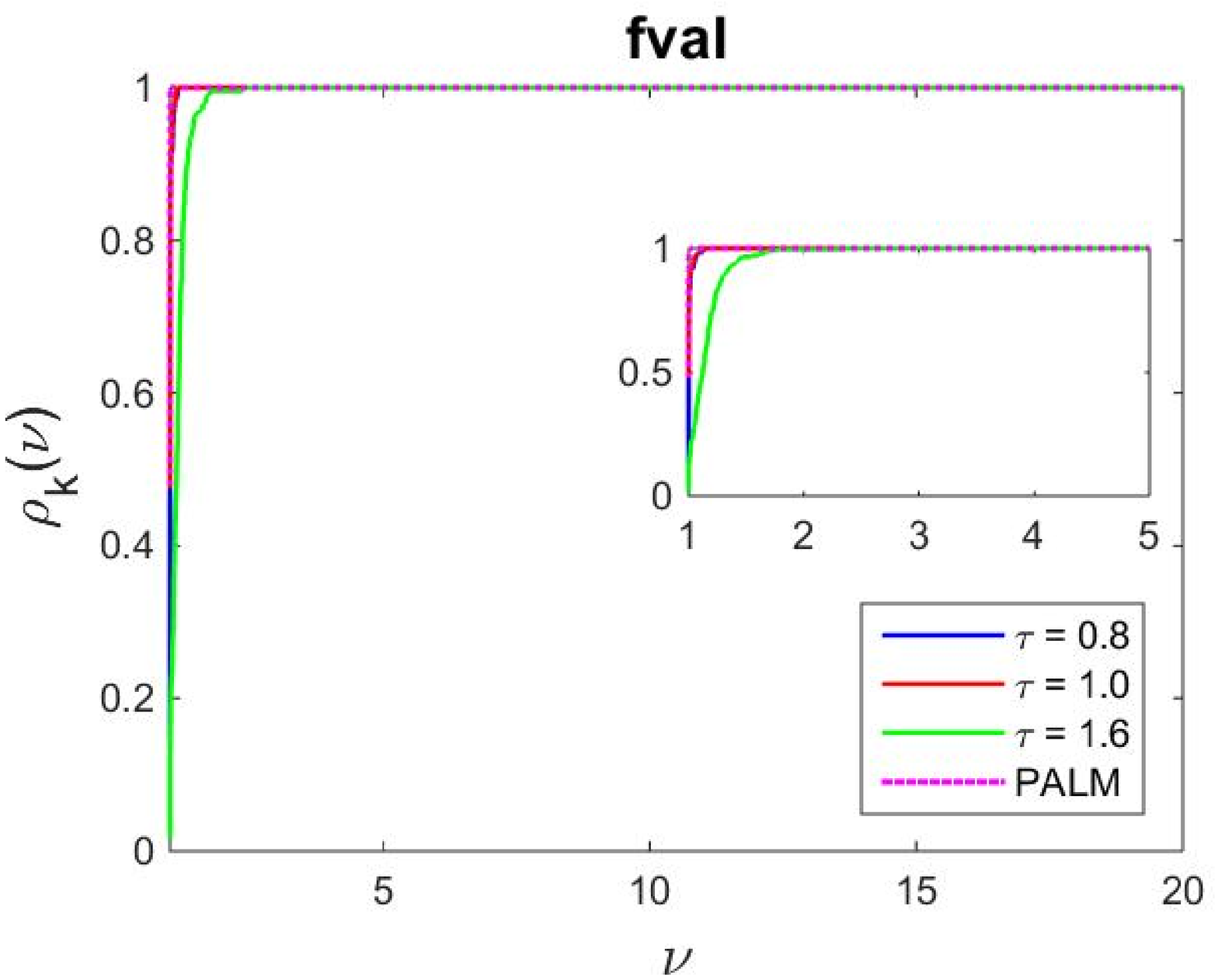}}
\subfigure[fraction regularizer]{\includegraphics[width=6.2cm,height=5cm]{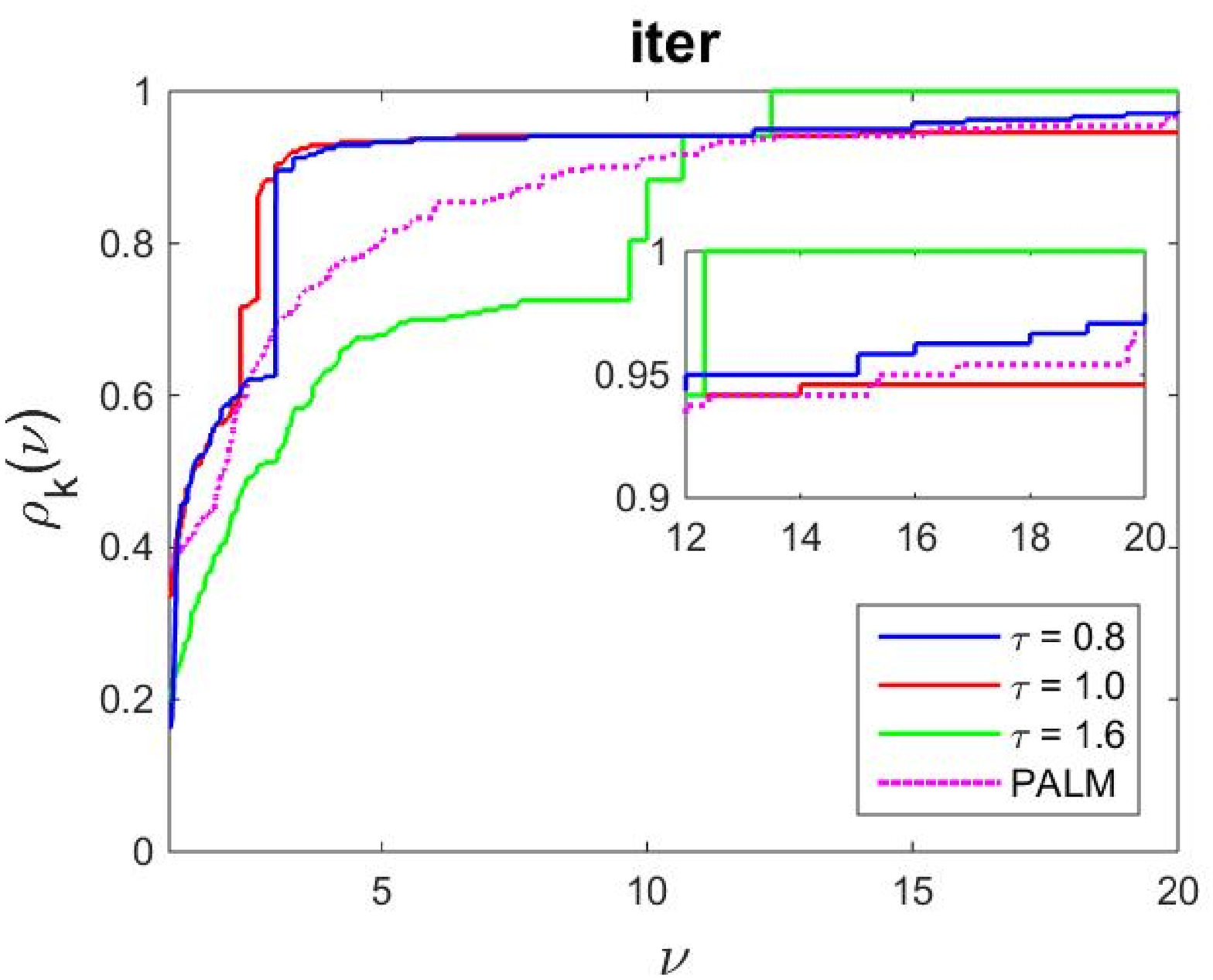}\includegraphics[width=6.2cm,height=5cm]{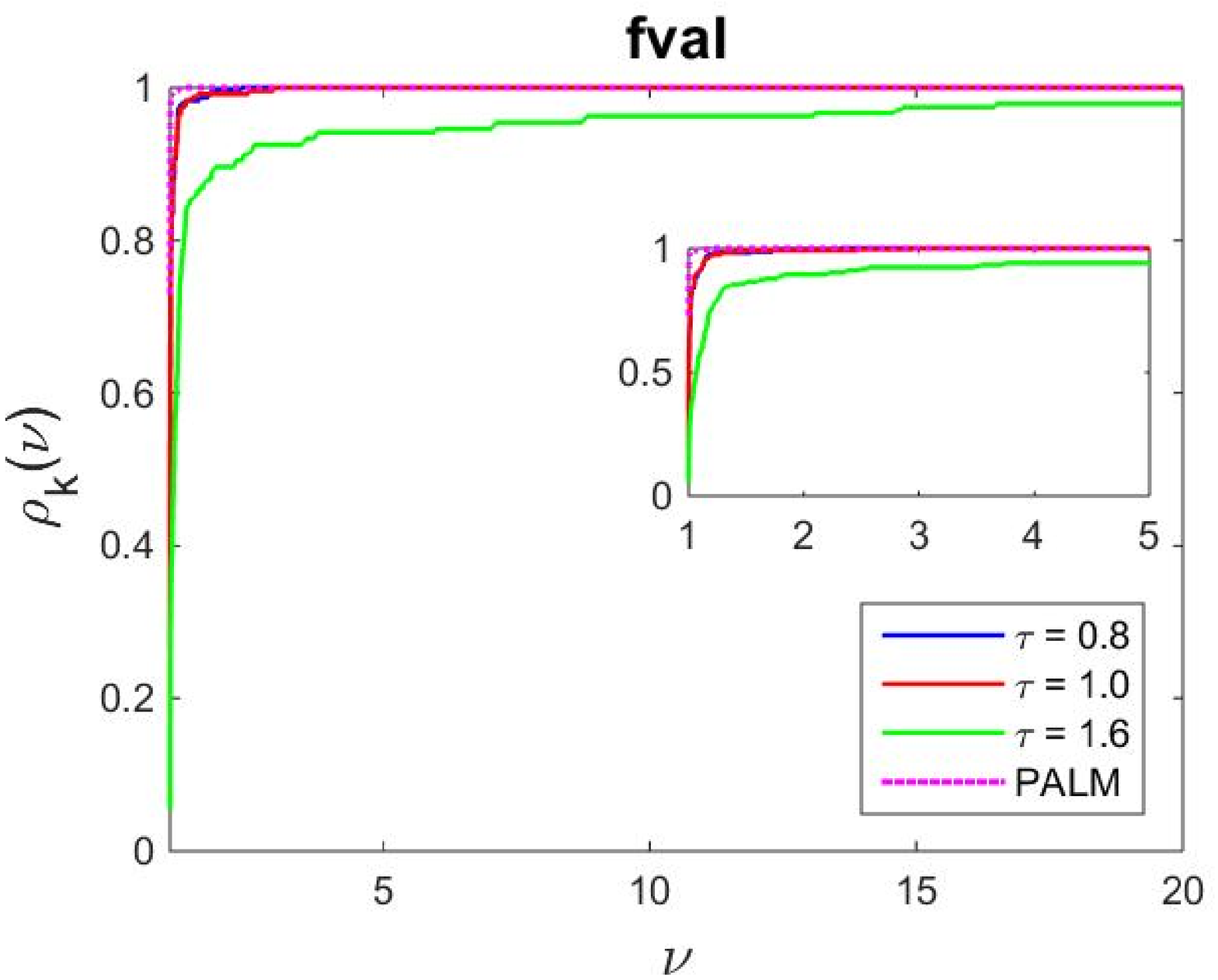}}
\subfigure[logistic regularizer]{\includegraphics[width=6.2cm,height=5cm]{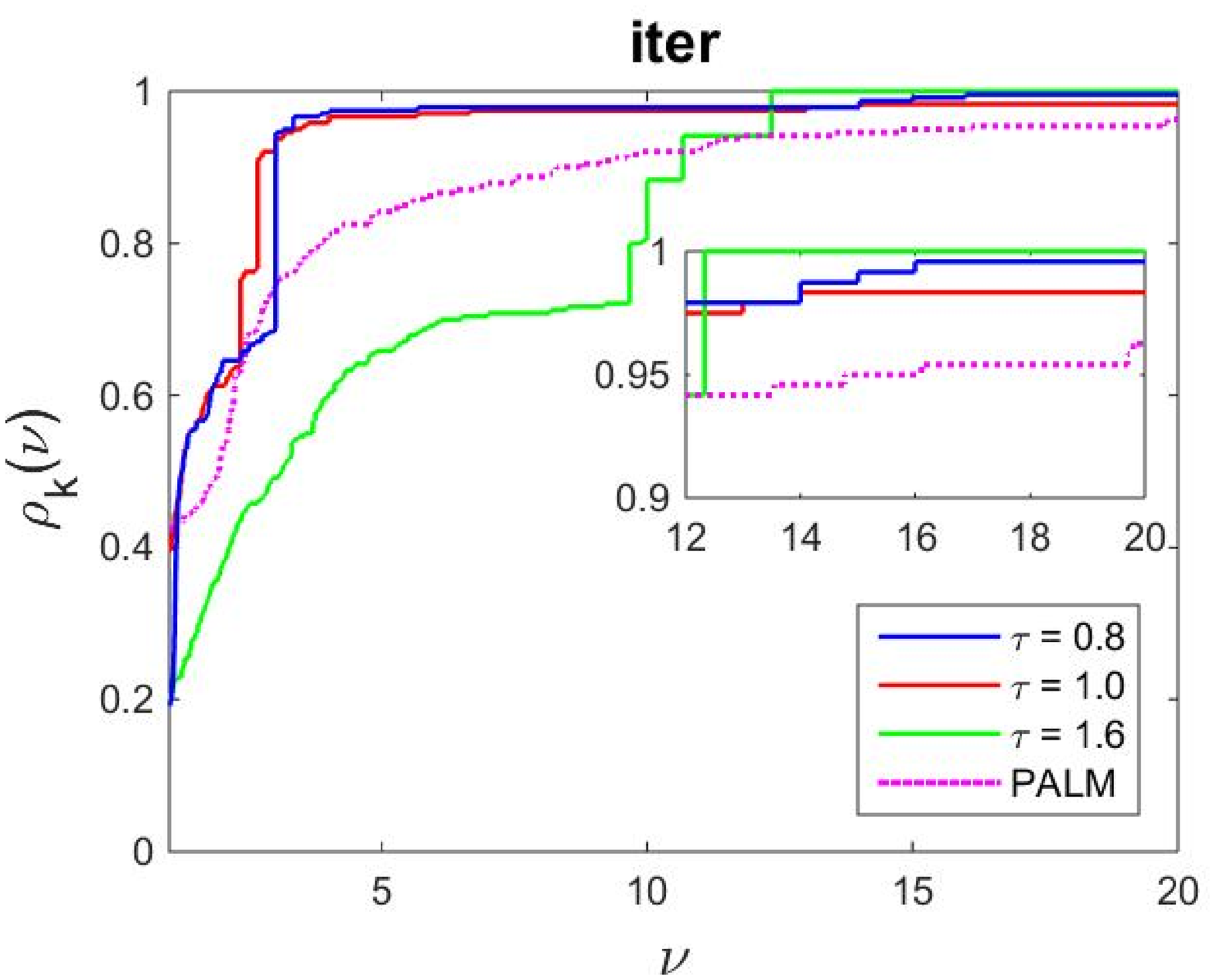}\includegraphics[width=6.2cm,height=5cm]{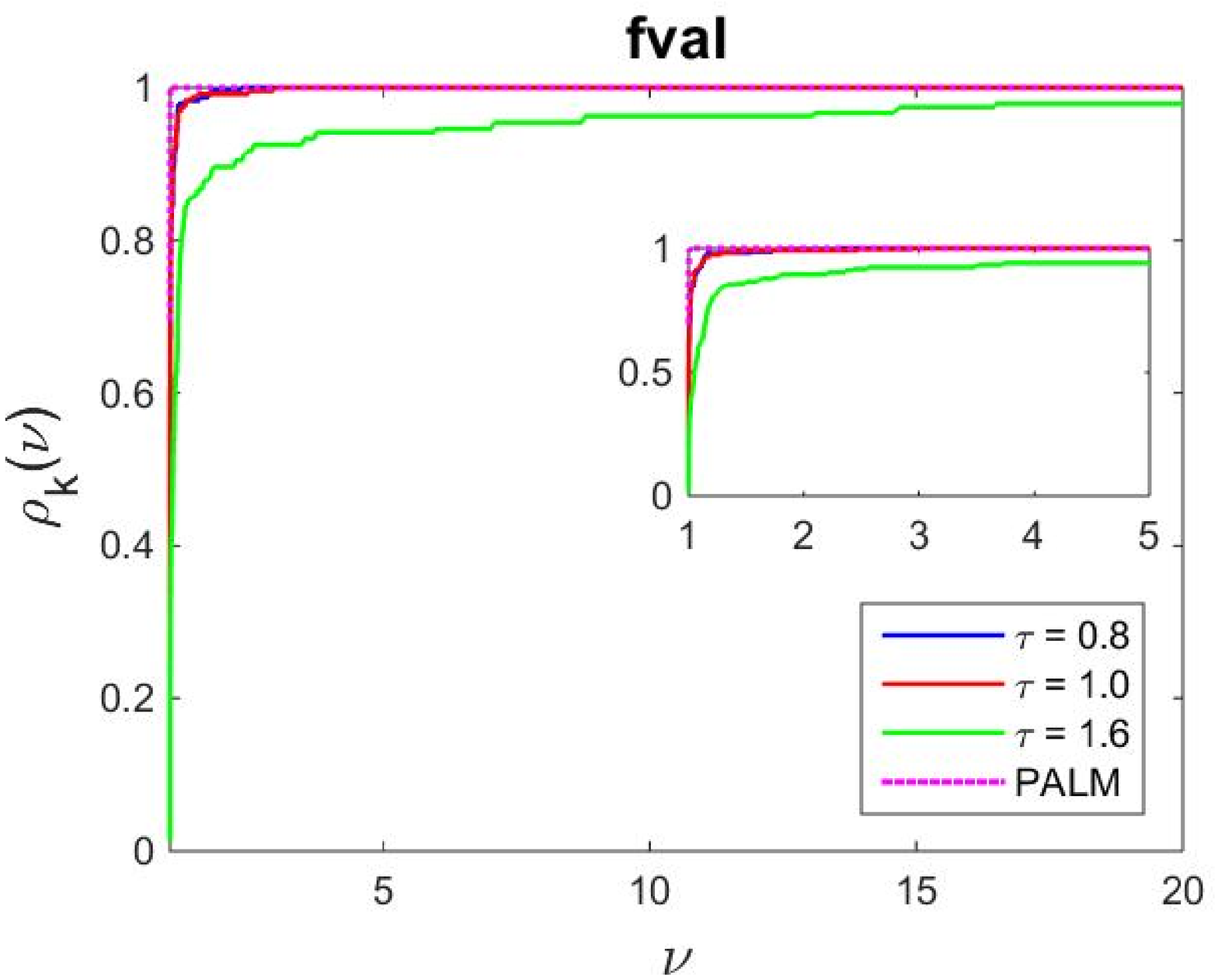}}
\caption{Performance profiles of iteration numbers (denoted by ``$\mathbf{iter}$" on the left) and function values (denoted by ``$\mathbf{fval}$" on the right) for each sparse regularizer with ${\cal A}(L+S) = L+S$. The blown-up subfigures are used to highlight the differences in a specific range of $\nu$.}\label{per_pro}
\end{figure}

For ${\cal A}(L+S) = L+S$, our test problems are described in Table \ref{proset}, where we use the four real videos introduced above as our input data in \eqref{testing_model}, with 3 choices of sparse regularizers, 10 choices of $\mu$, and 6 choices of $p$ and $\alpha$. Thus, we have 4 solvers and a total of $720$ test problems, with 240 test problems for each sparse regularizer. Moreover, we set $\mathrm{Tol}_{A,1} = 10^{-4}$, $\mathrm{Tol}_{A,2} = 5\times10^{-3}$ and $\mathrm{Tol}_{P} = 10^{-4}$. Fig.\,\ref{per_pro} shows the performance profiles of iteration numbers and function values for different regularizers under this scenario.

\begin{figure}[ht]
\centering
\subfigure[bridge regularizer]{\includegraphics[width=6.2cm,height=5cm]{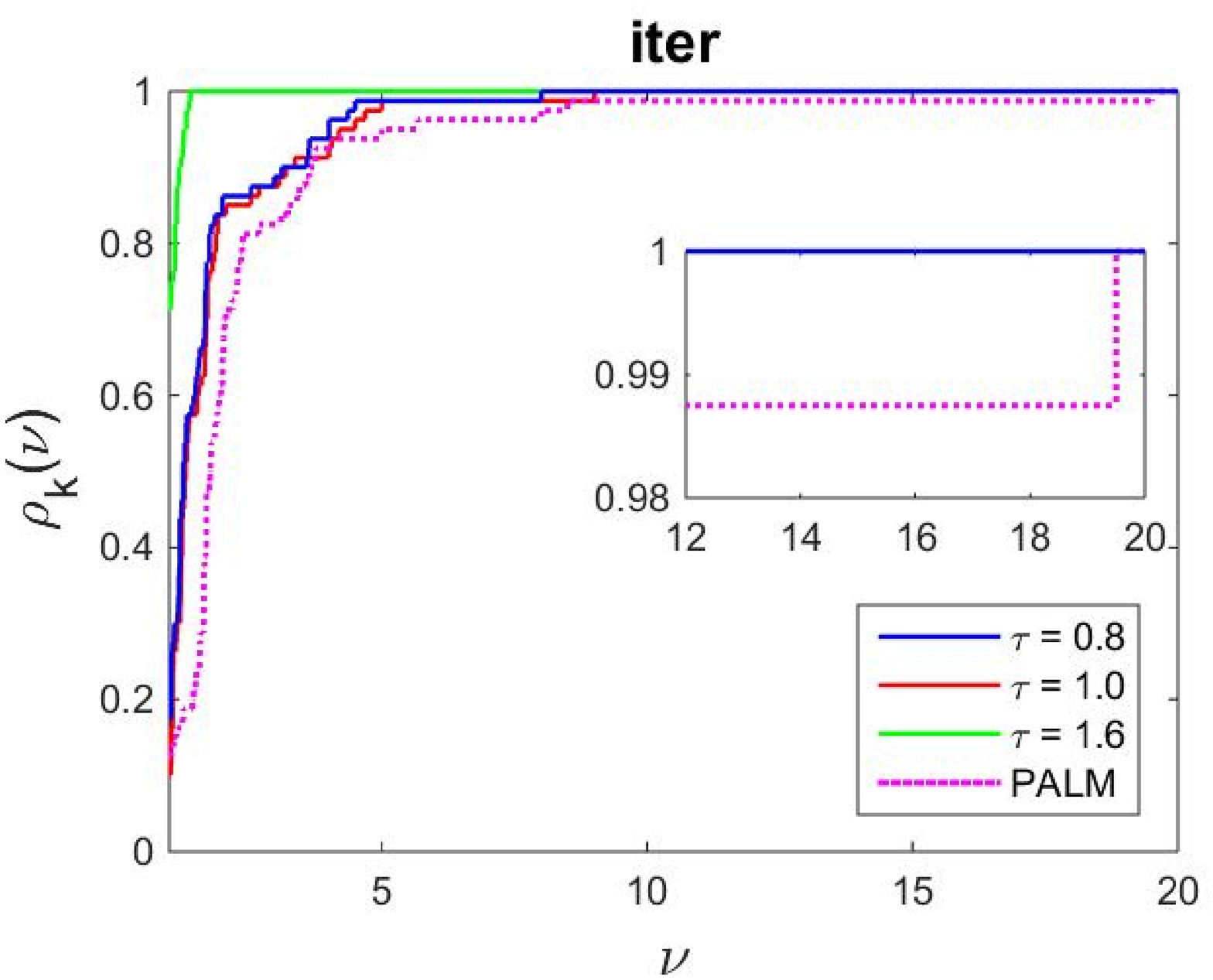}\includegraphics[width=6.2cm,height=5cm]{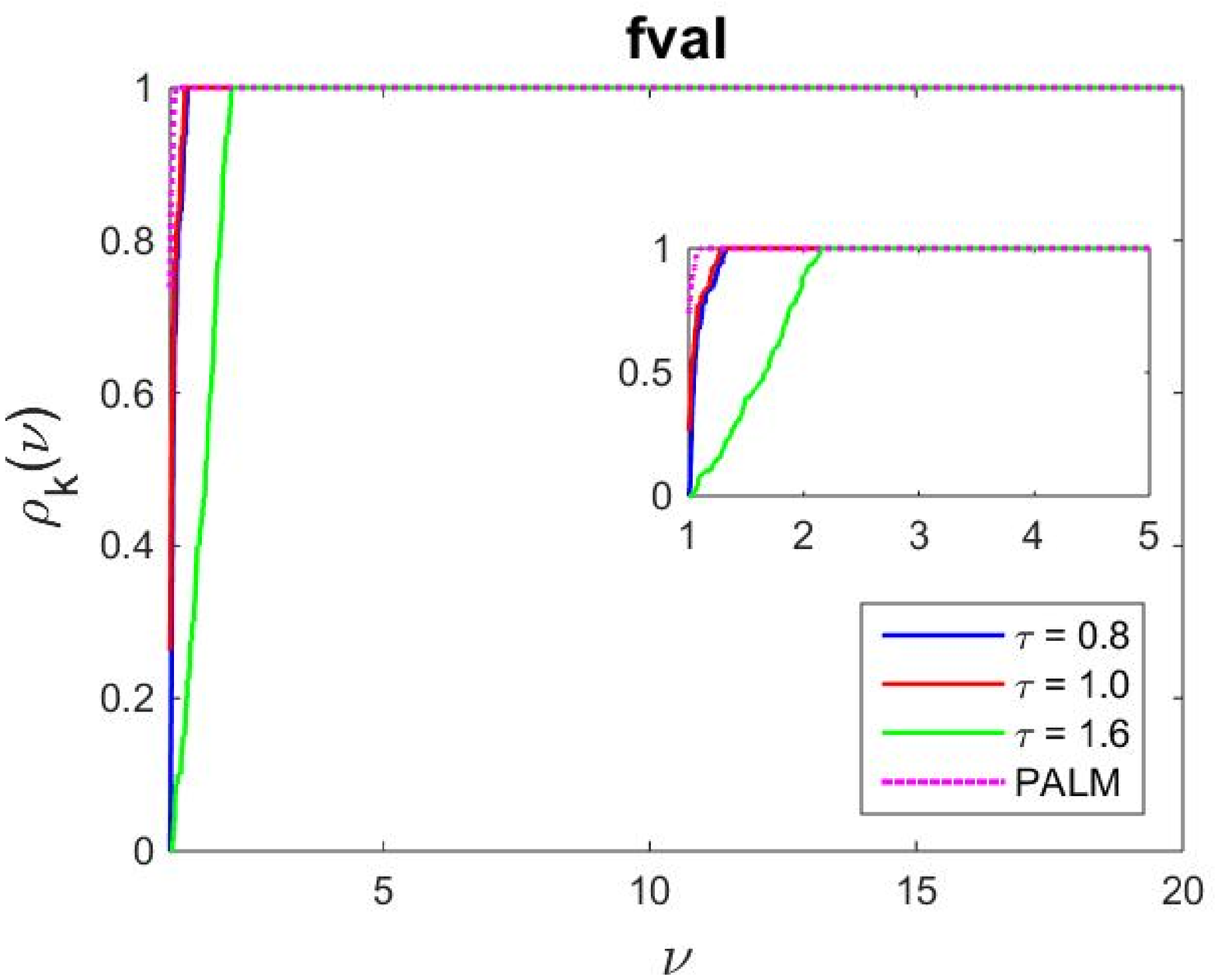}}
\subfigure[fraction regularizer]{\includegraphics[width=6.2cm,height=5cm]{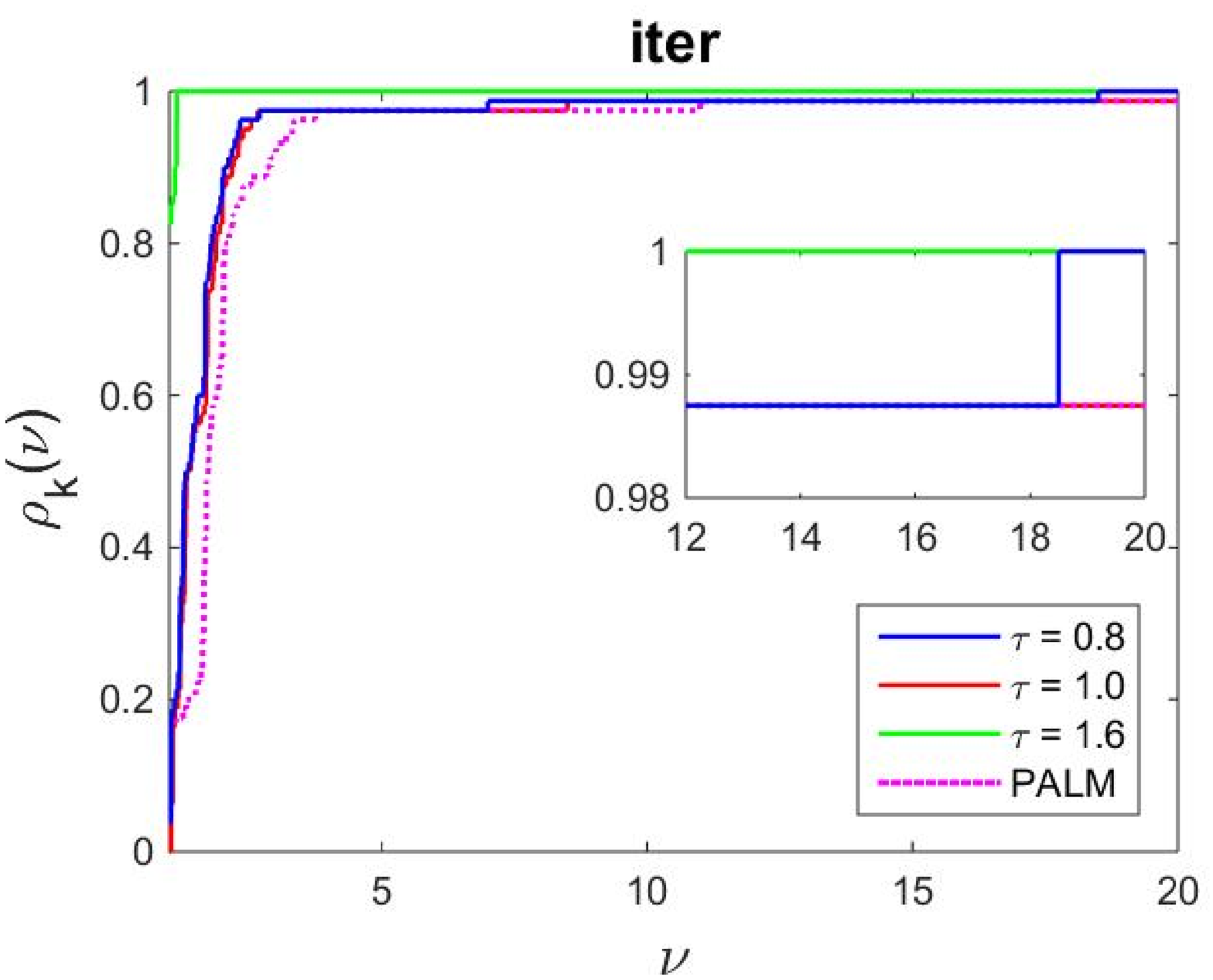}\includegraphics[width=6.2cm,height=5cm]{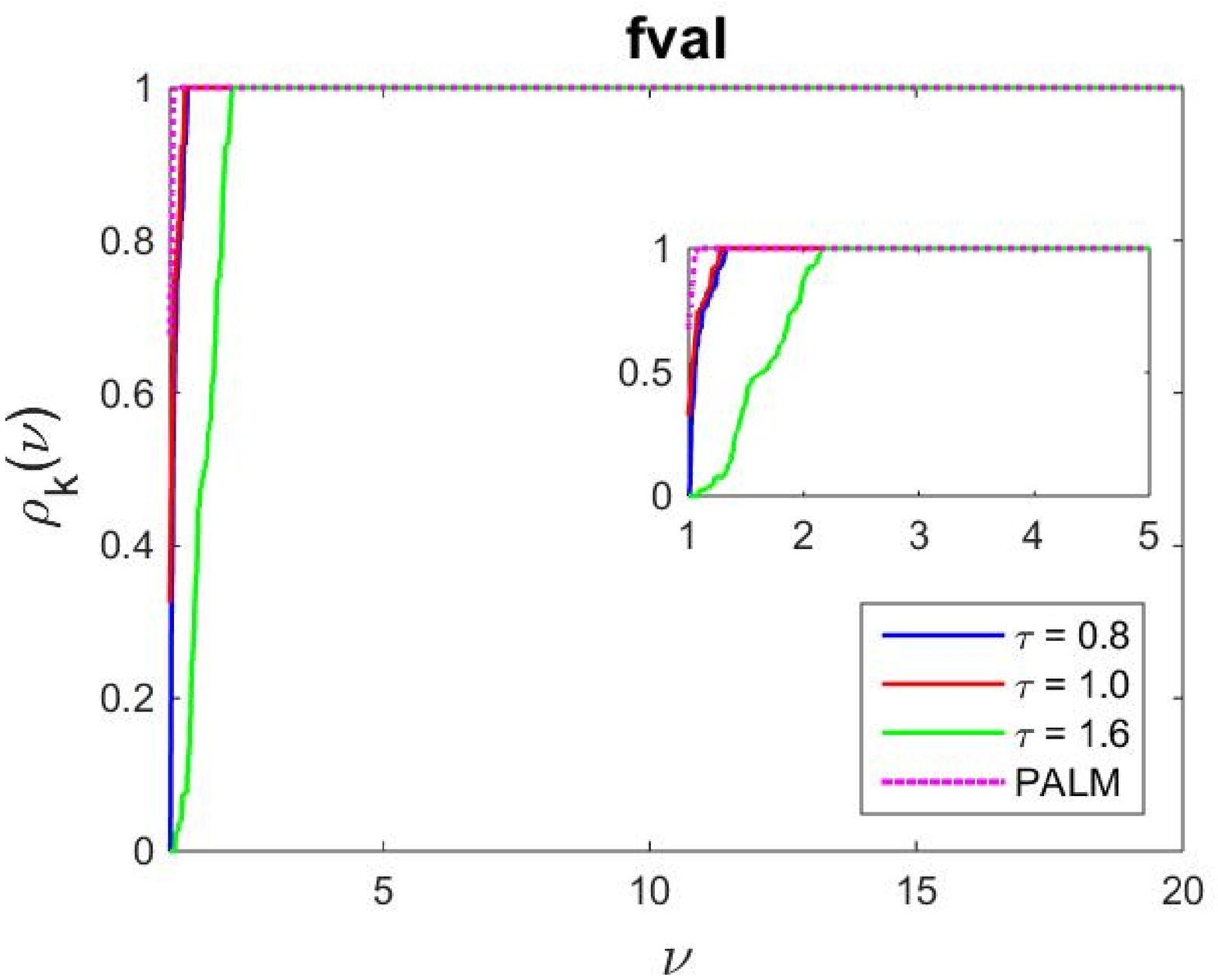}}
\subfigure[logistic regularizer]{\includegraphics[width=6.2cm,height=5cm]{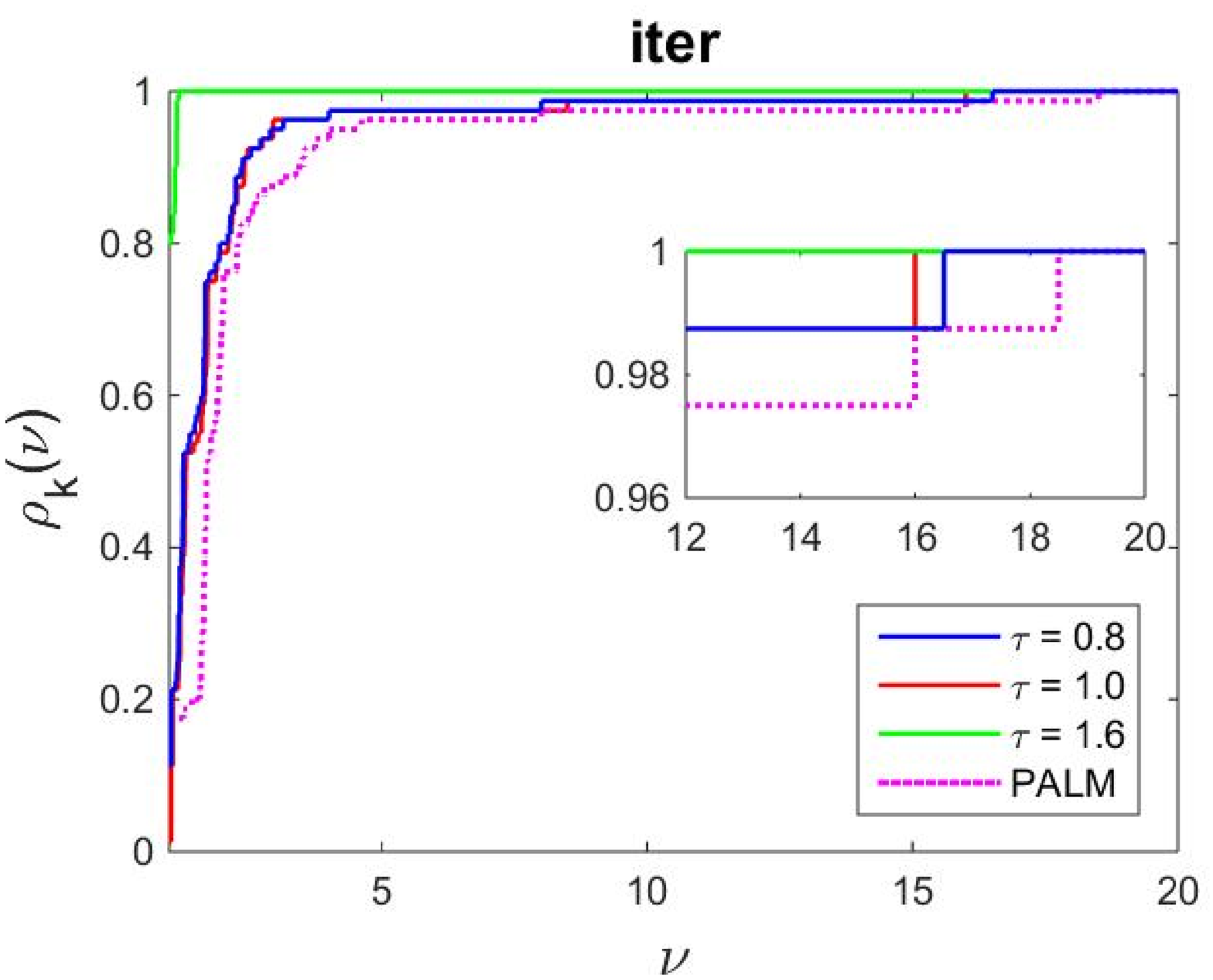}\includegraphics[width=6.2cm,height=5cm]{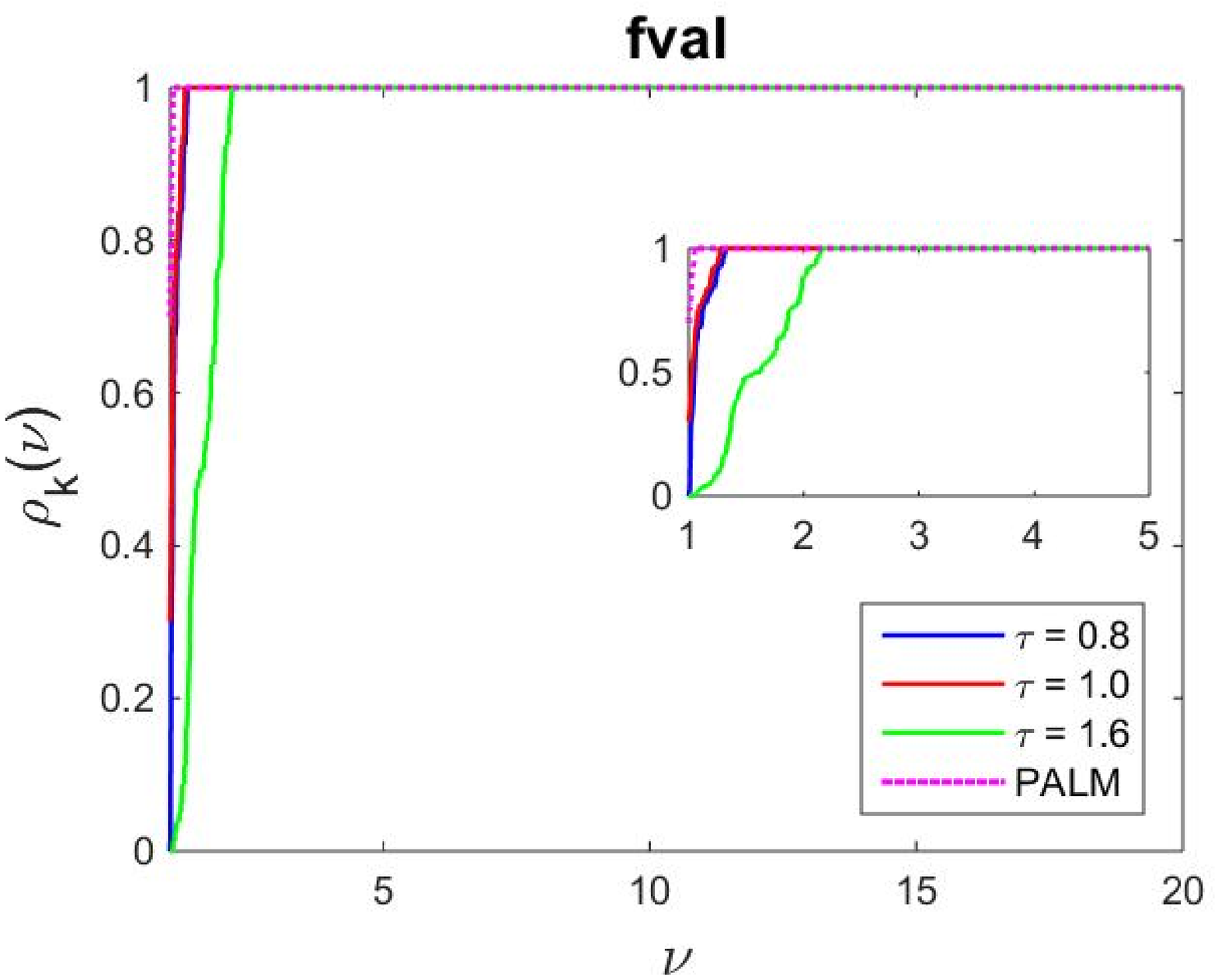}}
\caption{Performance profiles of iteration numbers (denoted by ``$\mathbf{iter}$" on the left) and function values (denoted by ``$\mathbf{fval}$" on the right) for each sparse regularizer with ${\cal A}(L+S) = H(L+S)$. The blown-up subfigures are used to highlight the differences in a specific range of $\nu$.}\label{per_pro_deblur}
\end{figure}

For ${\cal A}(L+S) = H(L+S)$, our test problems are described in Table \ref{proset_deblur}, where we use 2 choices of $p$ and $\alpha$. Thus, we have 4 solvers and a total of $240$ test problems, with 80 test problems for each sparse regularizer. In our experiments, we use the method described in \cite{hno2006} to generate the blurring matrix $H$, which can be represented as a Kronecker product $H=H_r \otimes H_c$ under the periodic boundary condition. The matlab codes\footnote{The codes are available at http://www.imm.dtu.dk/{\tiny$\sim$}pcha/HNO/ as a supplement to the book \cite{hno2006}.} that generate $H_r$ and $H_c$ are shown below, where ``frame\_size" is the size of each frame:
\vspace{2mm}
\begin{verbatim}
  [P, center] = psfGauss(frame_size, 1);
  [Hr, Hc] = kronDecomp(P, center, 'periodic');
\end{verbatim}
\vspace{2mm}
Moreover, we set $\mathrm{Tol}_{A,1} = 5\times10^{-3}$, $\mathrm{Tol}_{A,2} = 10^{-2}$ and $\mathrm{Tol}_{P} = 3\times10^{-3}$. Fig.\,\ref{per_pro_deblur} shows the performance profiles under this scenario.

It is not hard to see from Fig.\,\ref{per_pro} and Fig.\,\ref{per_pro_deblur} that the performance profiles of iteration numbers for the ADMM with $\tau=0.8$ and $\tau=1$ usually lie above those for the PALM; and their performance profiles of function values are almost the same. This shows that the ADMM with $\tau=0.8$ or $\tau=1$ takes less iterations for solving all the test problems while giving comparable function values. For bridge regularizer in the case where ${\cal A}(L+S) = L+S$ (see Fig.\,\ref{per_pro}(a)) and in the case where ${\cal A}(L+S) = H(L+S)$ (see Fig.\,\ref{per_pro_deblur}(a)), we can see that the ADMM with $\tau = 0.8$ sightly outperforms the ADMM with $\tau = 1$ in terms of the number of iterations. For other regularizers, their performances are comparable. Additionally, for the ADMM with $\tau=1.6$, we can see from Fig.\,\ref{per_pro} and Fig.\,\ref{per_pro_deblur} that it always terminates with the worst function value, although it is always fastest in the case where ${\cal A}(L+S) = H(L+S)$ (see Fig.\,\ref{per_pro_deblur}).

\begin{figure}[ht]
\centering
\subfigure[bridge regularizer]{\includegraphics[width=6.2cm,height=5cm]{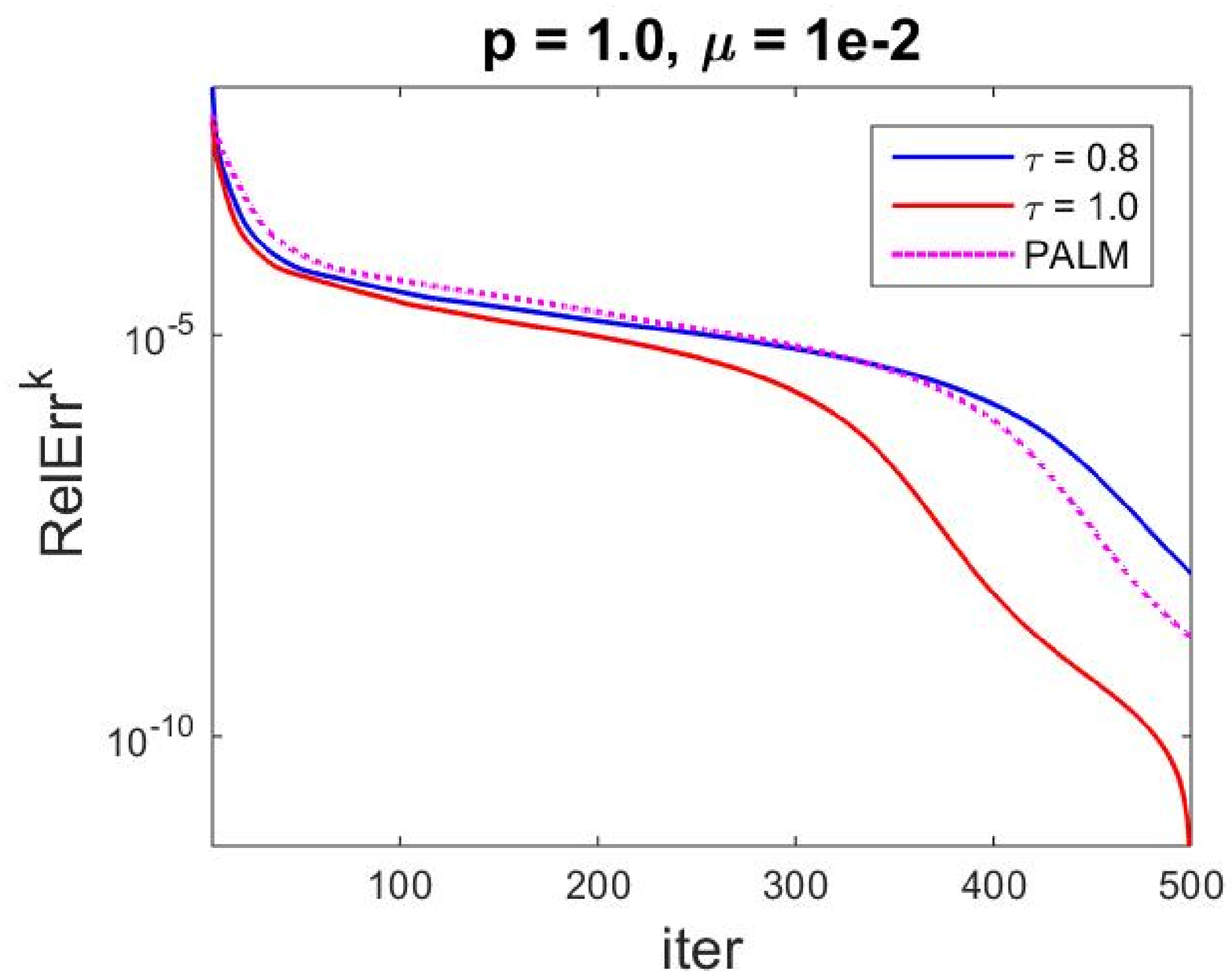}
\includegraphics[width=6.2cm,height=5cm]{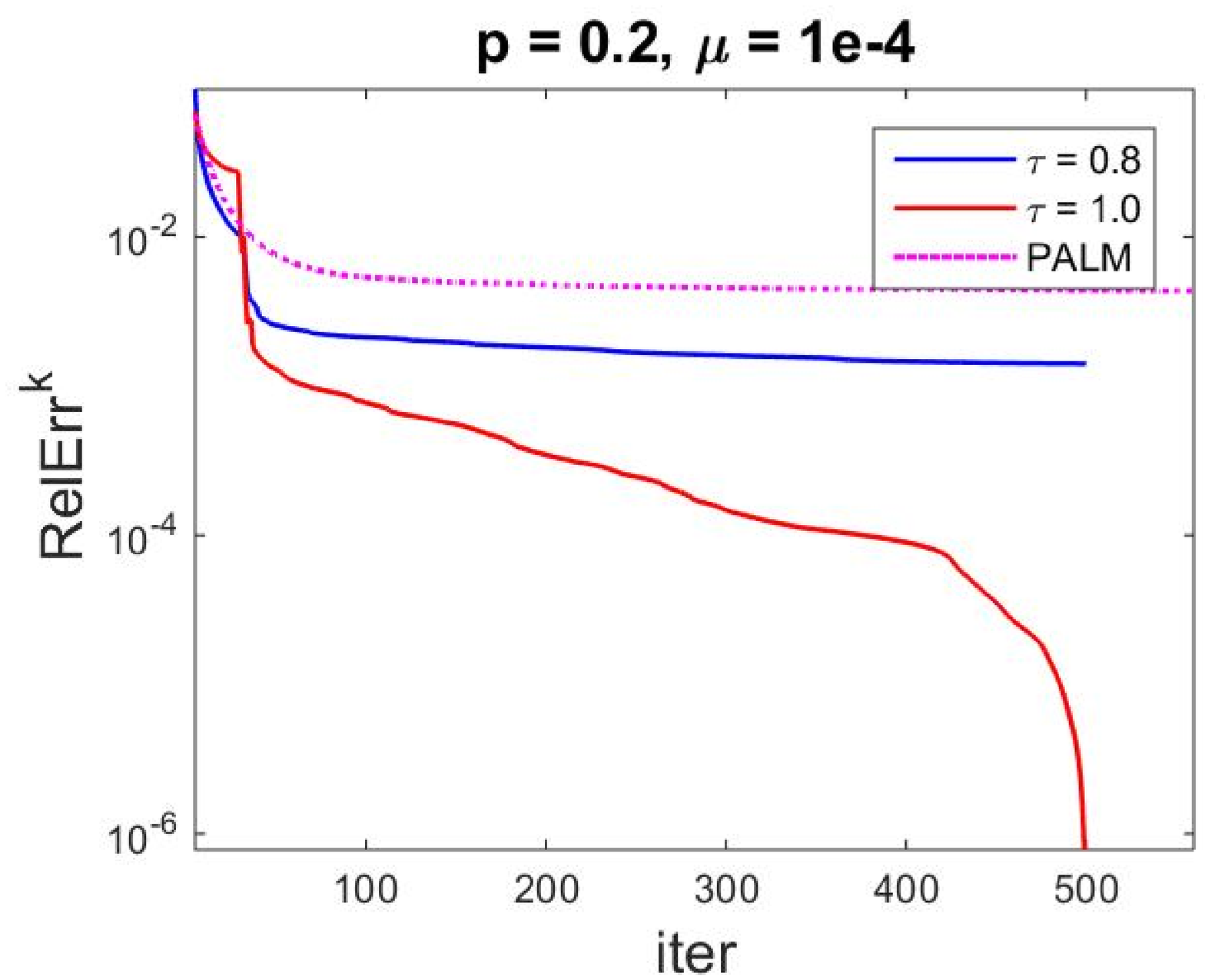}}
\subfigure[fraction regularizer]{\includegraphics[width=6.2cm,height=5cm]{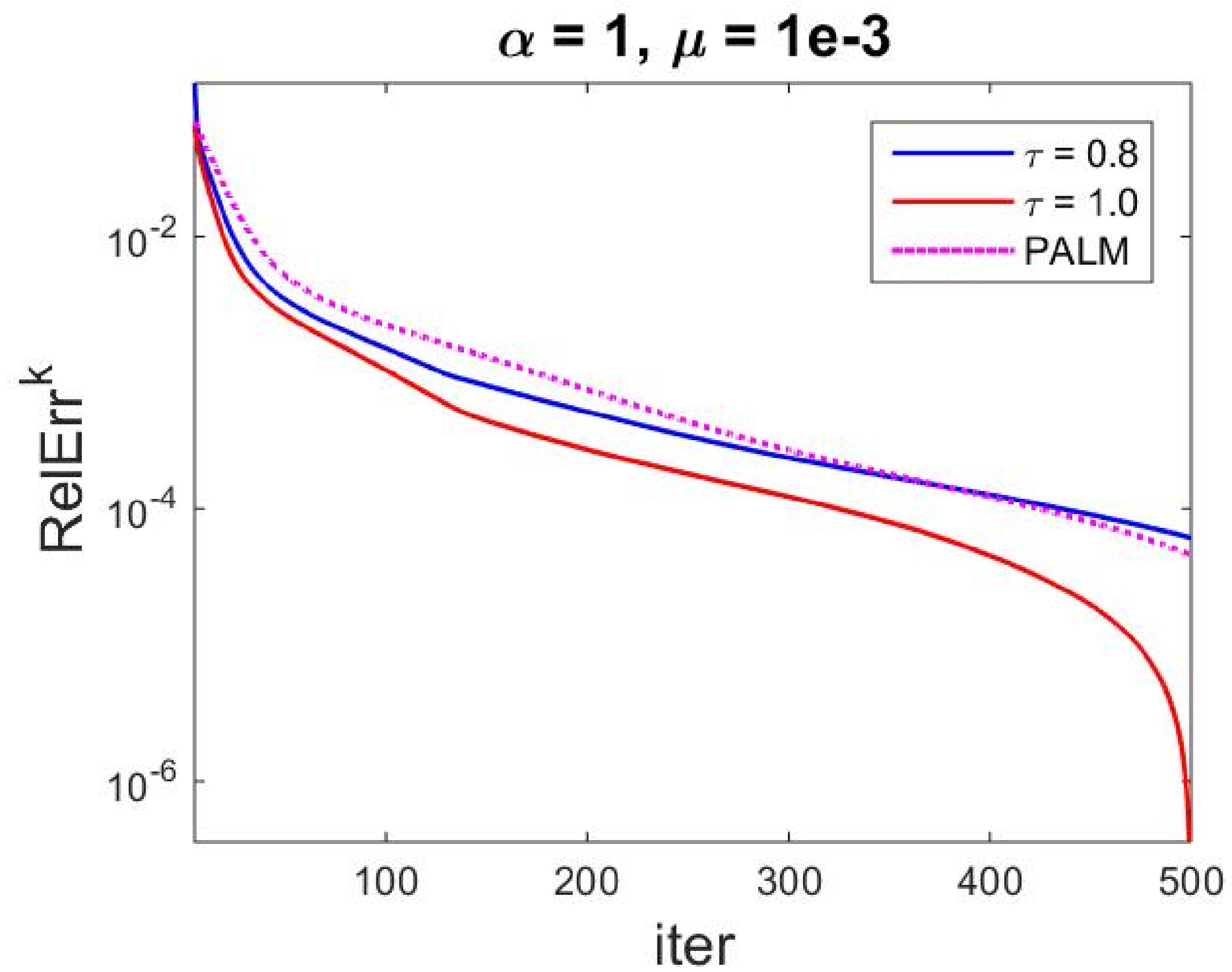}
\includegraphics[width=6.2cm,height=5cm]{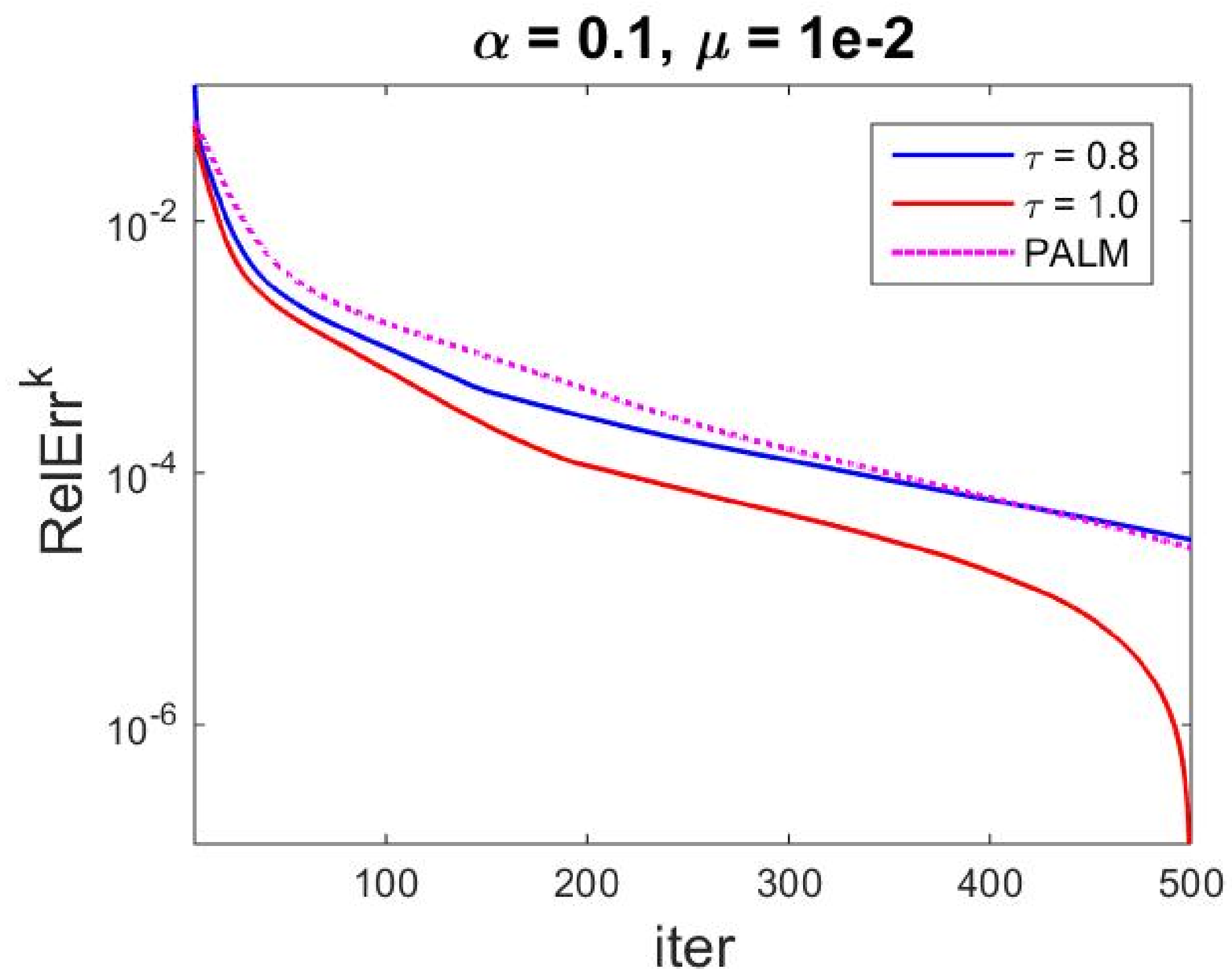}}
\subfigure[logistic regularizer]{\includegraphics[width=6.2cm,height=5cm]{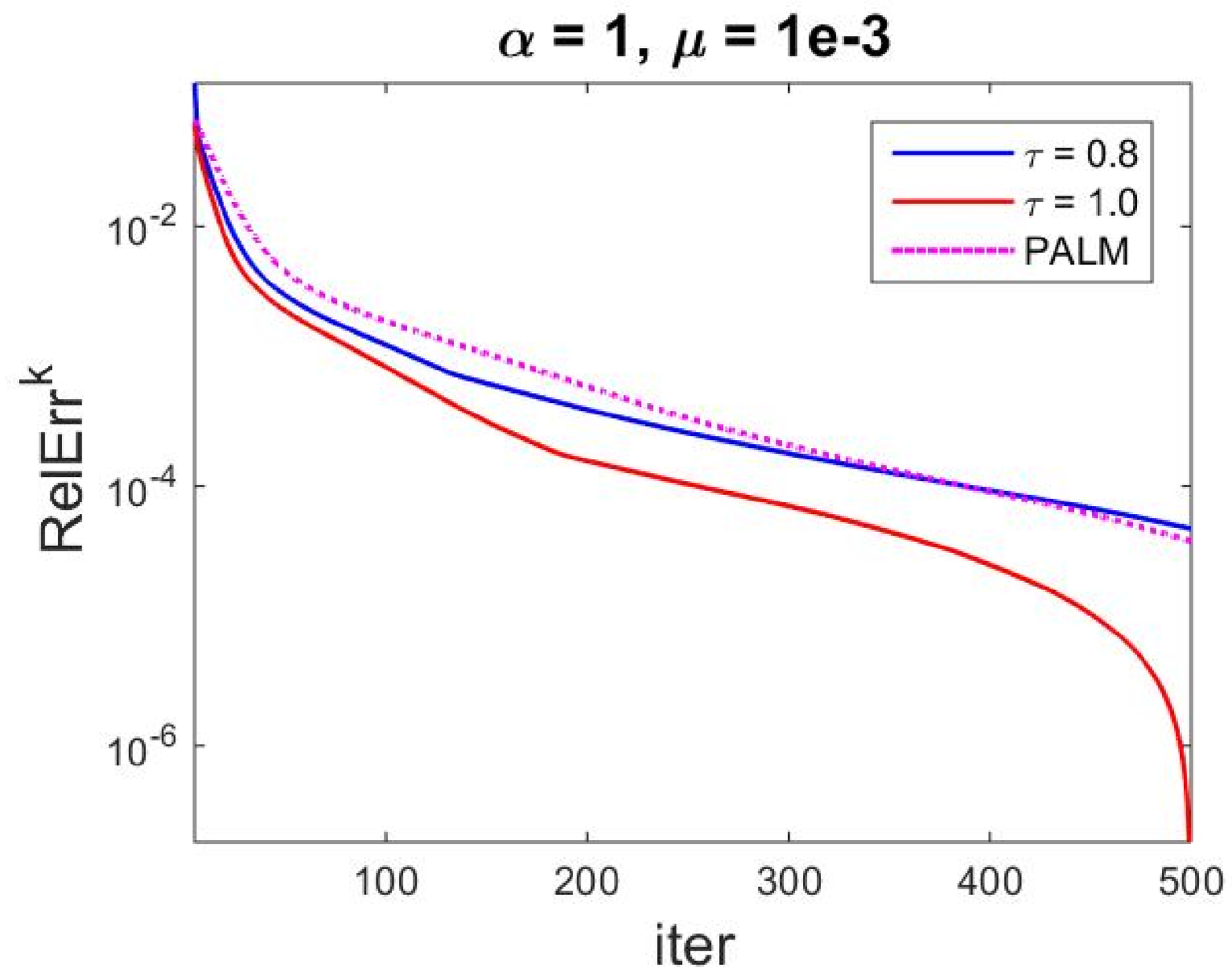}
\includegraphics[width=6.2cm,height=5cm]{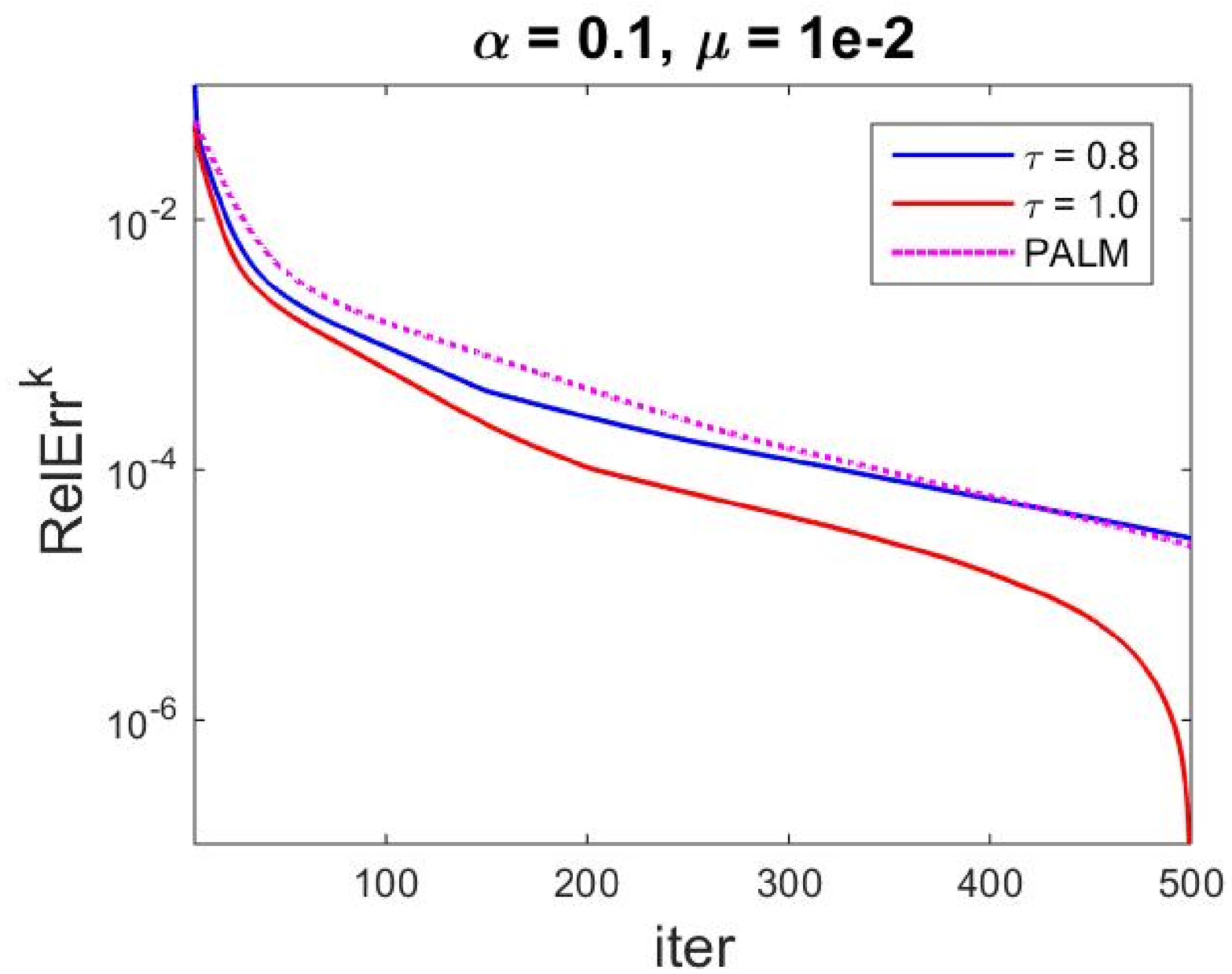}}
\caption{The $\mathrm{RelErr}^k$ vs the number of iterations for each sparse regularizer}\label{fval_vs_it}
\end{figure}

To better visualize the performance of the algorithms in terms of function values, we also plot $\mathrm{RelErr}^k := |\mathcal{F}(L^k, S^k) - \mathcal{F}_{\min}|/\mathcal{F}_{\min}$ against the number of iterations for each algorithm, where $\mathcal{F}(L^k, S^k)$ denotes the objective value obtained by each algorithm at $(L^k, S^k)$ and $\mathcal{F}_{\min}$ denotes the minimum of the objective values obtained from all algorithms. We only consider the ADMM with $\tau = 0.8$, the ADMM with $\tau = 1$ and the PALM, and terminate them only after {\em at least} 500 iterations {\em and} the termination criteria are satisfied with $\mathrm{Tol}_{A,1} = 10^{-5}$, $\mathrm{Tol}_{A,2} = 5\times10^{-4}$ and $\mathrm{Tol}_{P} = 10^{-5}$. For brevity, we focus on the scenario
${\cal A}(L+S) = L+S$ and use the ``Hall" video. The results are presented in Fig.\,\ref{fval_vs_it}, from which we can see that the ADMM with $\tau=1$ or $\tau=0.8$ performs better than PALM for those particular instances.

\begin{table}[ht]
\caption{Problem setting for ${\cal A}(L+S) = L+S$}\label{proset}
\centering \tabcolsep 3pt
\begin{tabular}{|c|c|c|}
\hline
data &  $\mu$  & regularizers   \\
\hline
\multirow{2}{*}{4 real videos} & 5e-1, 1e-1, 5e-2, 1e-2, 5e-3 & bridge: $p=0.2, 0.4, 0.5, 0.6, 0.8, 1$  \\
 &1e-3, 5e-4, 1e-4, 5e-5, 1e-5 & fraction/logistic: $\alpha=0.01, 0.1, 1, 2, 5, 10$  \\
\hline
\end{tabular}
\end{table}

\begin{table}[ht]
\caption{Problem setting for ${\cal A}(L+S) = H(L+S)$}\label{proset_deblur}
\centering \tabcolsep 3pt
\begin{tabular}{|c|c|c|}
\hline
data &  $\mu$  & regularizers   \\
\hline
\multirow{2}{*}{4 real videos} & 5e-1, 1e-1, 5e-2, 1e-2, 5e-3 & bridge: $p=0.5, 1$  \\
 &1e-3, 5e-4, 1e-4, 5e-5, 1e-5 & fraction/logistic: $\alpha=1, 2$  \\
\hline
\end{tabular}
\end{table}

\subsection{Simulation Results}

In this subsection, we present some simulation results for the background/foreground extraction problem. In order to evaluate the performance in background/foreground extraction, we compare the support of the recovered foreground $S^*$ with the support of the ground-truth $\tilde{S}$ by computing the following measurement:
\begin{eqnarray*}
\textrm{F-measure}:=2\times\frac{\mathrm{precision}\cdot\mathrm{recall}}{\mathrm{precision}+\mathrm{recall}},
\end{eqnarray*}
where precision and recall are defined as
\begin{eqnarray*}
\mathrm{precision}:=\frac{\mathrm{TP}}{\mathrm{TP}+\mathrm{FP}}, \quad \mathrm{recall}:=\frac{\mathrm{TP}}{\mathrm{TP}+\mathrm{FN}},
\end{eqnarray*}
in which, \vspace{1mm}
\begin{itemize}
\item TP stands for true positives: the number of true foreground pixels that are recovered; \vspace{1mm}
\item FP stands for false positives: the number of background pixels that are misdetected as foreground; \vspace{1mm}
\item FN stands for false negatives: the number of true foreground pixels that are missed. \vspace{1mm}
\end{itemize}
The support of the recovered foreground $S^*$ is obtained by thresholding $S^*$ entry-wise with a threshold value (we use 1e-3 in our numerical experiments). We would like to point out that F-measure varies between 0 and 1 according to the similarity of the support of $S^*$ and $\tilde{S}$. The higher the F-measure value, the better the recovery accuracy of the support of $\tilde{S}$. The F-measure approaches the maximum value 1 if the supports of $S^*$ and $\tilde{S}$ are the same, which means the foreground is recovered completely.

In our experiments below, we choose $\tau=0.8$ for the ADMM. We also use the aforementioned four real videos as input with 3 choices of sparse regularizers and 2 choices of $p$ and $\alpha$. For each fixed $p$ and $\alpha$, we experiment with different regularization parameters $\mu$ (5e-1, 1e-1, 5e-2, 1e-2, 5e-3, 1e-3, 5e-4, 1e-4, 5e-5, 1e-5) and present only the $\mu$ corresponding to the maximal F-measure.\footnote{If the F-measures are the same, we pick the $\mu$ that corresponds to the minimal number of iterations.}

\begin{figure}[ht]
\centering
\subfigure{\begin{minipage}[t]{0.3\textwidth}\centering bri.$p=1.0$\vspace{1mm} \\
\includegraphics[height=1.8cm,width=1.8cm]{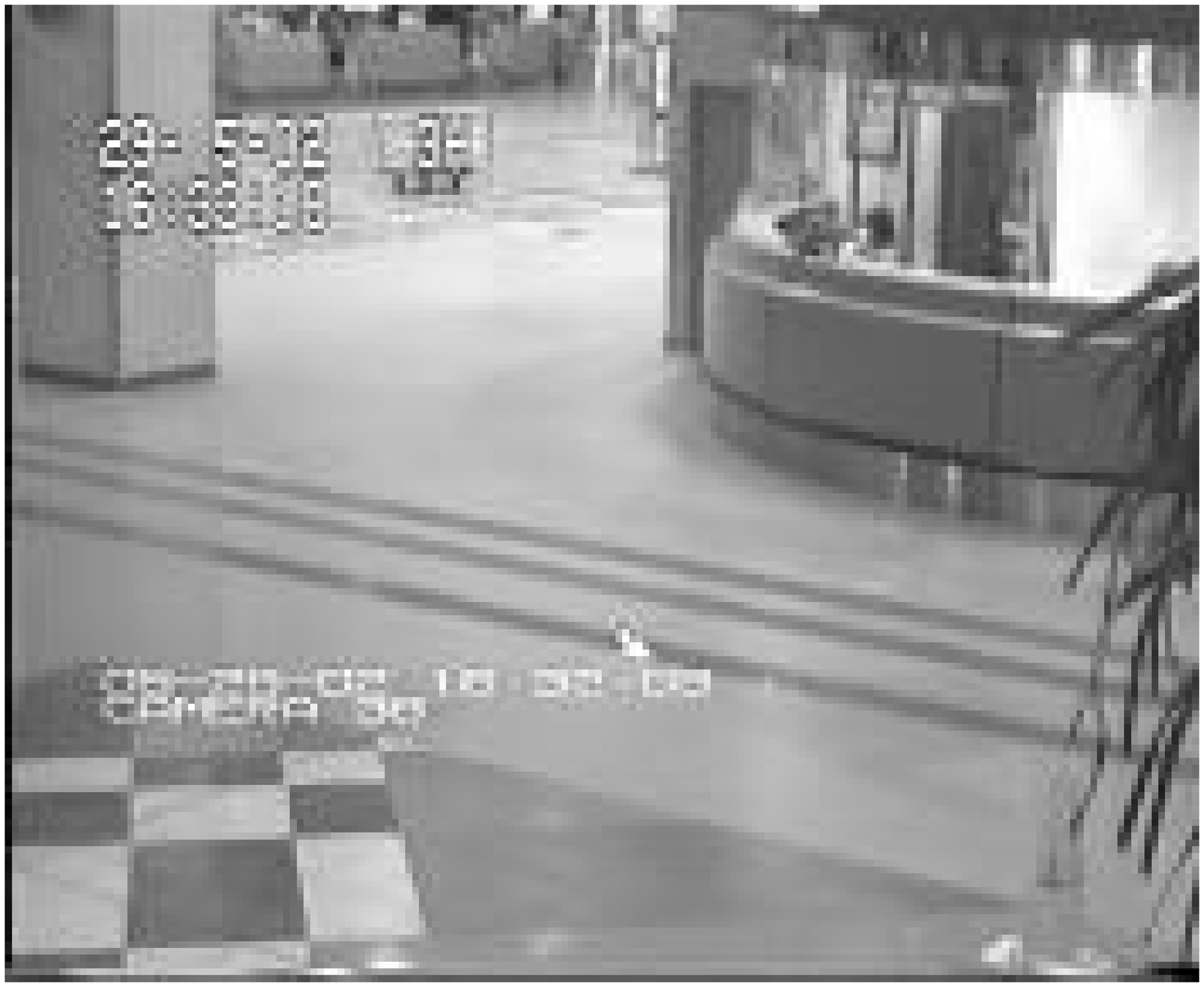}~\includegraphics[height=1.8cm,width=1.8cm]{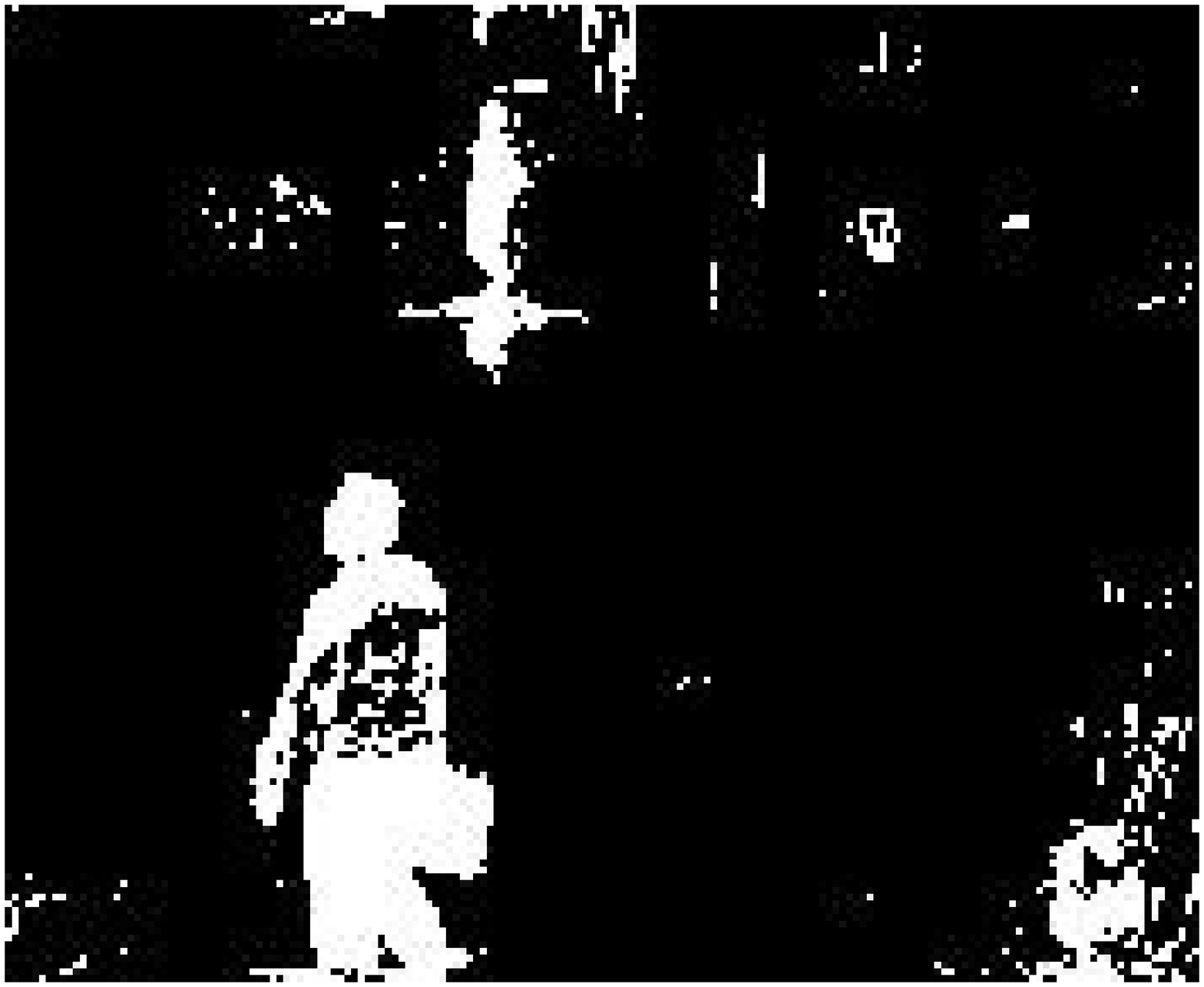}\vspace{0.5mm}\\
\includegraphics[height=1.8cm,width=1.8cm]{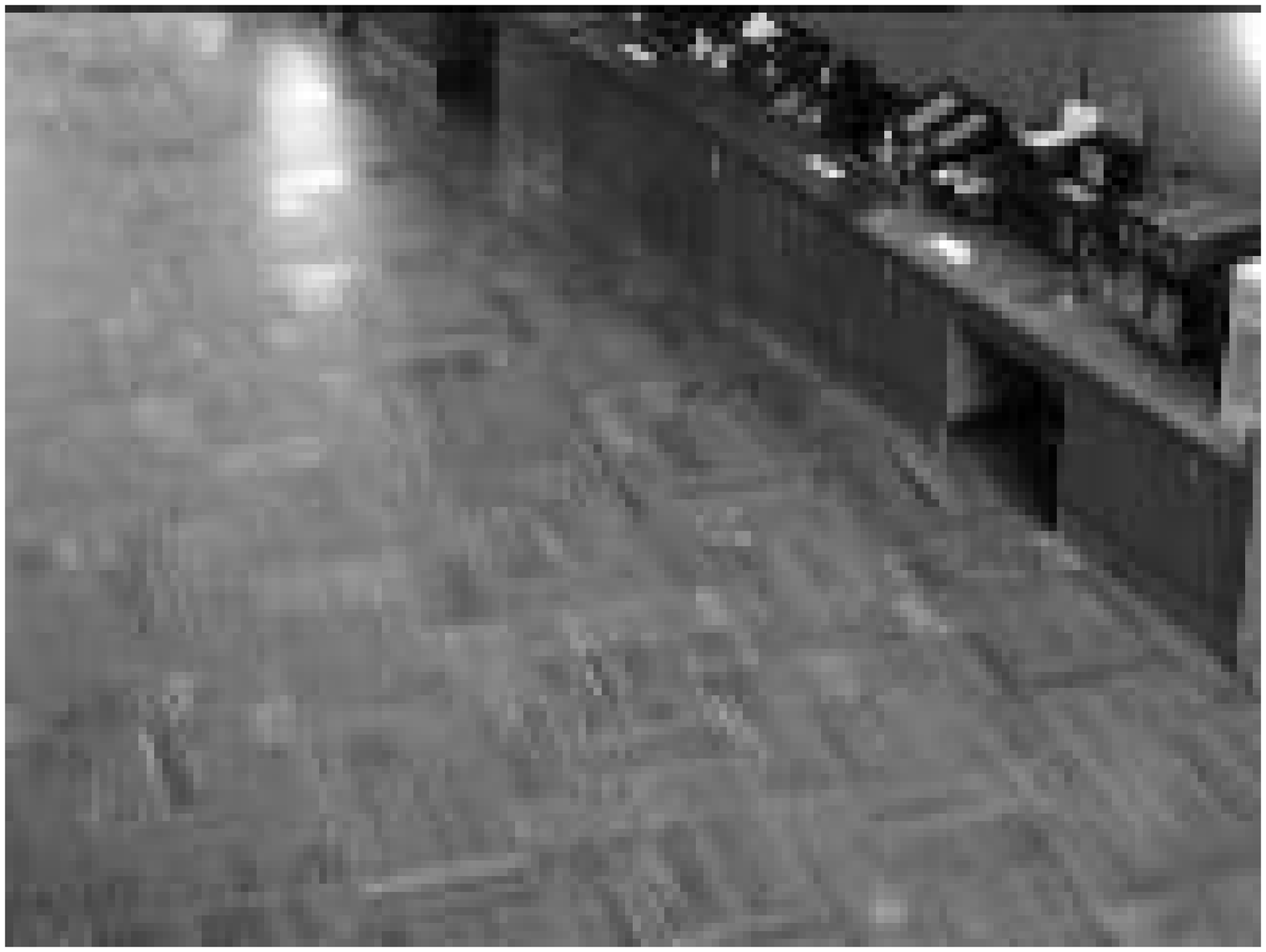}~\includegraphics[height=1.8cm,width=1.8cm]{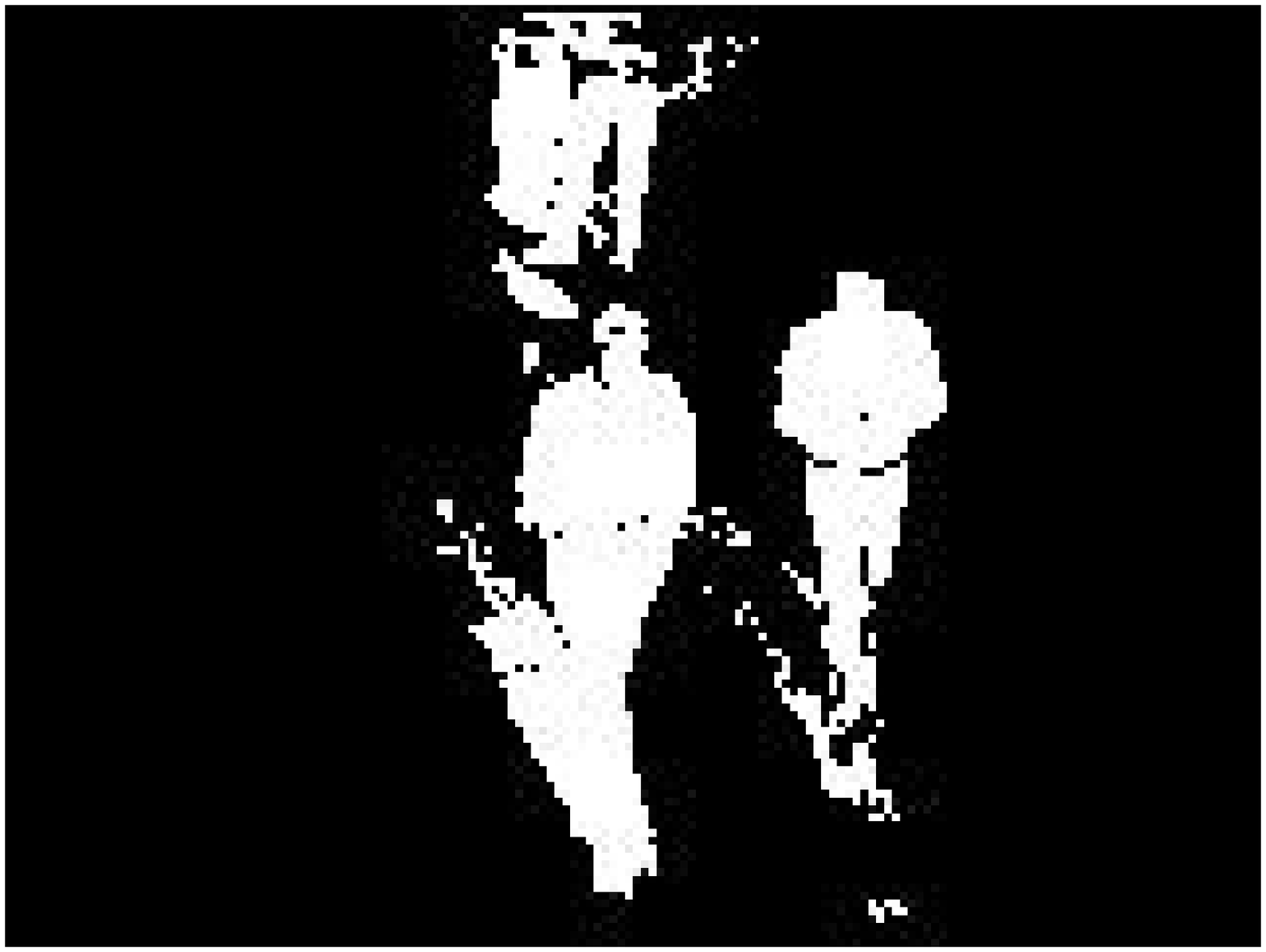}\vspace{0.5mm}\\
\includegraphics[height=1.8cm,width=1.8cm]{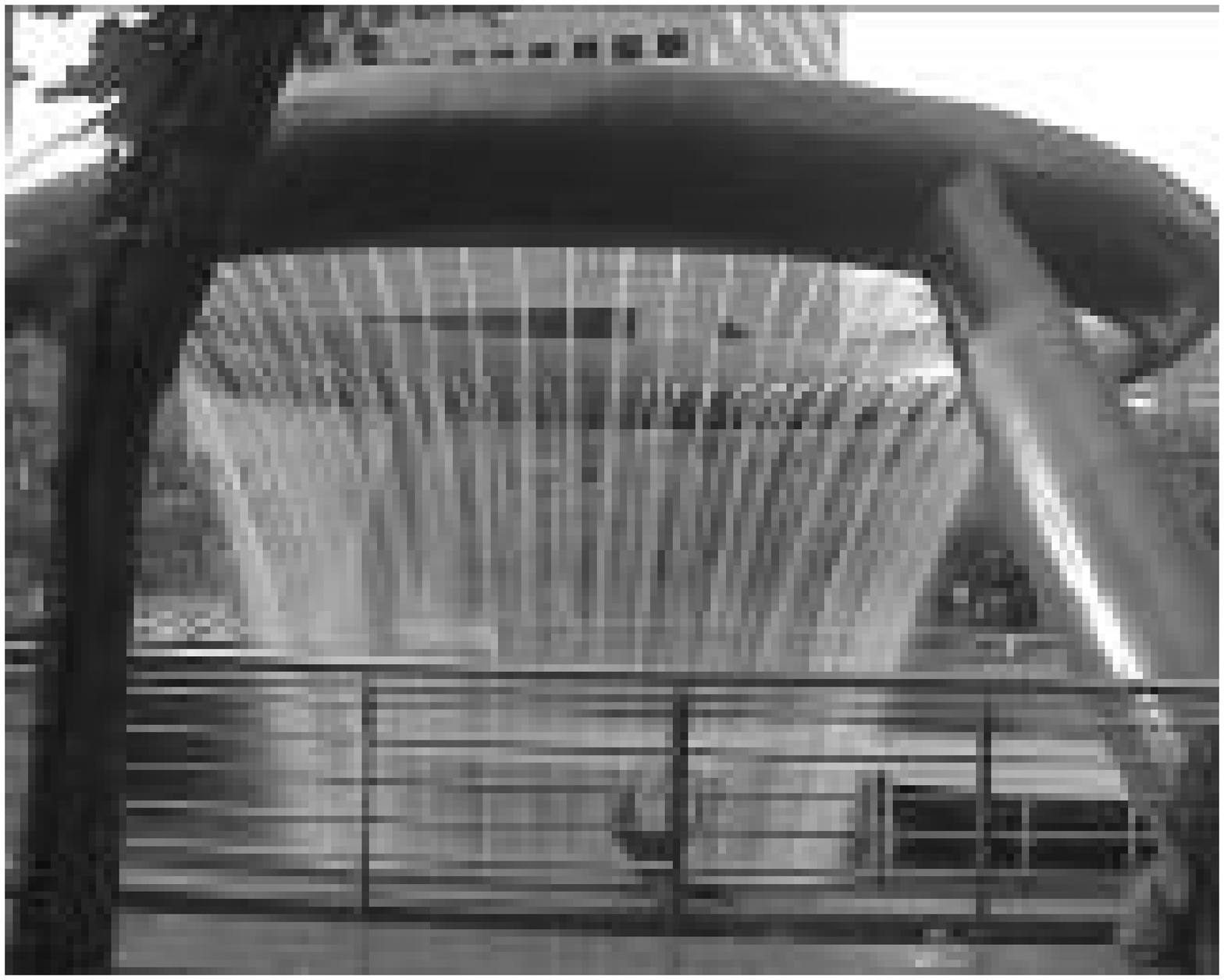}~\includegraphics[height=1.8cm,width=1.8cm]{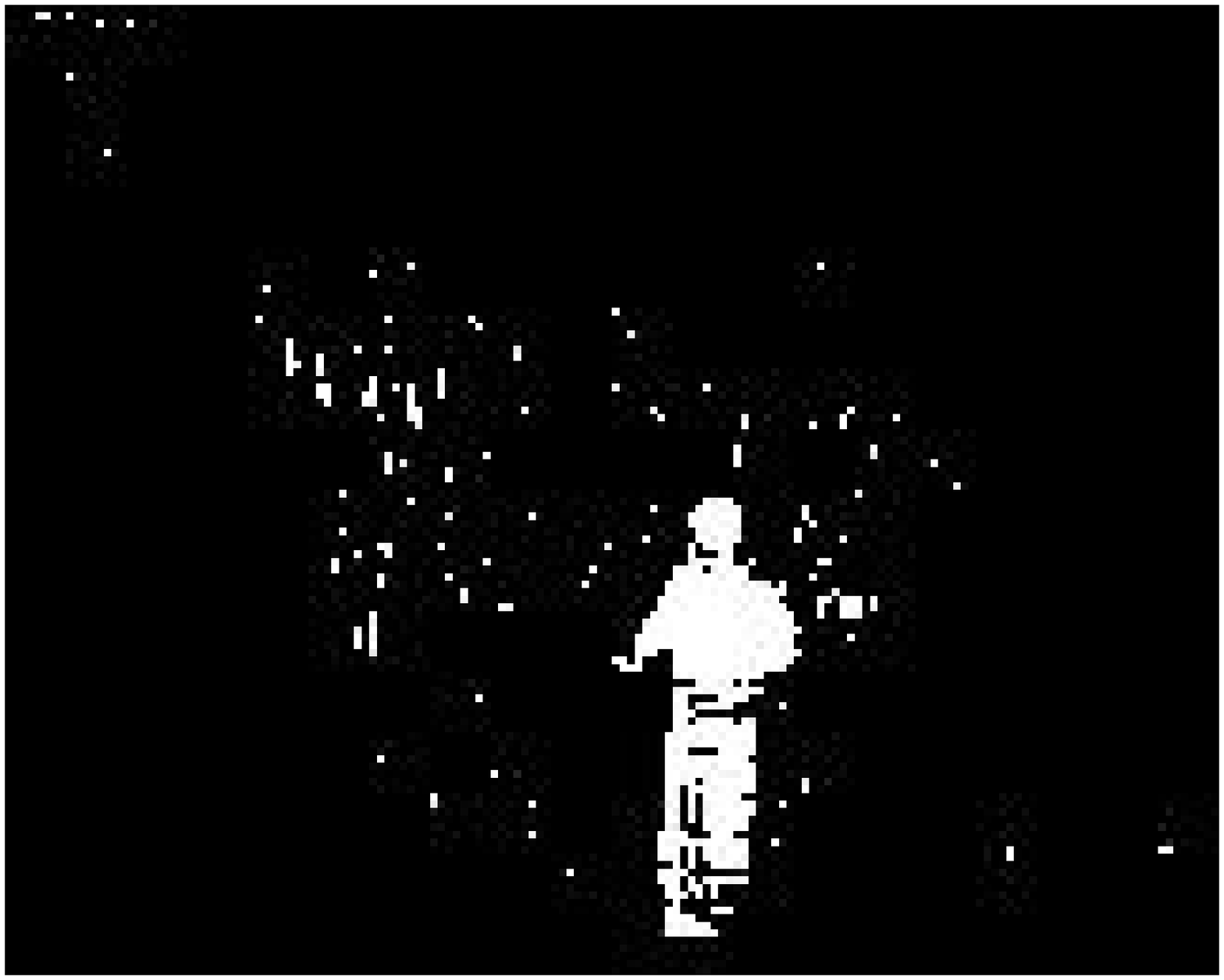}\vspace{0.5mm}\\
\includegraphics[height=1.8cm,width=1.8cm]{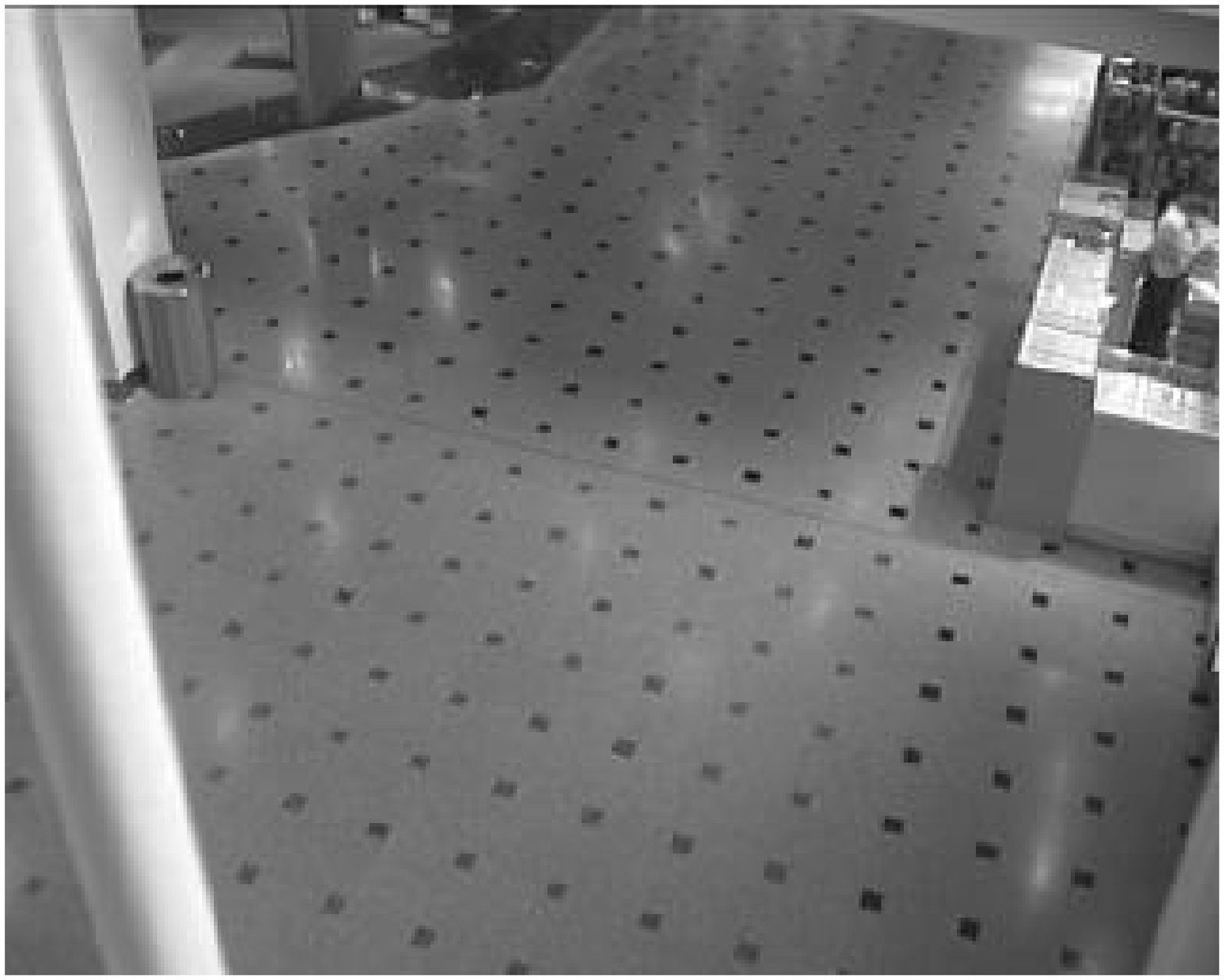}~\includegraphics[height=1.8cm,width=1.8cm]{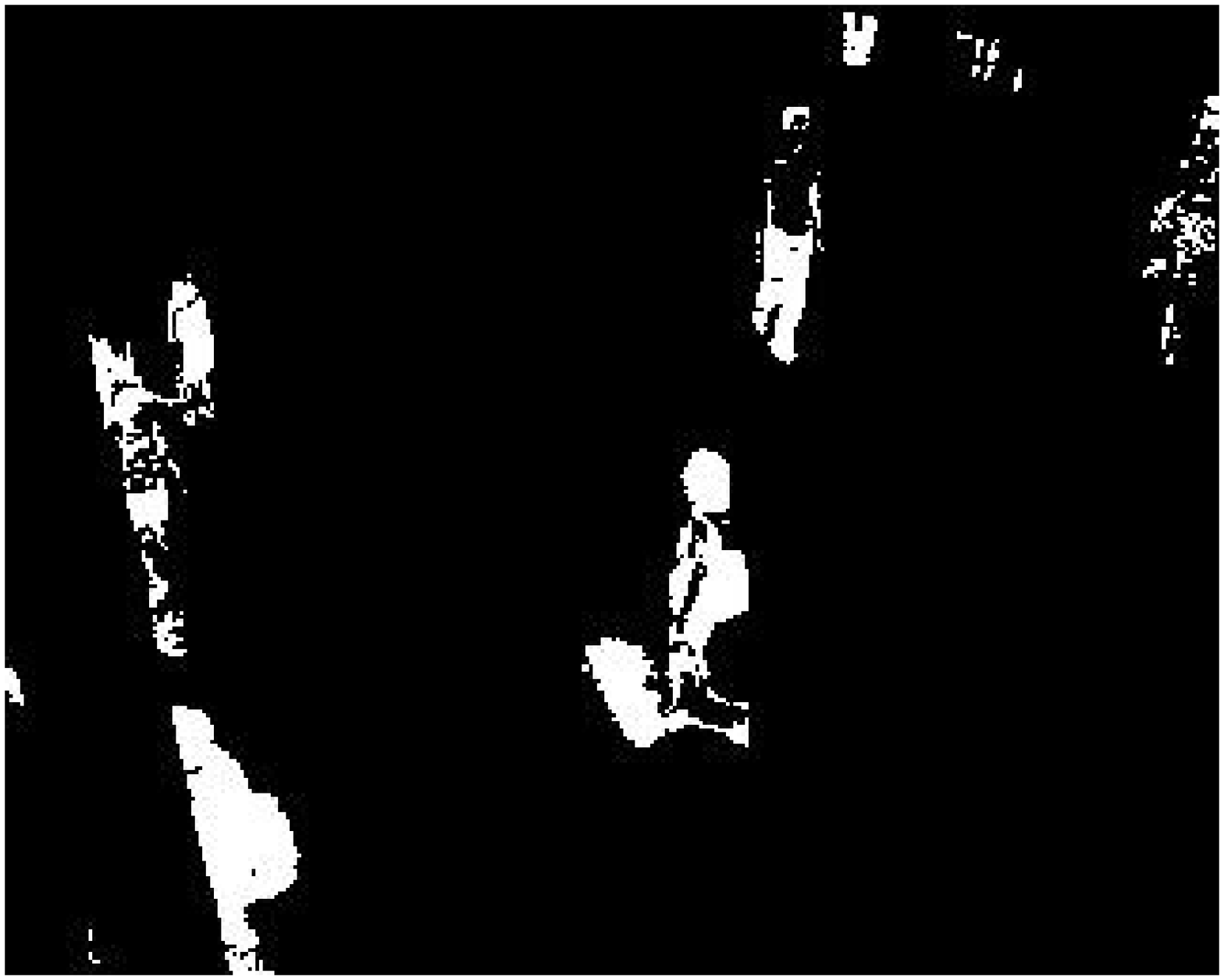}
\end{minipage}}
\subfigure{\begin{minipage}[t]{0.3\textwidth}\centering bri.$p=0.5$\vspace{1mm} \\
\includegraphics[height=1.8cm,width=1.8cm]{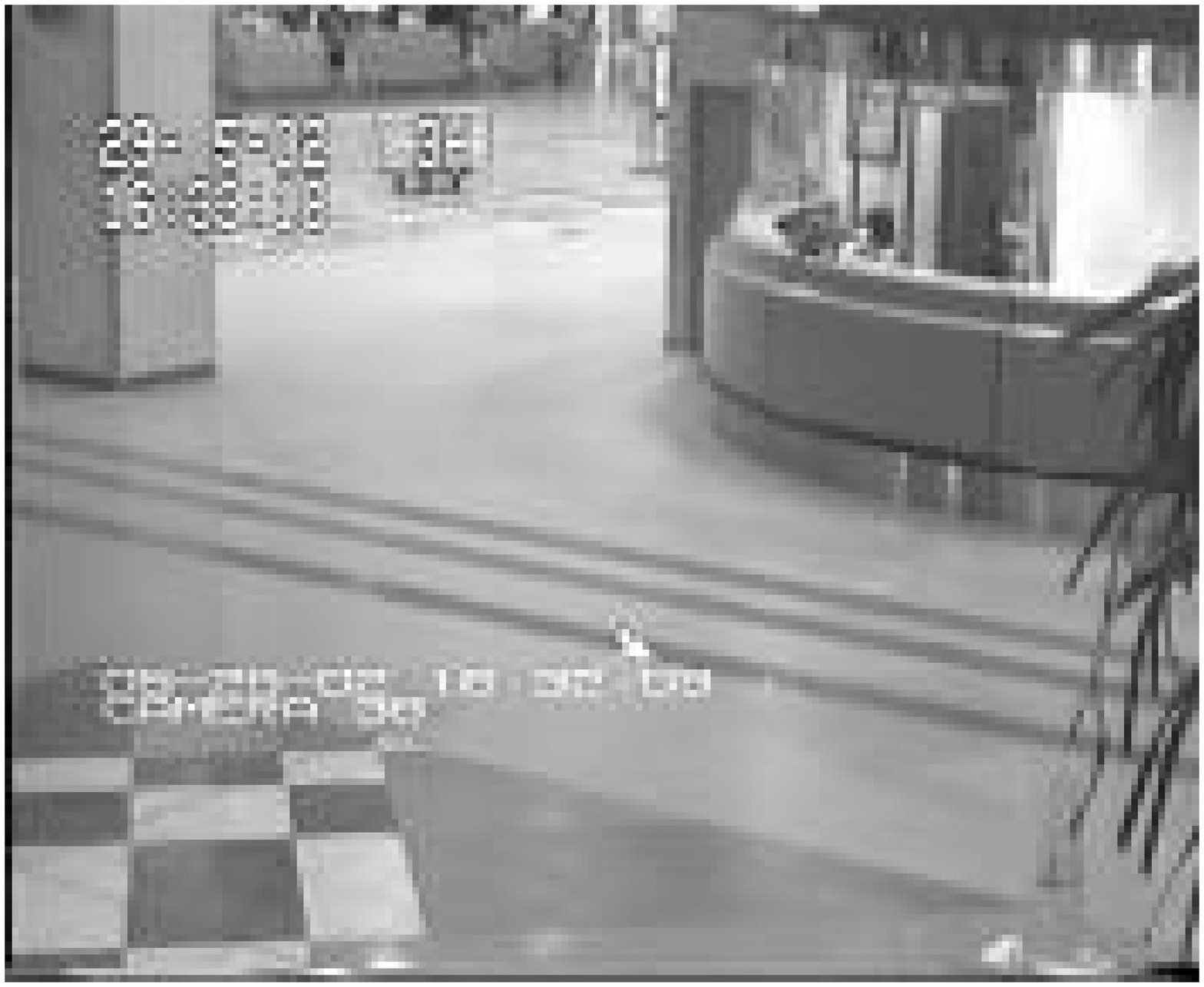}~\includegraphics[height=1.8cm,width=1.8cm]{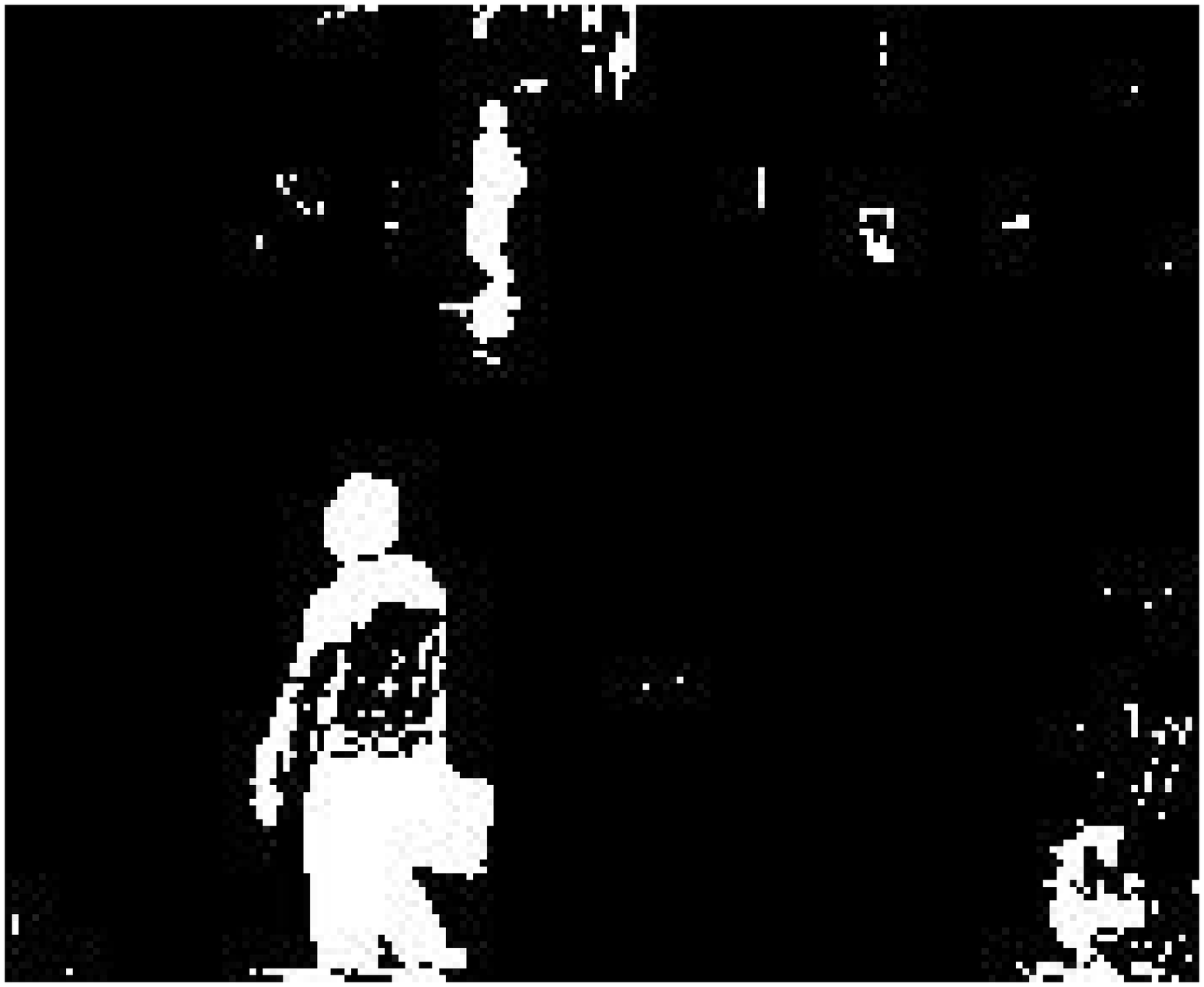}\vspace{0.5mm}\\
\includegraphics[height=1.8cm,width=1.8cm]{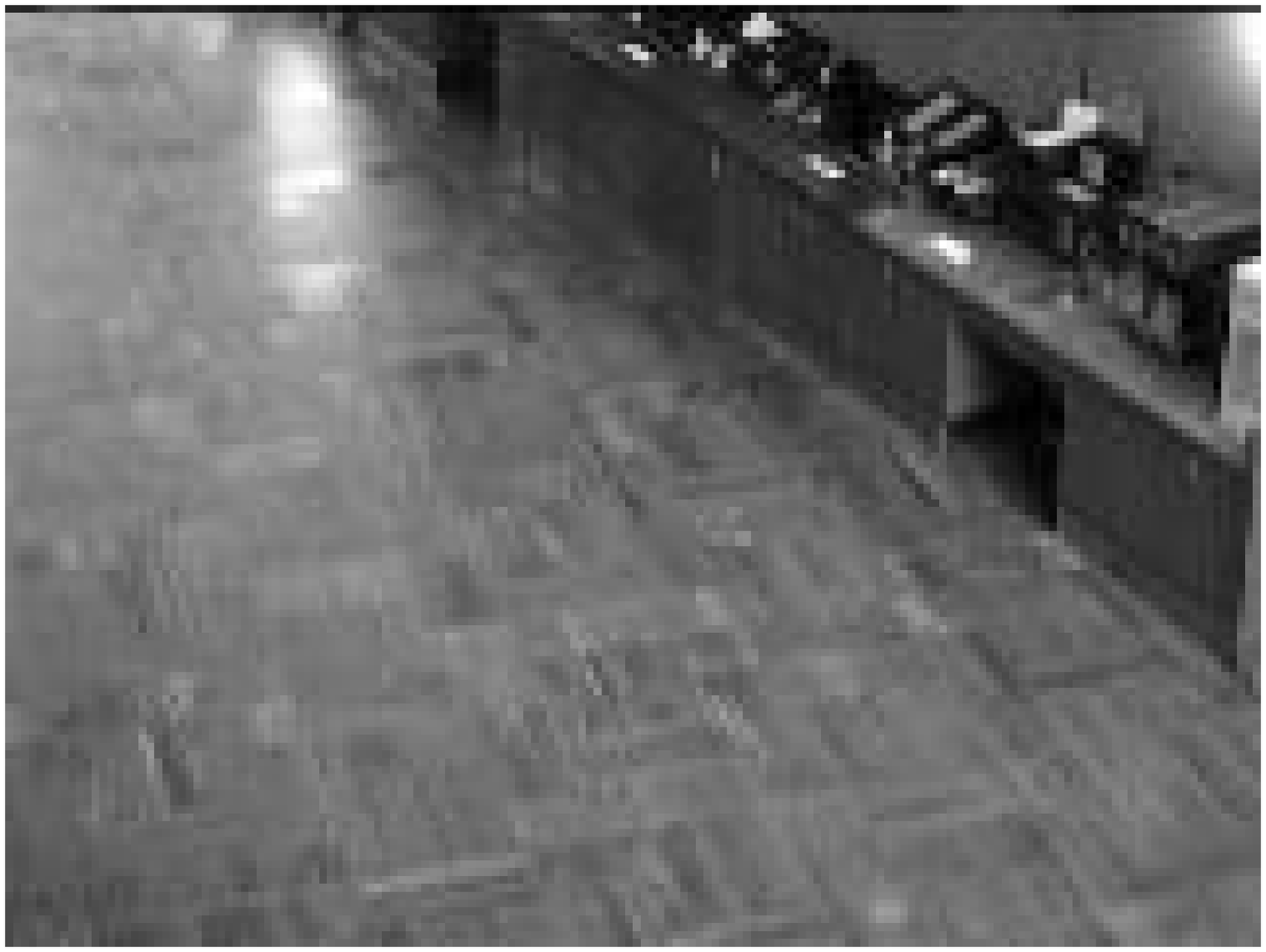}~\includegraphics[height=1.8cm,width=1.8cm]{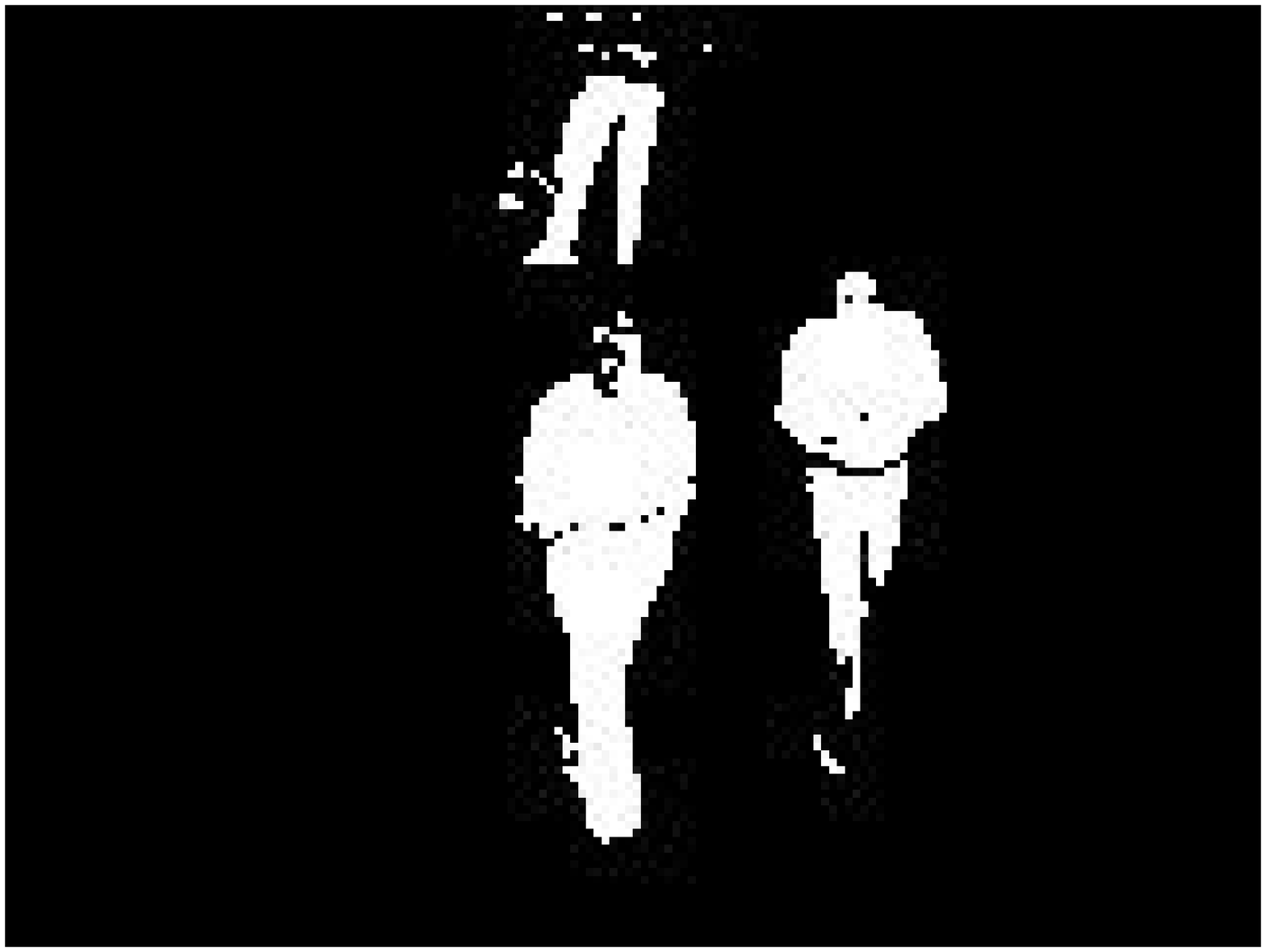}\vspace{0.5mm}\\
\includegraphics[height=1.8cm,width=1.8cm]{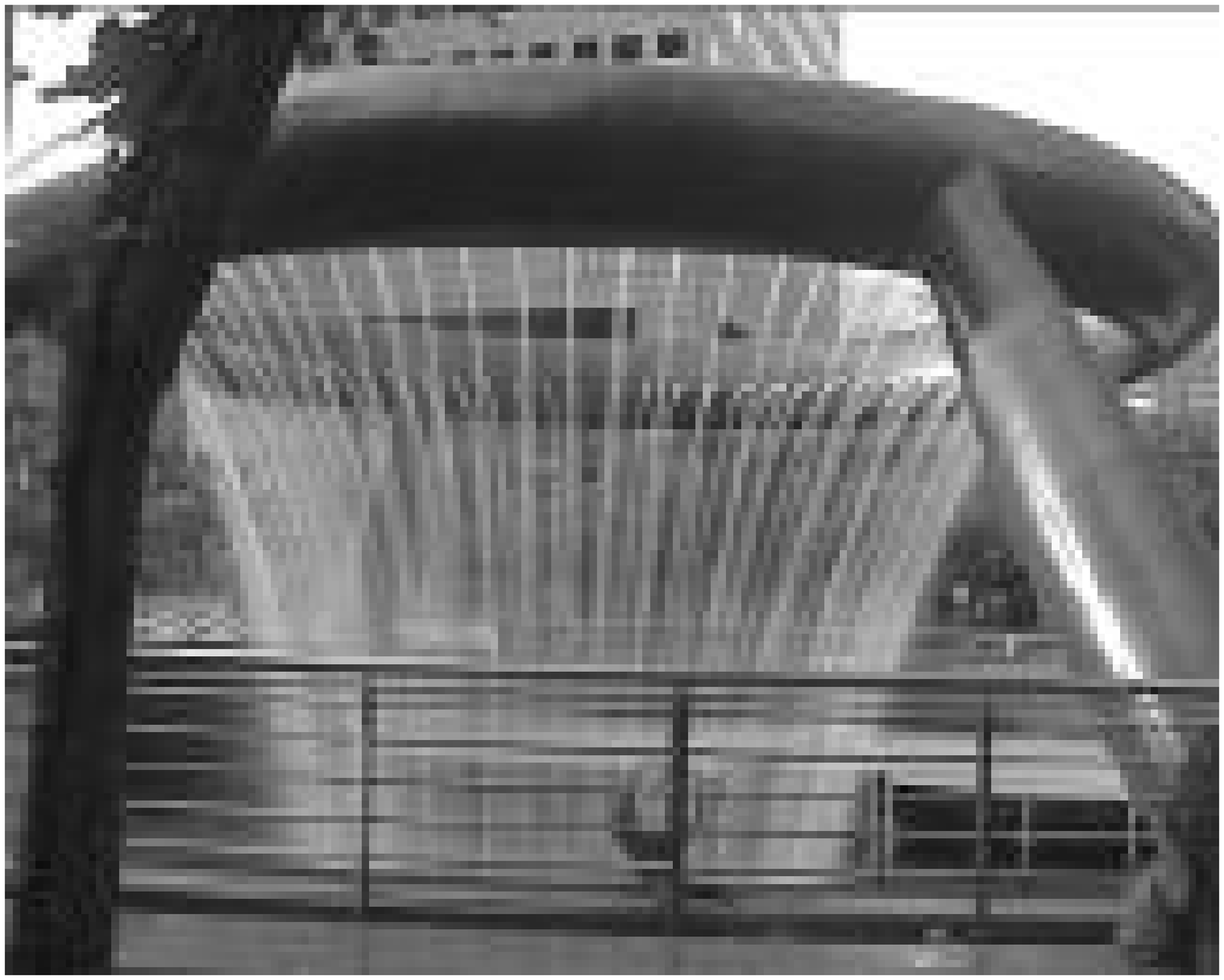}~\includegraphics[height=1.8cm,width=1.8cm]{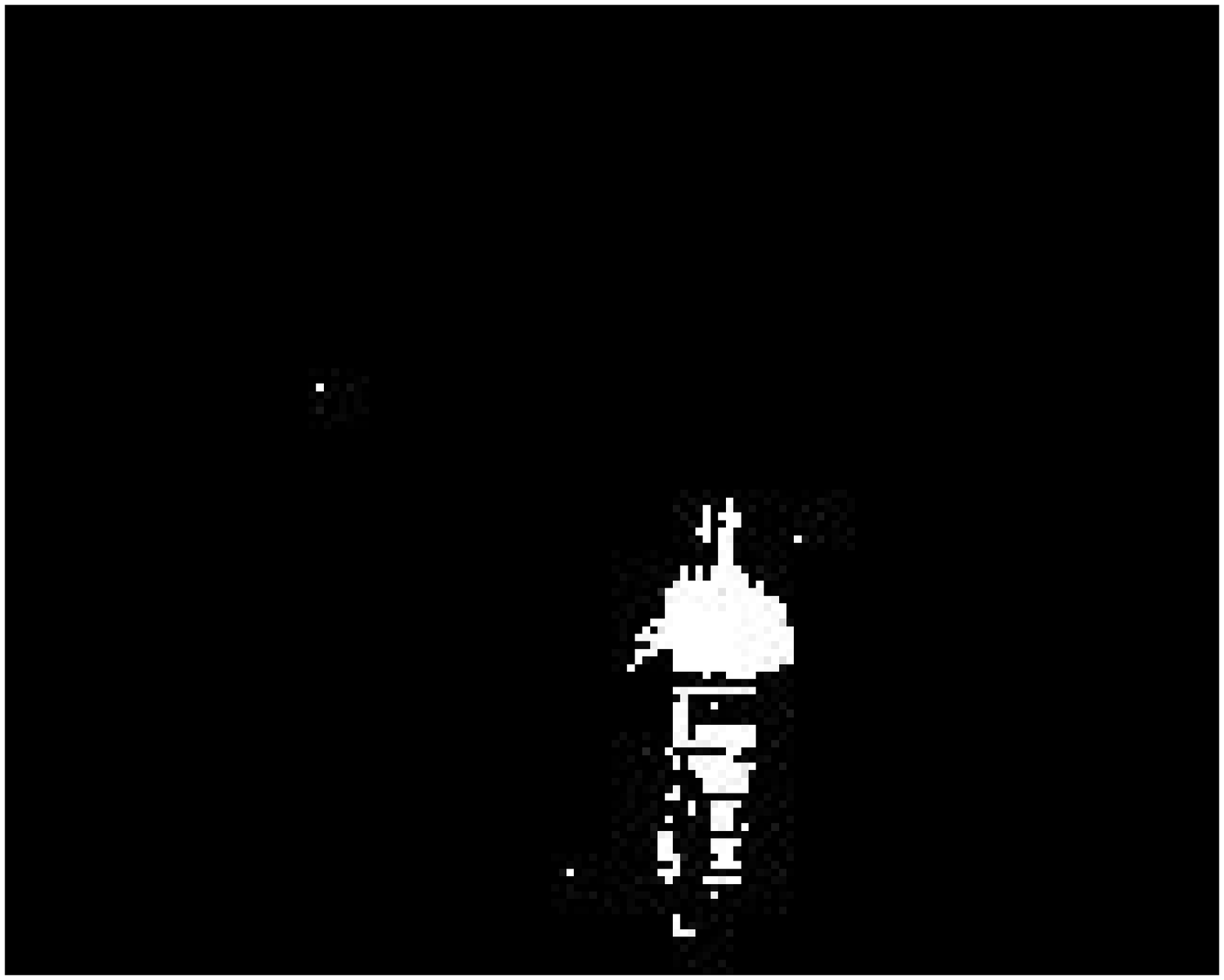}\vspace{0.5mm}\\
\includegraphics[height=1.8cm,width=1.8cm]{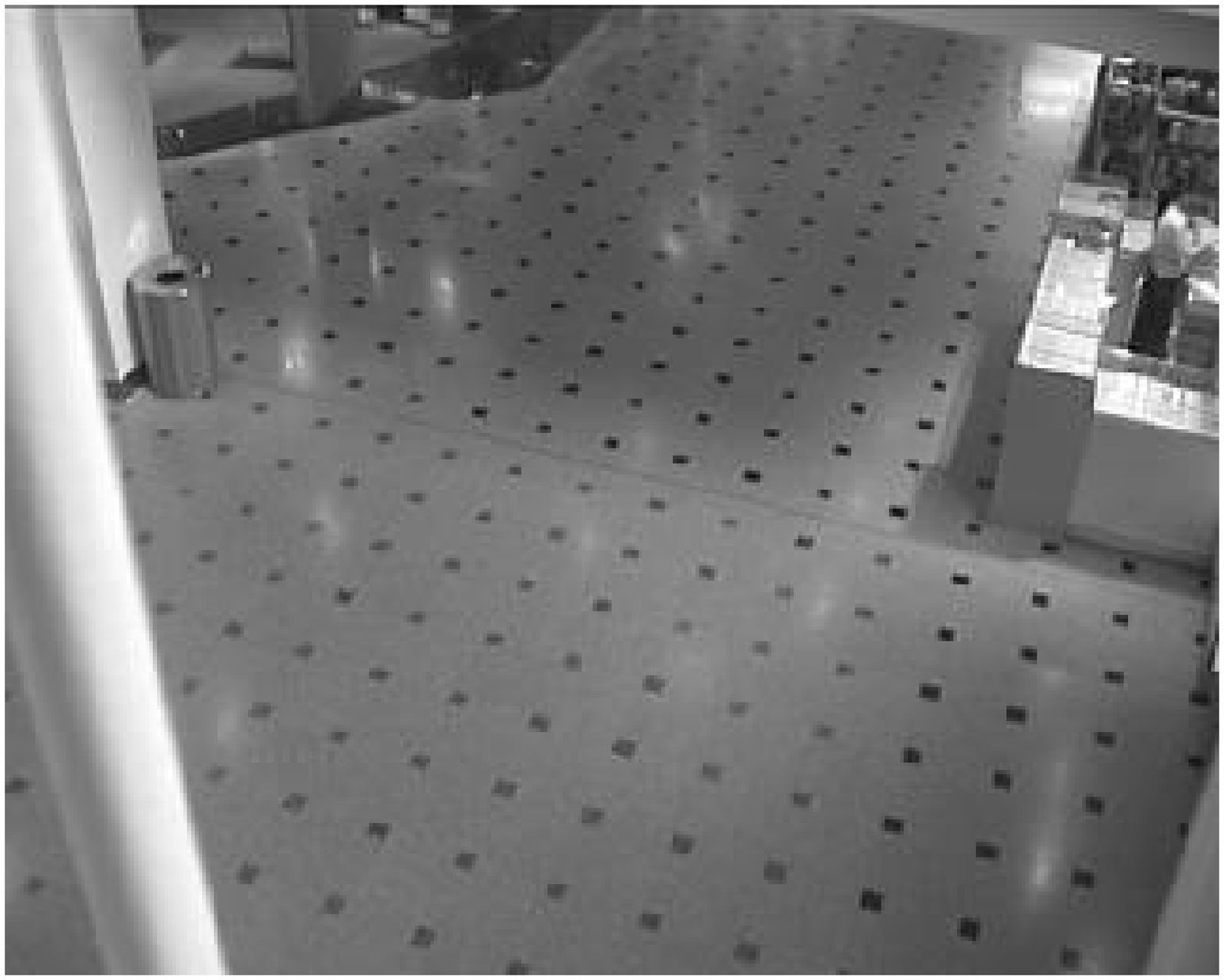}~\includegraphics[height=1.8cm,width=1.8cm]{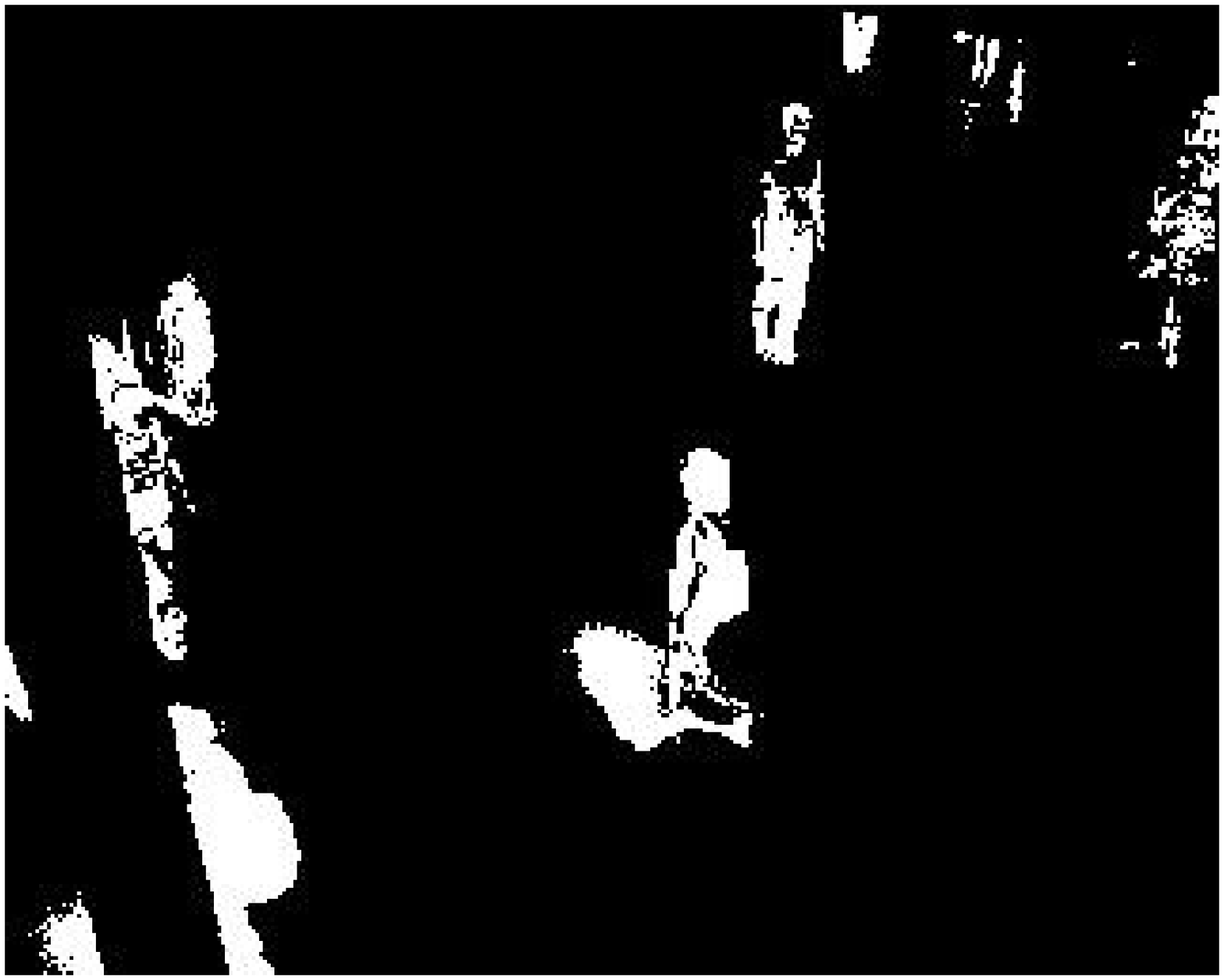}
\end{minipage}}
\subfigure{\begin{minipage}[t]{0.3\textwidth}\centering fra.$\alpha=1$\vspace{1.6mm} \\
\includegraphics[height=1.8cm,width=1.8cm]{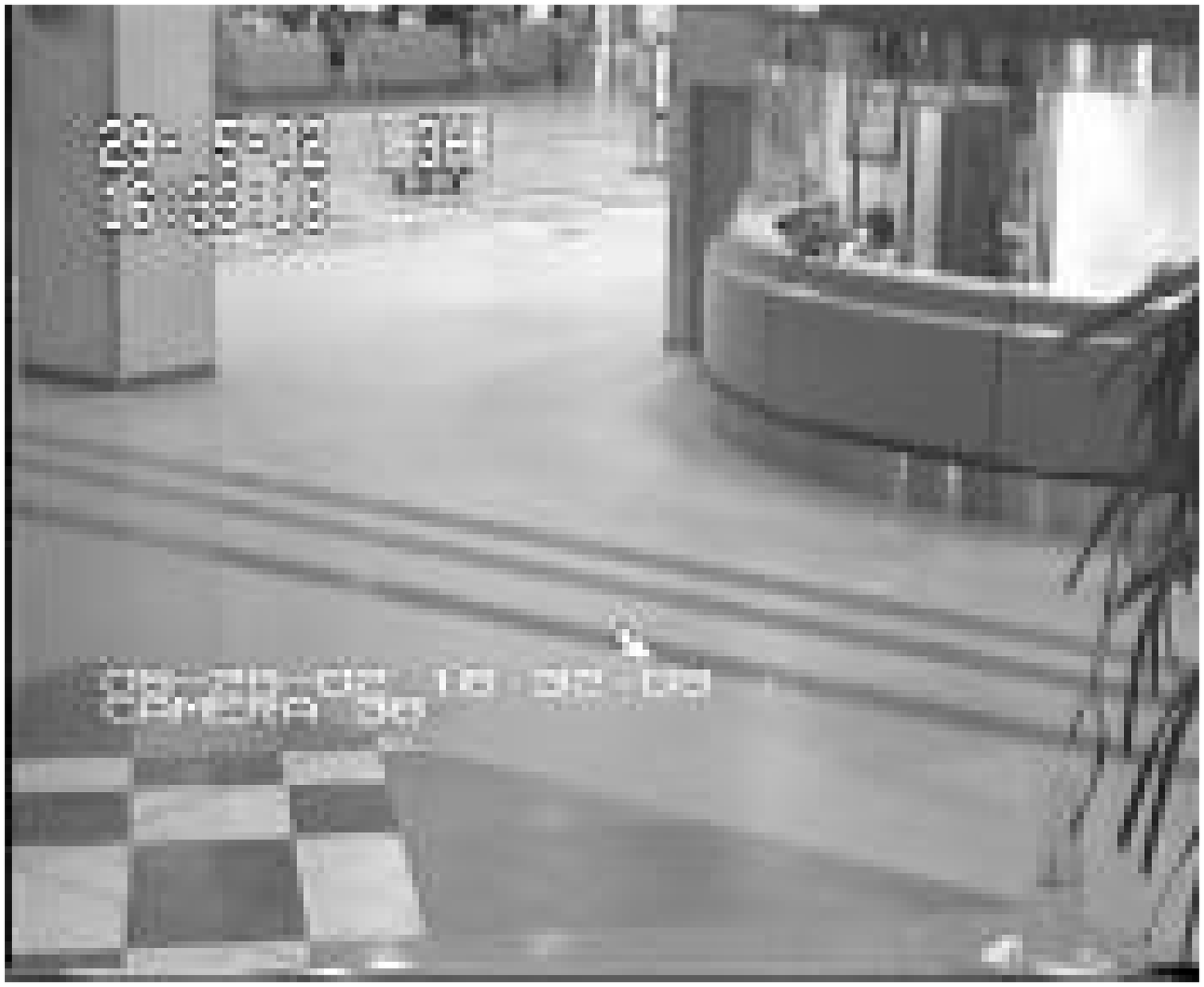}~\includegraphics[height=1.8cm,width=1.8cm]{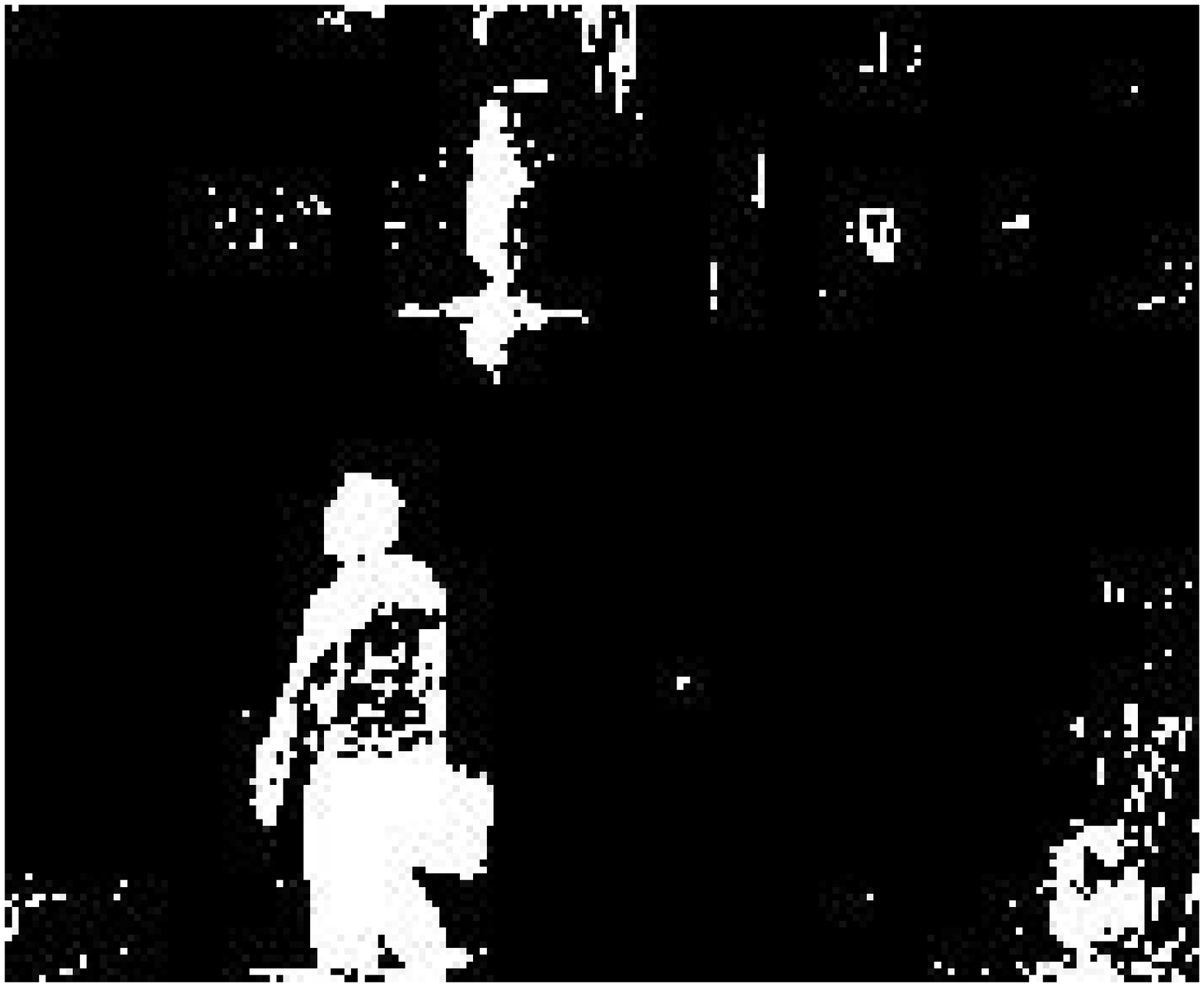}\vspace{0.5mm}\\
\includegraphics[height=1.8cm,width=1.8cm]{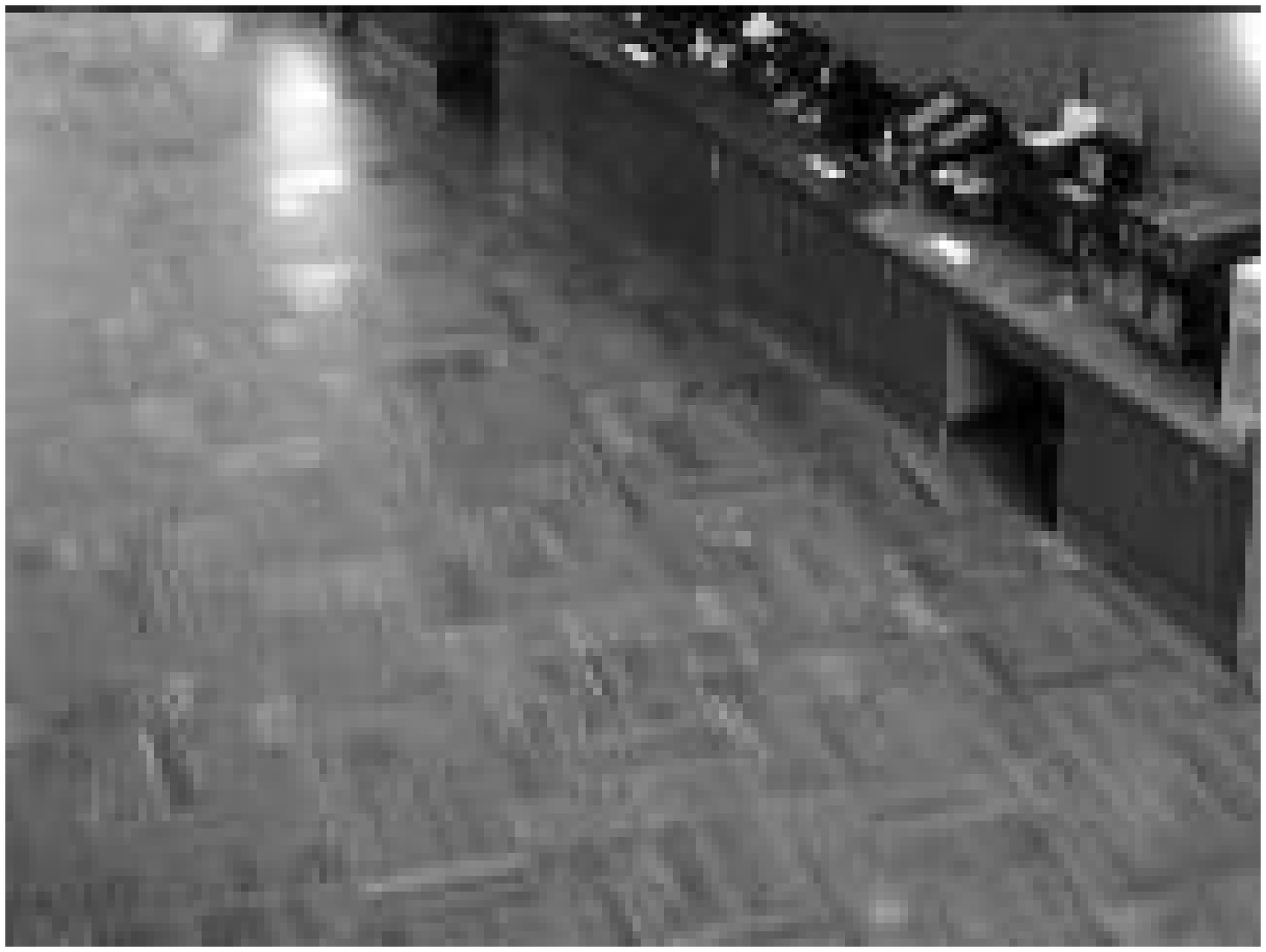}~\includegraphics[height=1.8cm,width=1.8cm]{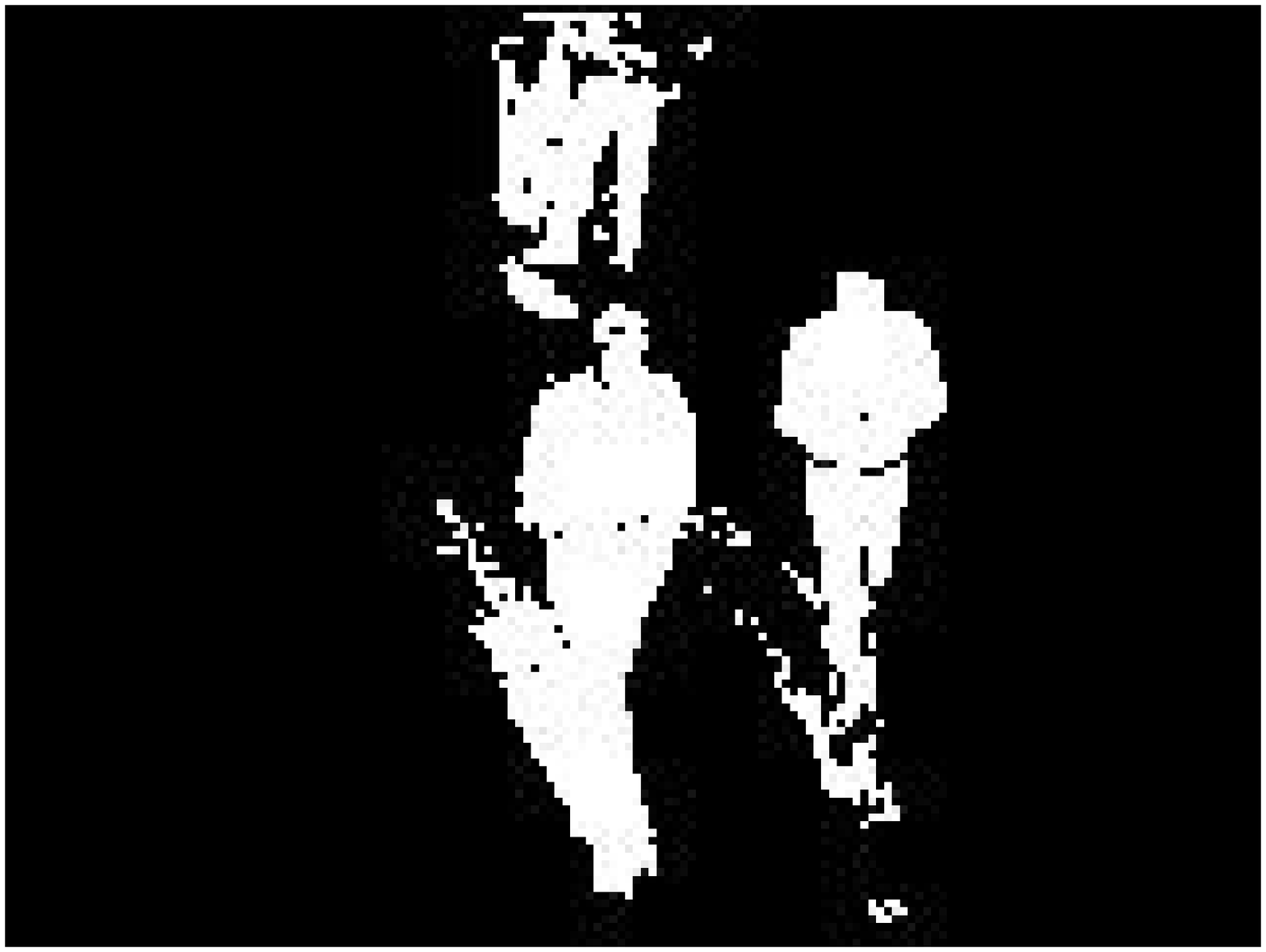}\vspace{0.5mm}\\
\includegraphics[height=1.8cm,width=1.8cm]{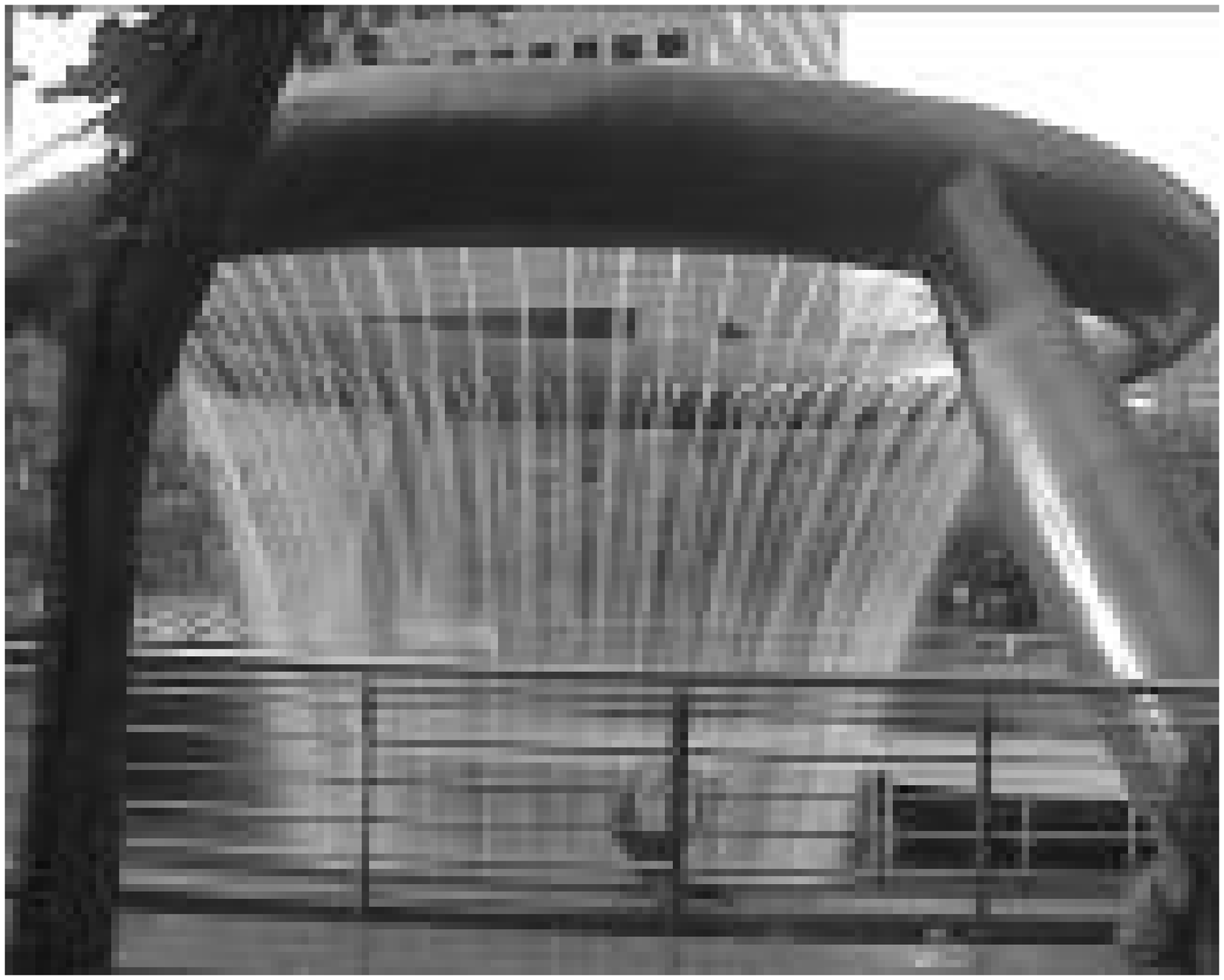}~\includegraphics[height=1.8cm,width=1.8cm]{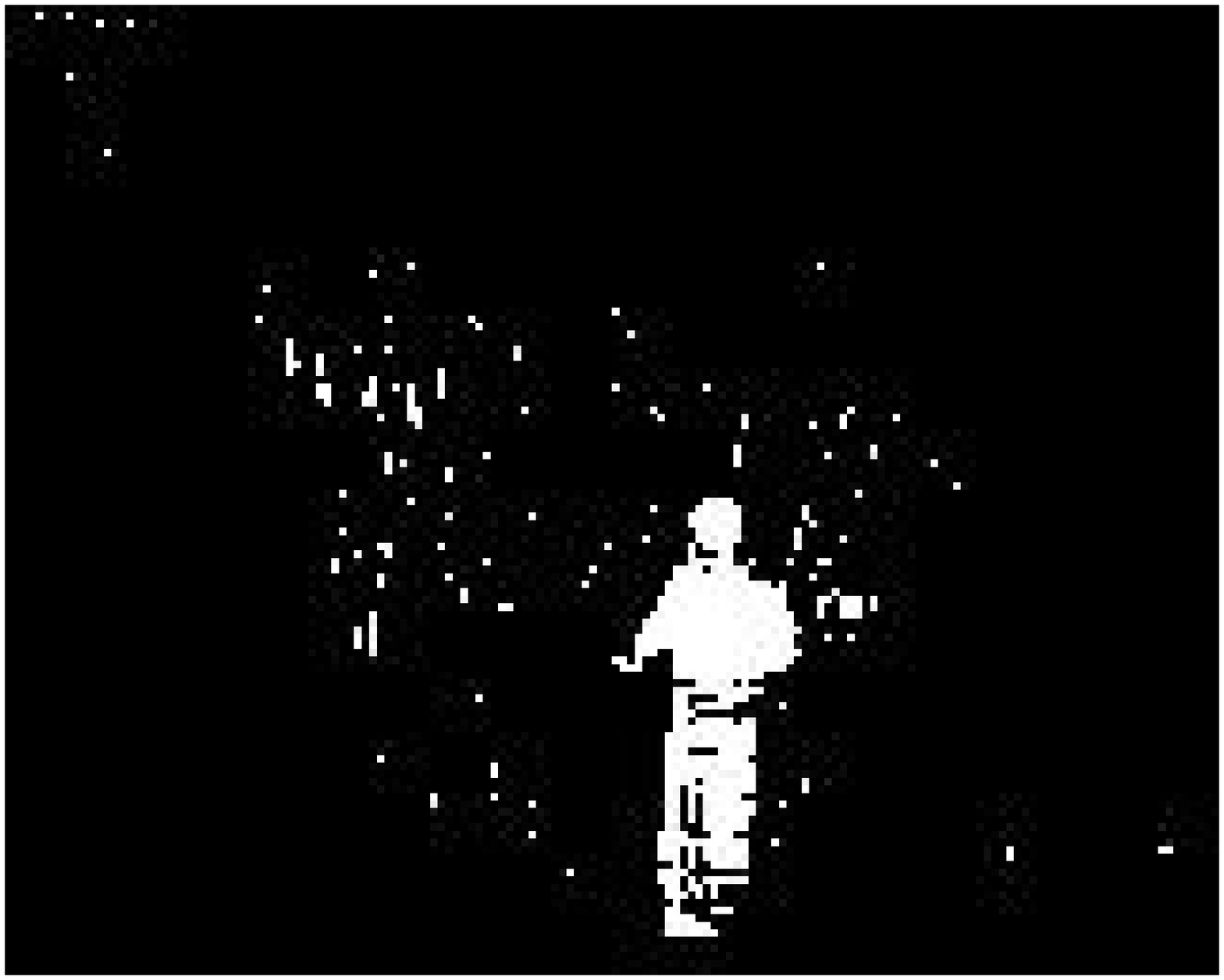}\vspace{0.5mm}\\
\includegraphics[height=1.8cm,width=1.8cm]{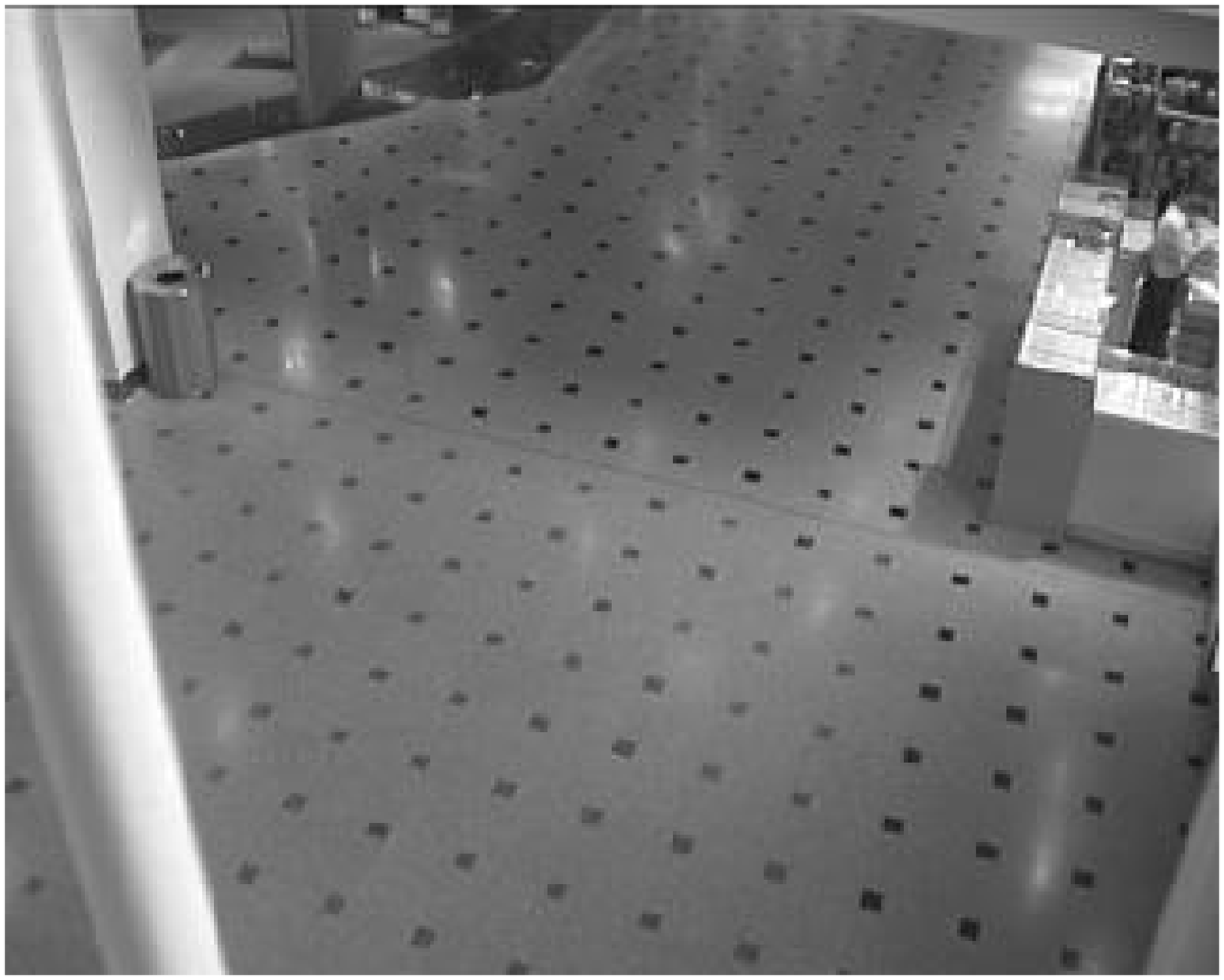}~\includegraphics[height=1.8cm,width=1.8cm]{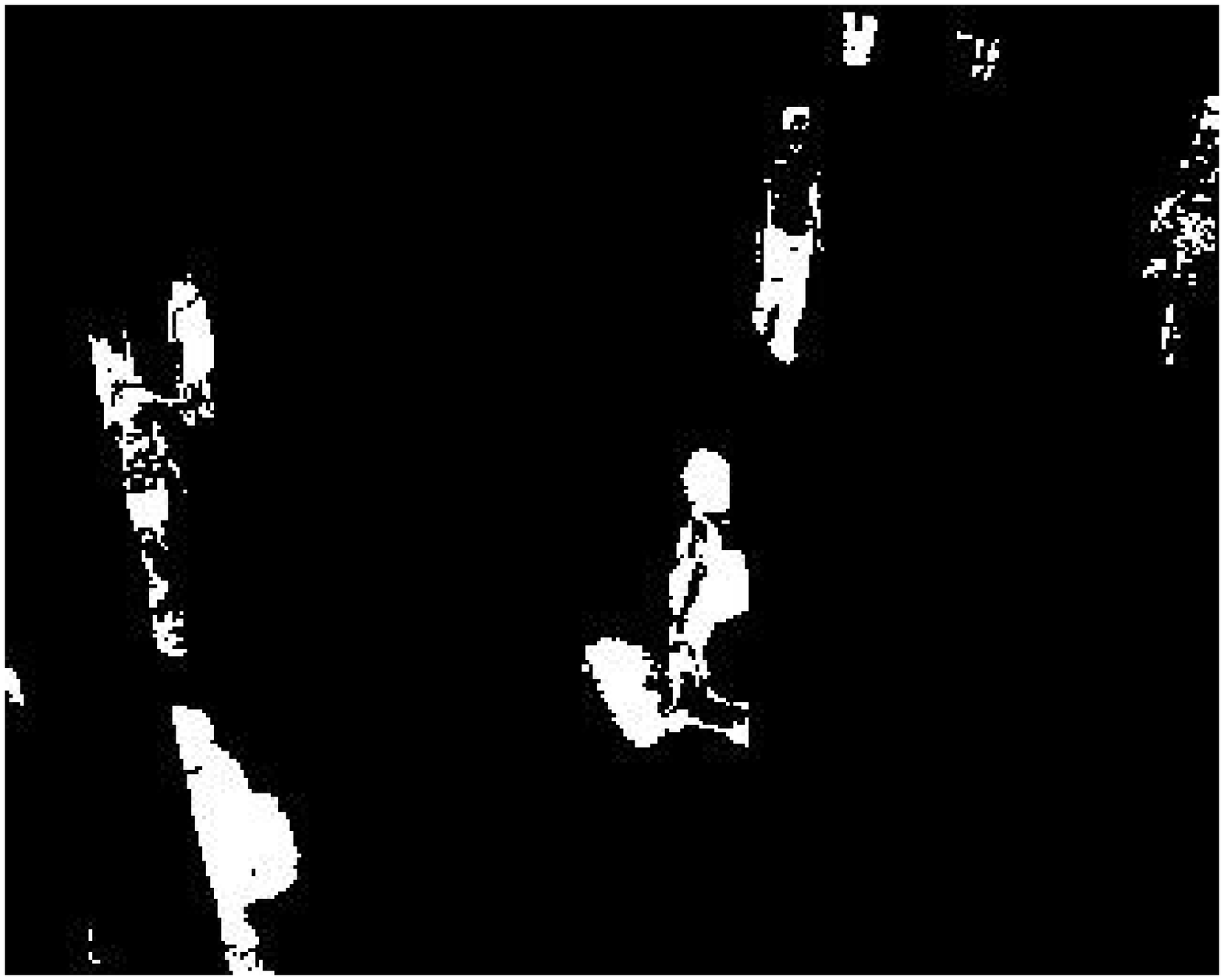}
\end{minipage}}

\subfigure{\begin{minipage}[t]{0.3\textwidth}\centering fra.$\alpha=2$\vspace{1.6mm} \\
\includegraphics[height=1.8cm,width=1.8cm]{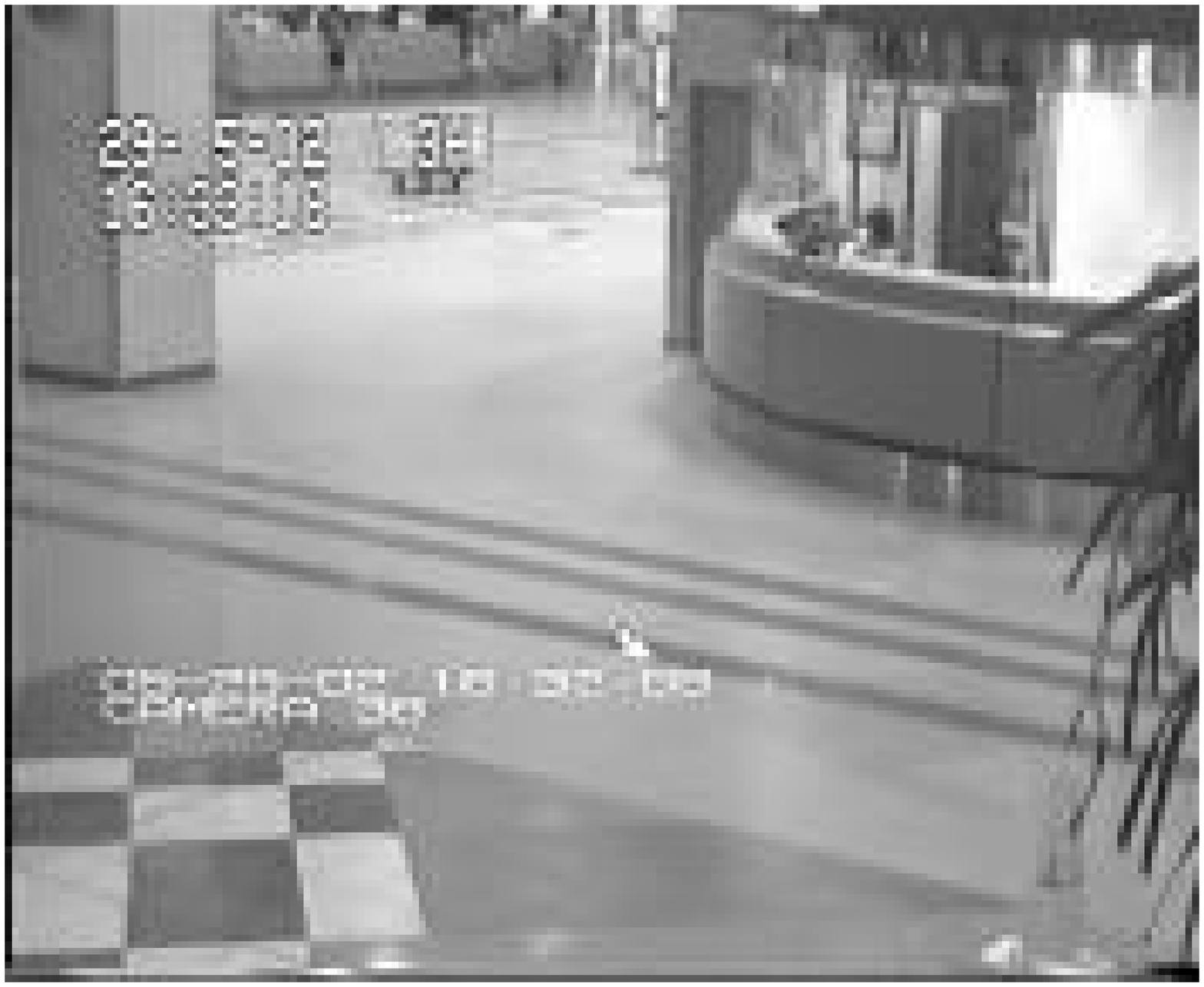}~\includegraphics[height=1.8cm,width=1.8cm]{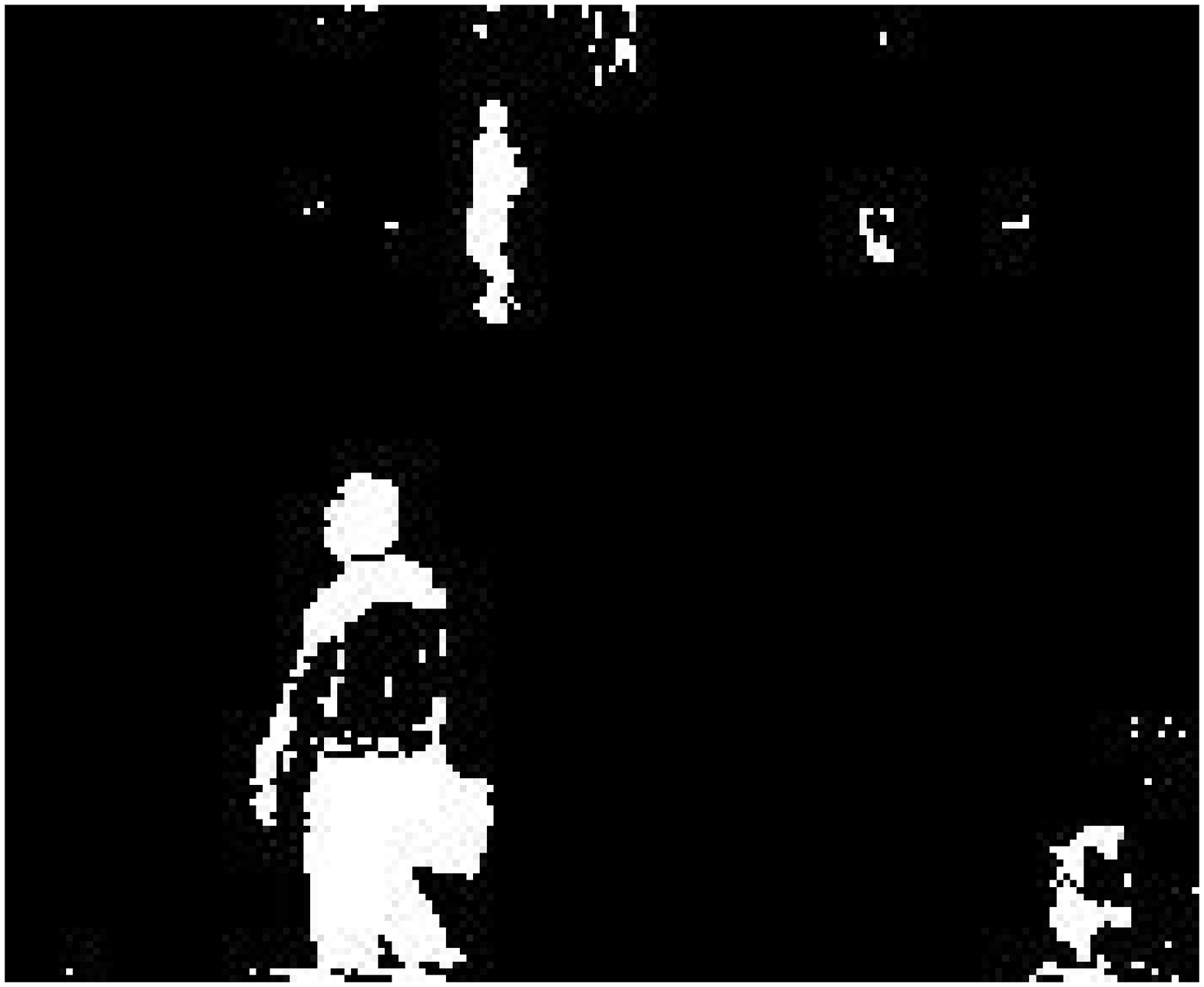}\vspace{0.5mm}\\
\includegraphics[height=1.8cm,width=1.8cm]{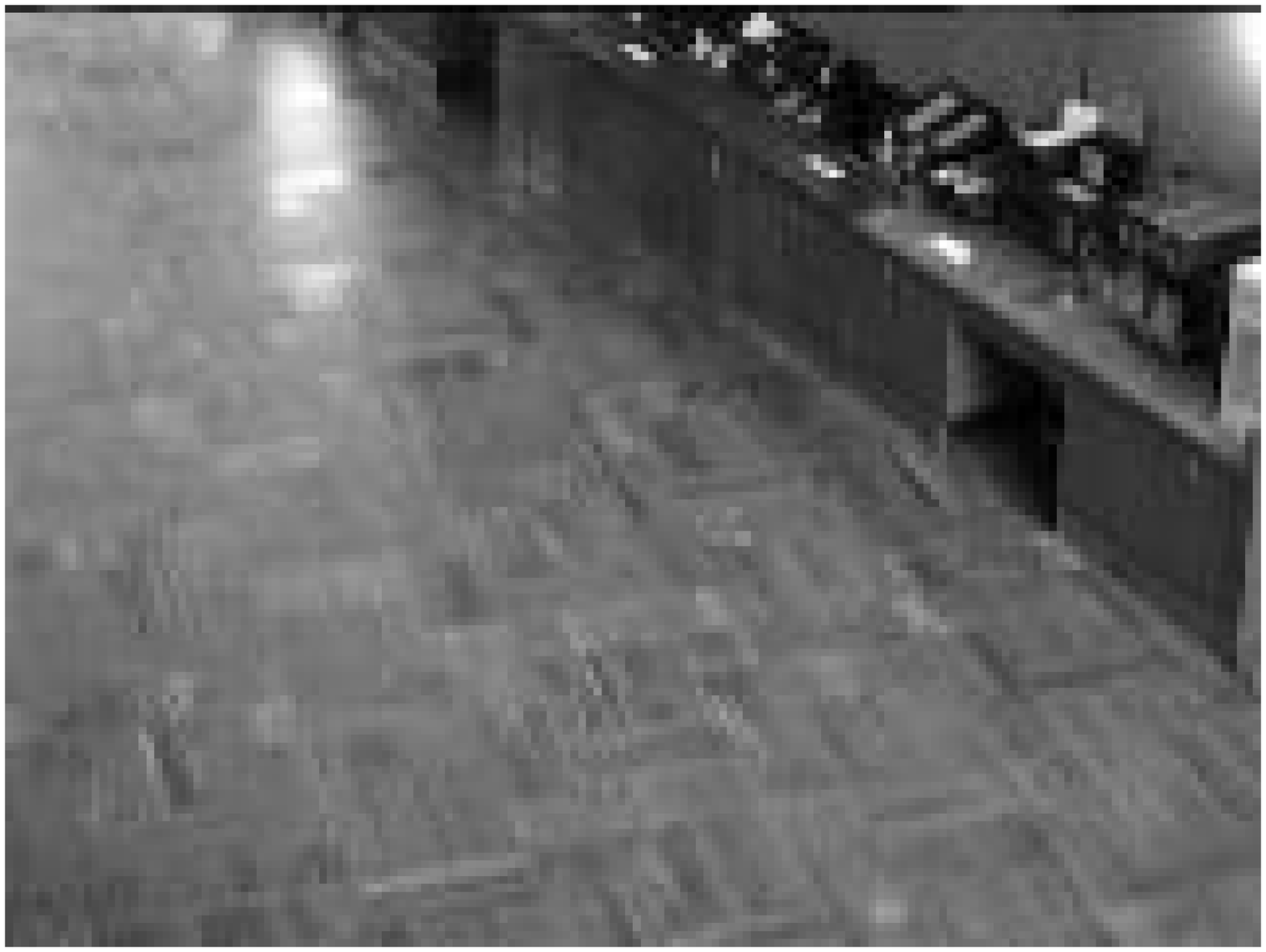}~\includegraphics[height=1.8cm,width=1.8cm]{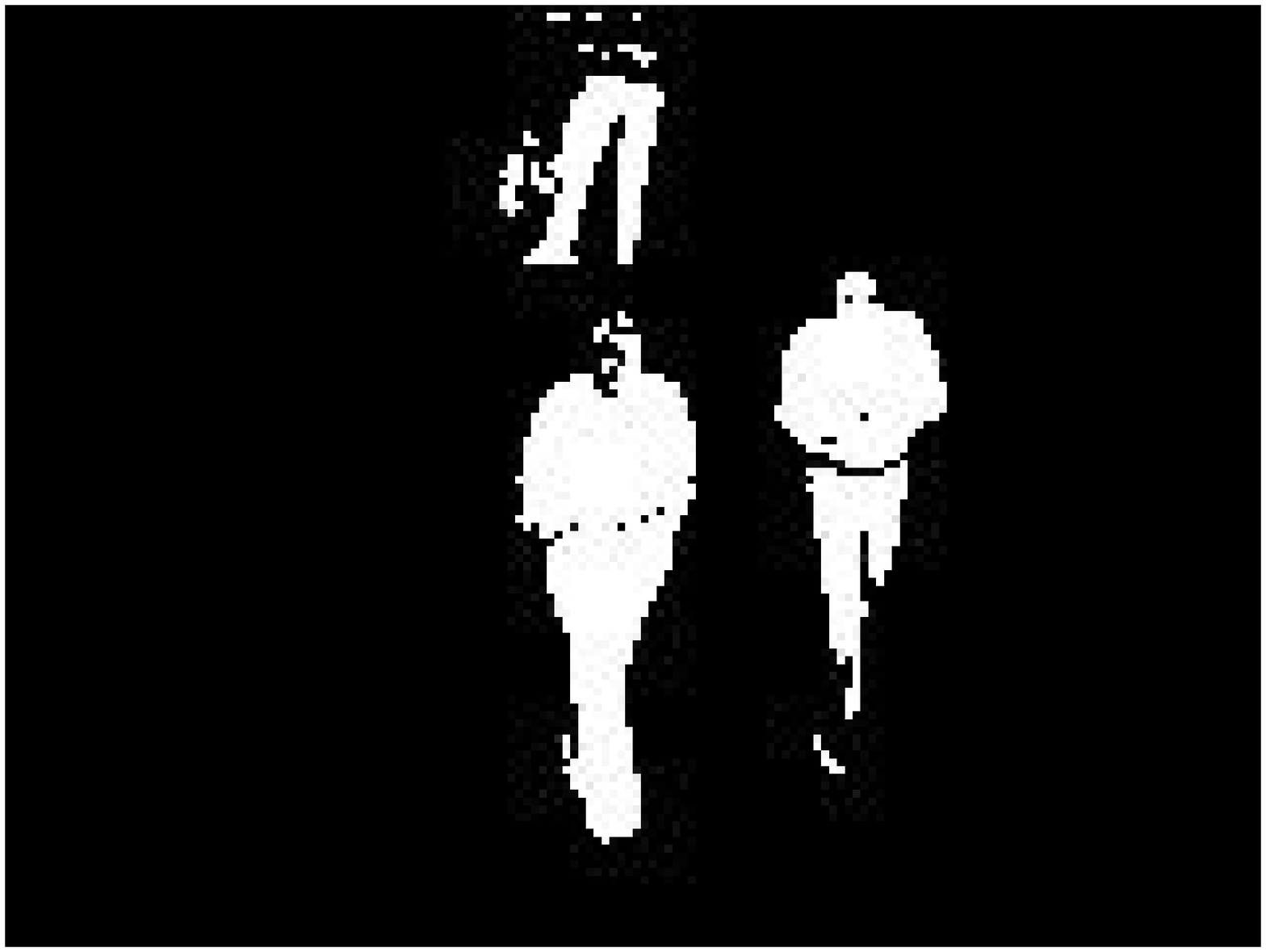}\vspace{0.5mm}\\
\includegraphics[height=1.8cm,width=1.8cm]{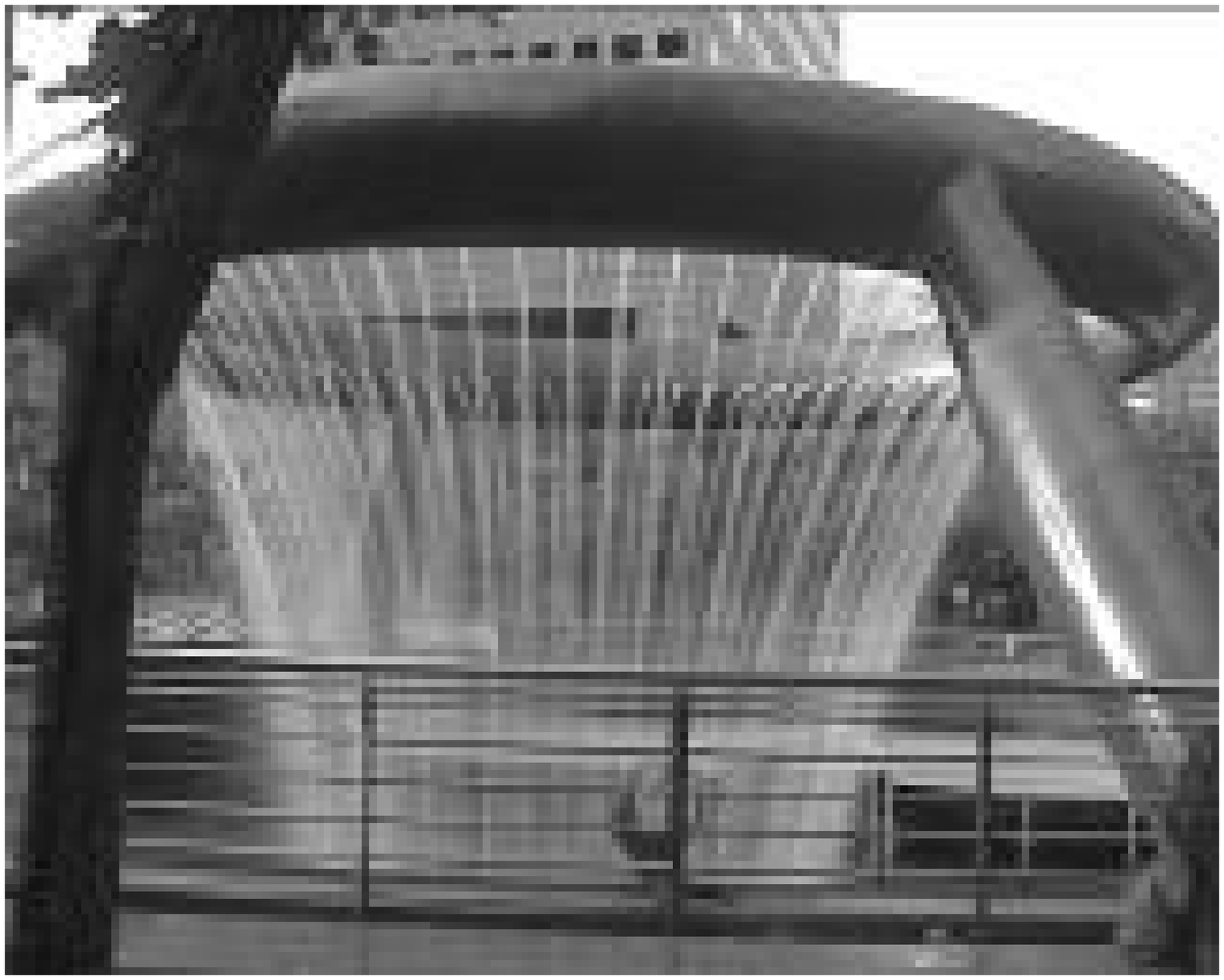}~\includegraphics[height=1.8cm,width=1.8cm]{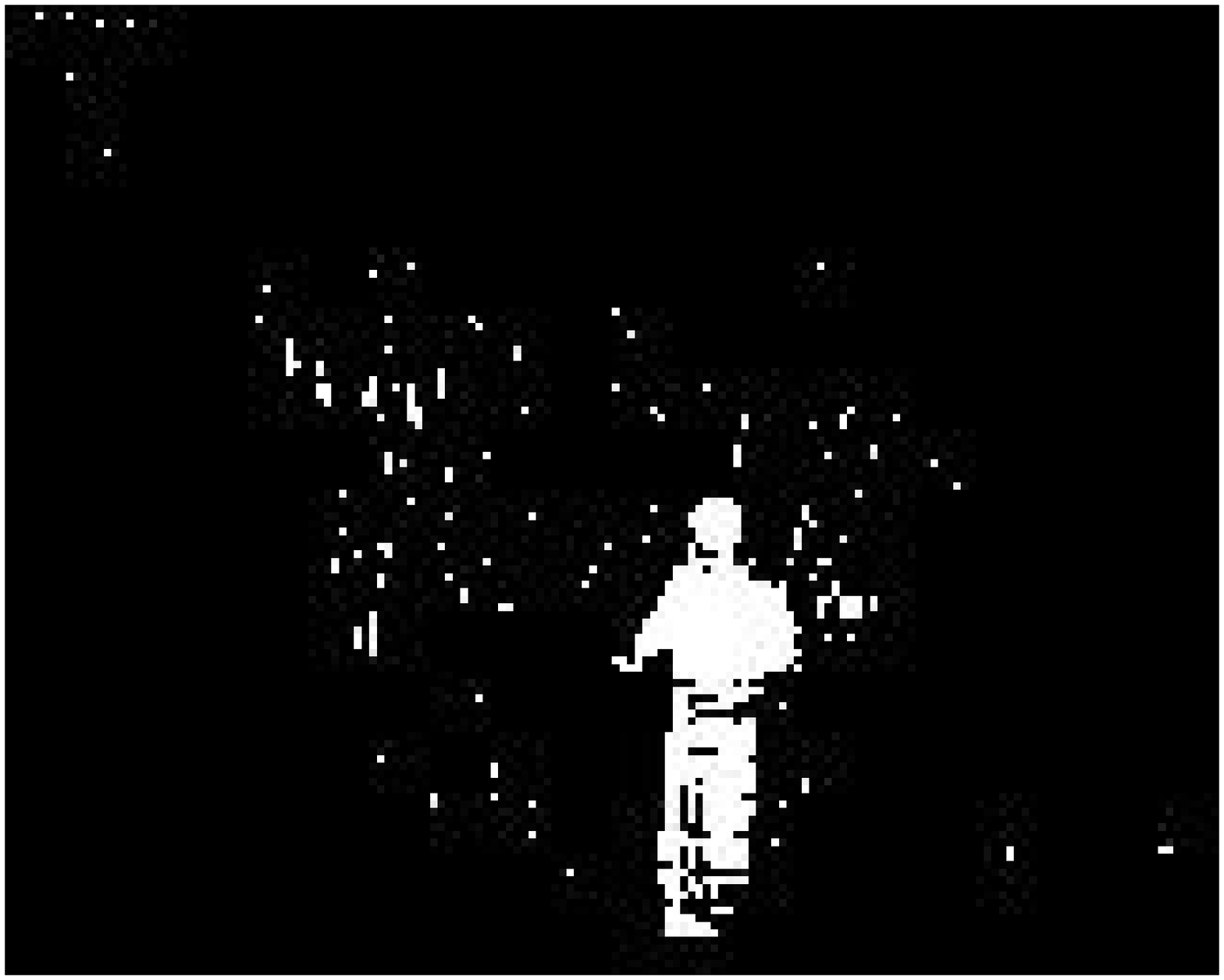}\vspace{0.5mm}\\
\includegraphics[height=1.8cm,width=1.8cm]{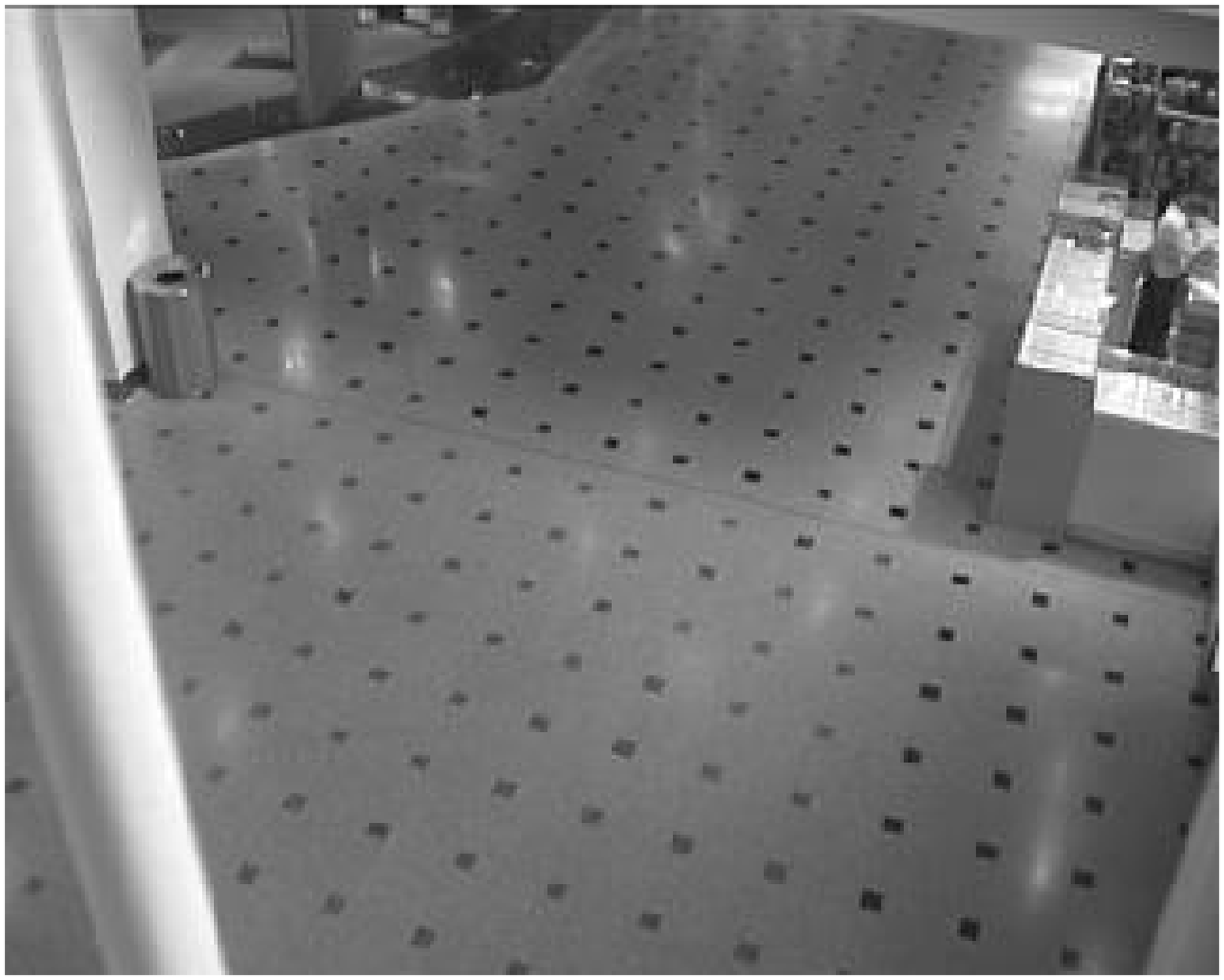}~\includegraphics[height=1.8cm,width=1.8cm]{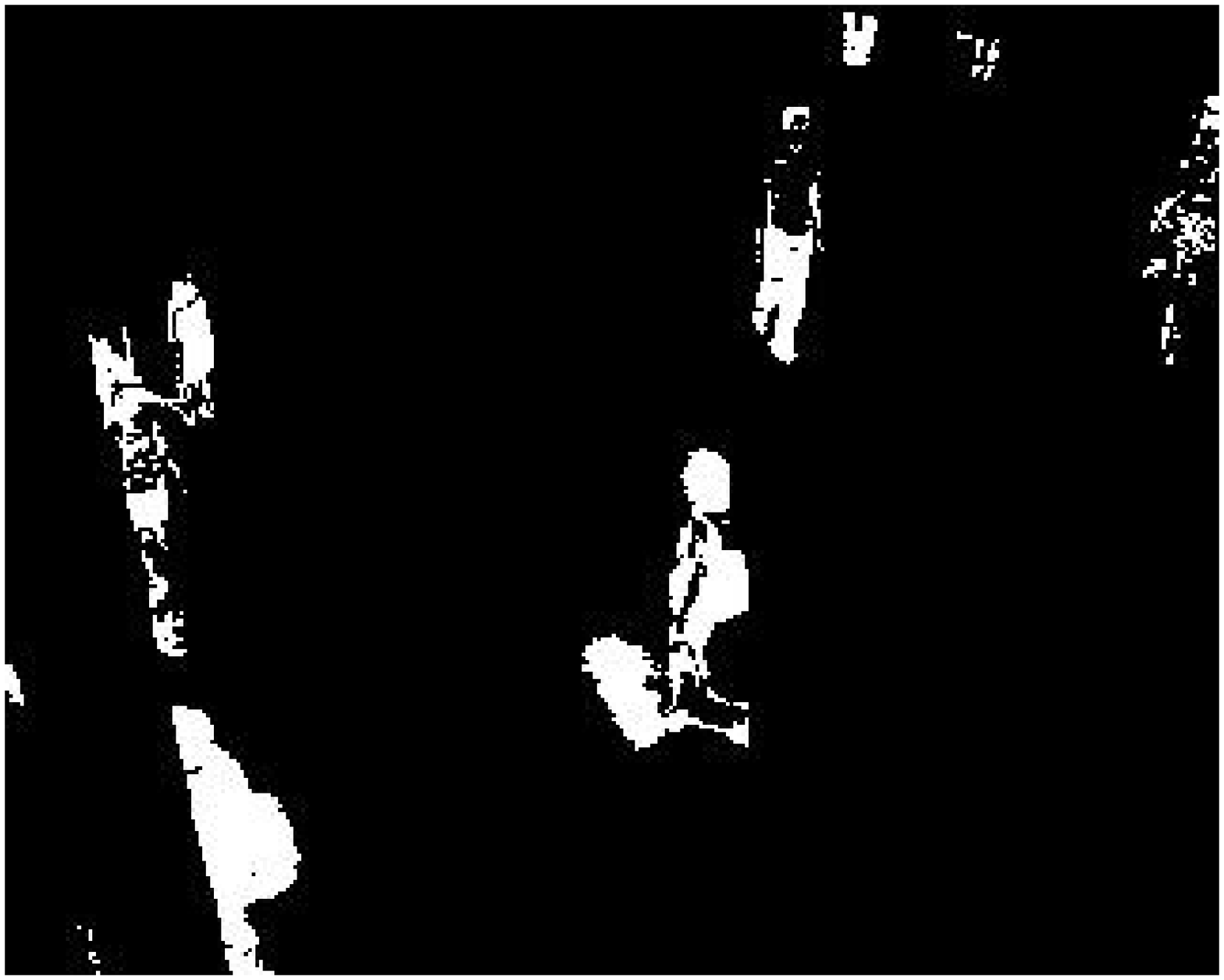}
\end{minipage}}
\subfigure{\begin{minipage}[t]{0.3\textwidth}\centering log.$\alpha=1$\vspace{1mm} \\
\includegraphics[height=1.8cm,width=1.8cm]{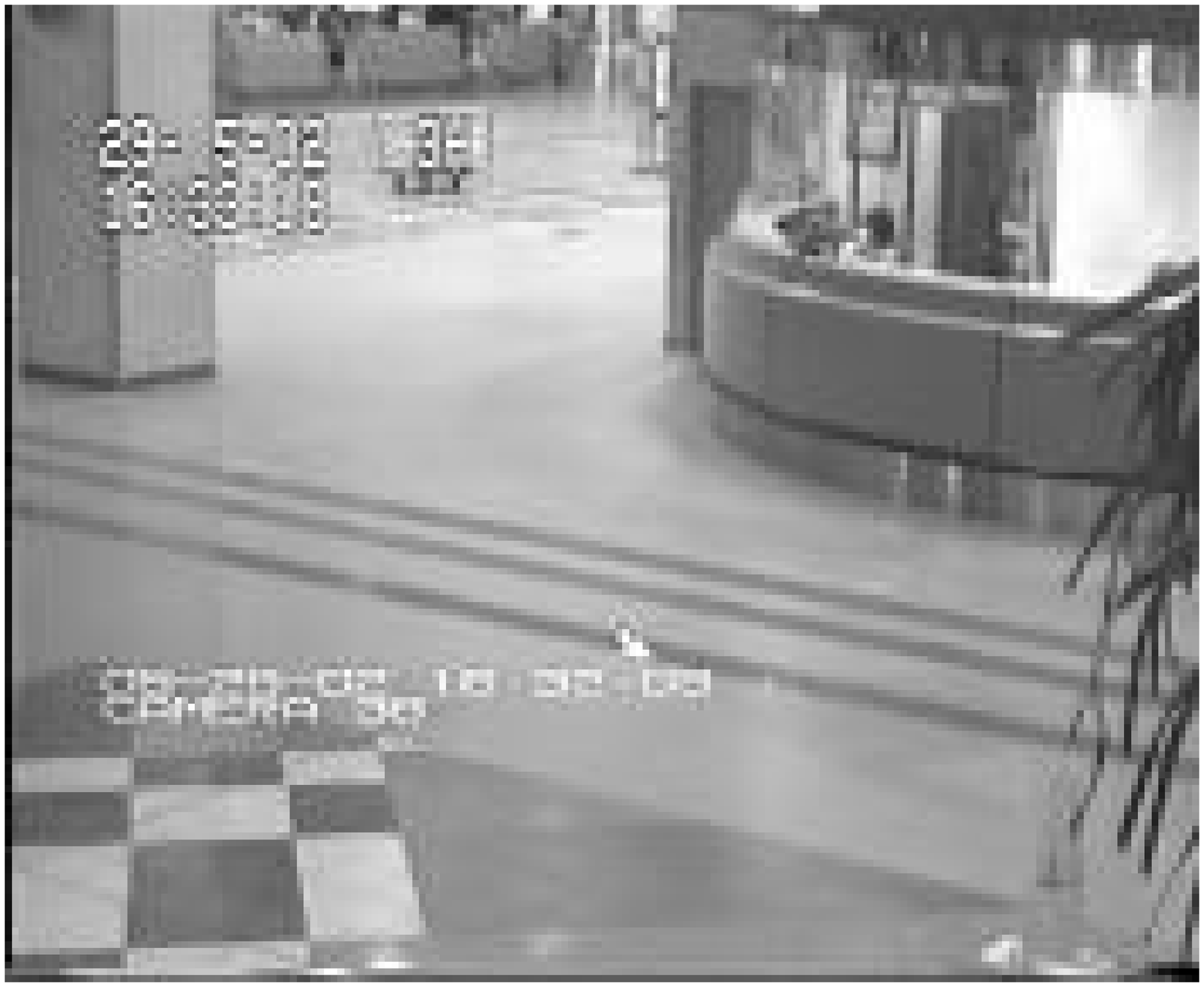}~\includegraphics[height=1.8cm,width=1.8cm]{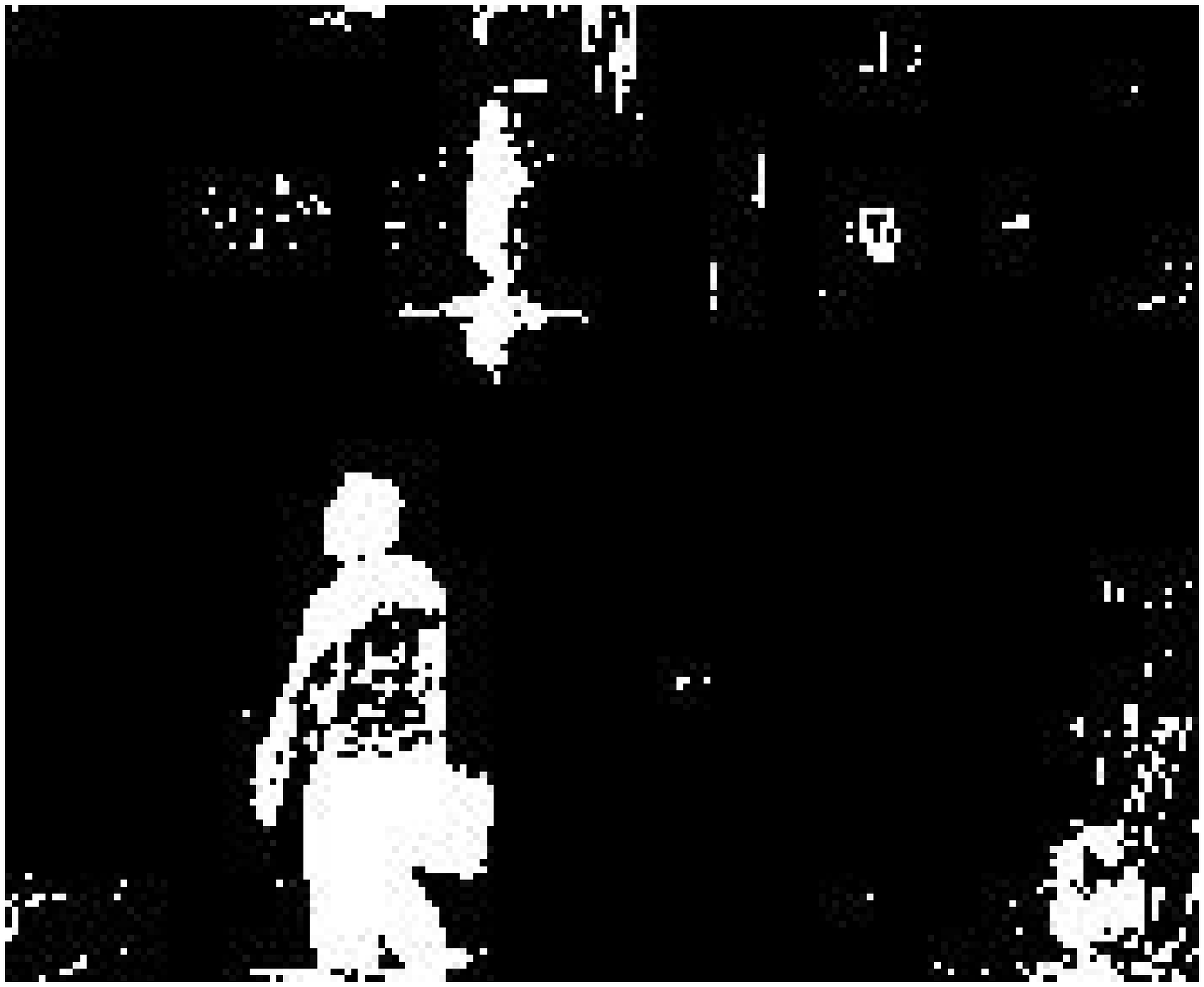}\vspace{0.5mm}\\
\includegraphics[height=1.8cm,width=1.8cm]{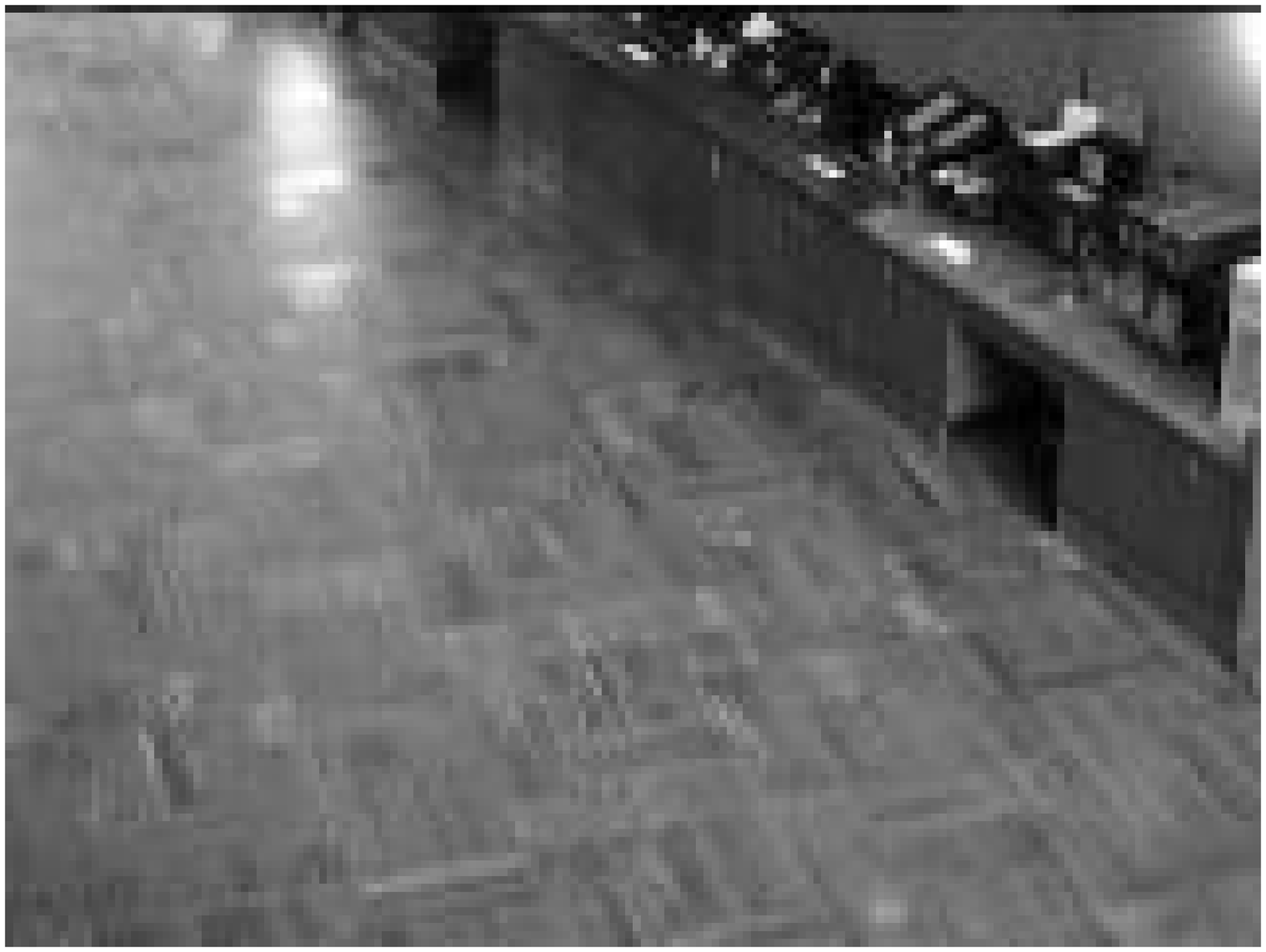}~\includegraphics[height=1.8cm,width=1.8cm]{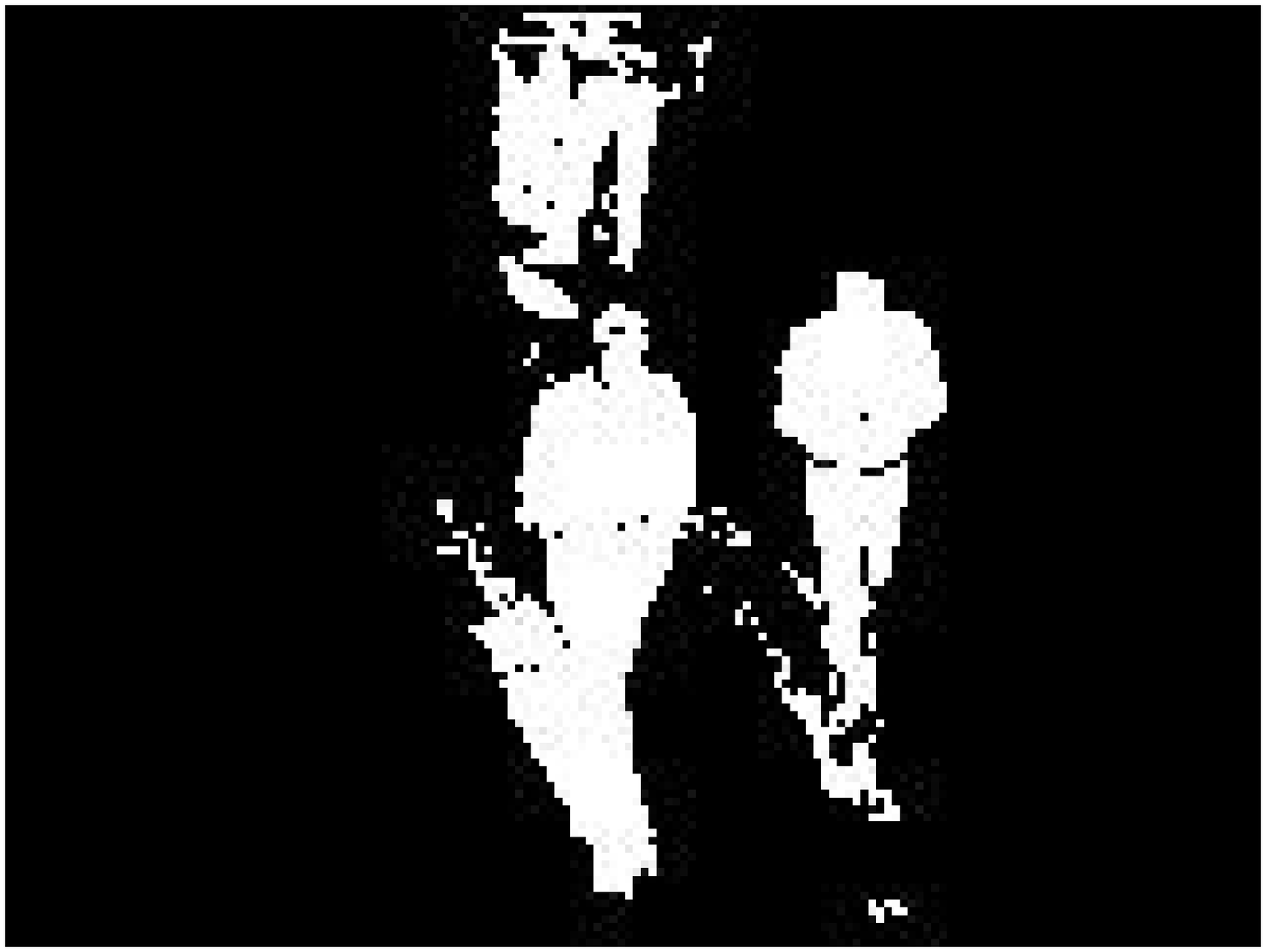}\vspace{0.5mm}\\
\includegraphics[height=1.8cm,width=1.8cm]{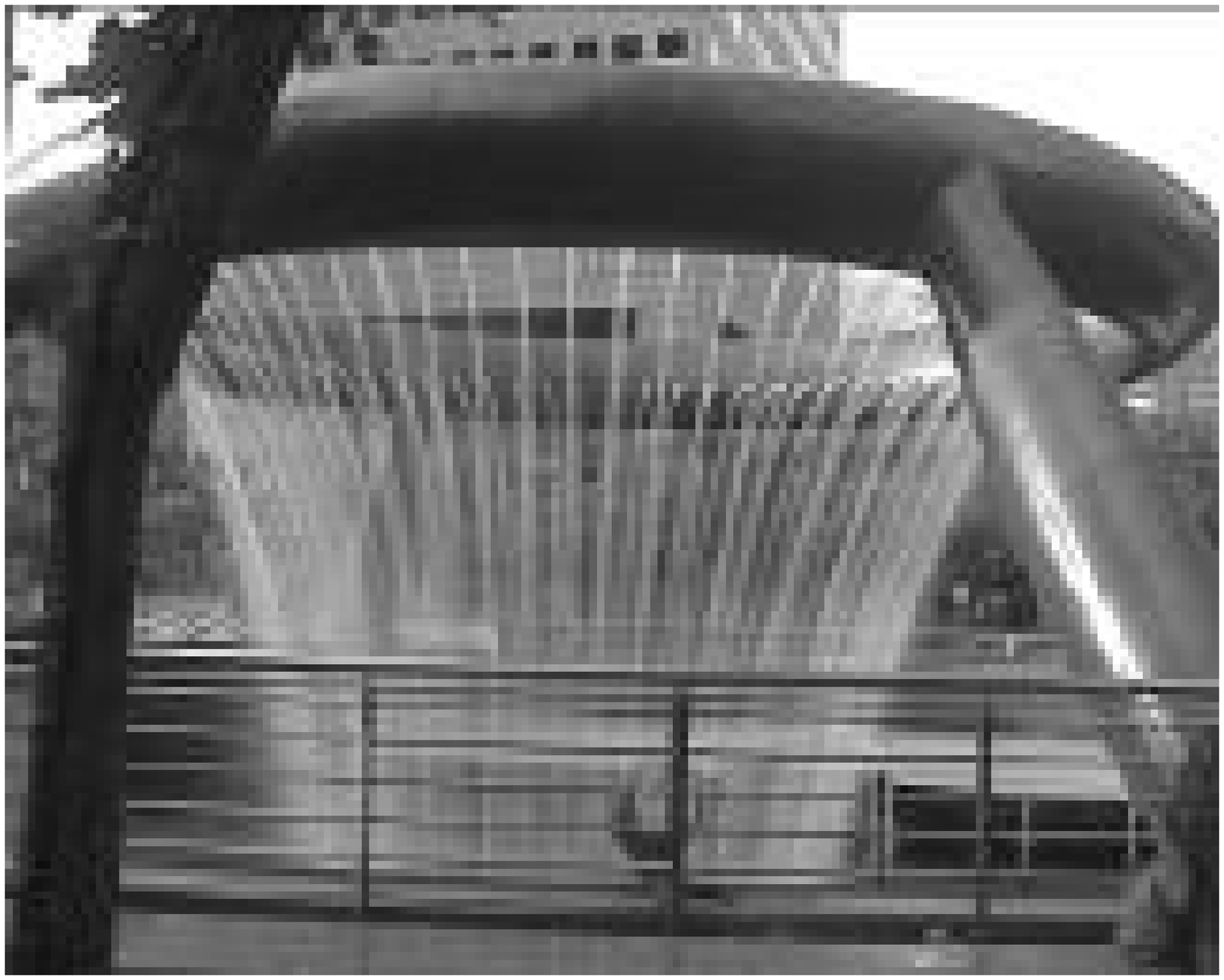}~\includegraphics[height=1.8cm,width=1.8cm]{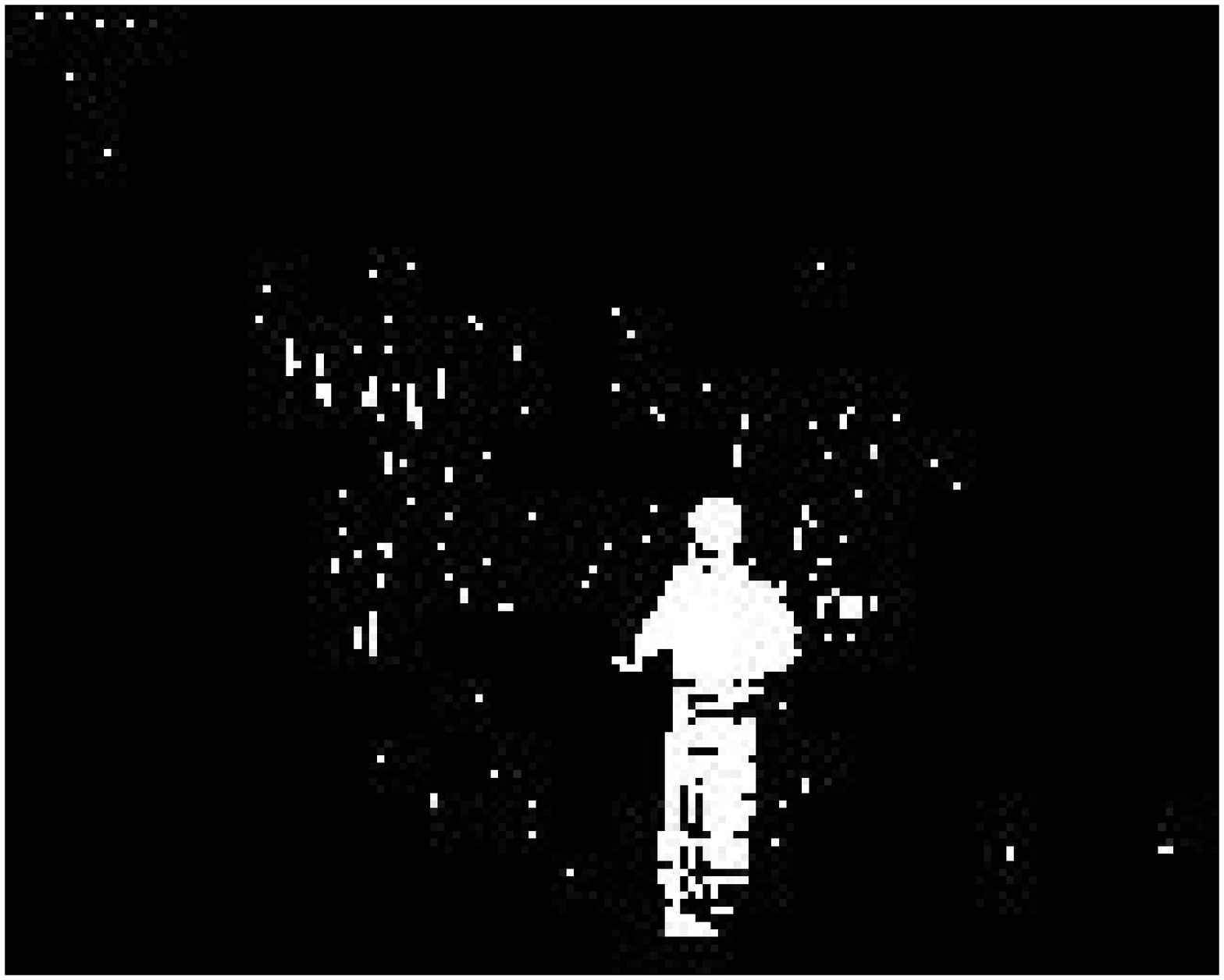}\vspace{0.5mm}\\
\includegraphics[height=1.8cm,width=1.8cm]{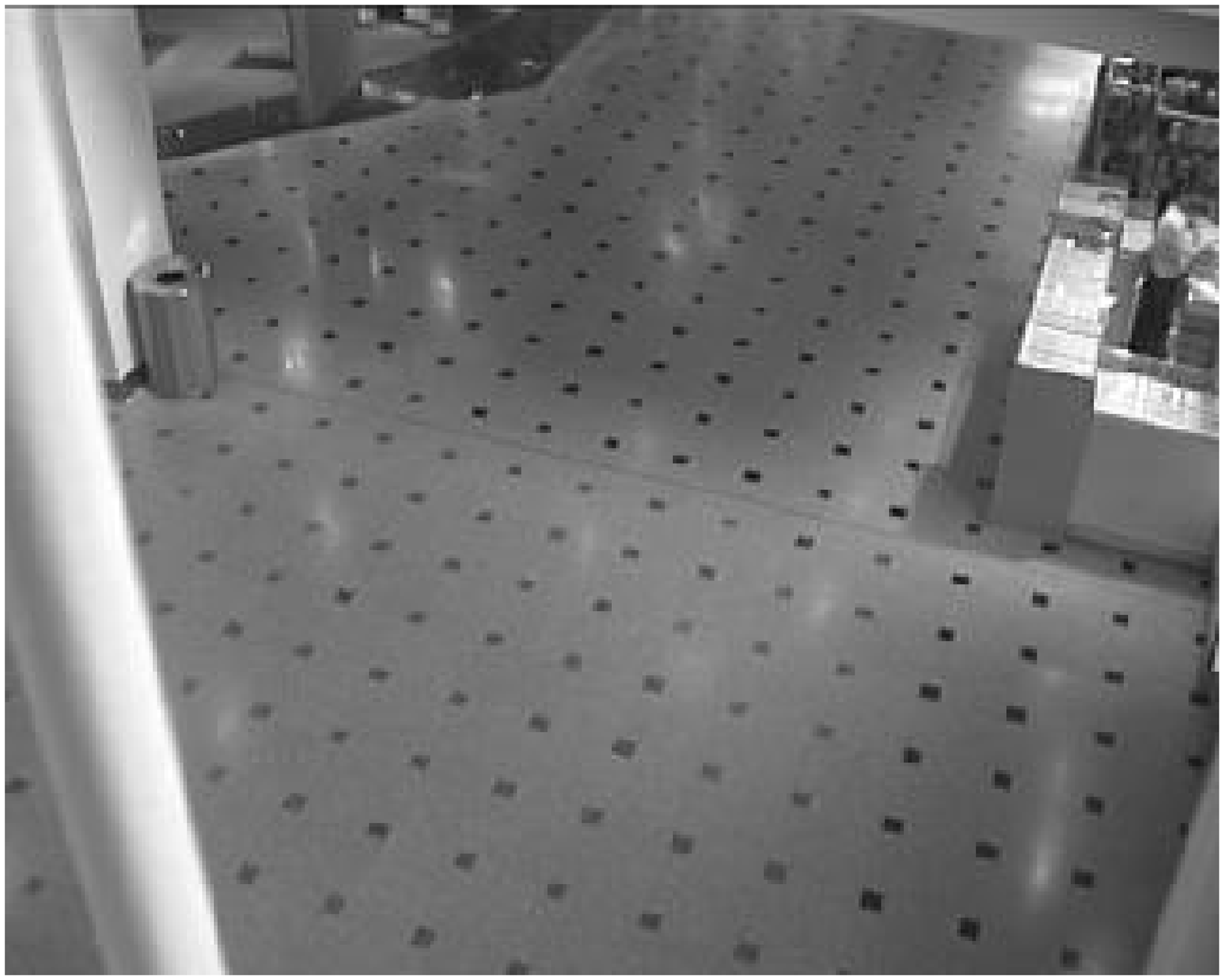}~\includegraphics[height=1.8cm,width=1.8cm]{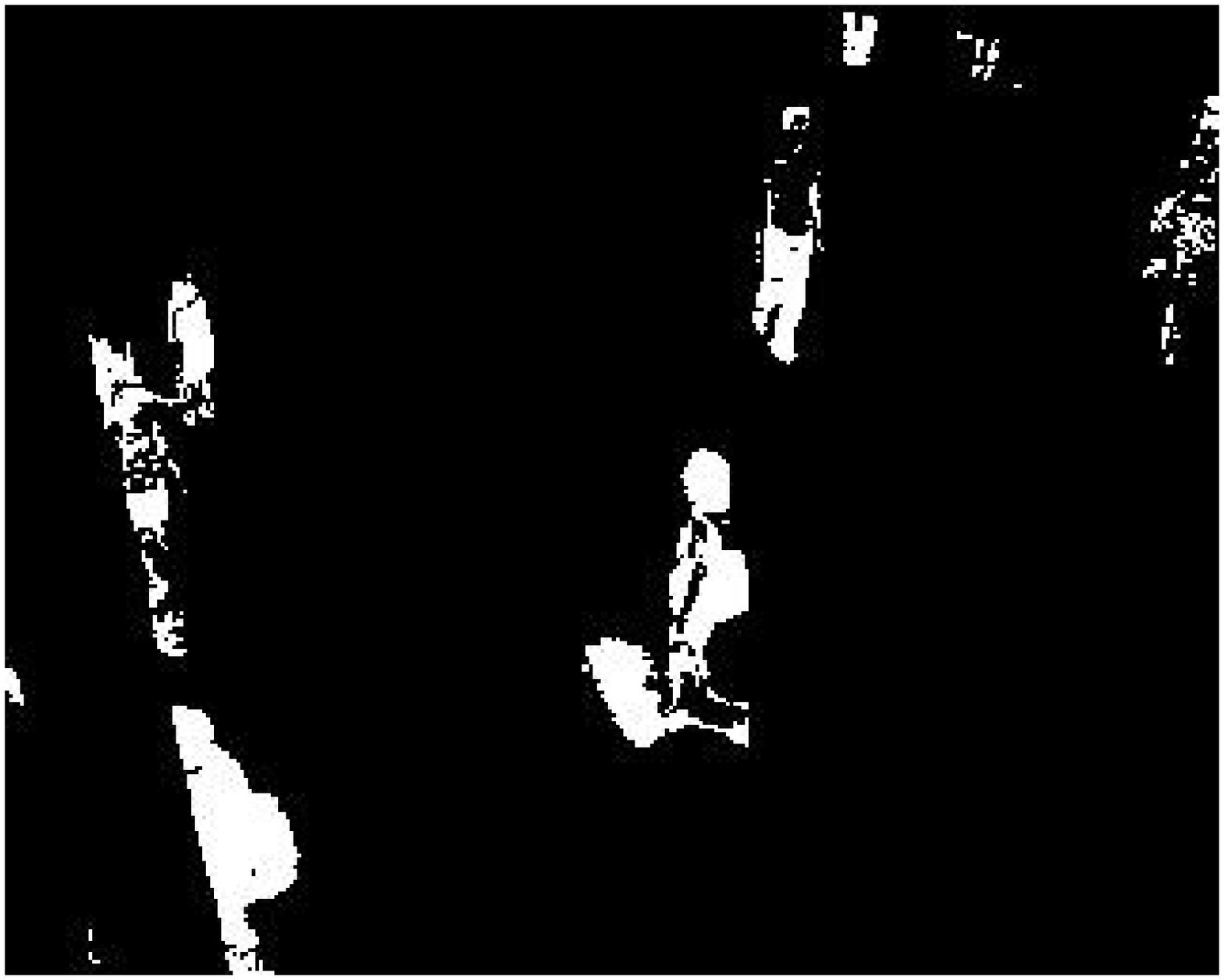}
\end{minipage}}
\subfigure{\begin{minipage}[t]{0.3\textwidth}\centering log.$\alpha=2$\vspace{1mm} \\
\includegraphics[height=1.8cm,width=1.8cm]{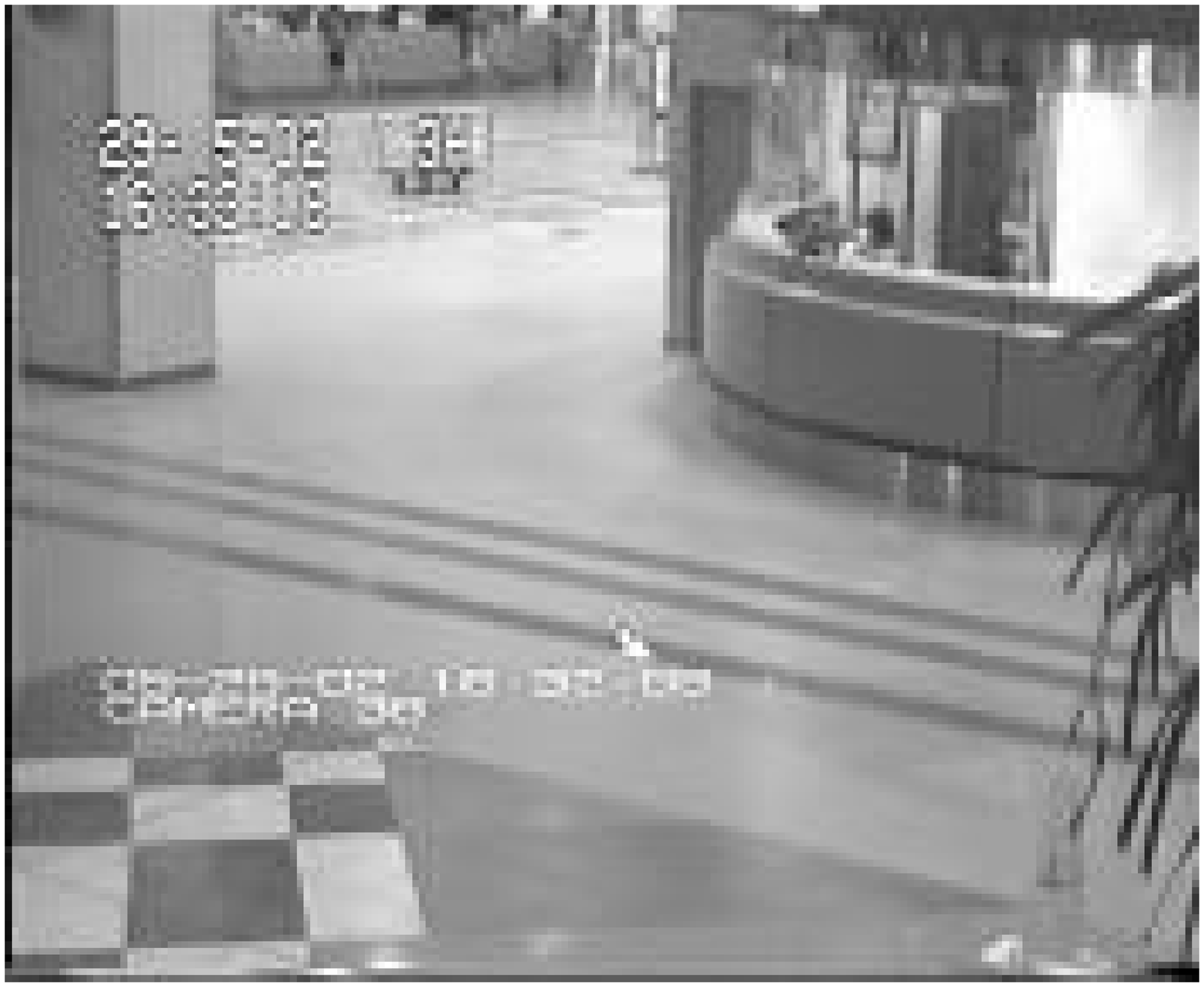}~\includegraphics[height=1.8cm,width=1.8cm]{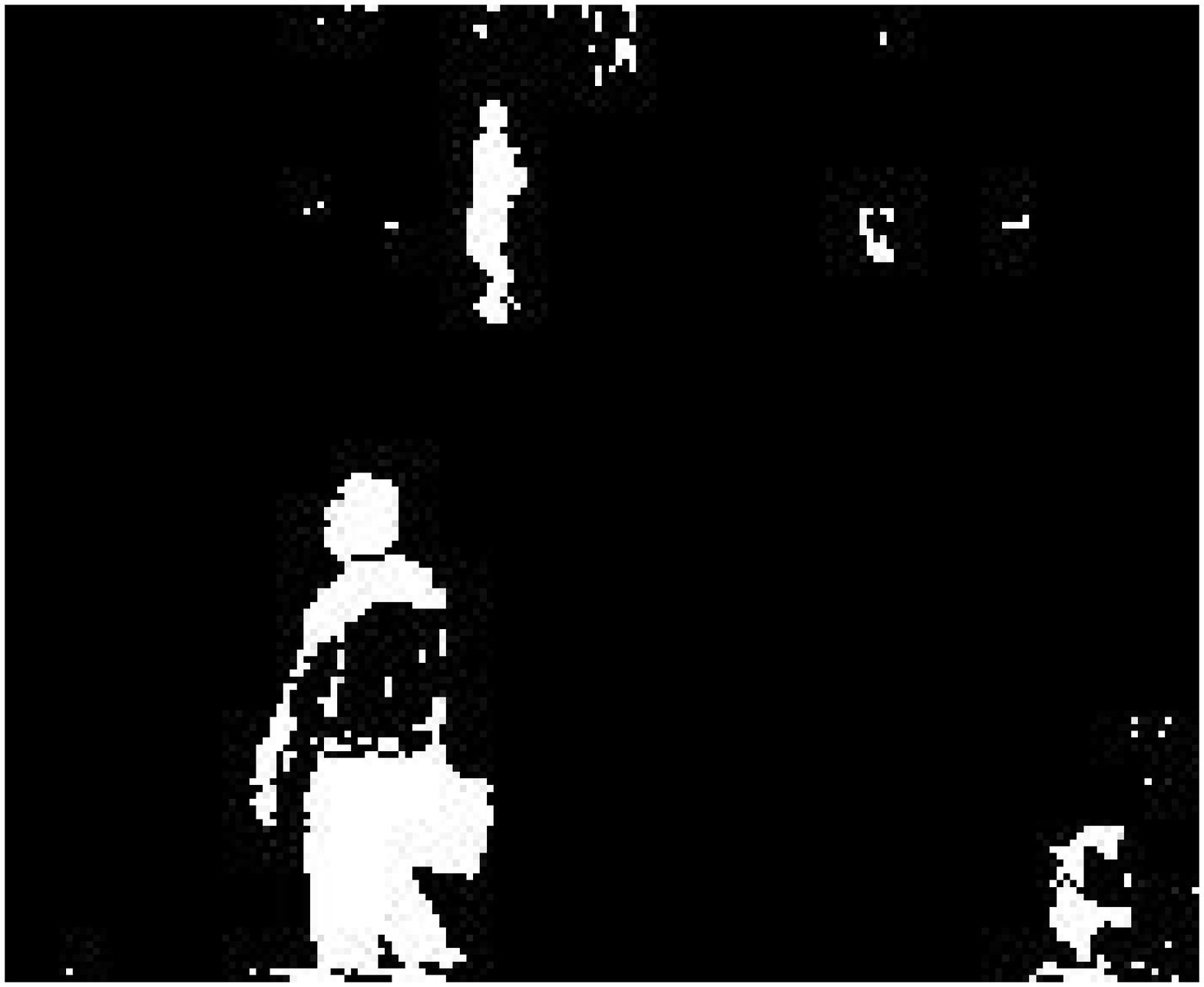}\vspace{0.5mm}\\
\includegraphics[height=1.8cm,width=1.8cm]{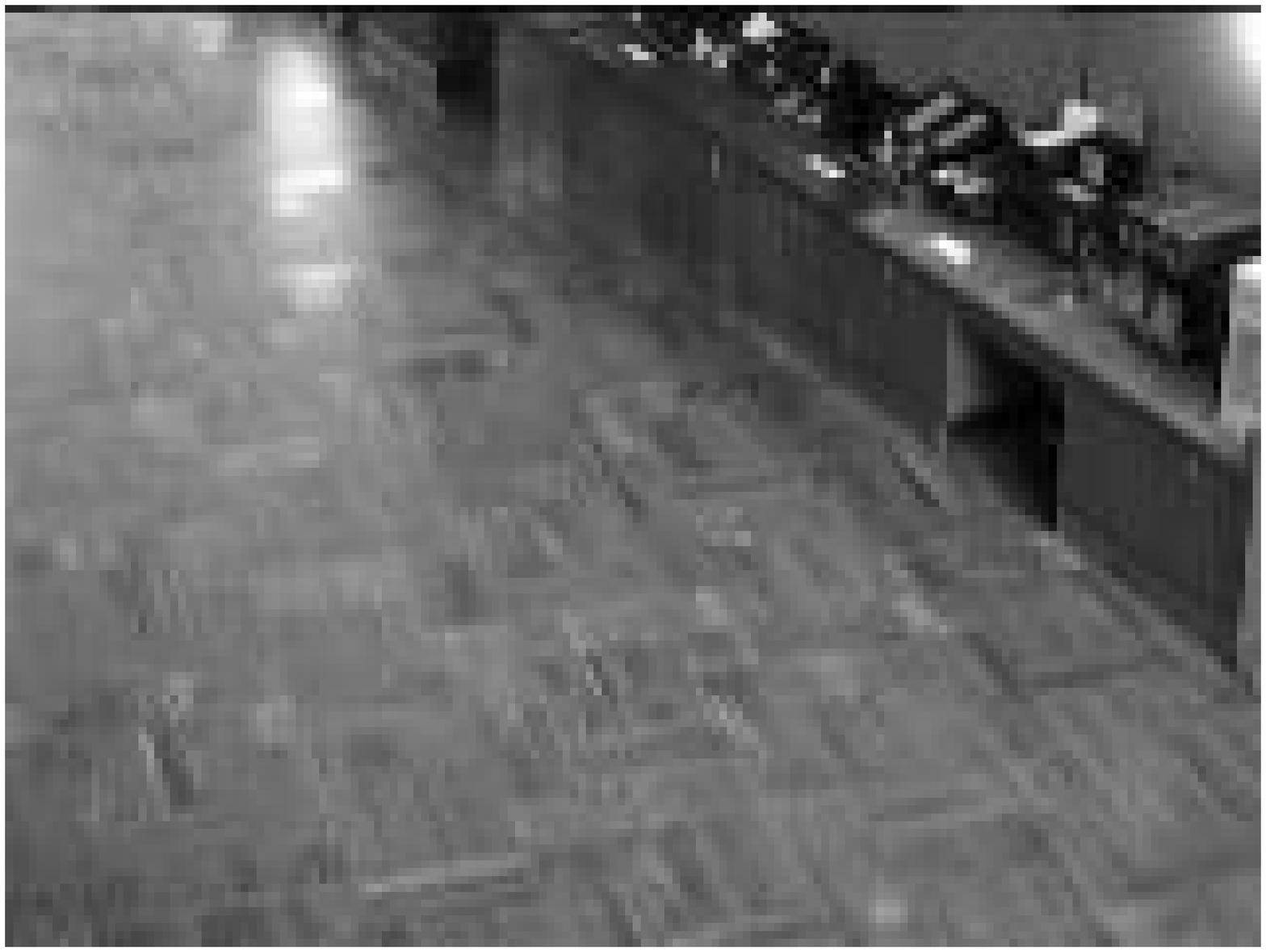}~\includegraphics[height=1.8cm,width=1.8cm]{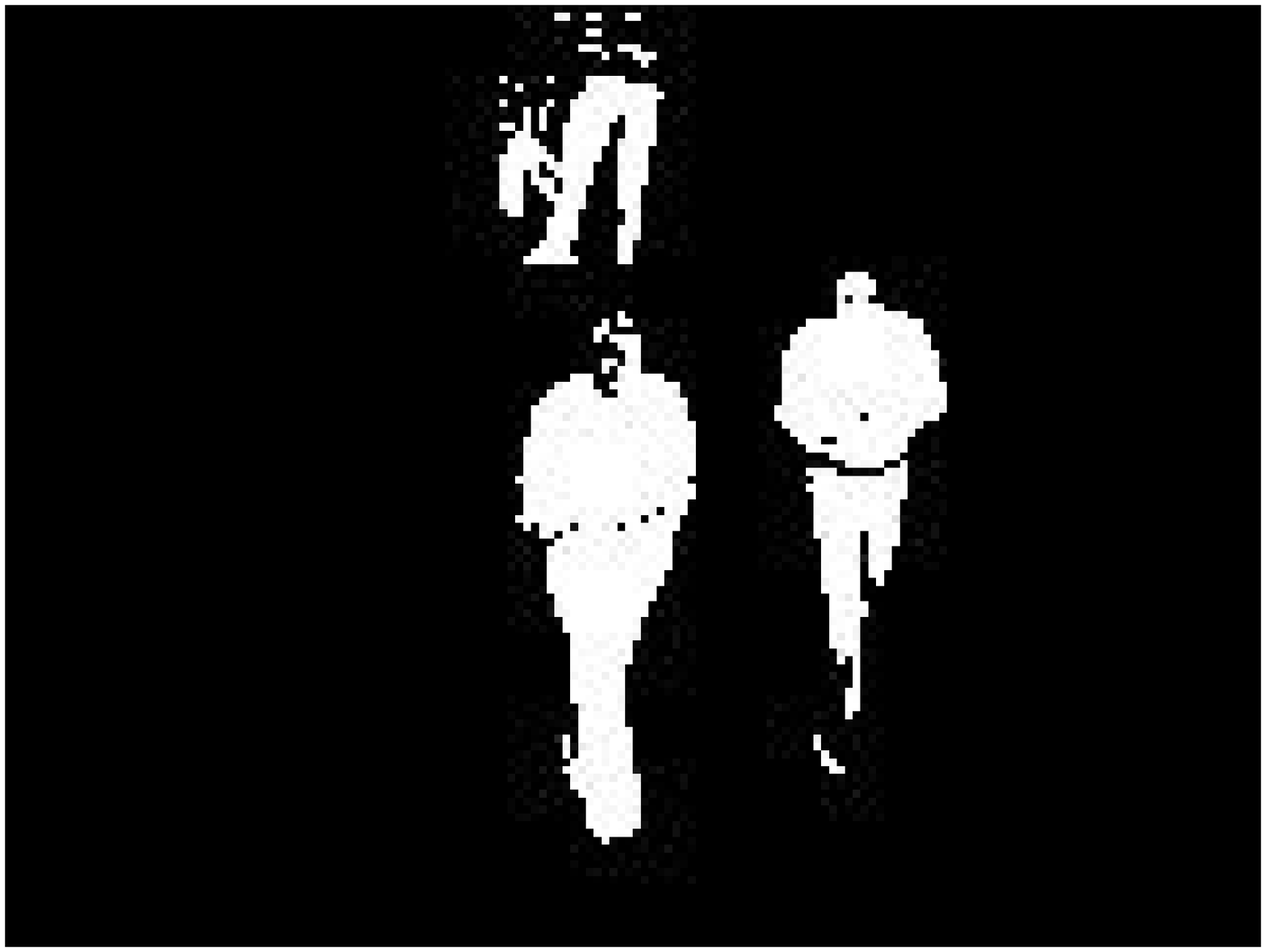}\vspace{0.5mm}\\
\includegraphics[height=1.8cm,width=1.8cm]{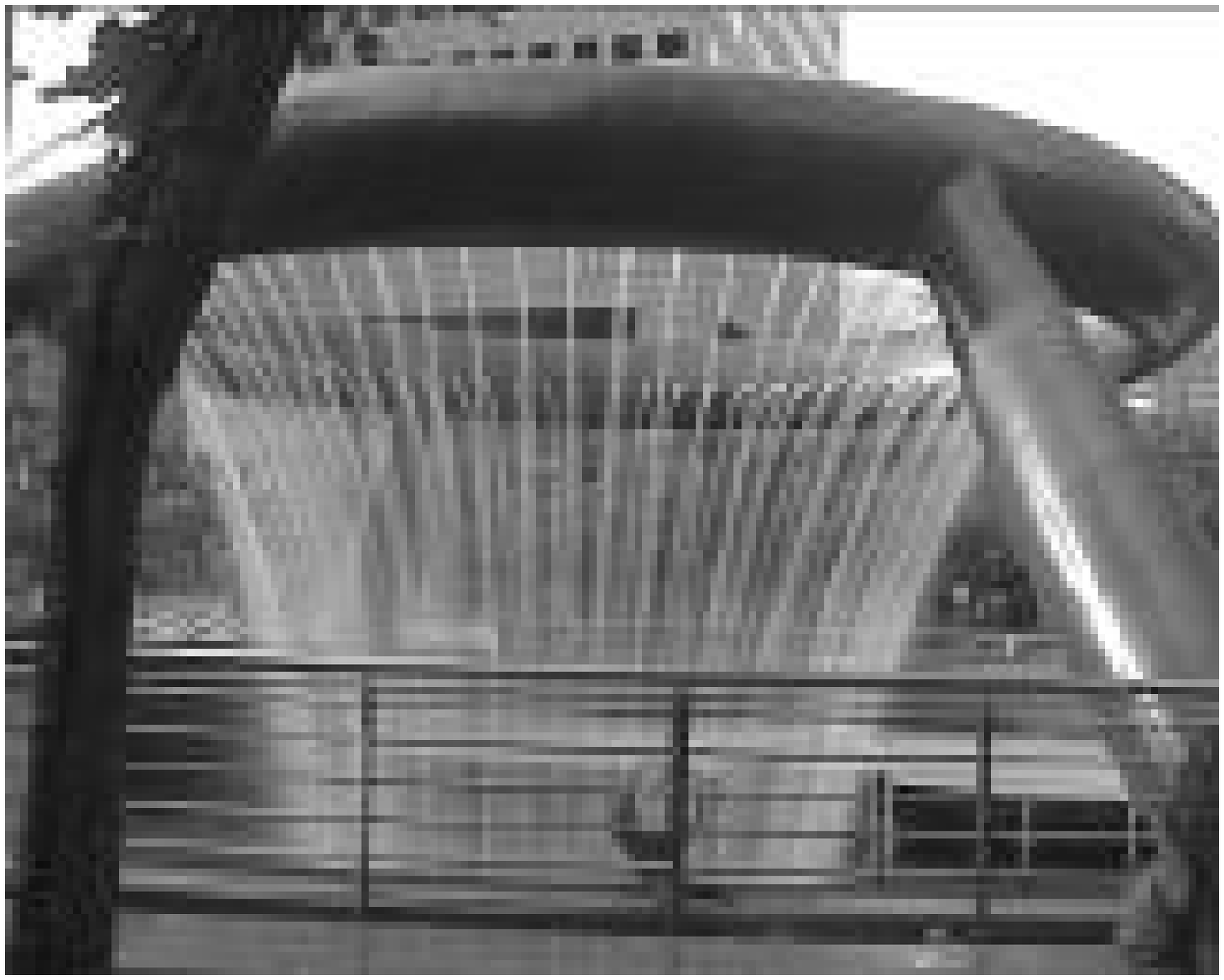}~\includegraphics[height=1.8cm,width=1.8cm]{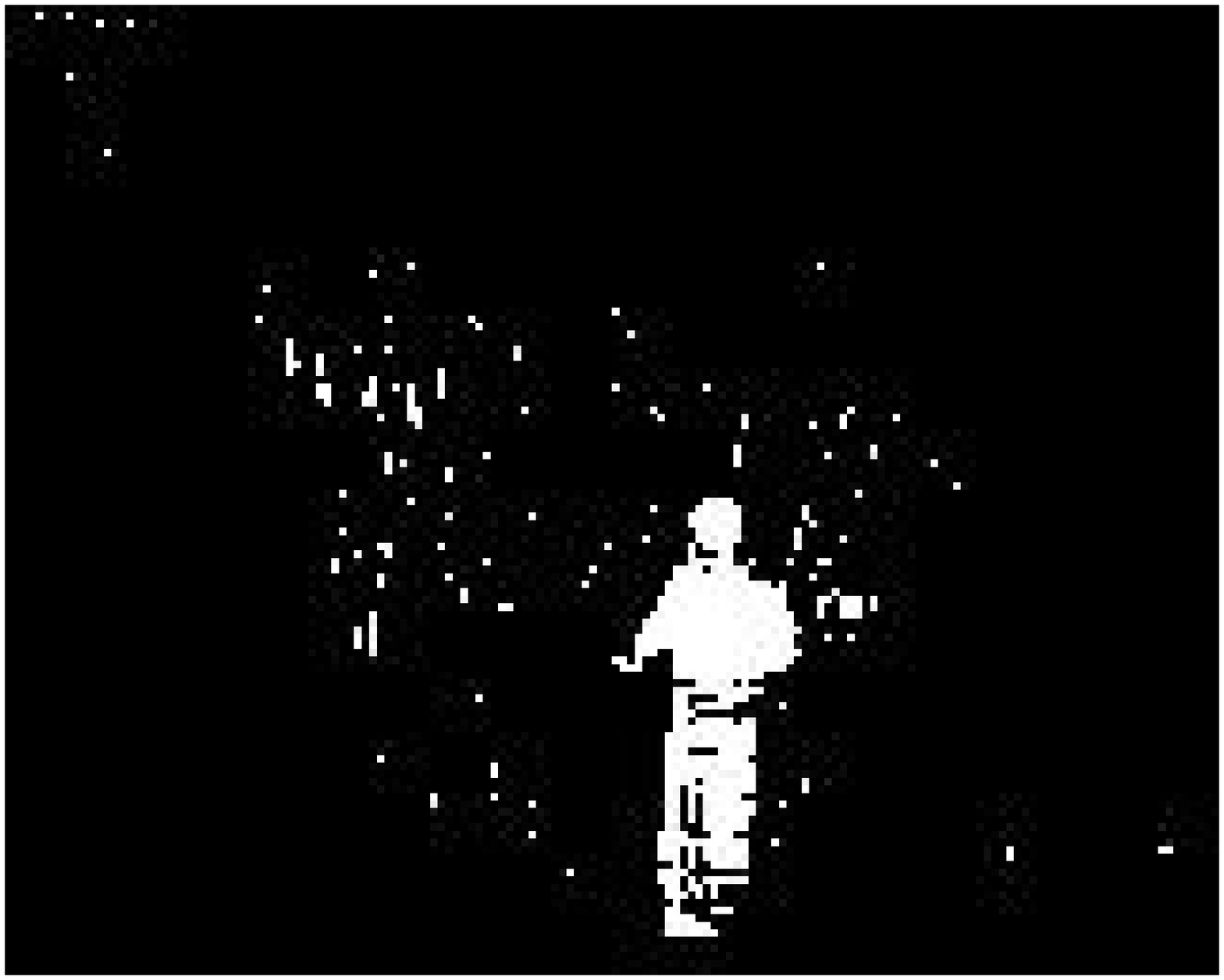}\vspace{0.5mm}\\
\includegraphics[height=1.8cm,width=1.8cm]{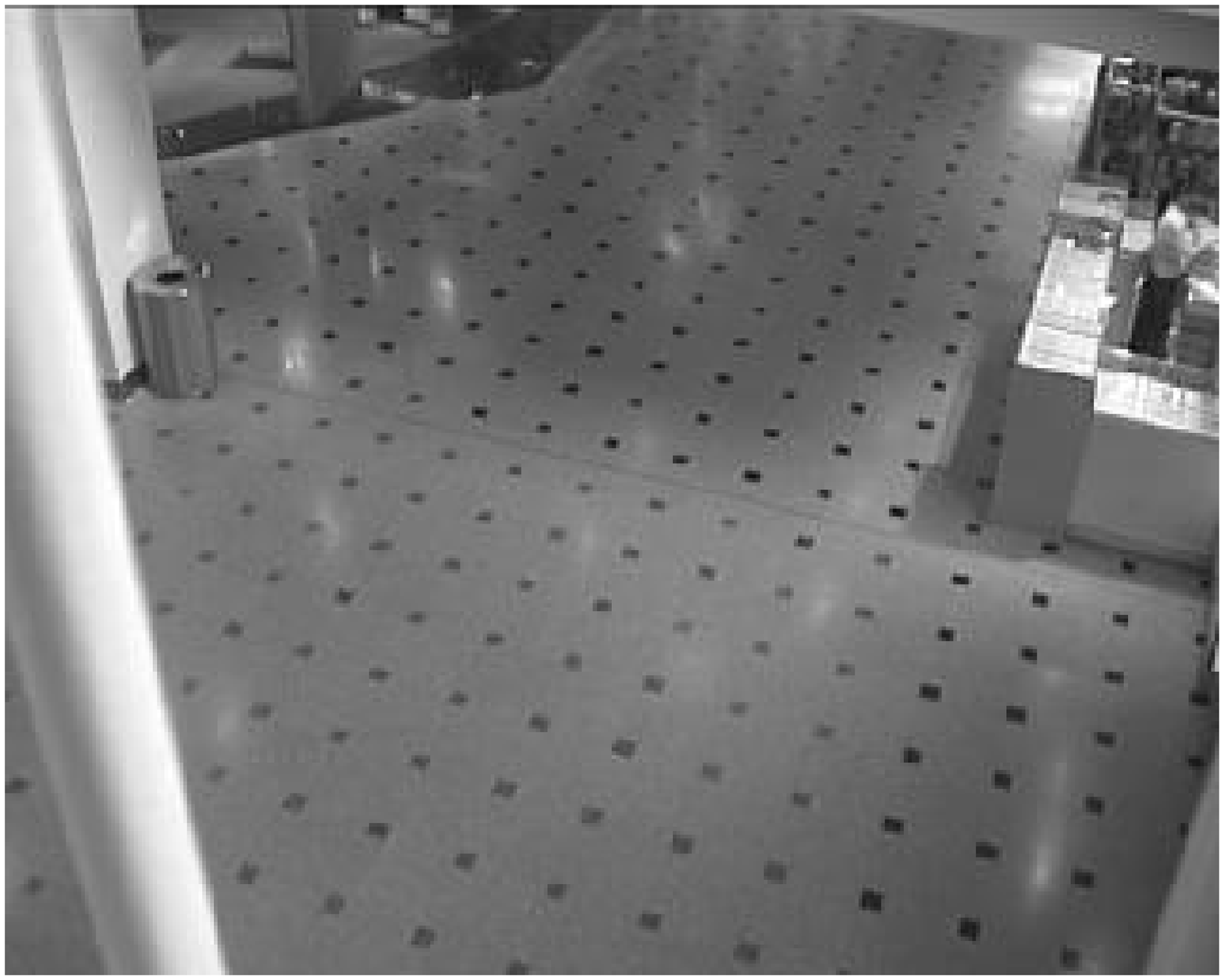}~\includegraphics[height=1.8cm,width=1.8cm]{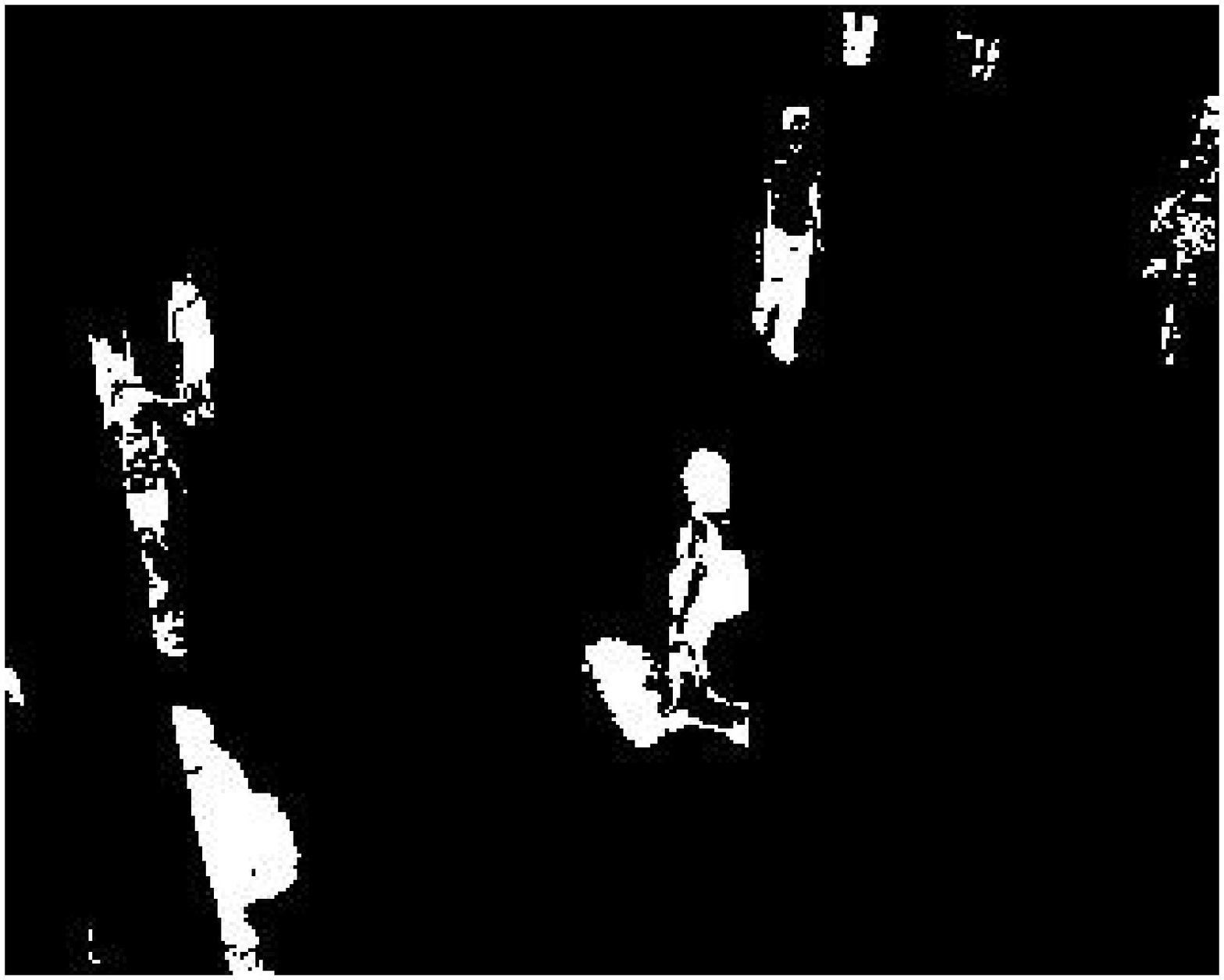}
\end{minipage}}
\caption{Extracted backgrounds and foregrounds given by the ADMM for noisy surveillance videos.}\label{fig_noisy}
\end{figure}

\paragraph{Extraction from noisy surveillance videos} In this case, ${\cal A}(L+S)=L+S$, $\lambda_{\max}=\lambda_{\min}=1$ and we set $\mathrm{Tol}_{A,1} = 10^{-4}$, $\mathrm{Tol}_{A,2} = 5\times10^{-3}$ and $\mathrm{Tol}_{P} = 10^{-4}$. The computational results are reported in Table \ref{noisyresults}, where we report $p$ and $\alpha$, the optimal $\mu$, the number of iterations, the CPU time (seconds) and F-measure. We also show the extracted backgrounds and foregrounds given by the ADMM in Fig.\,\ref{fig_noisy}.

\begin{table}[ht]
\caption{Numerical results for extraction from noisy surveillance videos}\label{noisyresults}
\centering \tabcolsep 3pt
\begin{tabular}{|clc|cccc|cccc|}
\hline
&&& \multicolumn{4}{c|}{ADMM} & \multicolumn{4}{c|}{PALM}               \\
Data & \multicolumn{2}{c|}{regularizer} & $\mu$ & iter & time & \footnotesize{F-measure} & $\mu$ & iter & time & \footnotesize{F-measure} \\
\hline
\multirow{6}{*}{\small{Hall}}
& bri.\,$p$       &   1.0  &  5e-02  & 10  & 3.21  & 0.7562 & 5e-02  & 19  & 3.96 & 0.7560 \\
&                 &   0.5  &  1e-02  & 32  & 11.26  & 0.7634 & 1e-02  & 36  & 9.29 & 0.7624 \\
& fra.\,$\alpha$  &   1.0  &  5e-02  & 23  & 8.26  & 0.7578 & 5e-02  & 33  & 8.53 & 0.7578 \\
&                 &   2.0  &  5e-02  & 12  & 4.17  & 0.7368 & 5e-02  & 15  & 3.69 & 0.7371 \\
& log.\,$\alpha$  &   1.0  &  5e-02  & 12  & 21.12  & 0.7566 & 5e-02  & 39  & 68.70 & 0.7576 \\
&                 &   2.0  &  5e-02  & 12  & 16.00  & 0.7368 & 5e-02  & 16  & 29.04 & 0.7368 \\
\hline

\multirow{6}{*}{\small{Bootstrap}}
& bri.\,$p$      &   1.0  &  1e-01  & 14  & 3.30  & 0.8180 & 1e-01  & 19  & 3.15 & 0.8180 \\
&                &   0.5  &  5e-02  & 23  & 6.77  & 0.8206 & 5e-02  & 22  & 4.93 & 0.8209 \\
& fra.\,$\alpha$ &   1.0  &  1e-01  & 15  & 4.91  & 0.8163 & 1e-01  & 20  & 5.32 & 0.8165 \\
&                &   2.0  &  1e-01  & 14  & 4.18  & 0.8264 & 1e-01  & 16  & 3.72 & 0.8261 \\
& log.\,$\alpha$ &   1.0  &  1e-01  & 16  & 21.92  & 0.8195 & 1e-01  & 22  & 28.62 & 0.8195 \\
&                &   2.0  &  1e-01  & 12  & 8.91  & 0.8363 & 1e-01  & 10  & 6.50 & 0.8363 \\
\hline

\multirow{6}{*}{\small{Fountain}}
& bri.\,$p$      &   1.0  &  1e-01  & 9  & 2.19  & 0.7749 & 1e-01  & 7  & 1.10 & 0.7749 \\
&                &   0.5  &  5e-02  & 13  & 3.54  & 0.7000 & 5e-02  & 11  & 2.13 & 0.6922 \\
& fra.\,$\alpha$ &   1.0  &  1e-01  & 9  & 2.39  & 0.7717 & 1e-01  & 8  & 1.63 & 0.7717 \\
&                &   2.0  &  5e-02  & 10  & 2.82  & 0.7717 & 5e-02  & 9  & 1.87 & 0.7717 \\
& log.\,$\alpha$ &   1.0  &  1e-01  & 9  & 13.41  & 0.7738 & 1e-01  & 7  & 9.65 & 0.7738 \\
&                &   2.0  &  5e-02  & 9  & 12.46  & 0.7717 & 5e-02  & 8  & 11.51 & 0.7717 \\
\hline

\multirow{6}{*}{\small{ShoppingMall}}
& bri.\,$p$       &   1.0  &  1e-01  & 10  & 9.66  & 0.7046 & 1e-01  & 13  & 8.73 & 0.7043 \\
&                 &   0.5  &  1e-02  & 39  & 52.39  & 0.7087 & 1e-02  & 79  & 83.62 & 0.7078 \\
& fra.\,$\alpha$  &   1.0  &  1e-01  & 12  & 14.33  & 0.7055 & 1e-01  & 18  & 16.95 & 0.7055 \\
&                 &   2.0  &  5e-02  & 15  & 18.46  & 0.7062 & 5e-02  & 26  & 25.34 & 0.7064 \\
& log.\,$\alpha$  &   1.0  &  1e-01  & 11  & 66.96  & 0.7055 & 1e-01  & 16  & 94.06 & 0.7055 \\
&                 &   2.0  &  5e-02  & 12  & 40.23  & 0.7057 & 5e-02  & 18  & 74.83 & 0.7057 \\
\hline
\end{tabular}
\end{table}

\begin{figure}[ht]
\centering
\subfigure{\begin{minipage}[t]{0.3\textwidth}\centering bri.$p=1.0$\vspace{1mm} \\
\includegraphics[height=1.8cm,width=1.8cm]{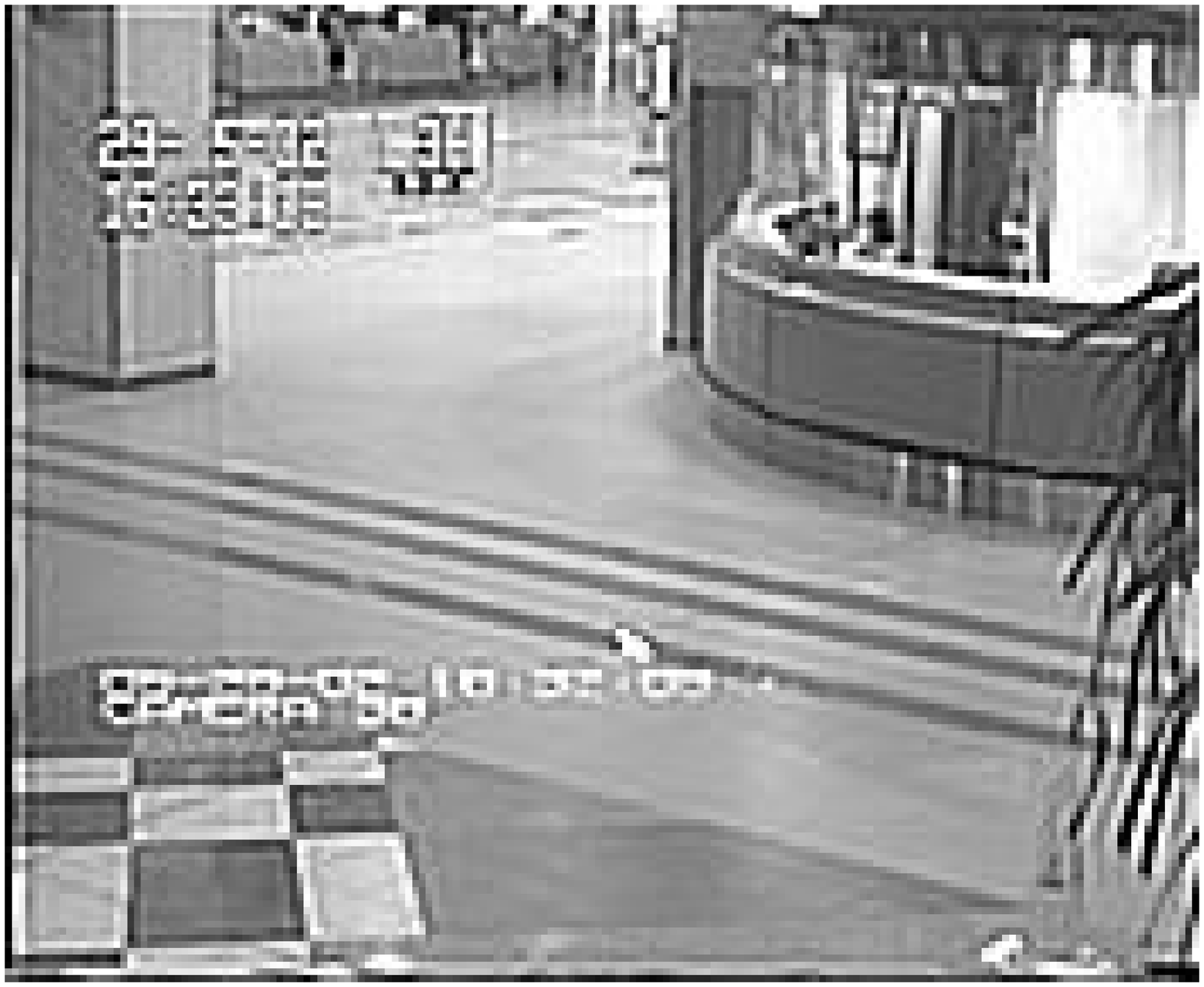}~\includegraphics[height=1.8cm,width=1.8cm]{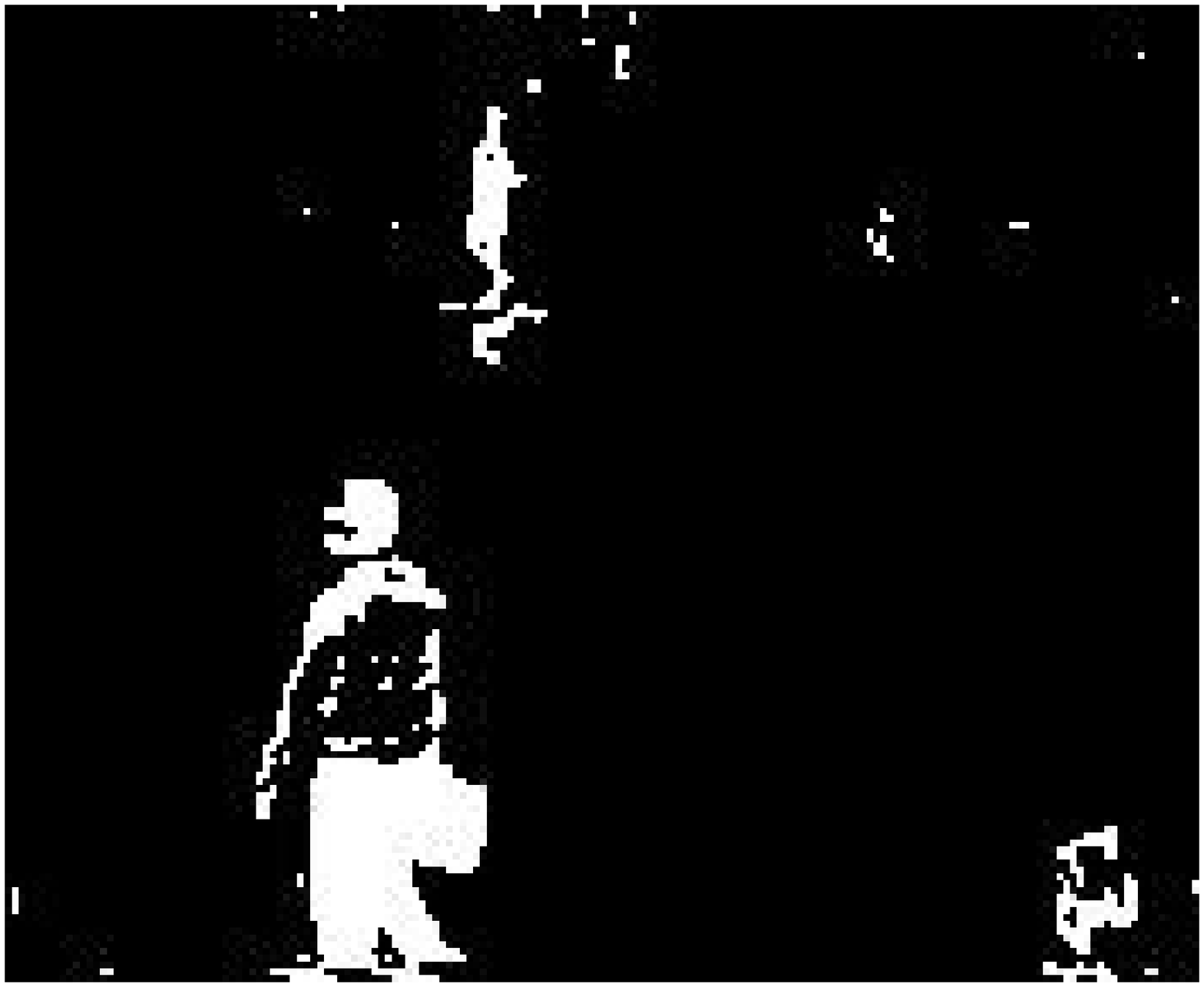}\vspace{0.5mm}\\
\includegraphics[height=1.8cm,width=1.8cm]{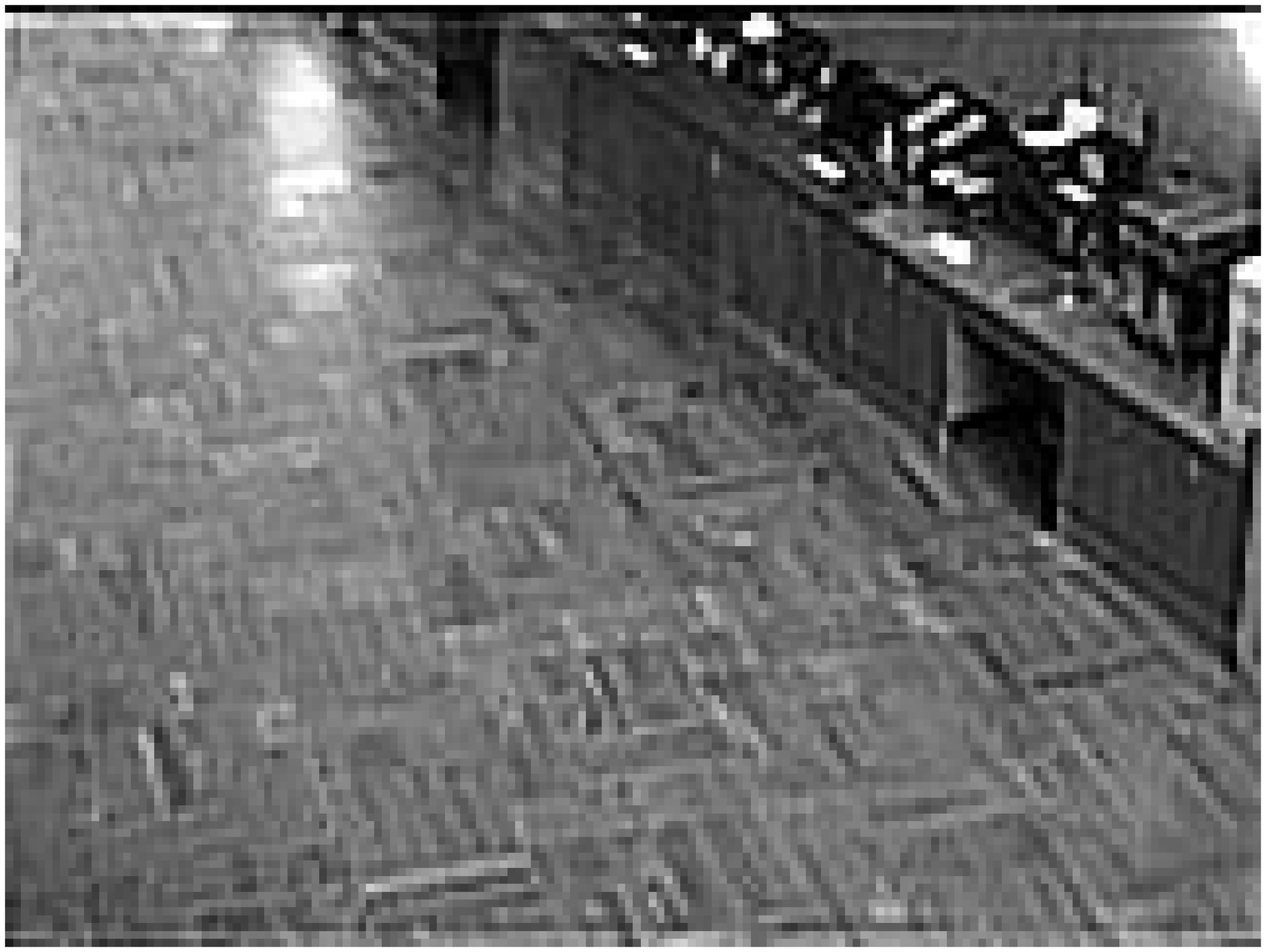}~\includegraphics[height=1.8cm,width=1.8cm]{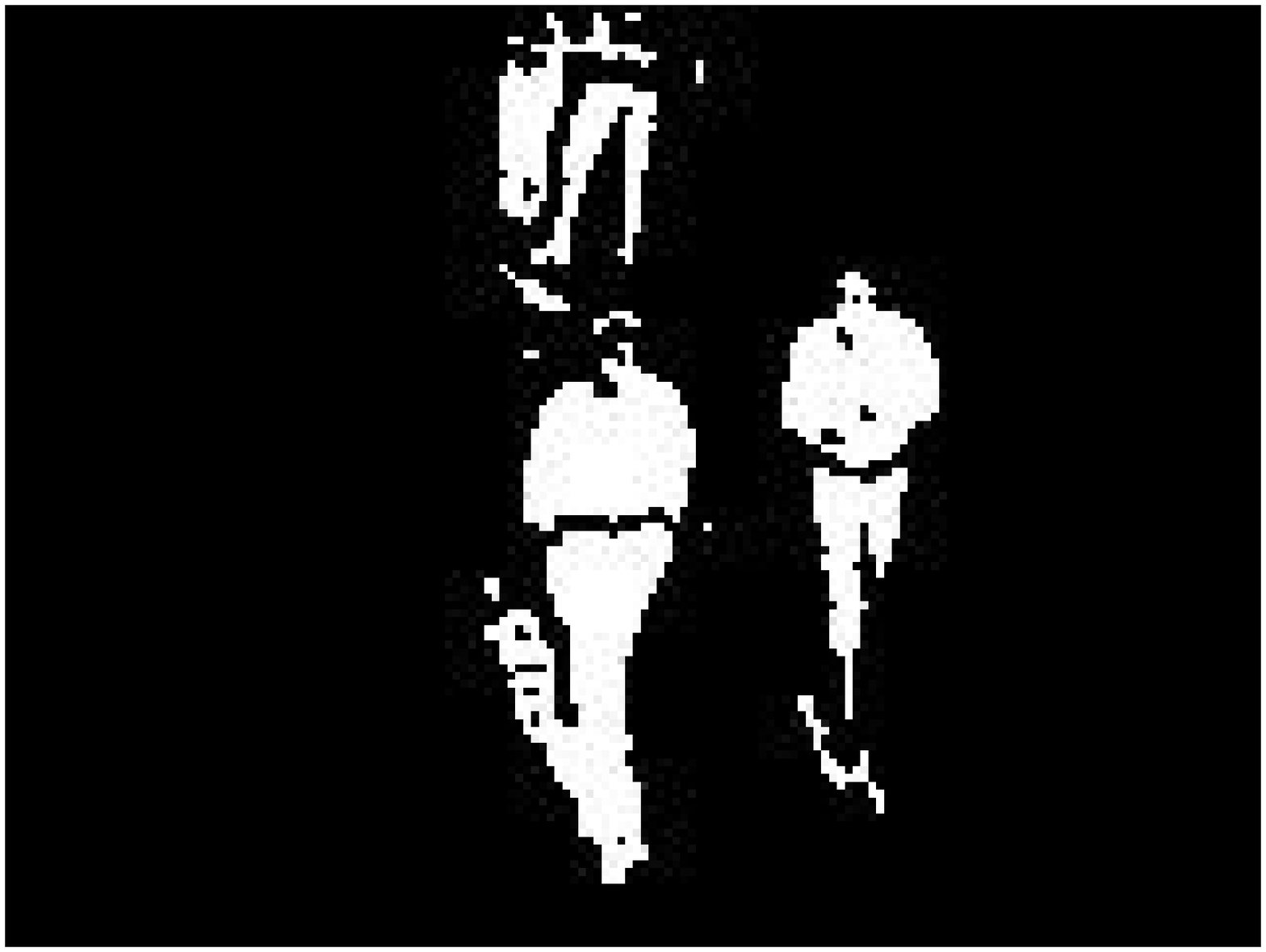}\vspace{0.5mm}\\
\includegraphics[height=1.8cm,width=1.8cm]{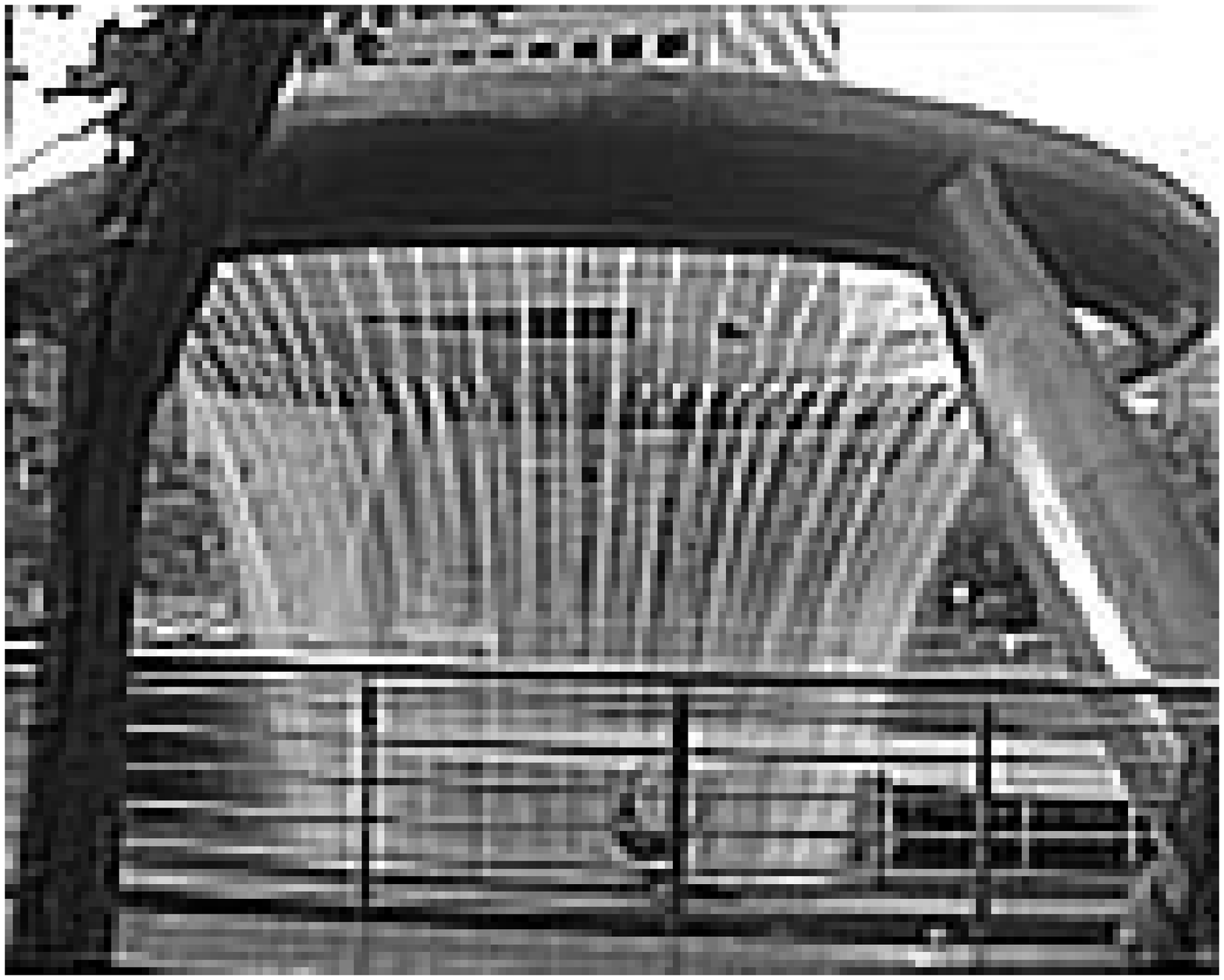}~\includegraphics[height=1.8cm,width=1.8cm]{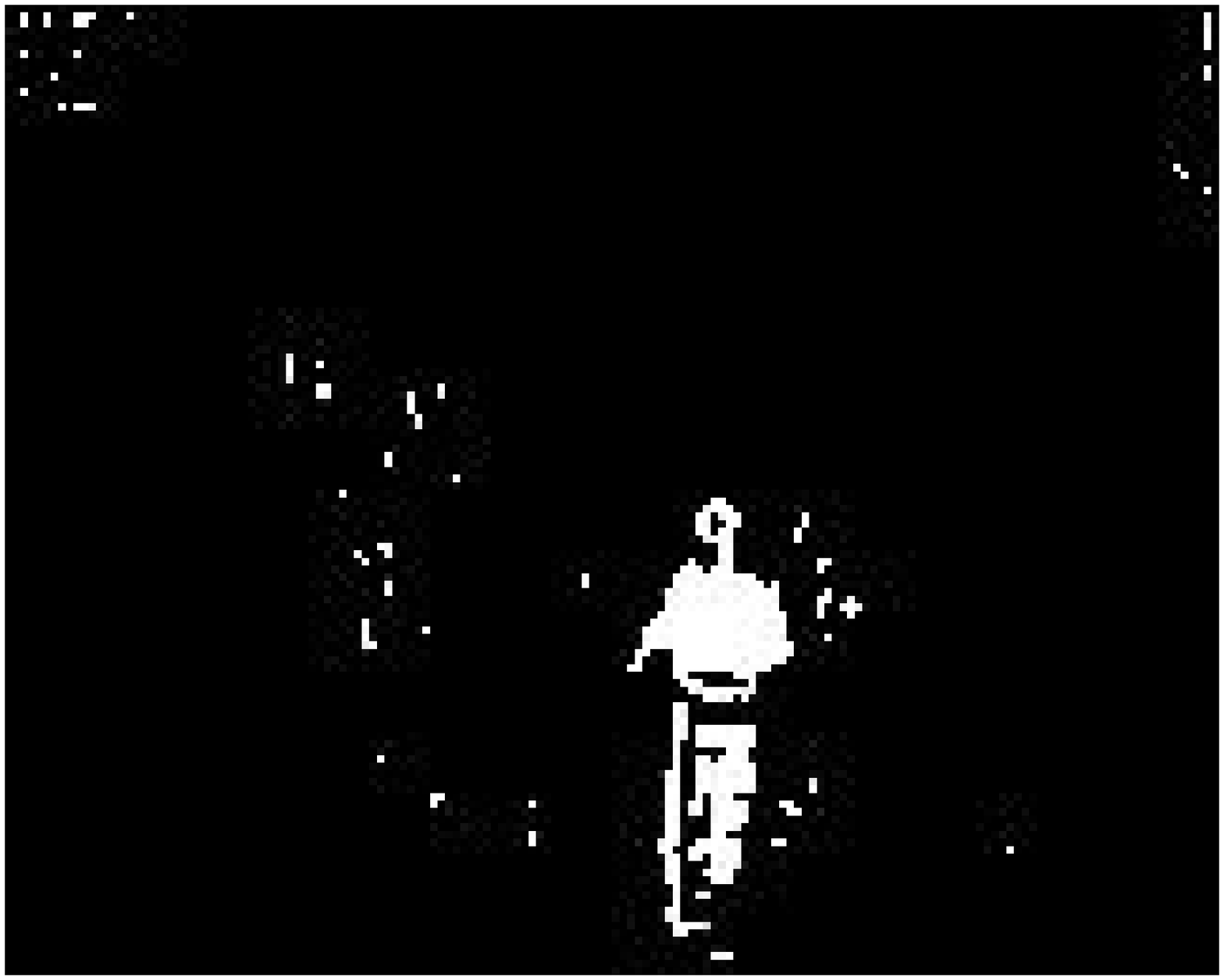}\vspace{0.5mm}\\
\includegraphics[height=1.8cm,width=1.8cm]{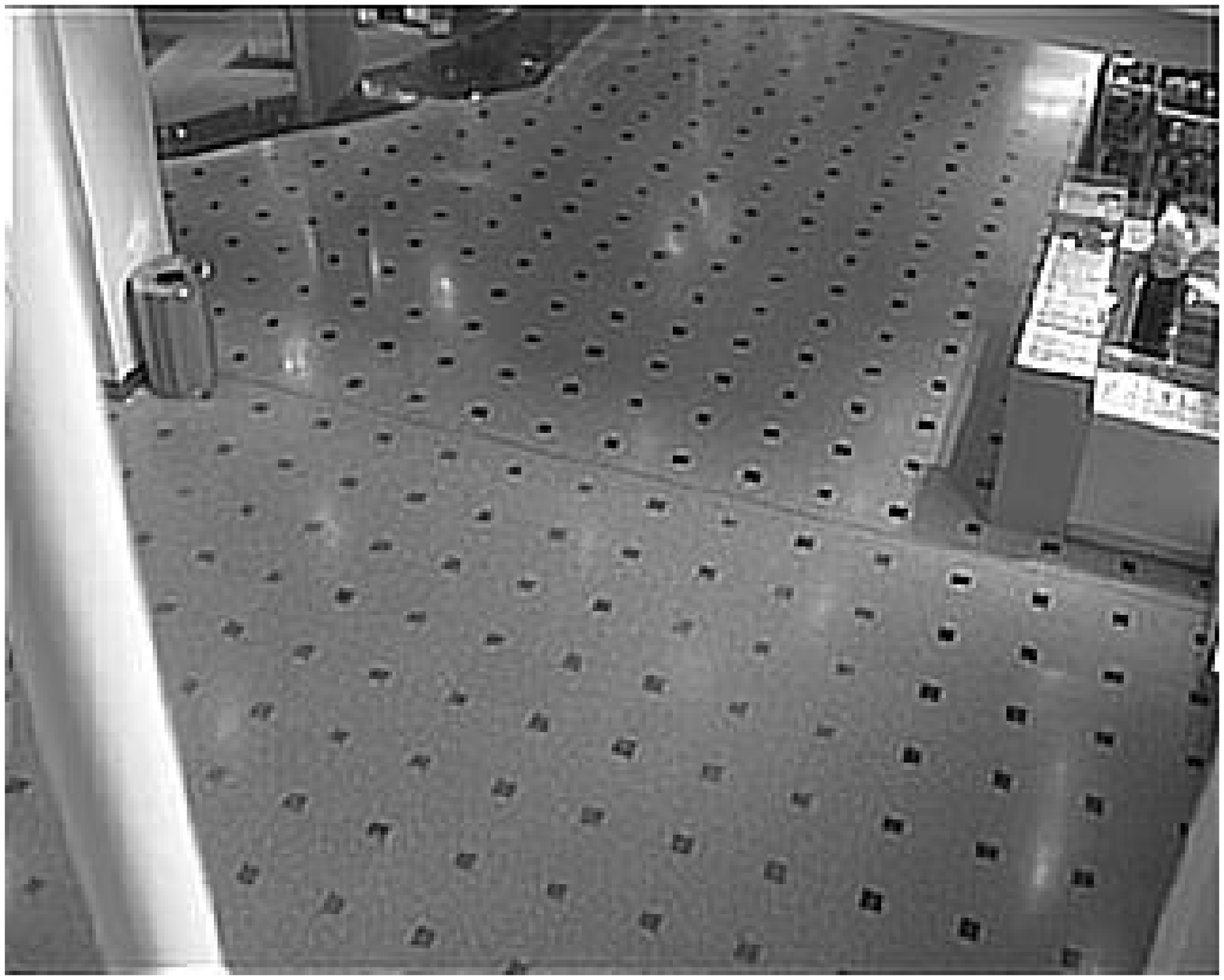}~\includegraphics[height=1.8cm,width=1.8cm]{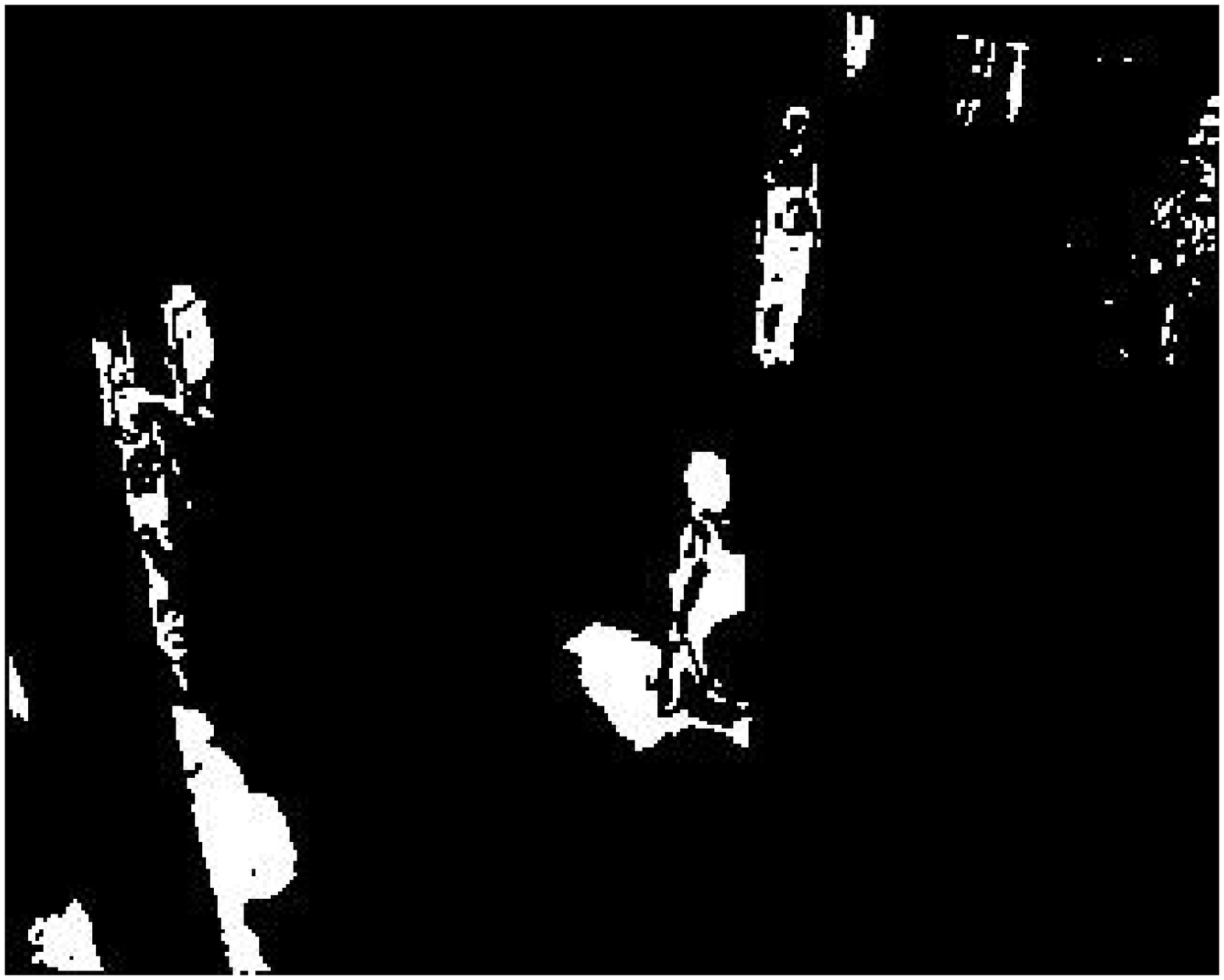}
\end{minipage}}
\subfigure{\begin{minipage}[t]{0.3\textwidth}\centering bri.$p=0.5$\vspace{1mm} \\
\includegraphics[height=1.8cm,width=1.8cm]{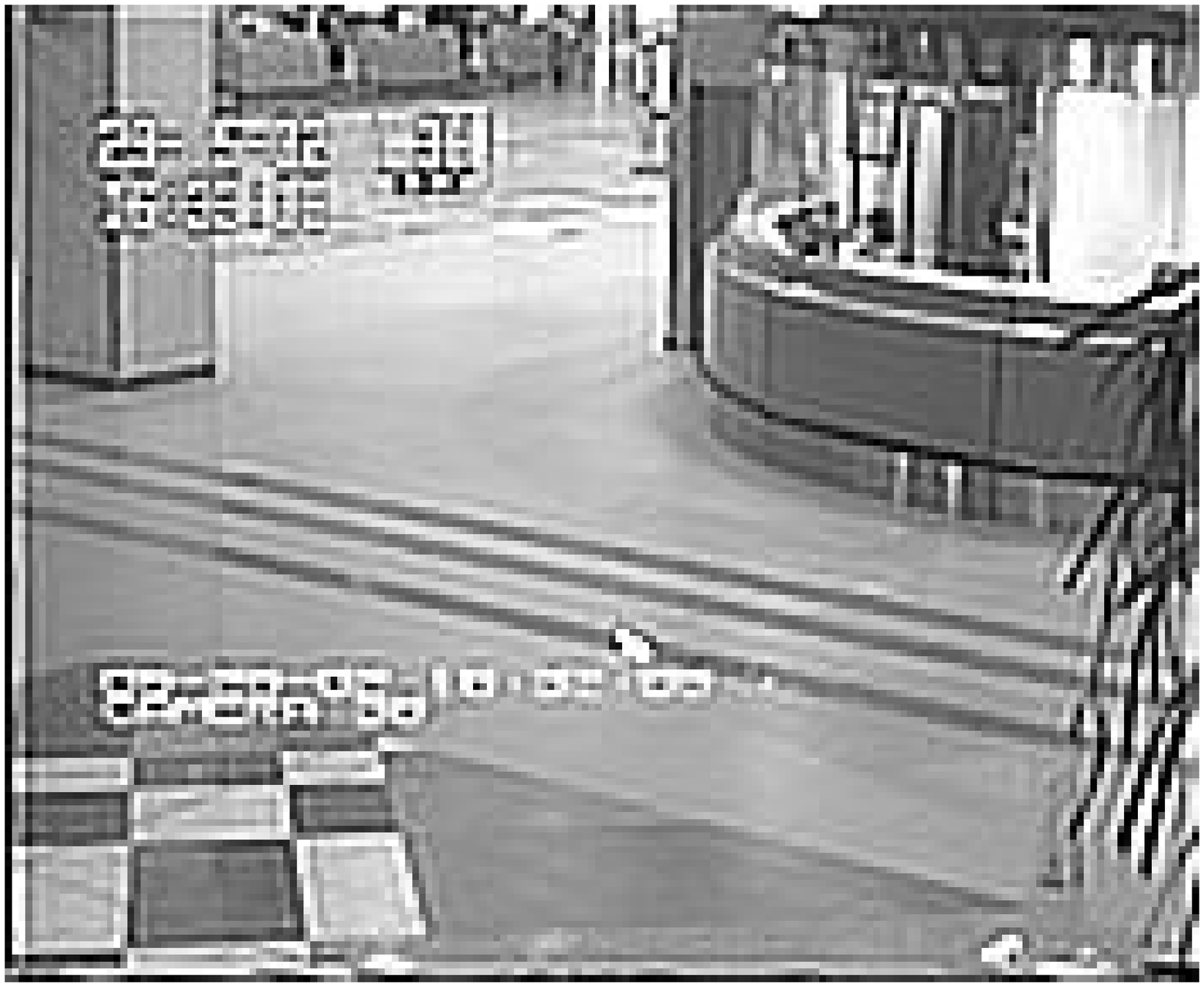}~\includegraphics[height=1.8cm,width=1.8cm]{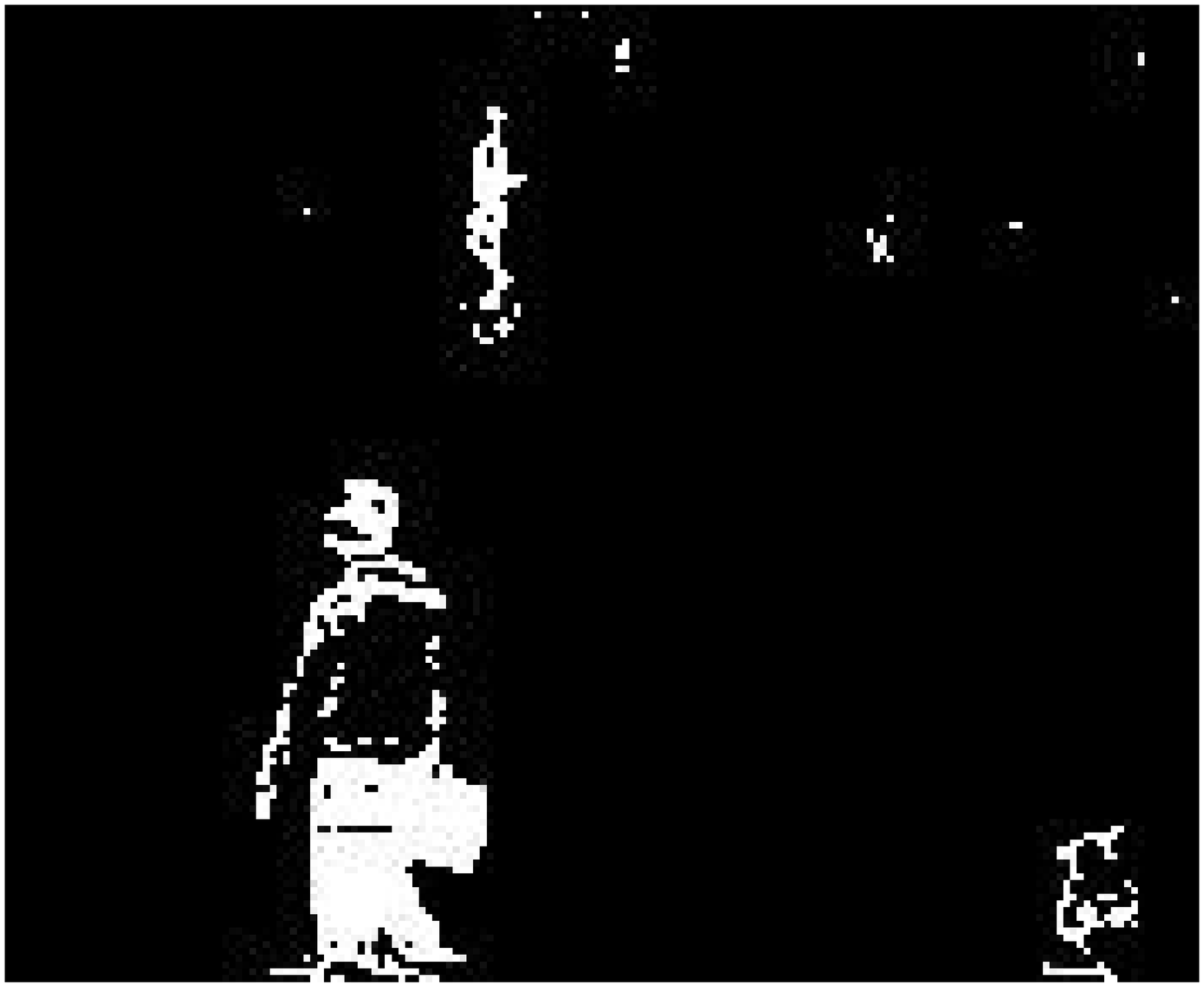}\vspace{0.5mm}\\
\includegraphics[height=1.8cm,width=1.8cm]{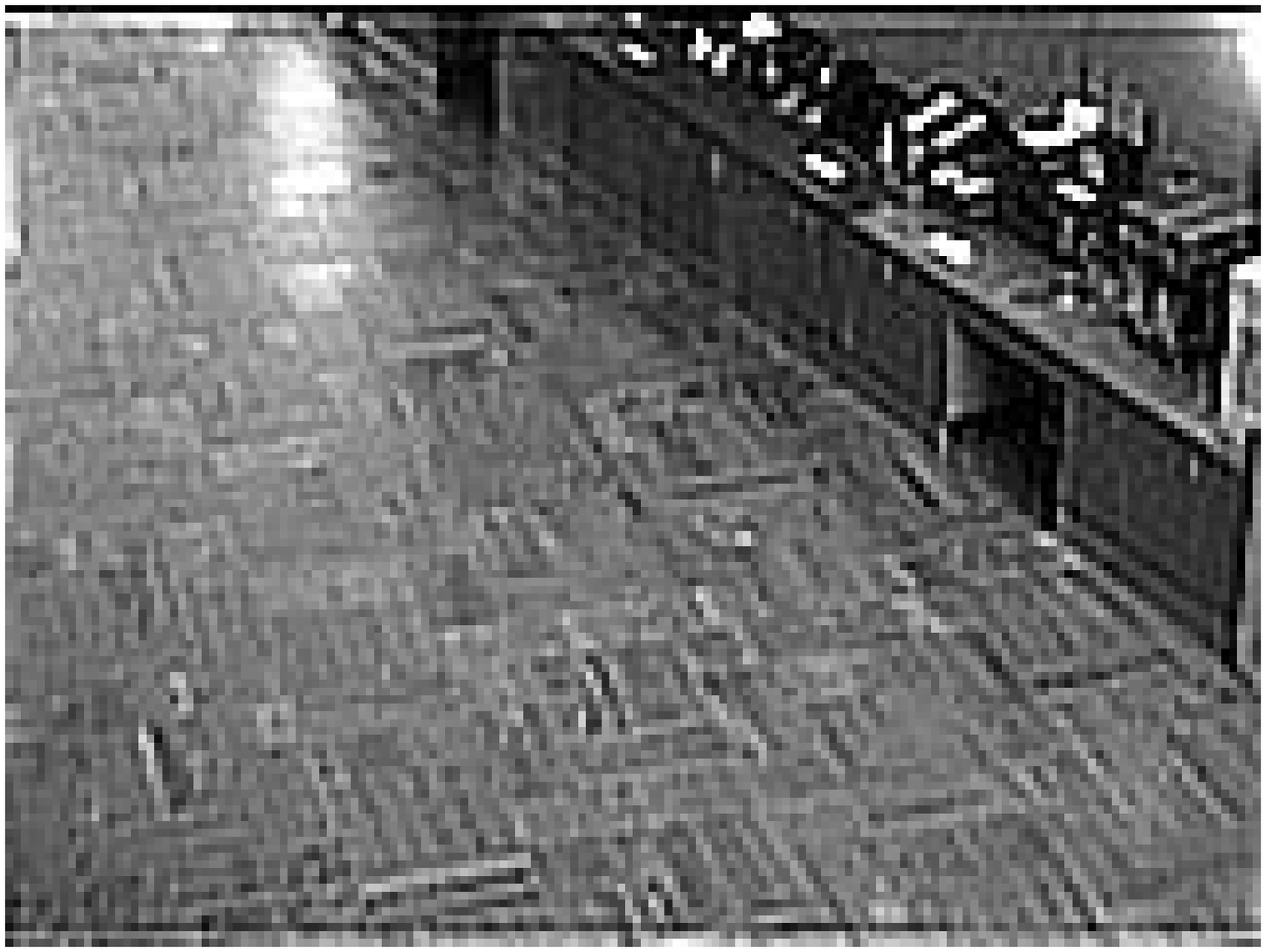}~\includegraphics[height=1.8cm,width=1.8cm]{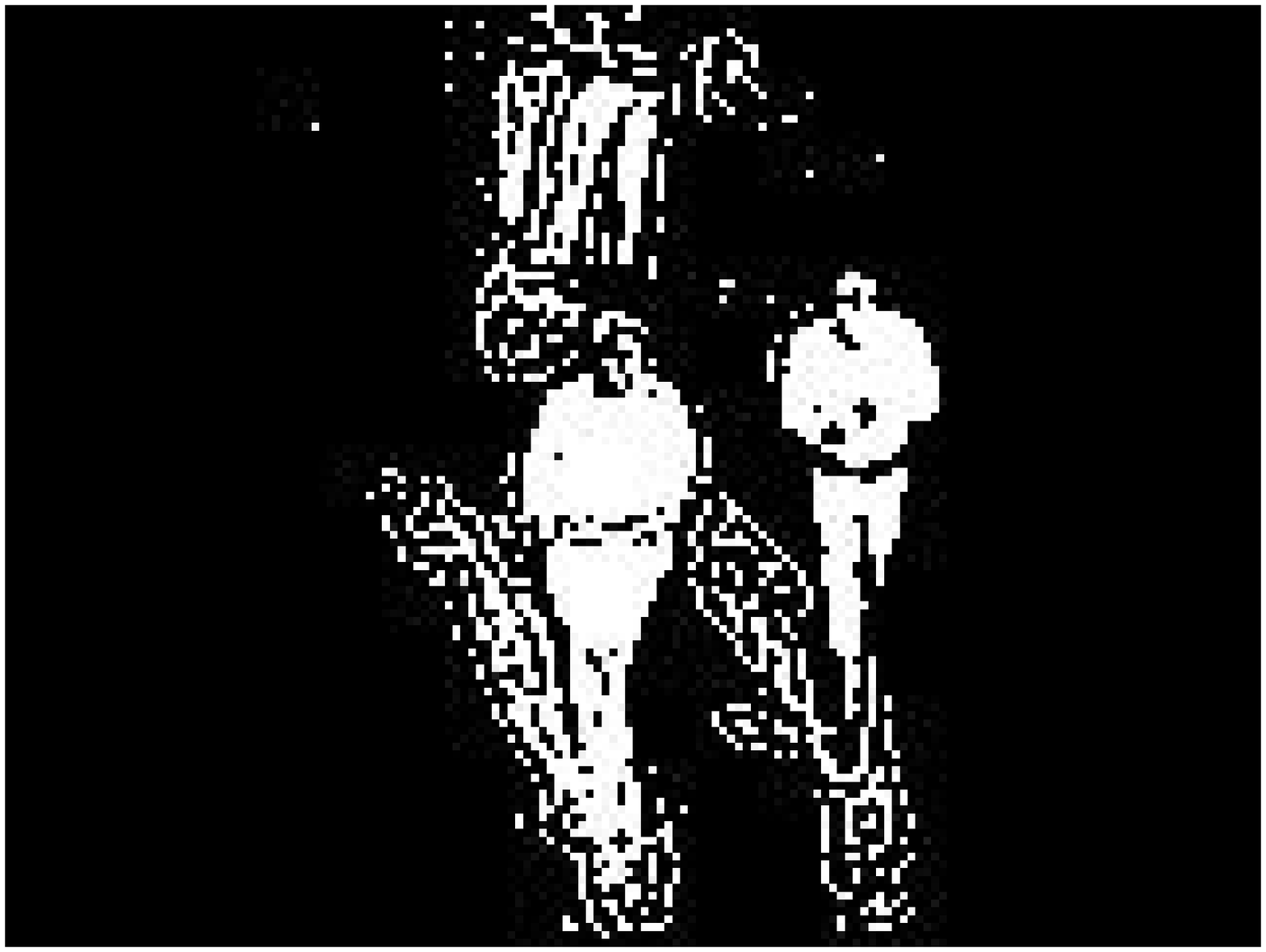}\vspace{0.5mm}\\
\includegraphics[height=1.8cm,width=1.8cm]{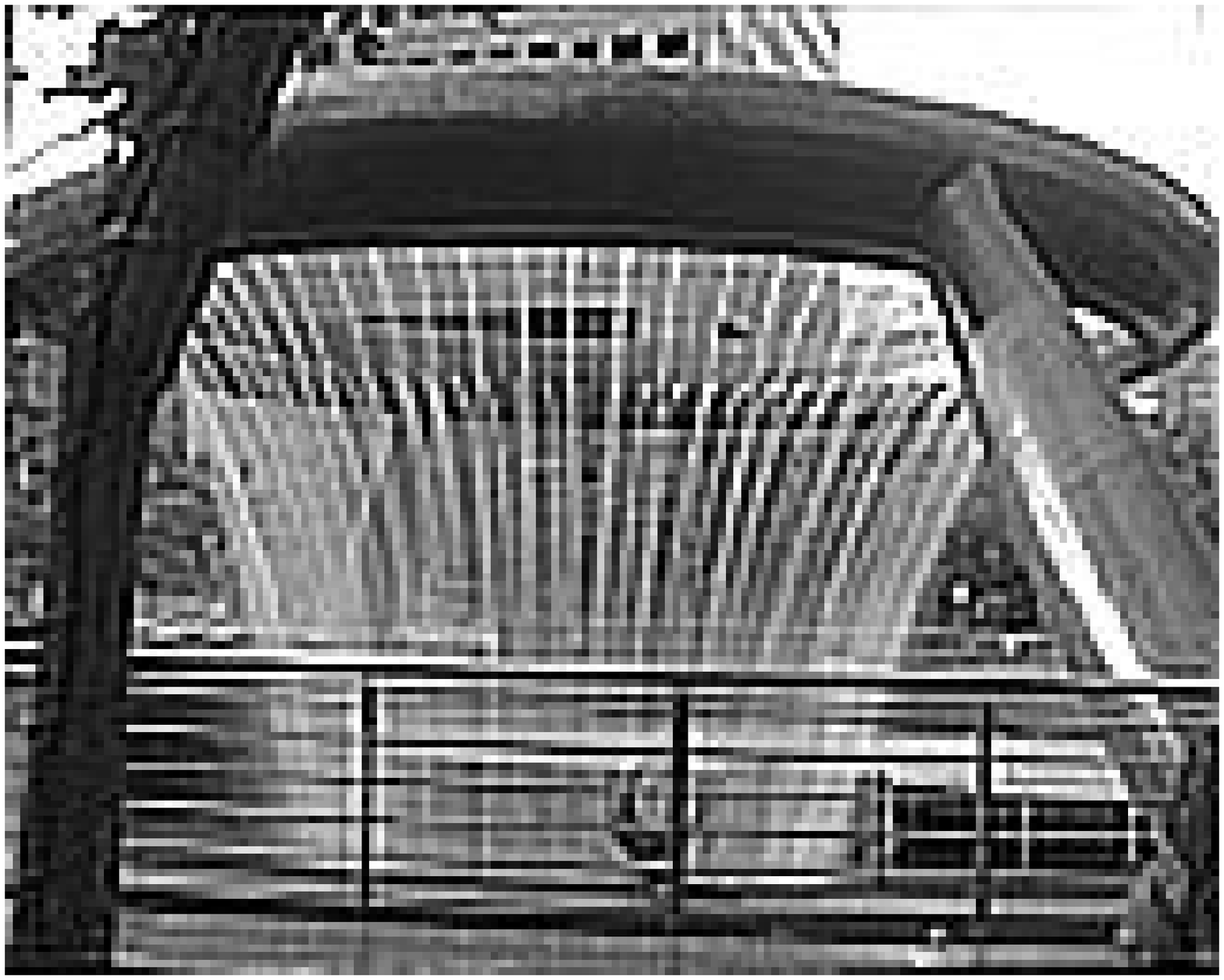}~\includegraphics[height=1.8cm,width=1.8cm]{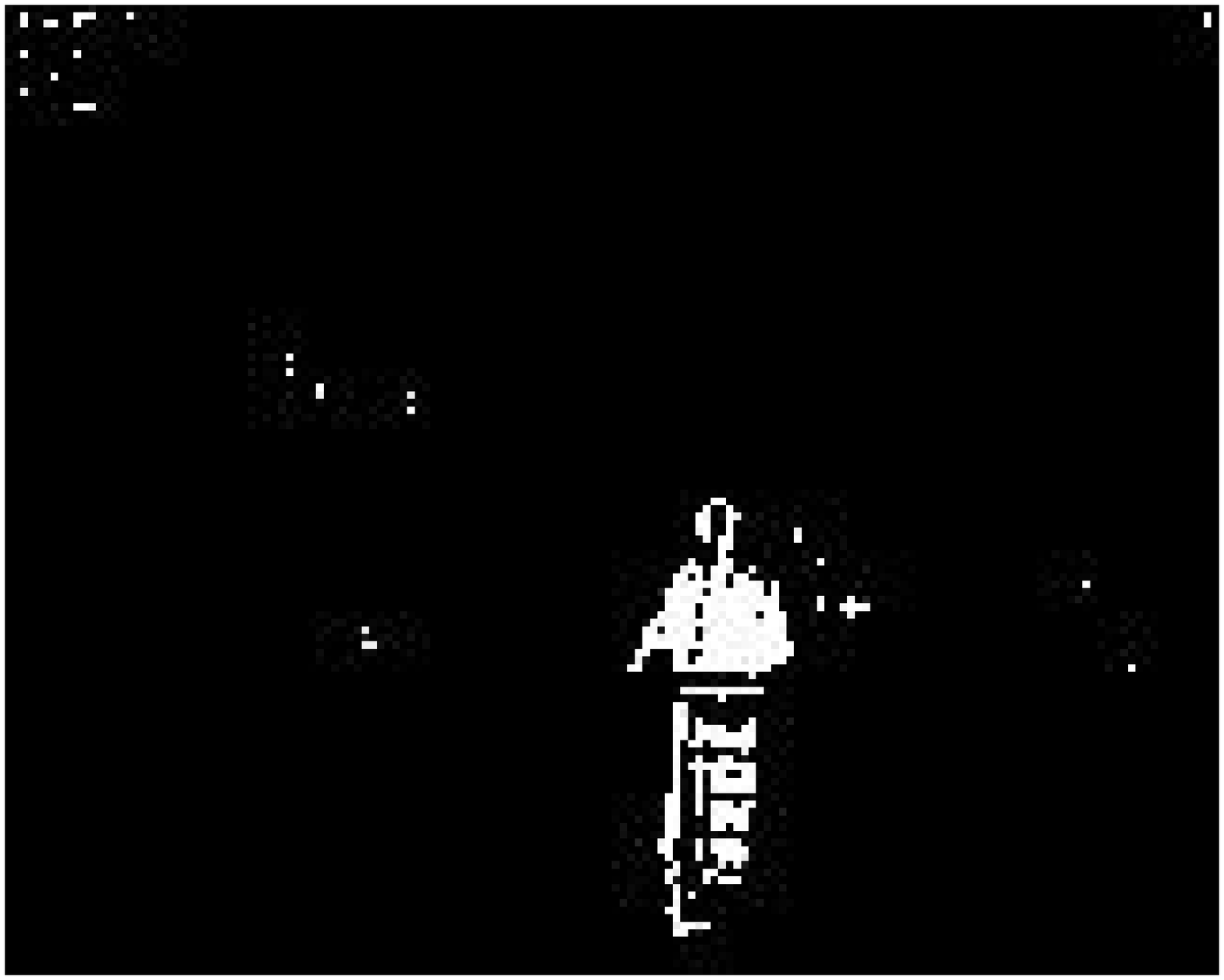}\vspace{0.5mm}\\
\includegraphics[height=1.8cm,width=1.8cm]{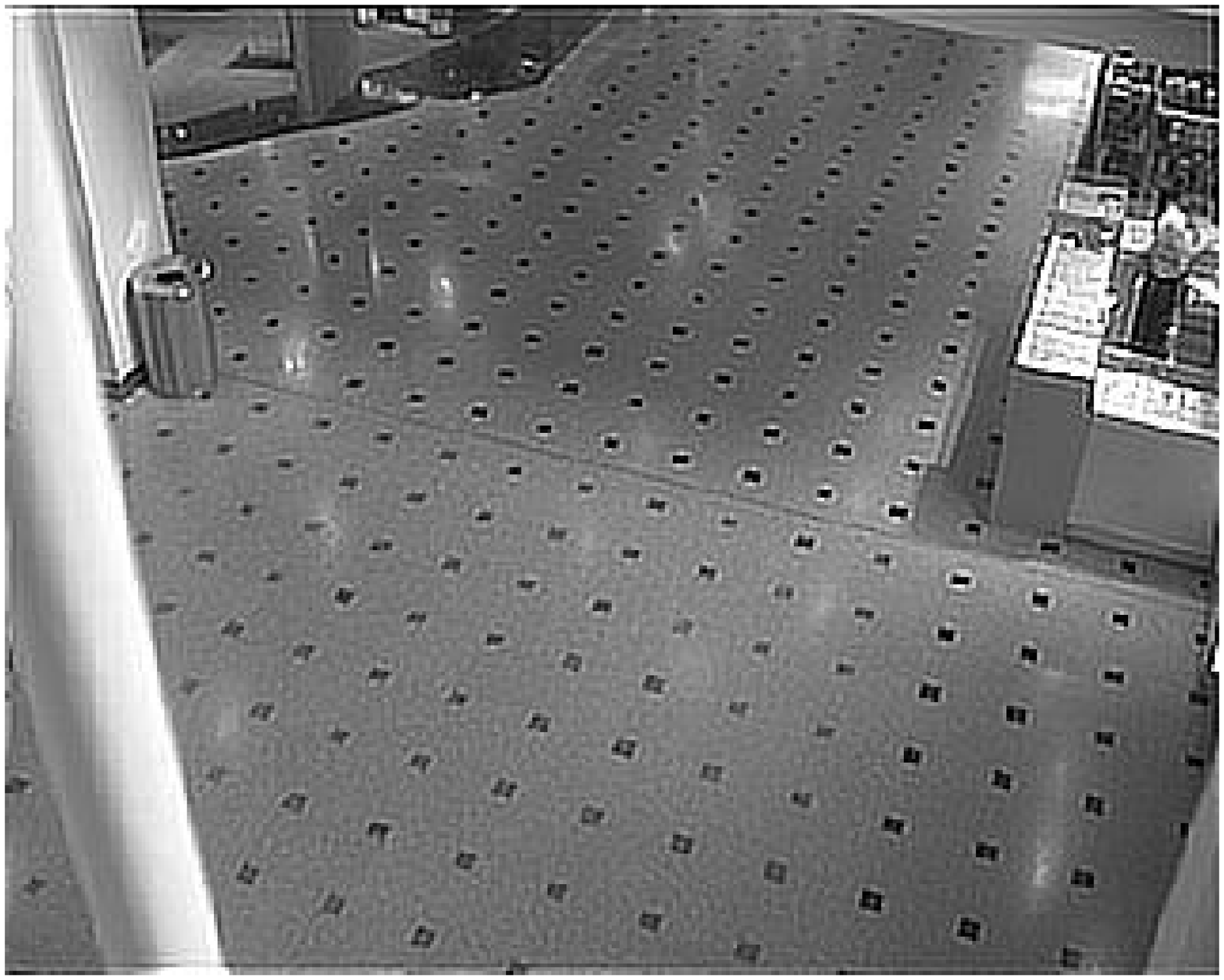}~\includegraphics[height=1.8cm,width=1.8cm]{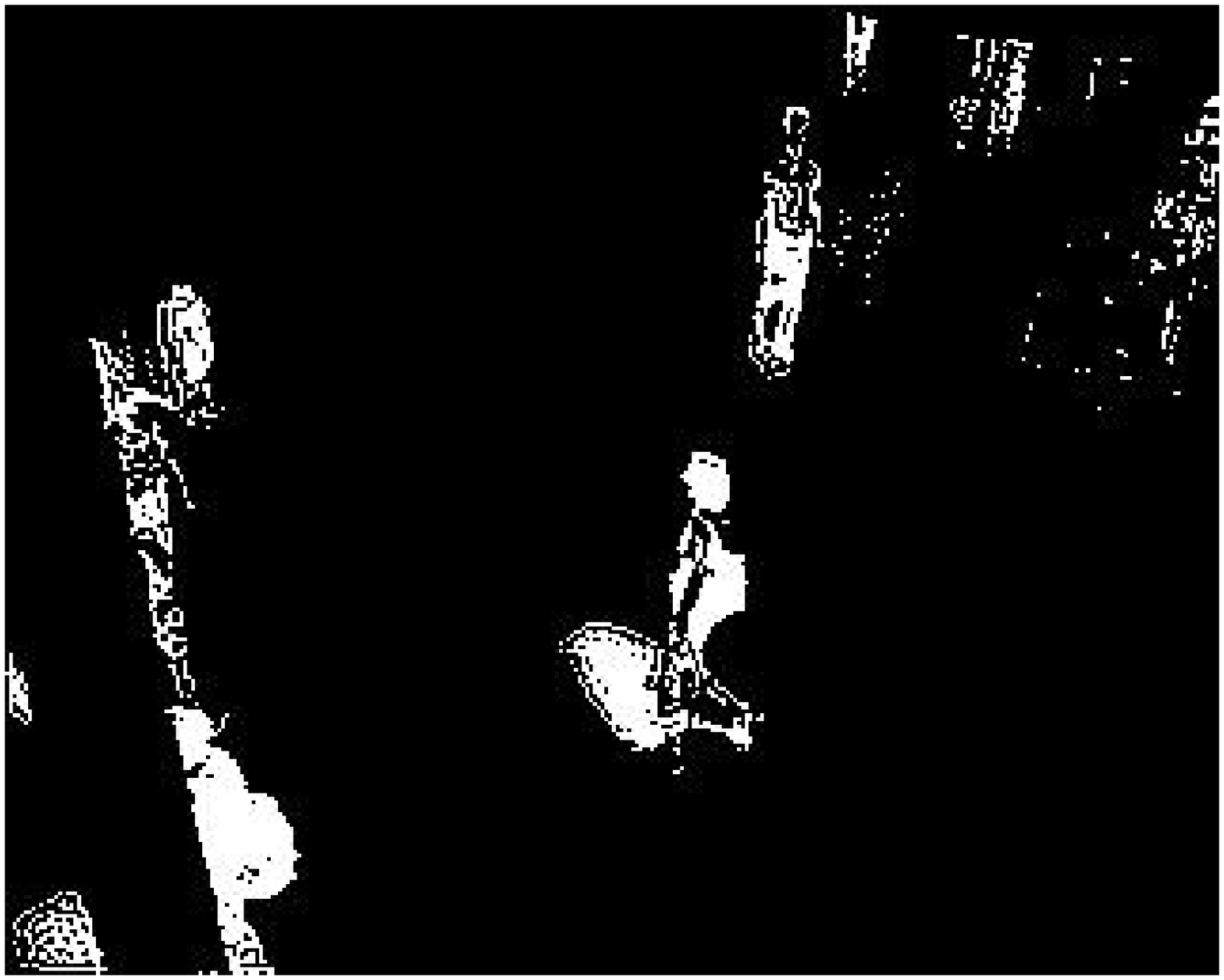}
\end{minipage}}
\subfigure{\begin{minipage}[t]{0.3\textwidth}\centering fra.$\alpha=1$\vspace{1.6mm} \\
\includegraphics[height=1.8cm,width=1.8cm]{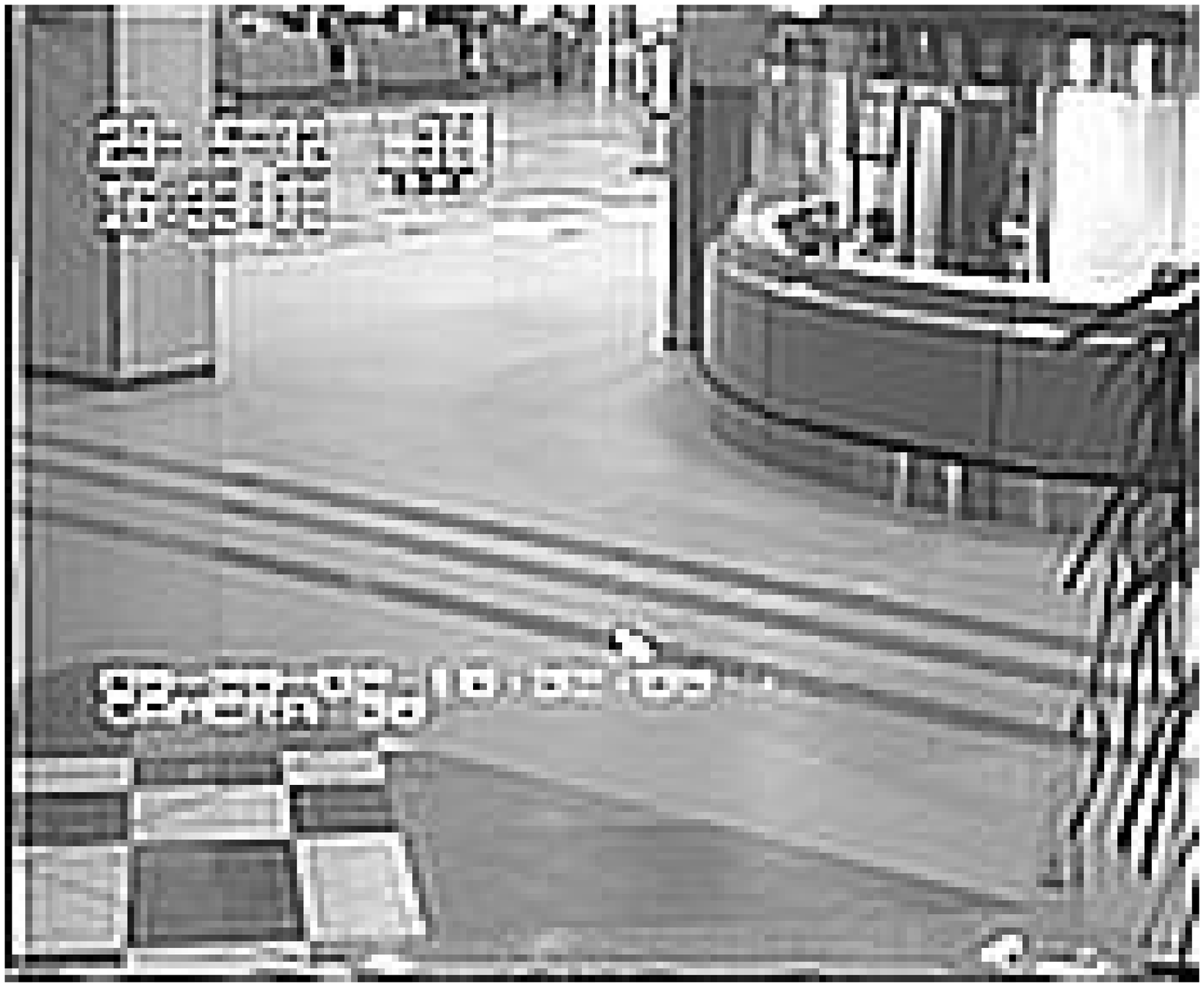}~\includegraphics[height=1.8cm,width=1.8cm]{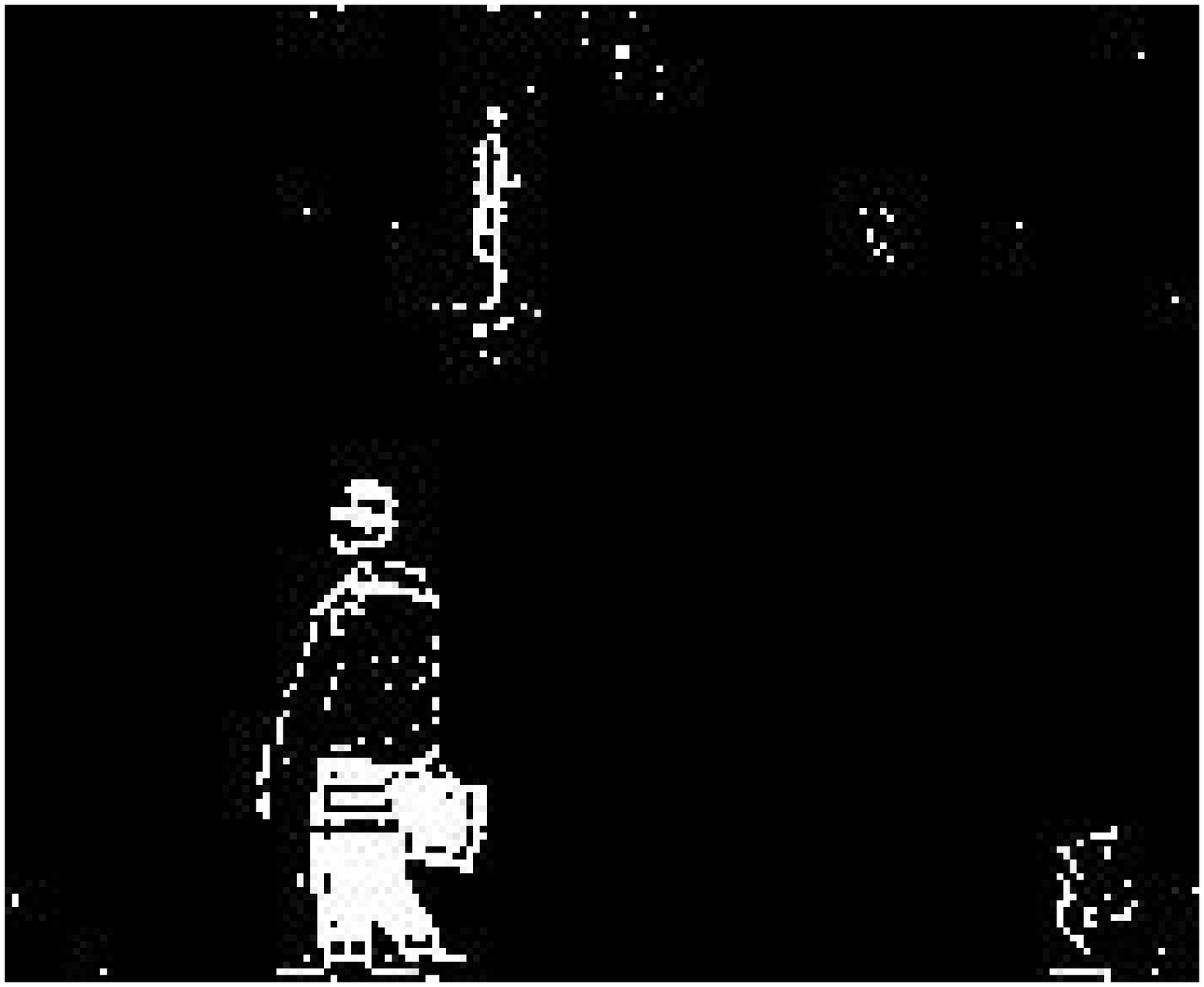}\vspace{0.5mm}\\
\includegraphics[height=1.8cm,width=1.8cm]{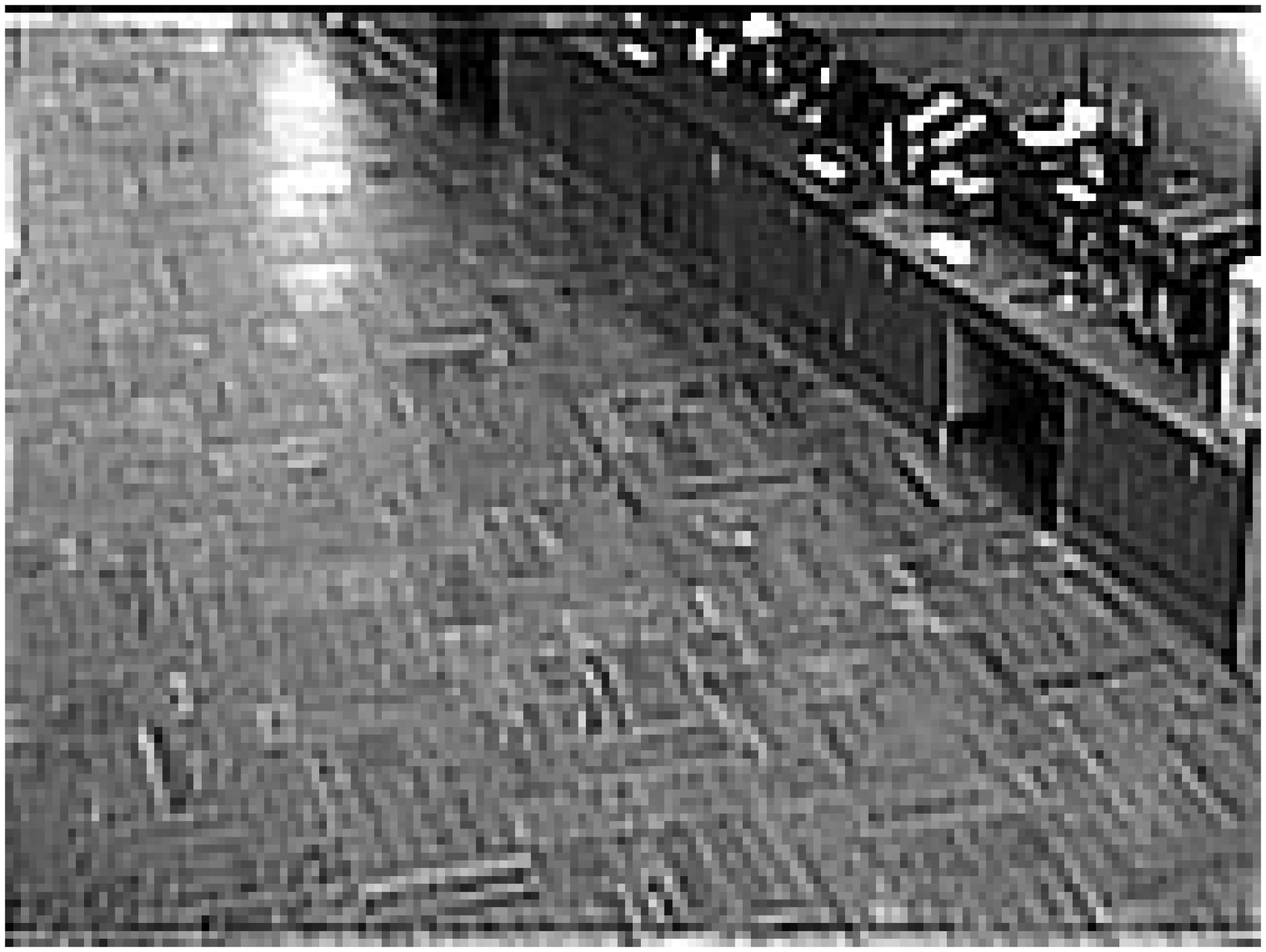}~\includegraphics[height=1.8cm,width=1.8cm]{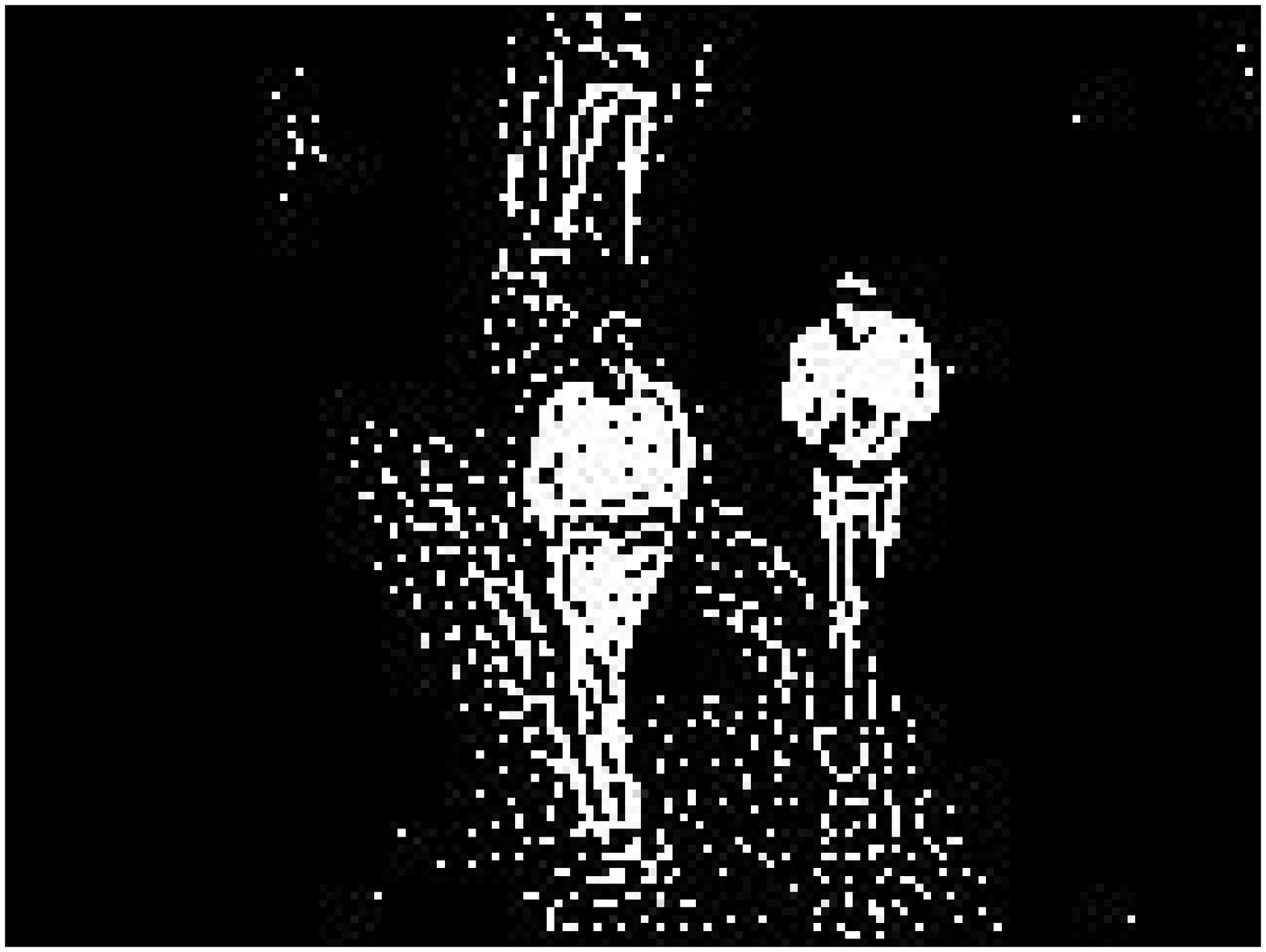}\vspace{0.5mm}\\
\includegraphics[height=1.8cm,width=1.8cm]{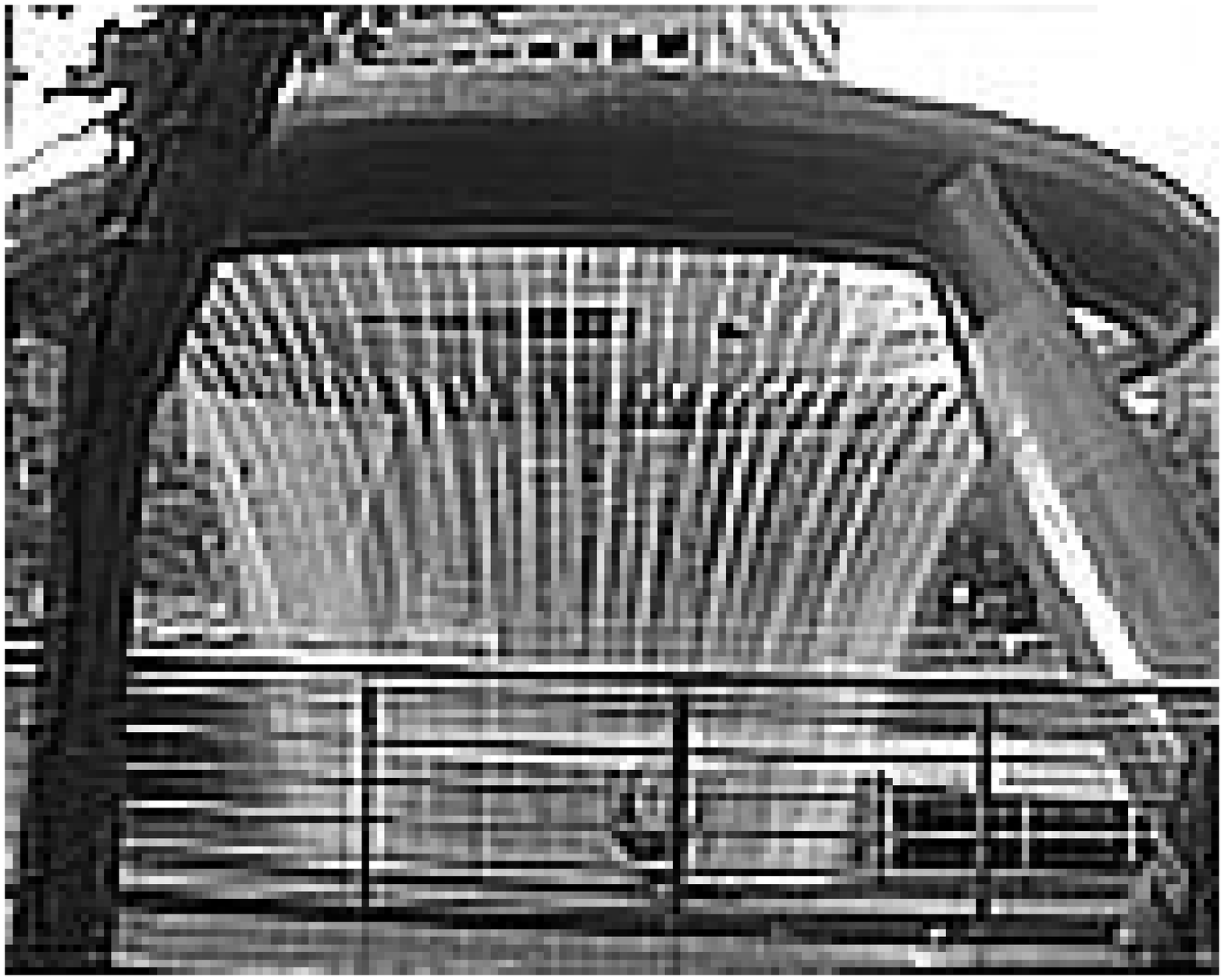}~\includegraphics[height=1.8cm,width=1.8cm]{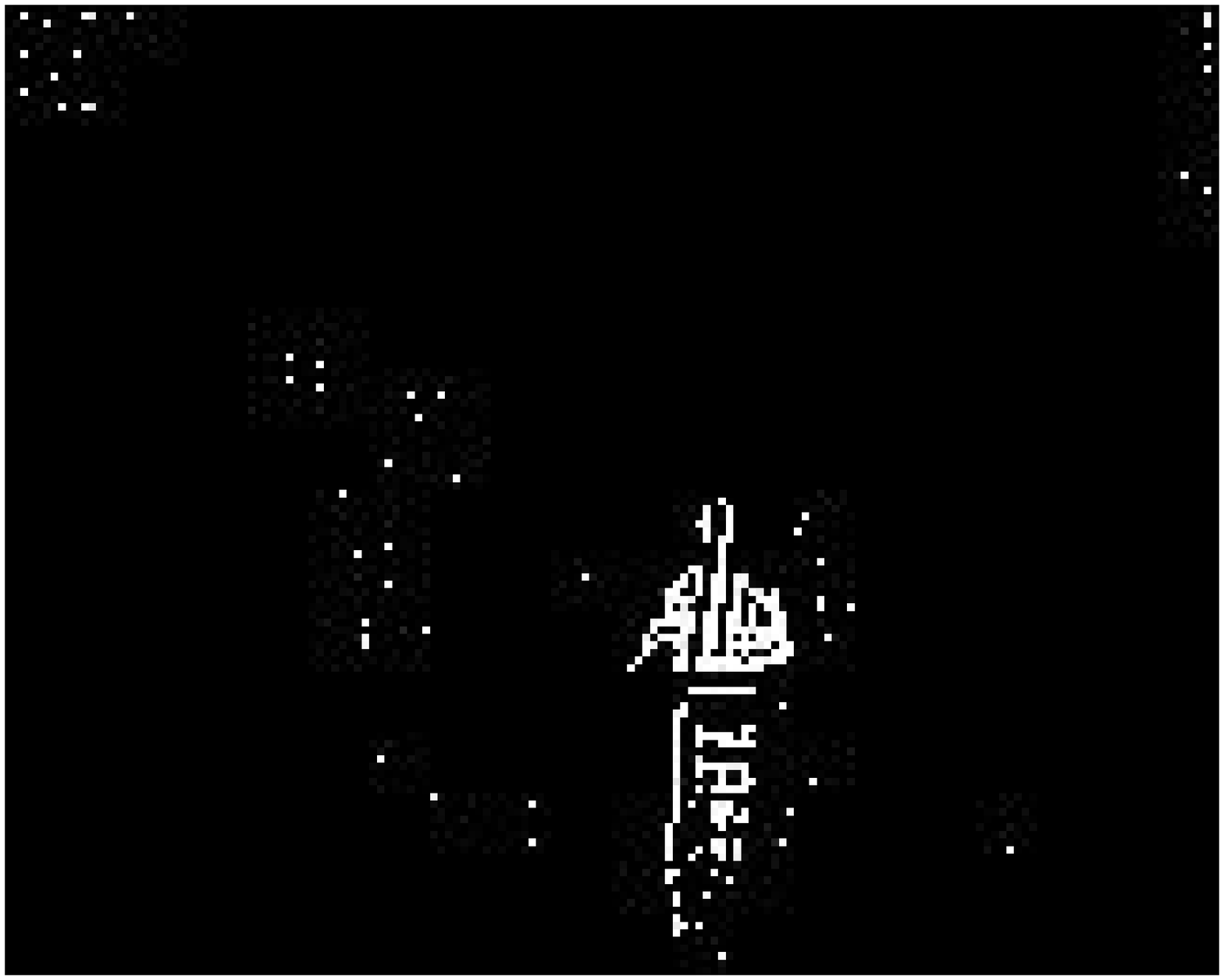}\vspace{0.5mm}\\
\includegraphics[height=1.8cm,width=1.8cm]{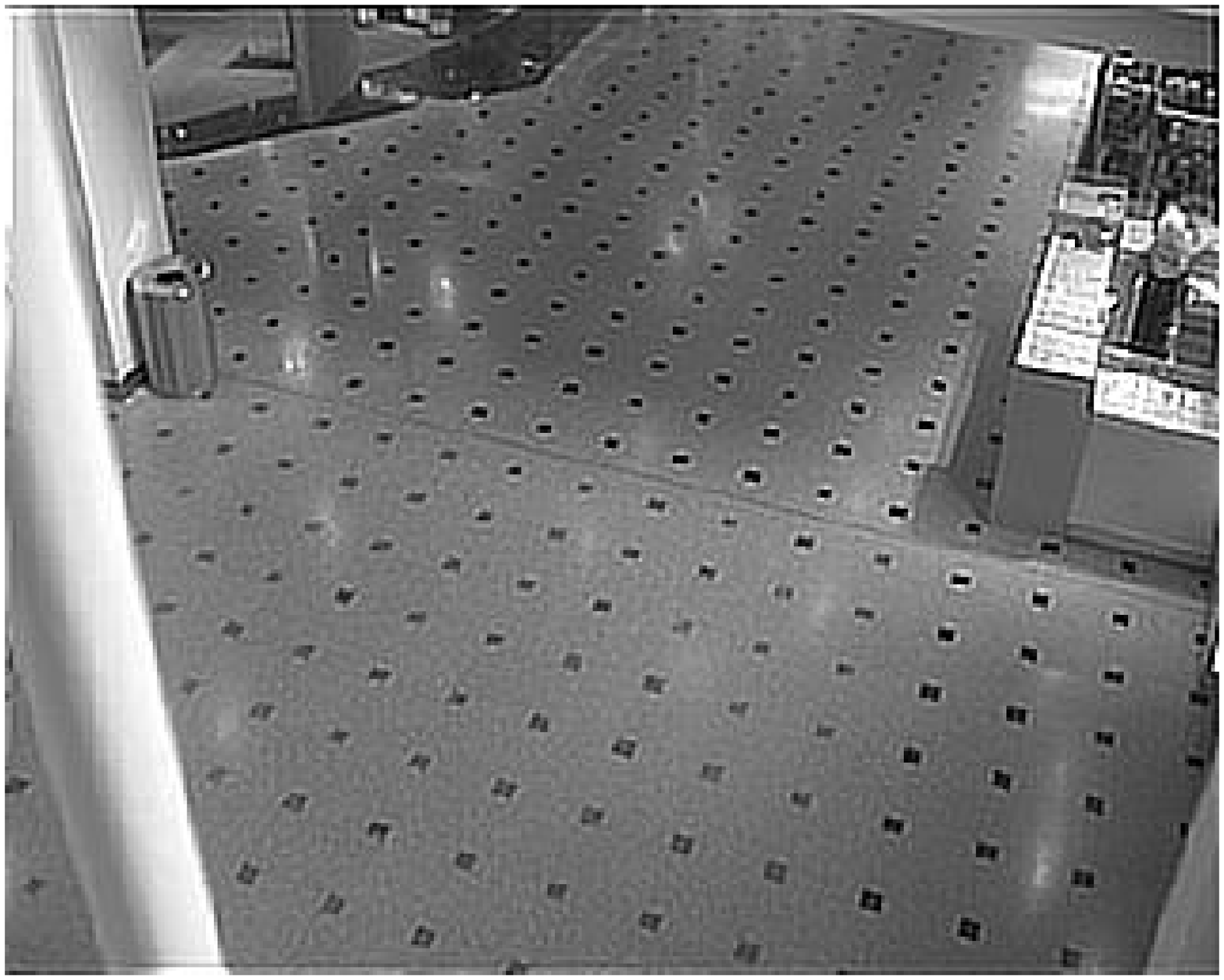}~\includegraphics[height=1.8cm,width=1.8cm]{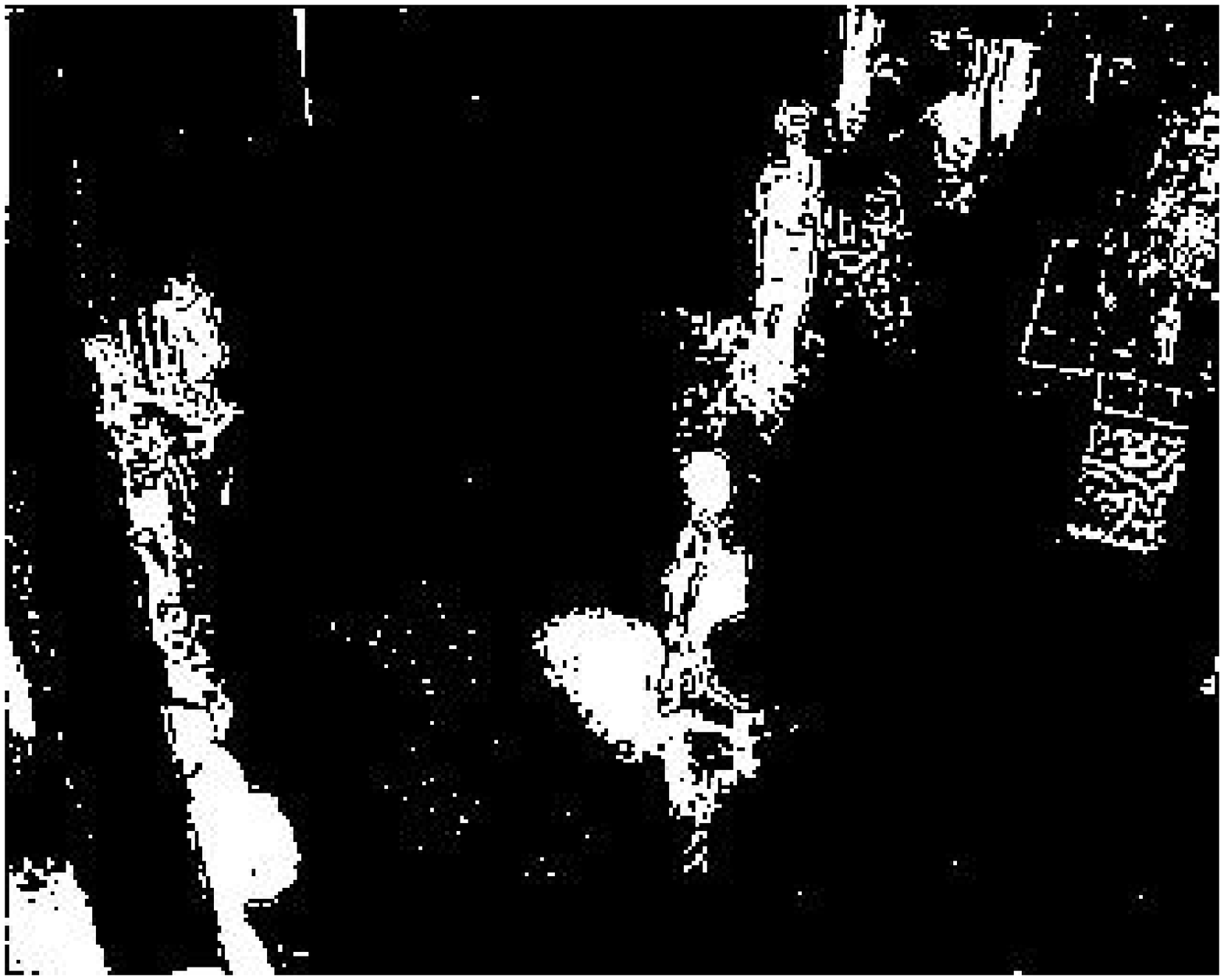}
\end{minipage}}  \\

\subfigure{\begin{minipage}[t]{0.3\textwidth}\centering fra.$\alpha=2$\vspace{1.6mm} \\
\includegraphics[height=1.8cm,width=1.8cm]{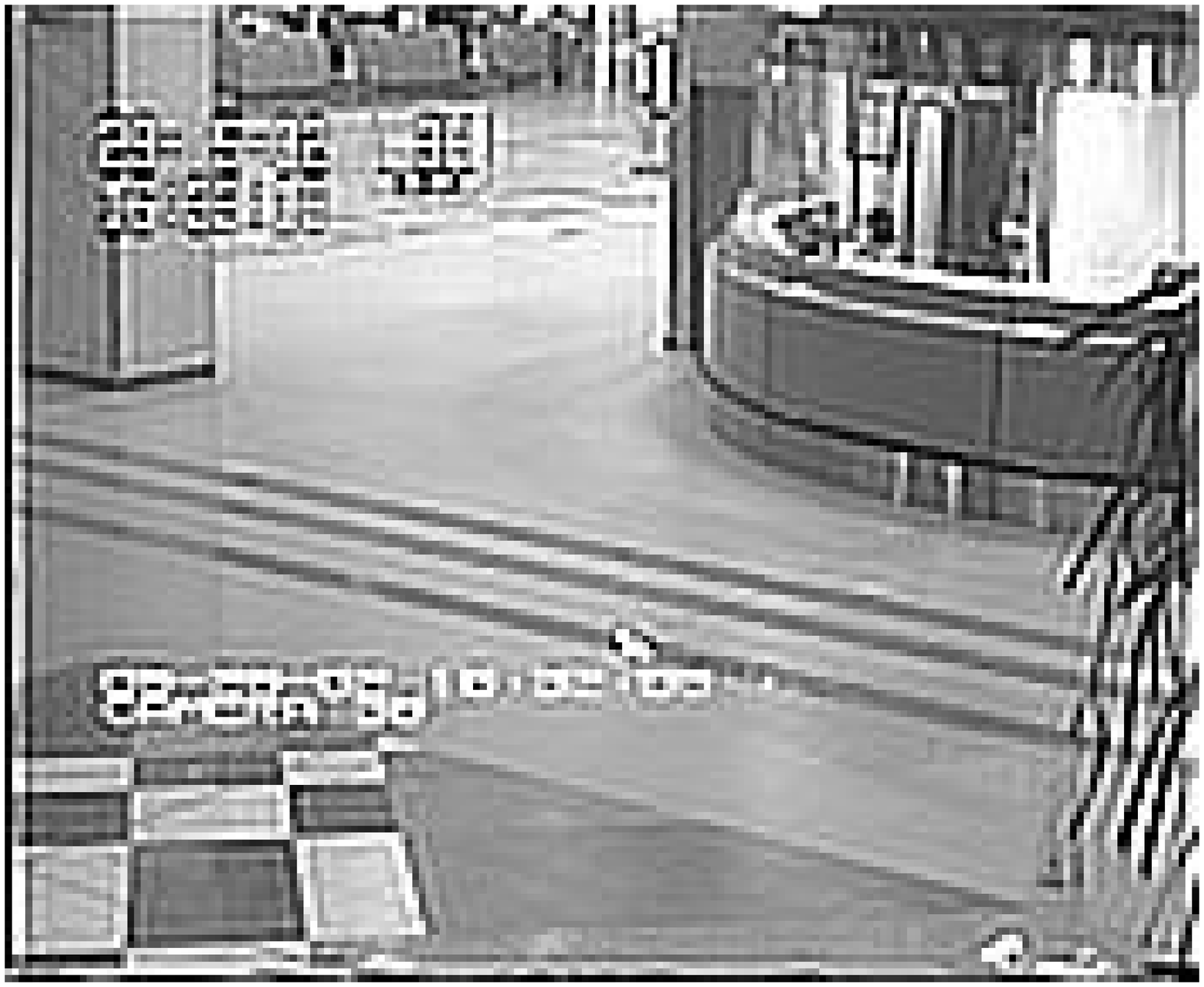}~\includegraphics[height=1.8cm,width=1.8cm]{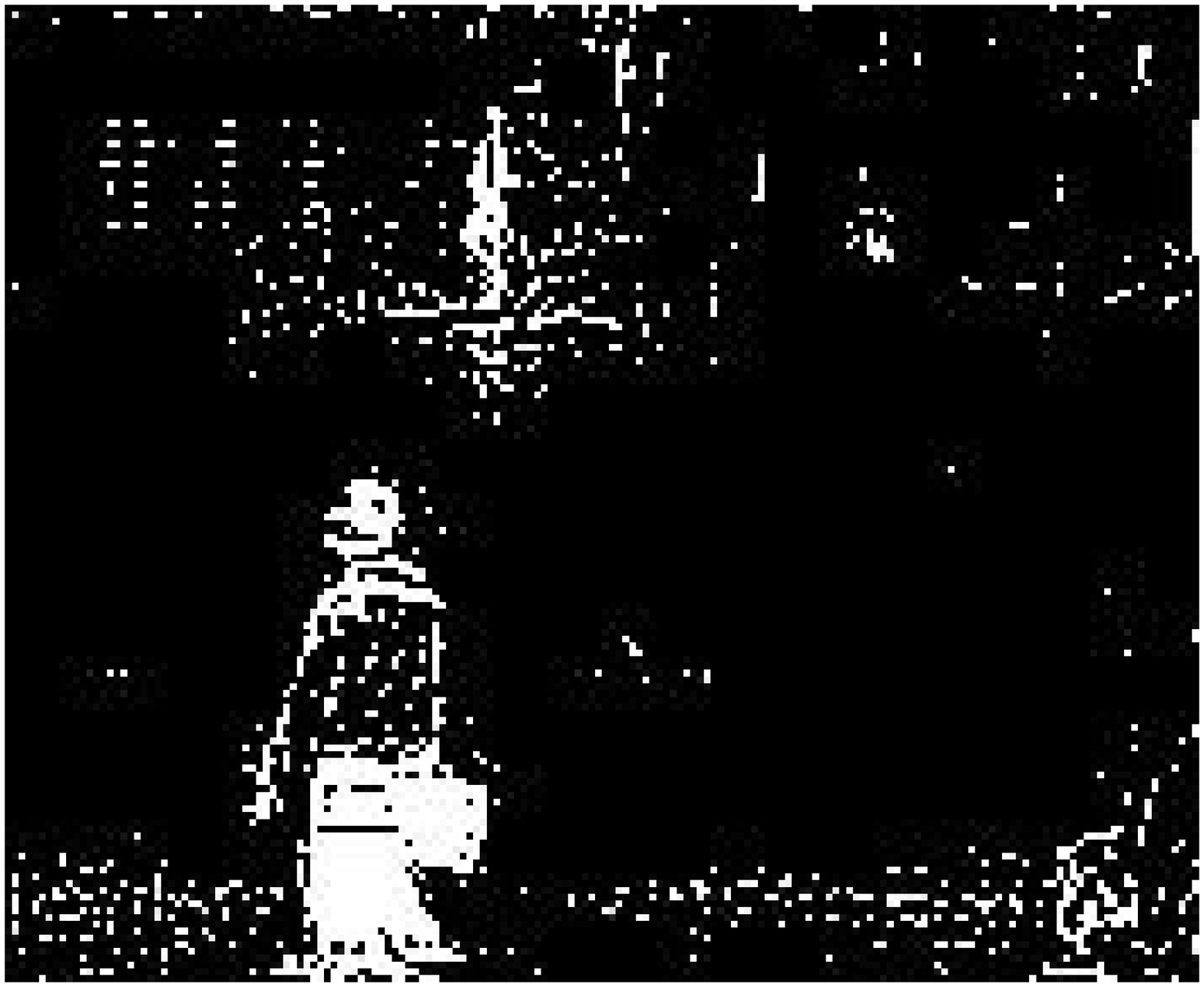}\vspace{0.5mm}\\
\includegraphics[height=1.8cm,width=1.8cm]{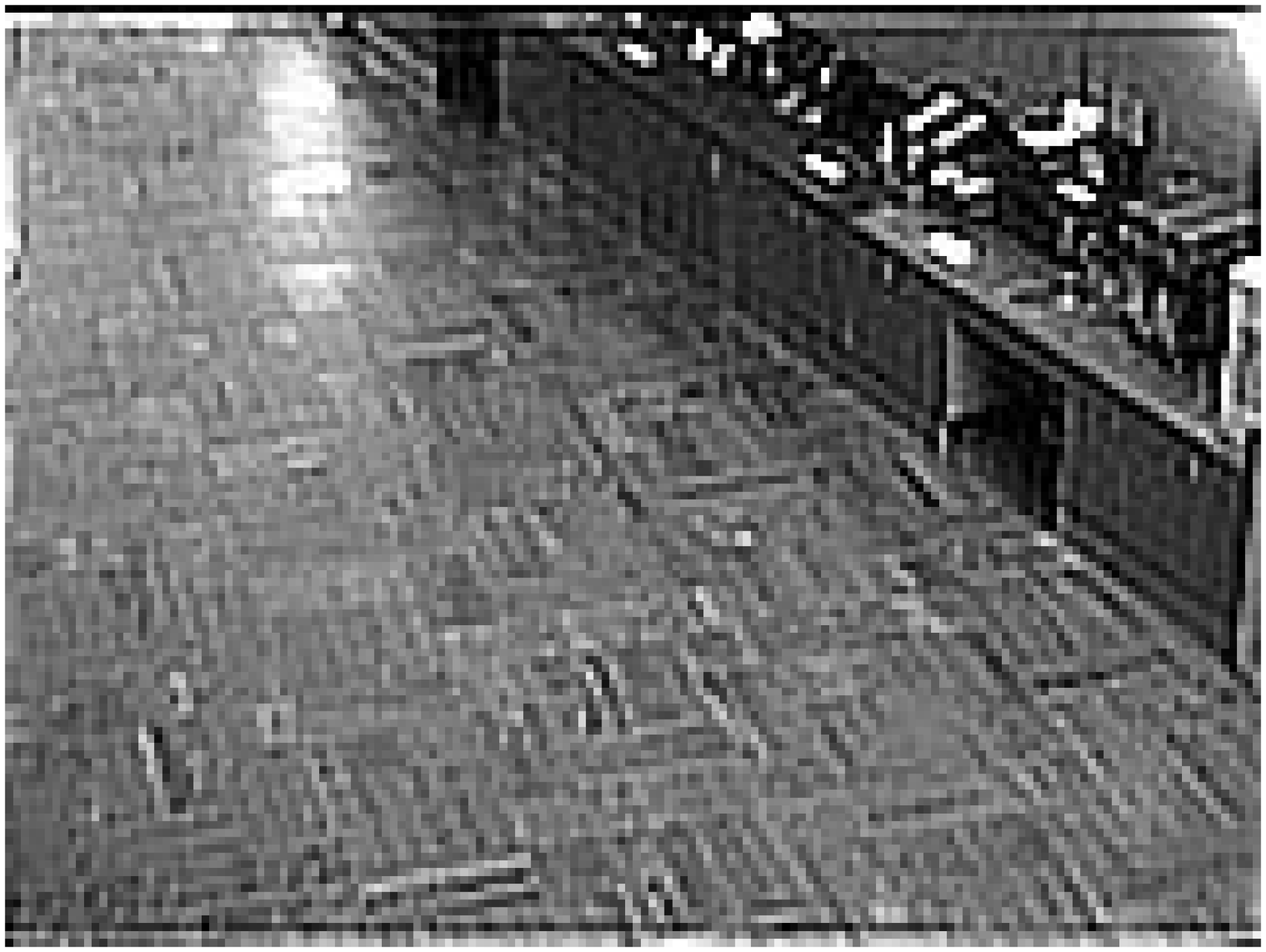}~\includegraphics[height=1.8cm,width=1.8cm]{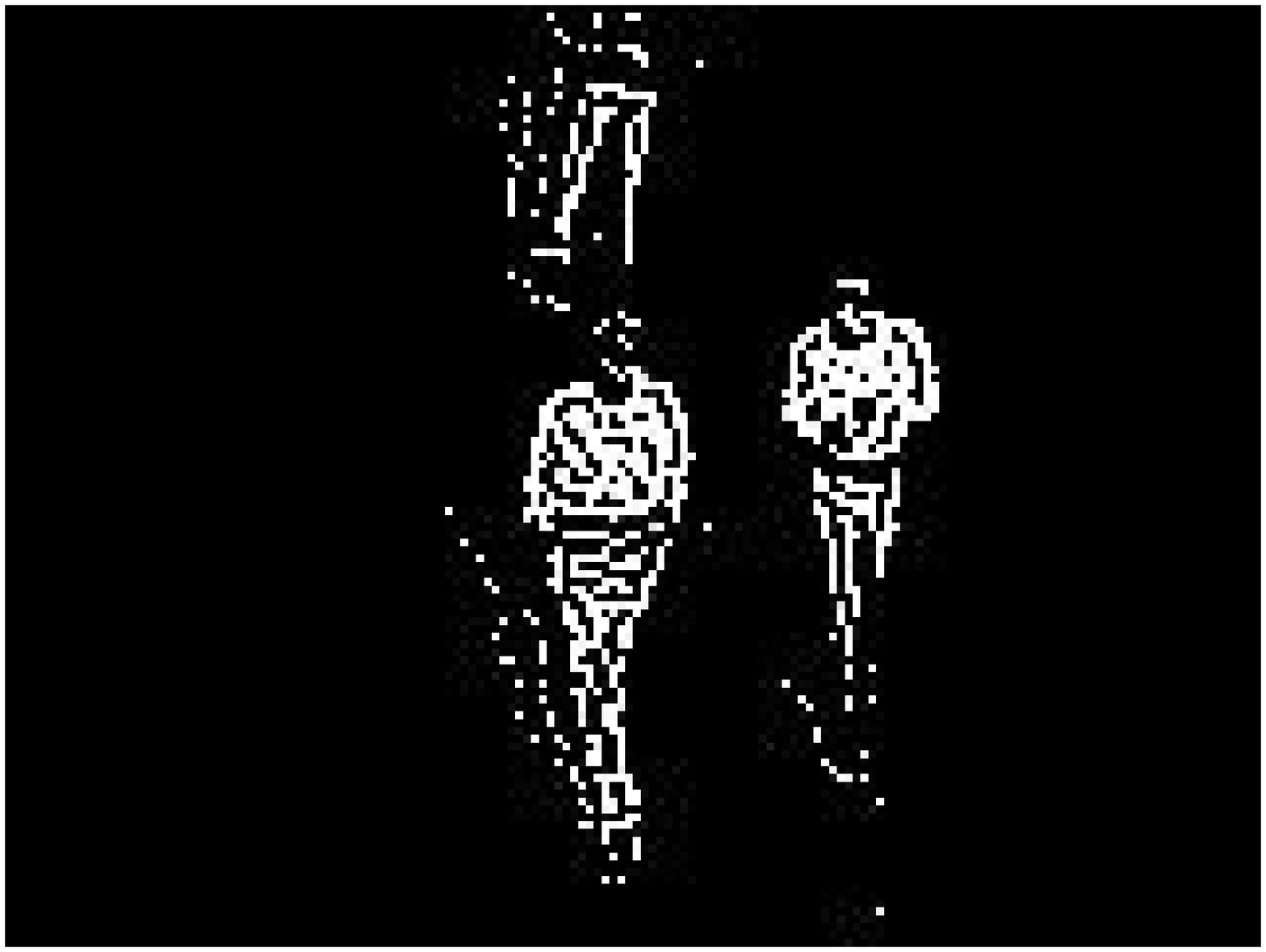}\vspace{0.5mm}\\
\includegraphics[height=1.8cm,width=1.8cm]{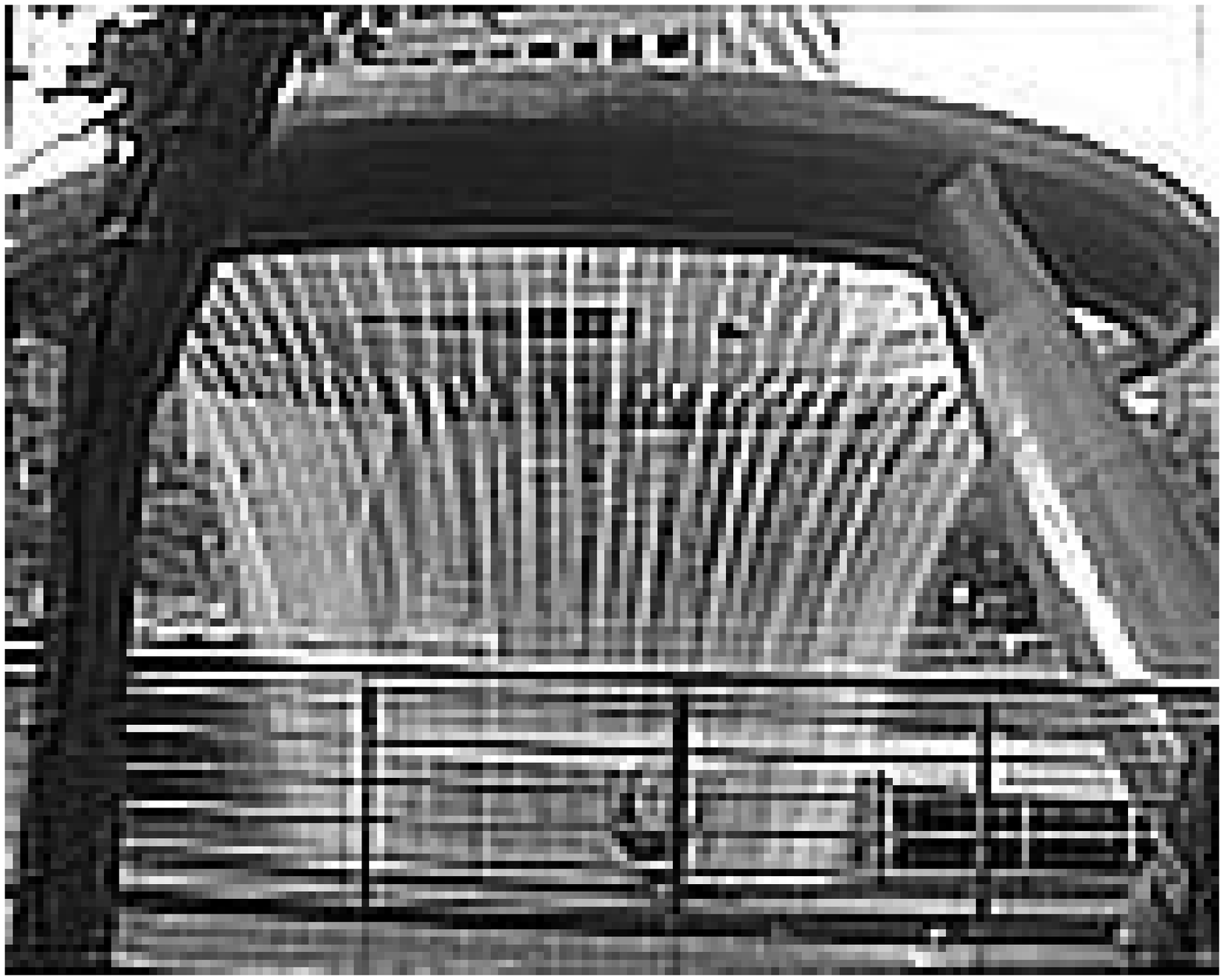}~\includegraphics[height=1.8cm,width=1.8cm]{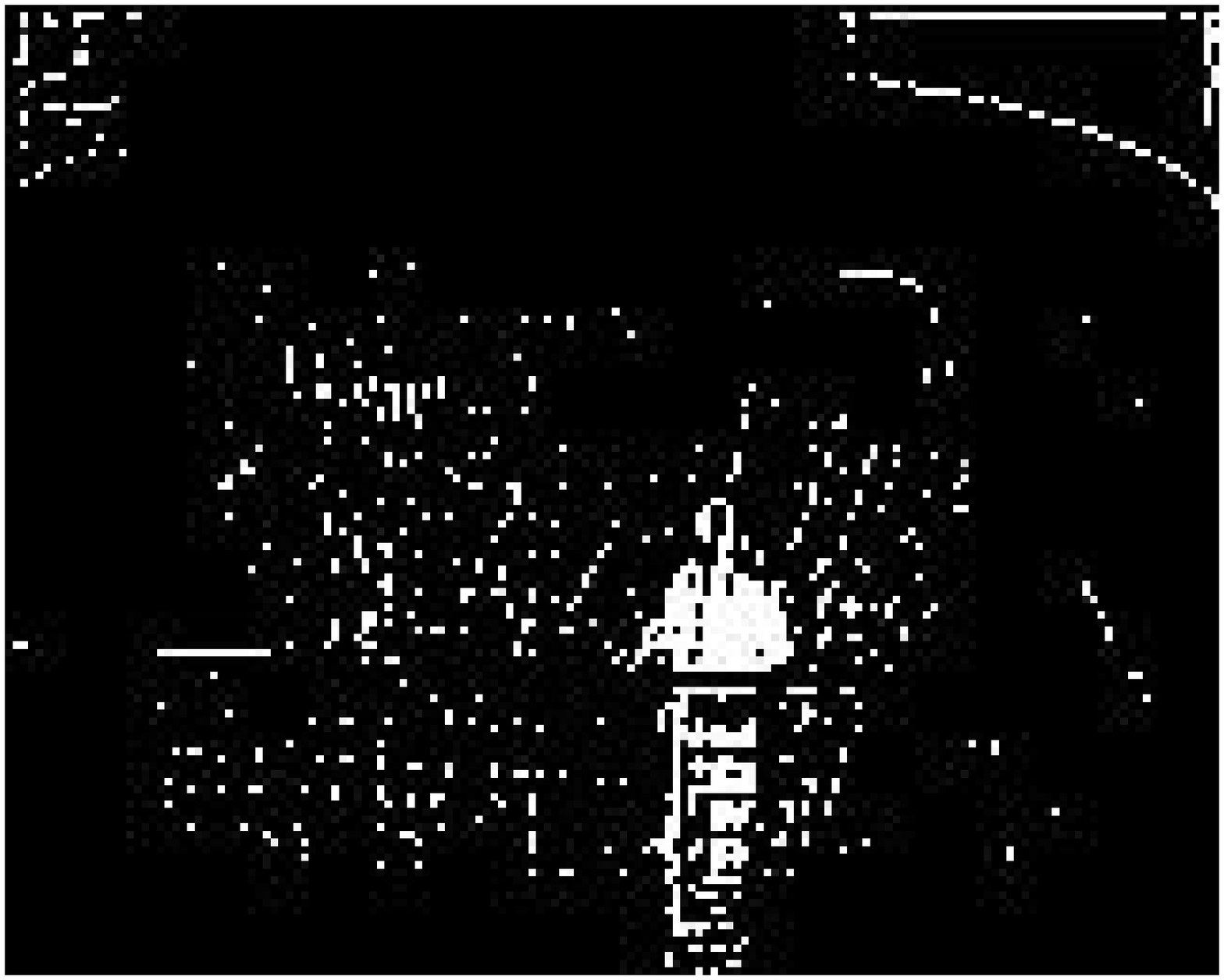}\vspace{0.5mm}\\
\includegraphics[height=1.8cm,width=1.8cm]{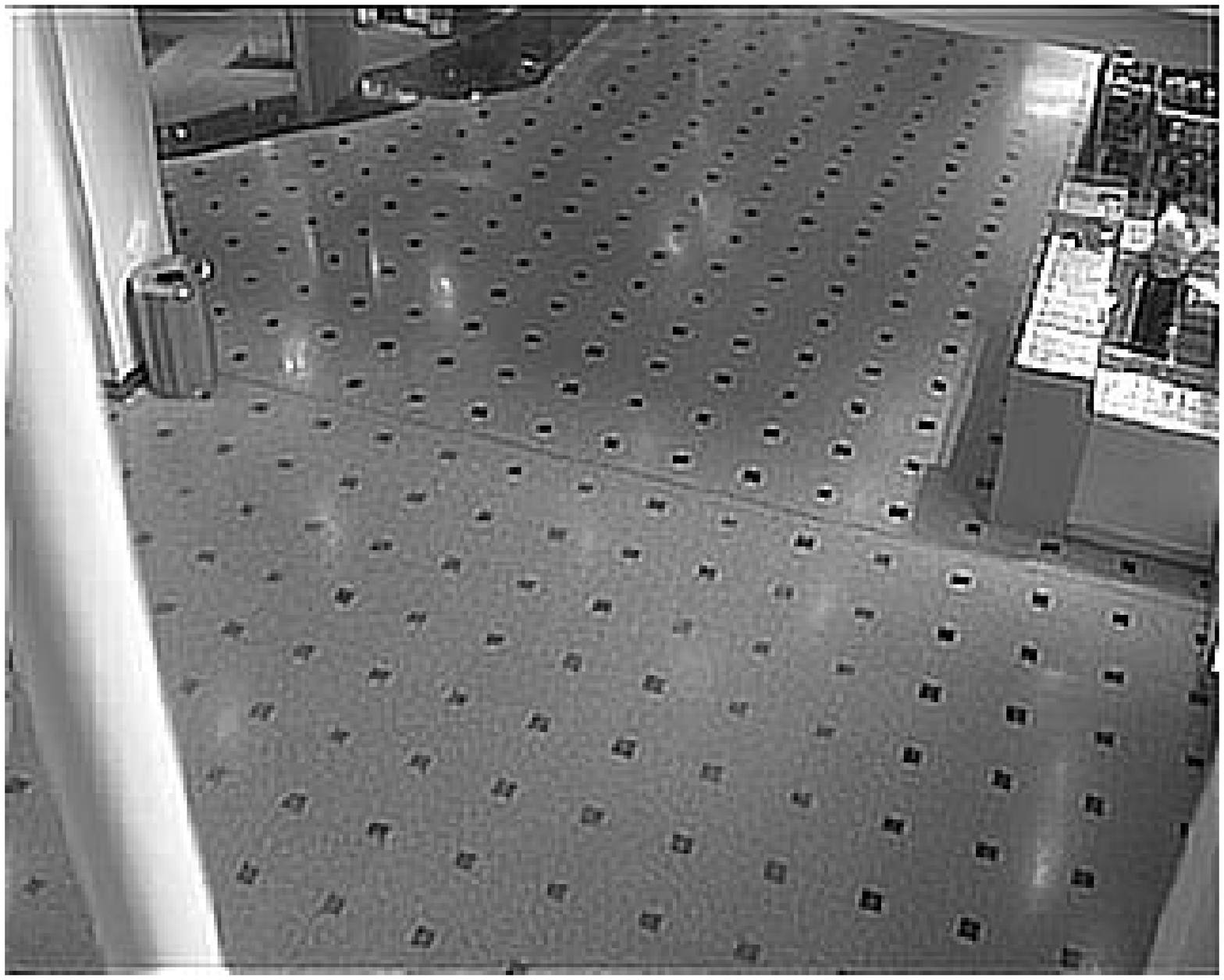}~\includegraphics[height=1.8cm,width=1.8cm]{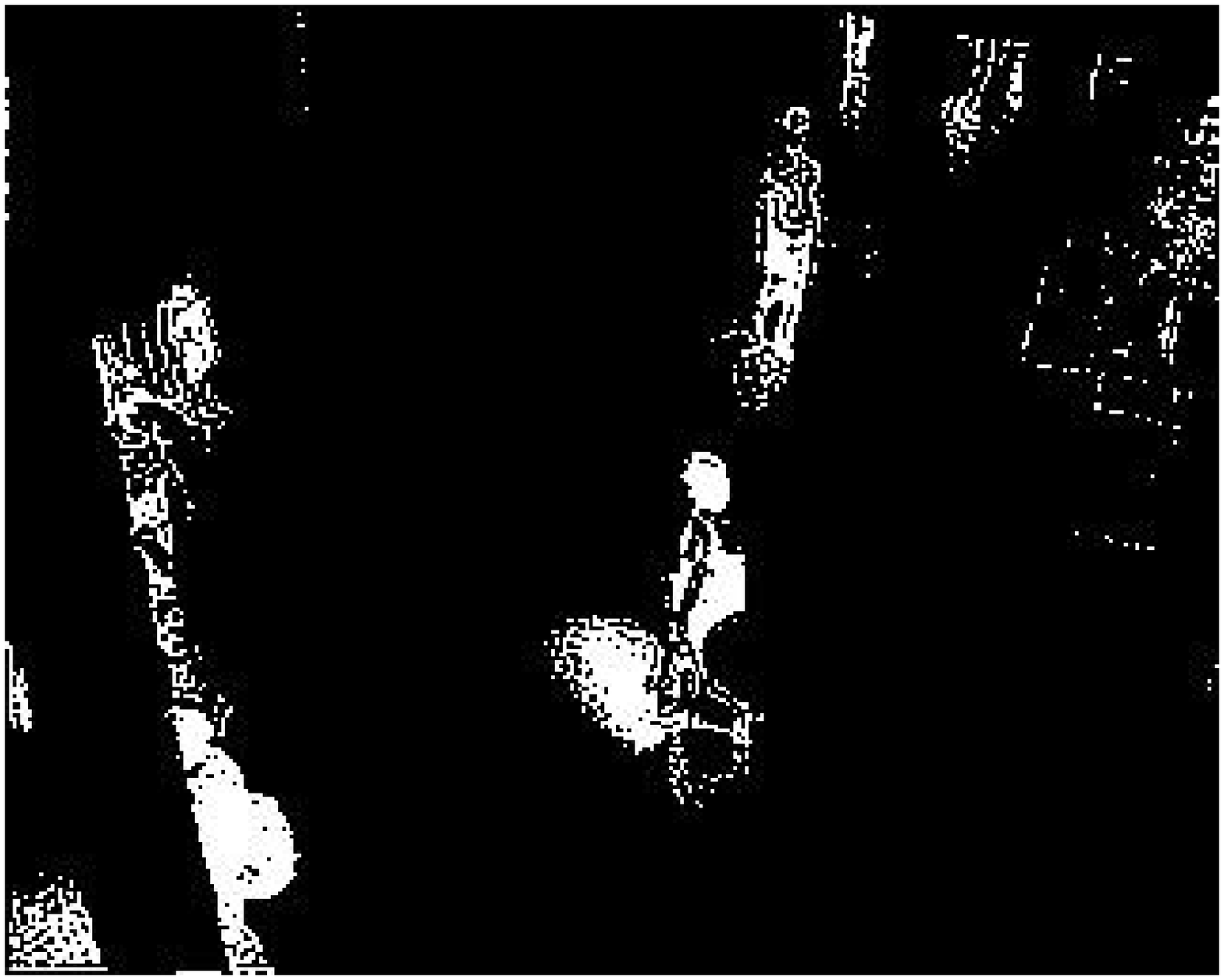}
\end{minipage}}
\subfigure{\begin{minipage}[t]{0.3\textwidth}\centering log.$\alpha=1$\vspace{1mm} \\
\includegraphics[height=1.8cm,width=1.8cm]{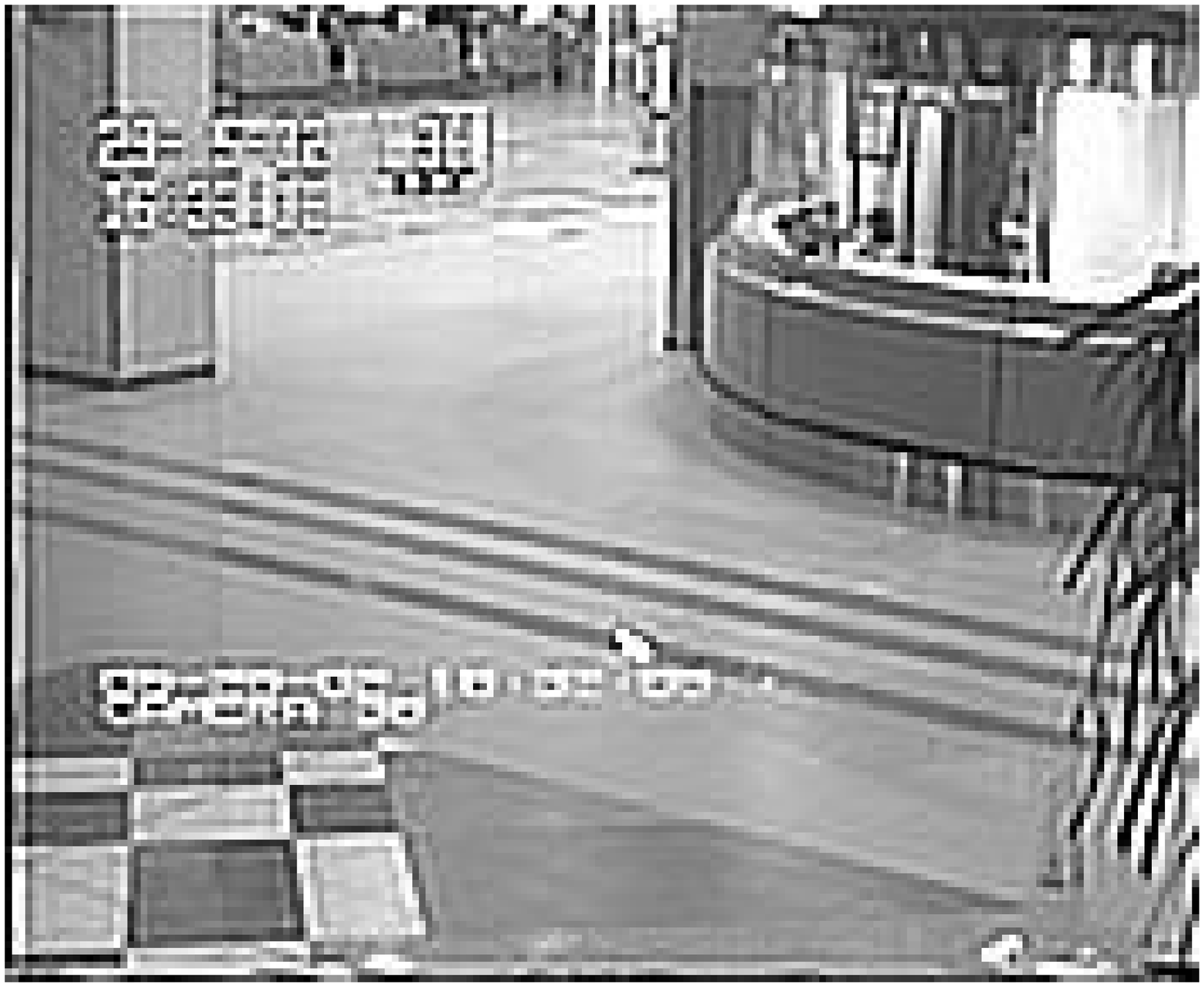}~\includegraphics[height=1.8cm,width=1.8cm]{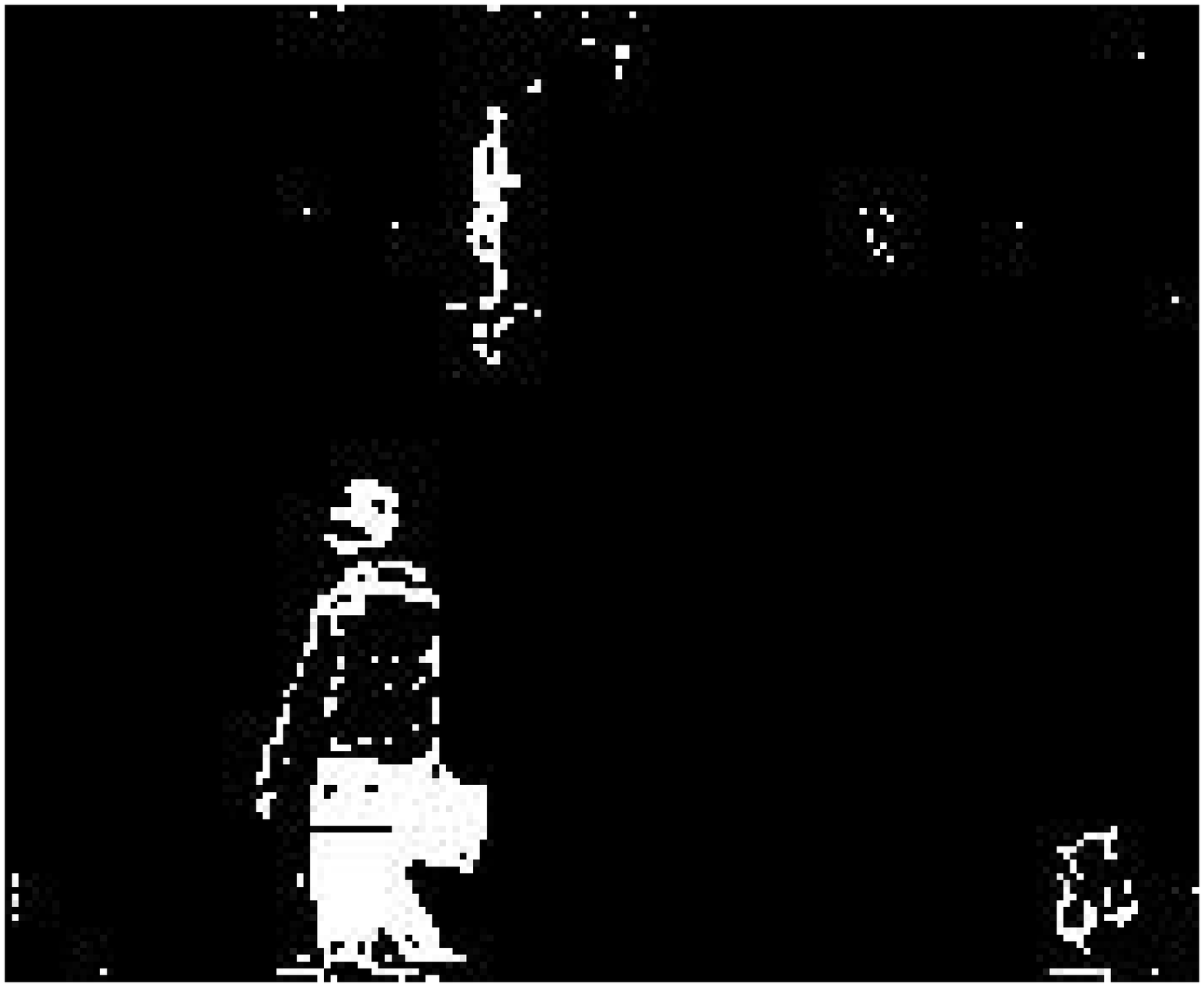}\vspace{0.5mm}\\
\includegraphics[height=1.8cm,width=1.8cm]{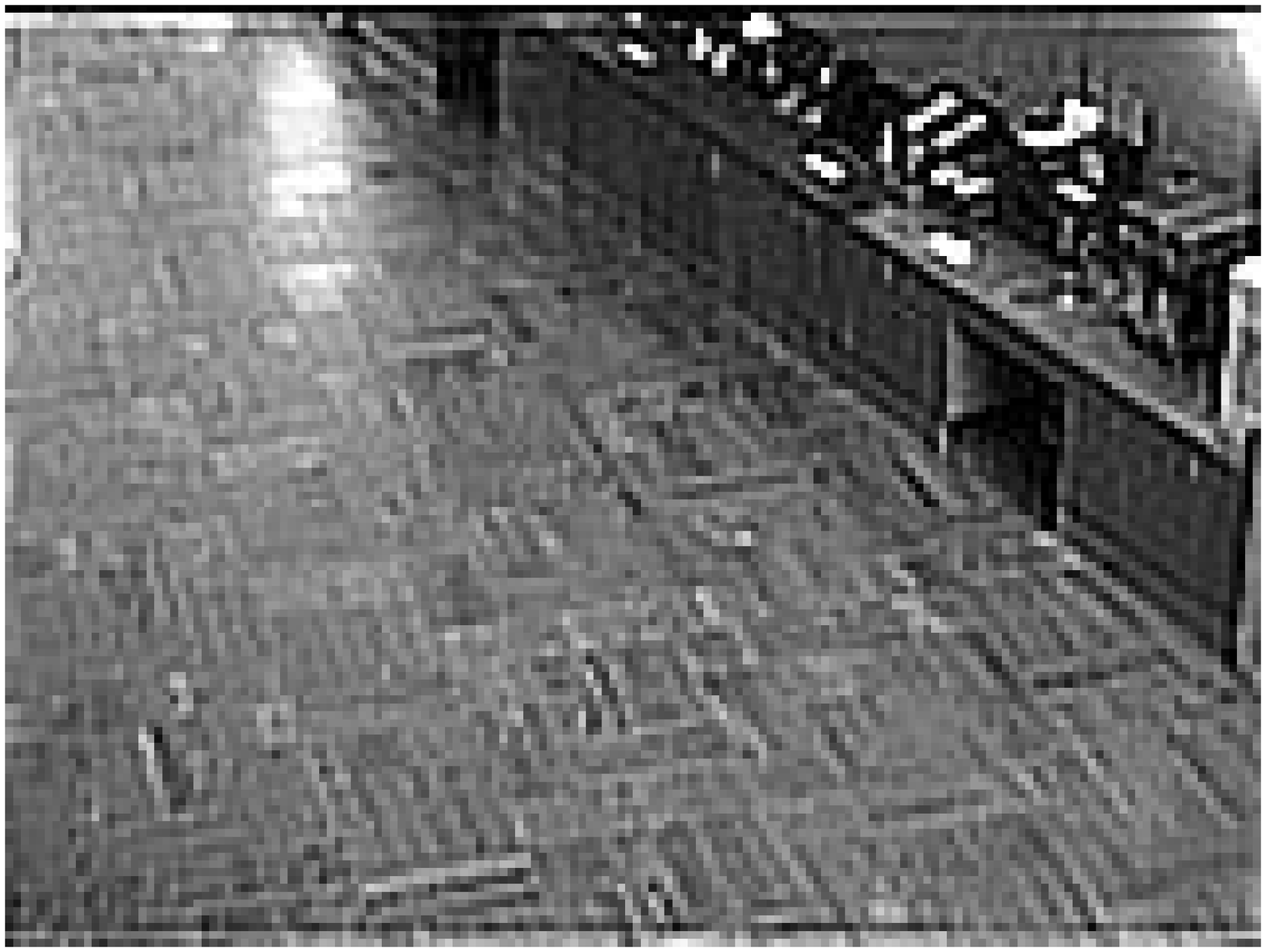}~\includegraphics[height=1.8cm,width=1.8cm]{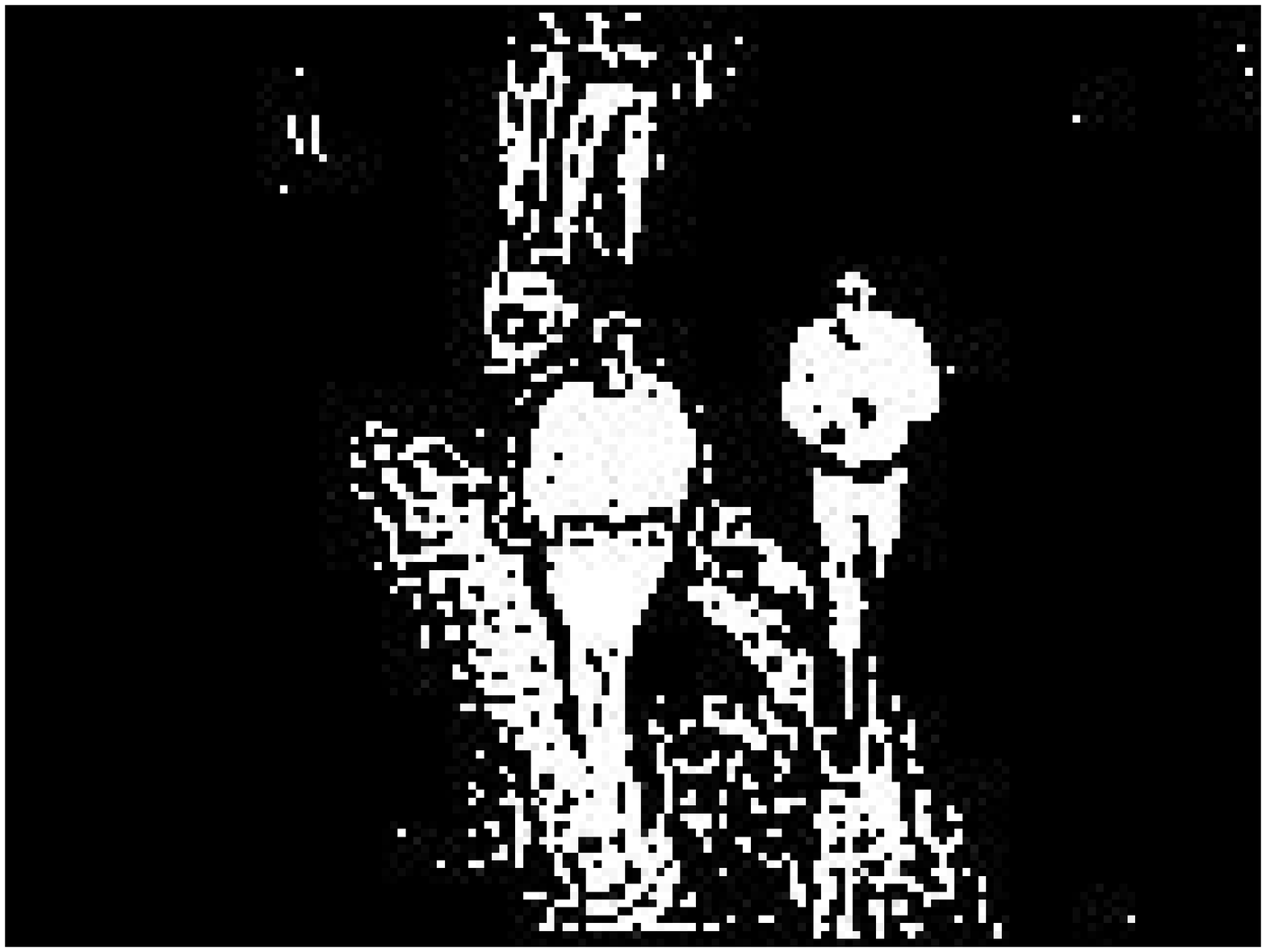}\vspace{0.5mm}\\
\includegraphics[height=1.8cm,width=1.8cm]{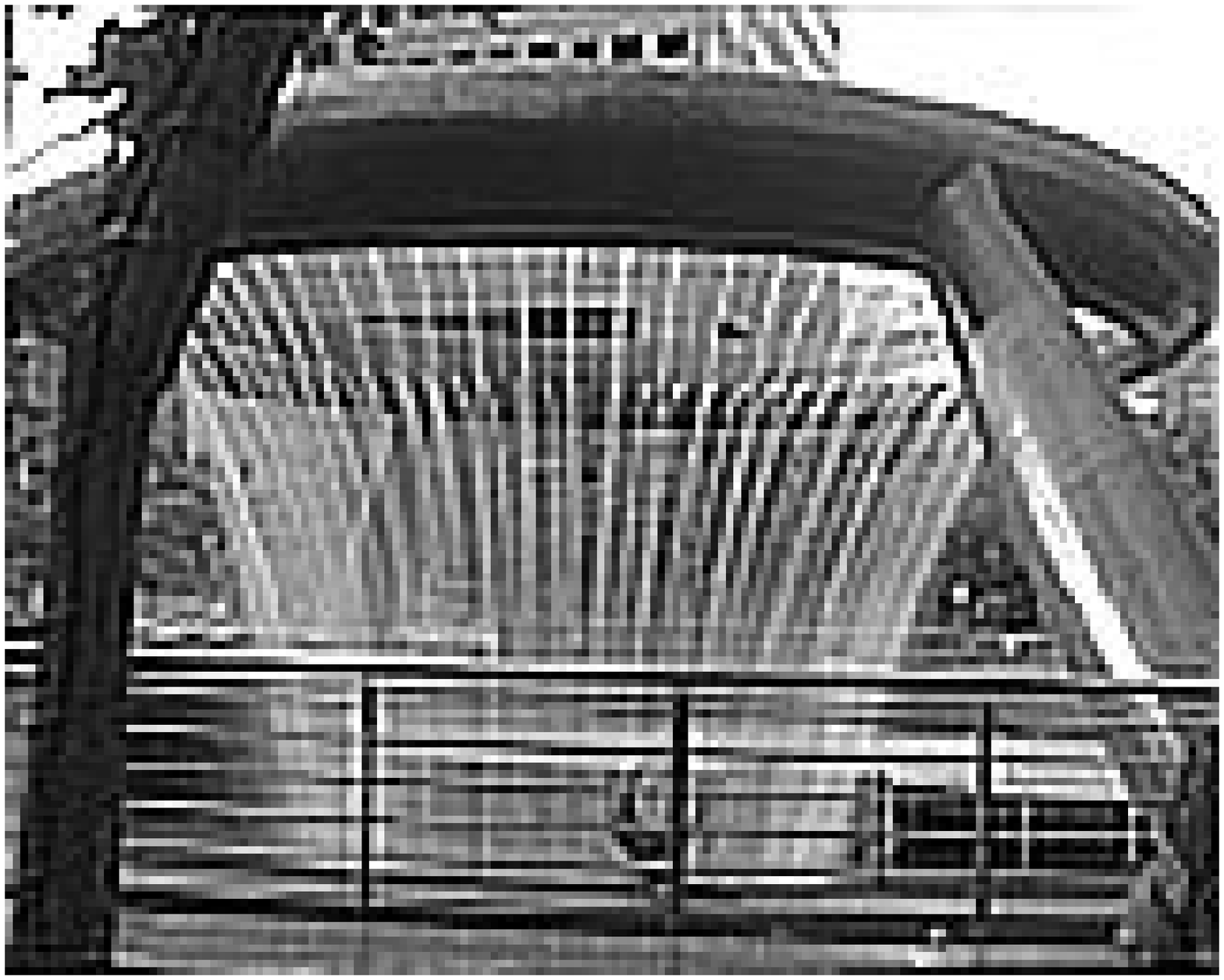}~\includegraphics[height=1.8cm,width=1.8cm]{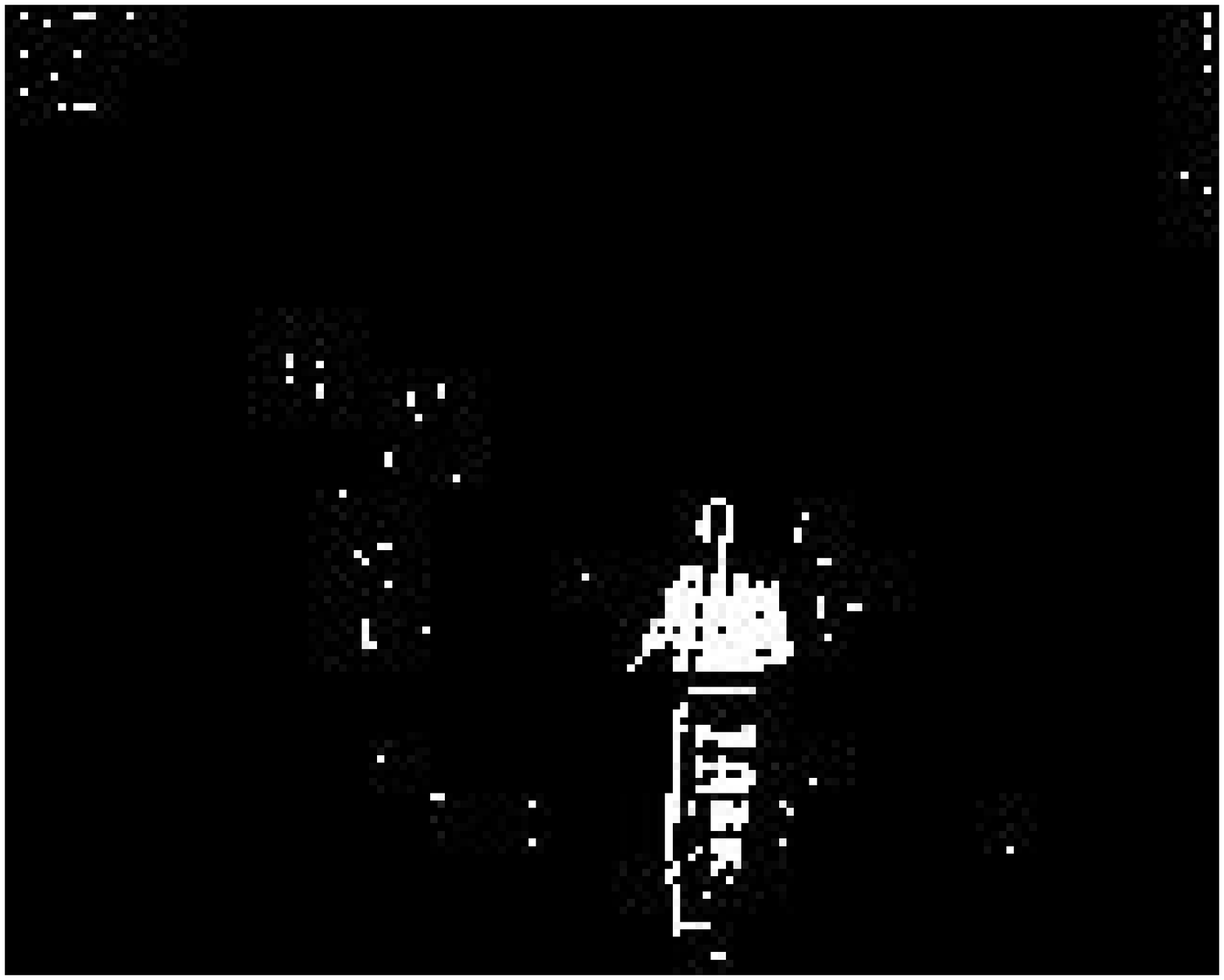}\vspace{0.5mm}\\
\includegraphics[height=1.8cm,width=1.8cm]{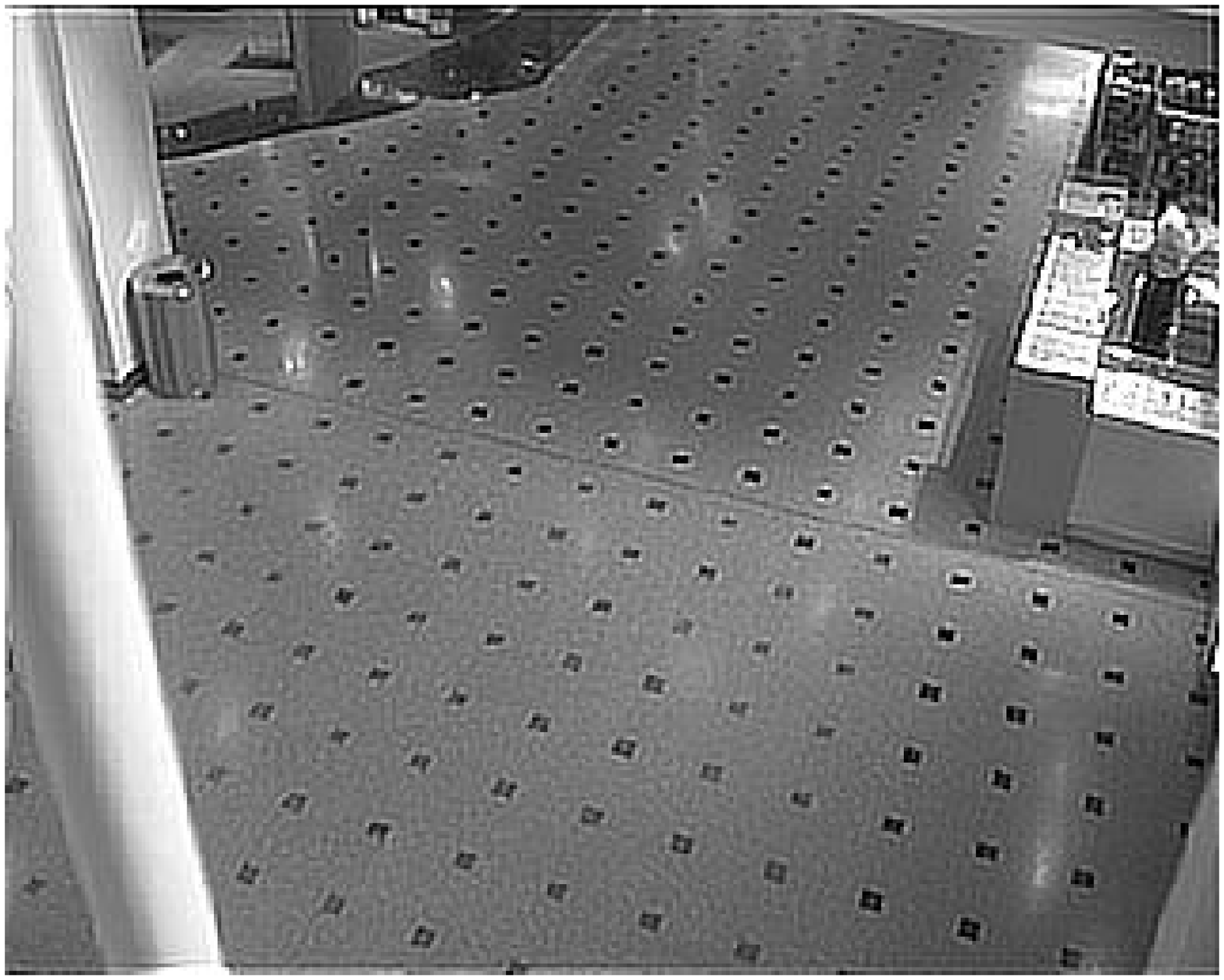}~\includegraphics[height=1.8cm,width=1.8cm]{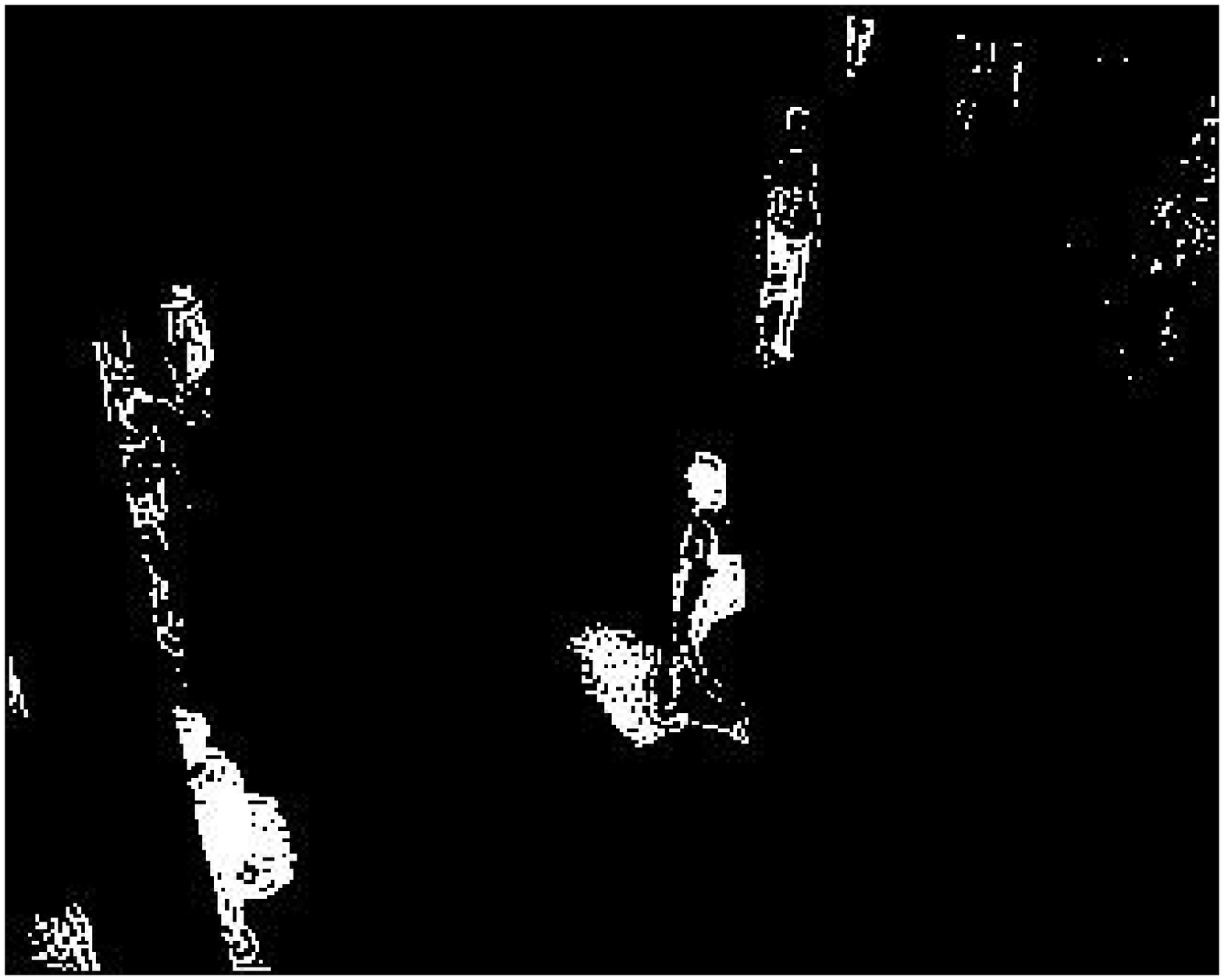}
\end{minipage}}
\subfigure{\begin{minipage}[t]{0.3\textwidth}\centering log.$\alpha=2$\vspace{1mm} \\
\includegraphics[height=1.8cm,width=1.8cm]{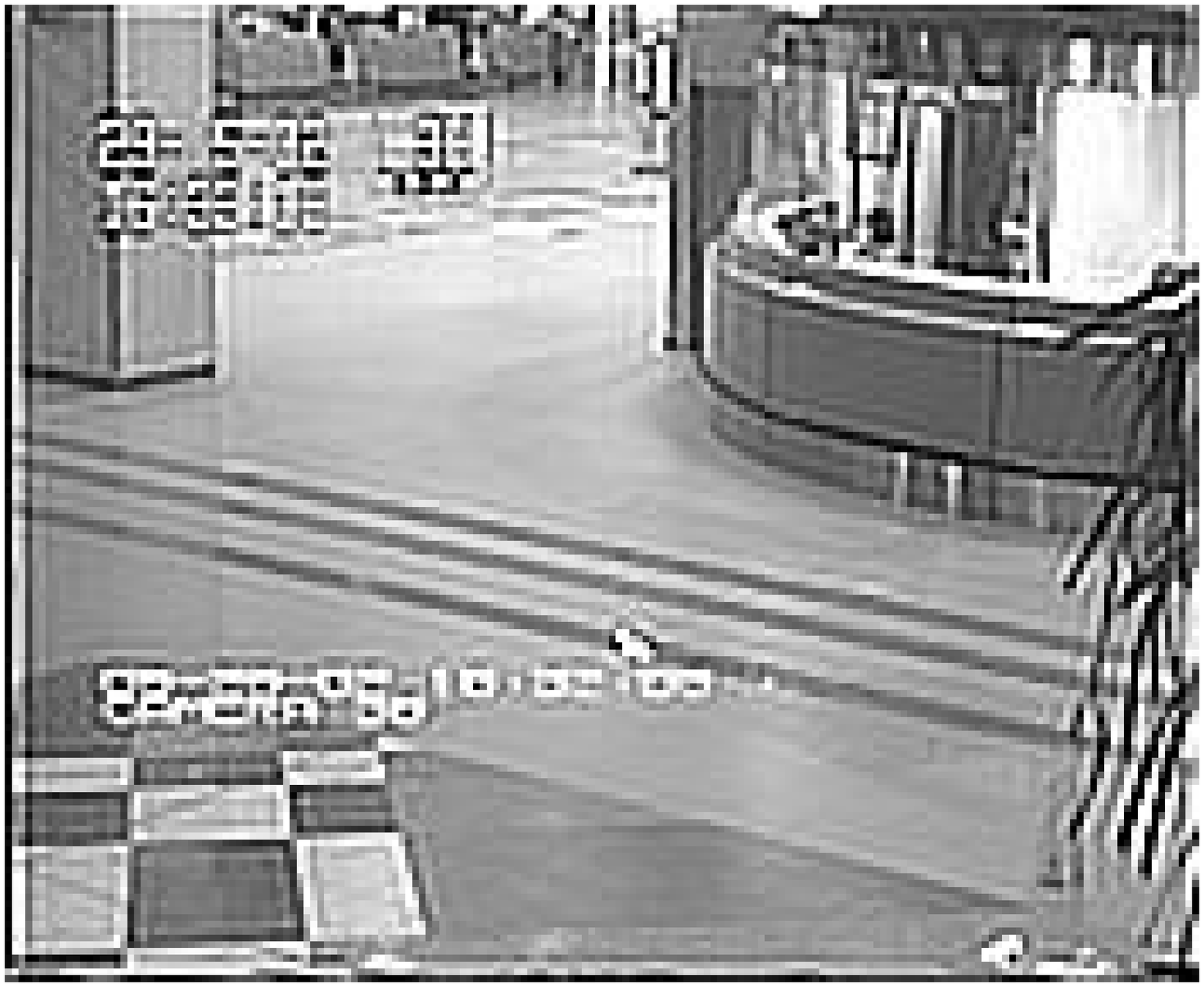}~\includegraphics[height=1.8cm,width=1.8cm]{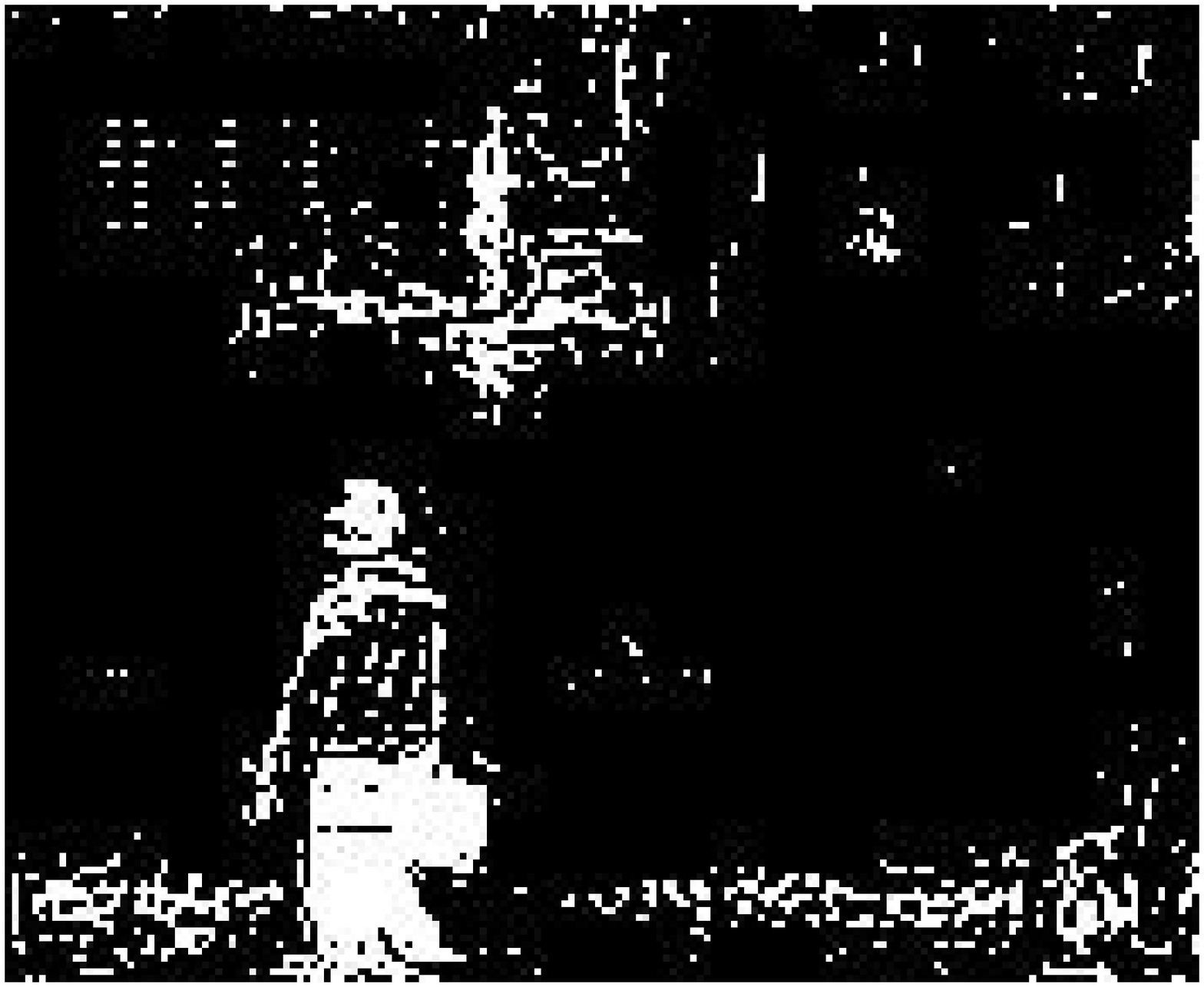}\vspace{0.5mm}\\
\includegraphics[height=1.8cm,width=1.8cm]{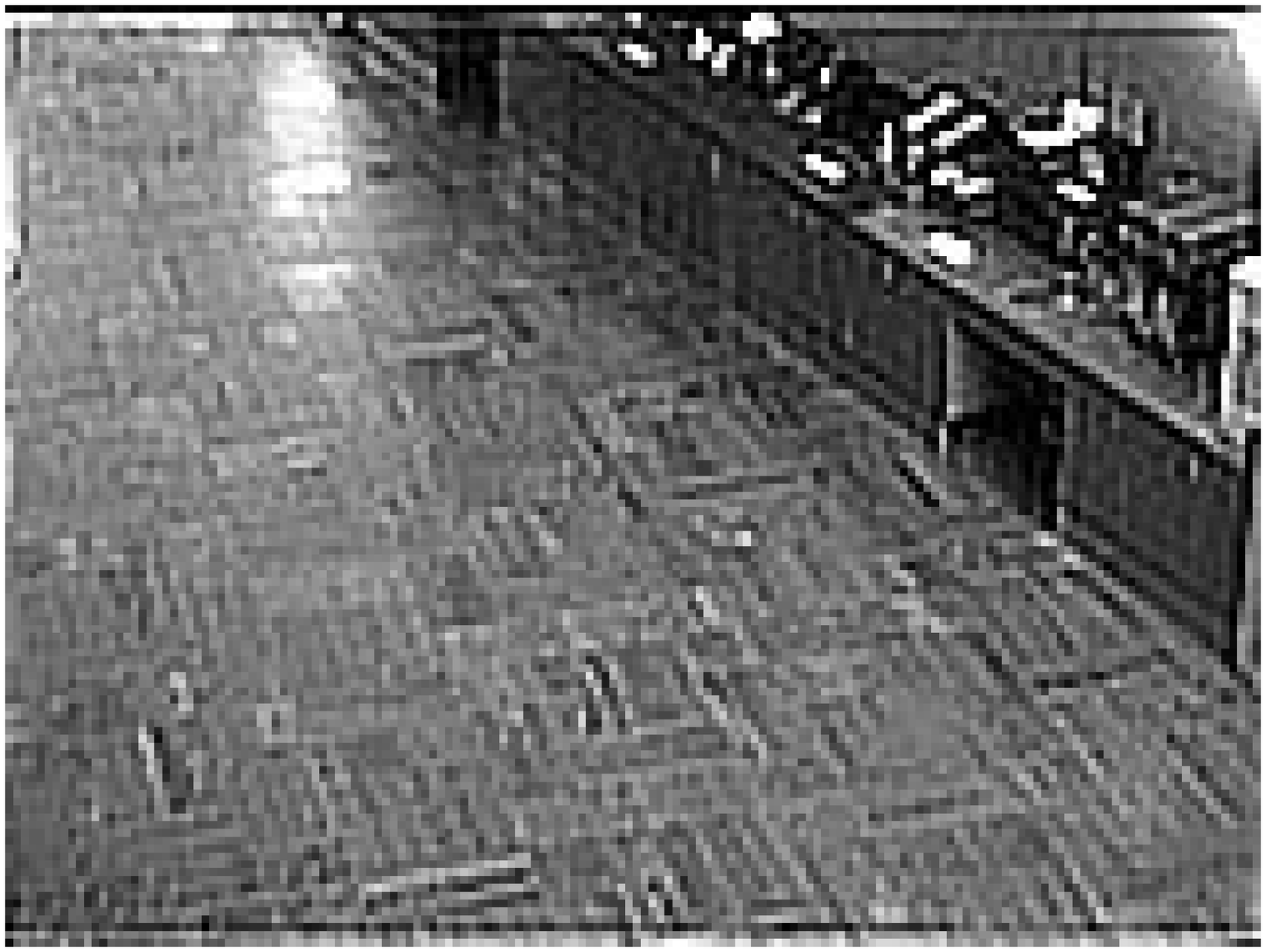}~\includegraphics[height=1.8cm,width=1.8cm]{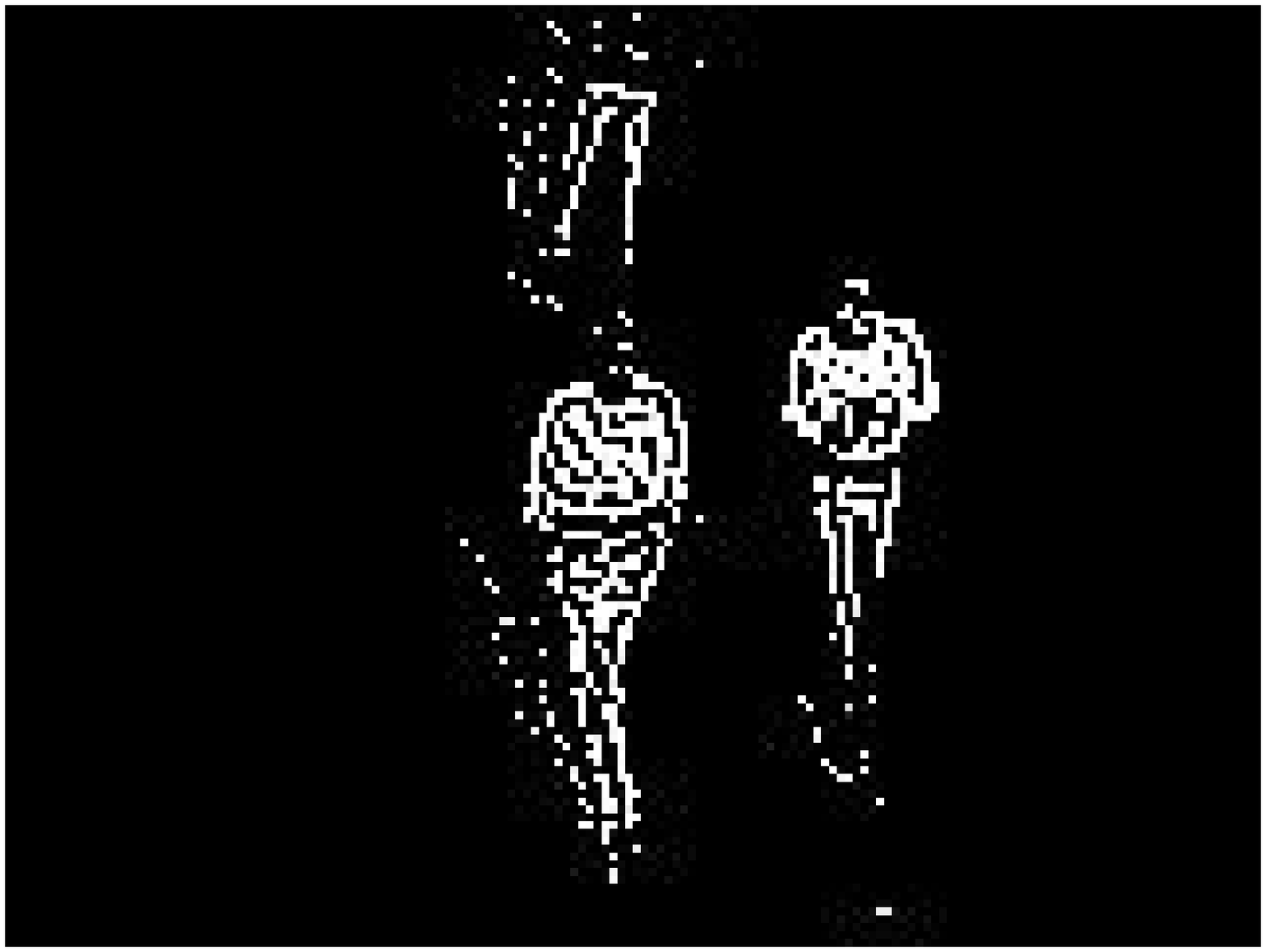}\vspace{0.5mm}\\
\includegraphics[height=1.8cm,width=1.8cm]{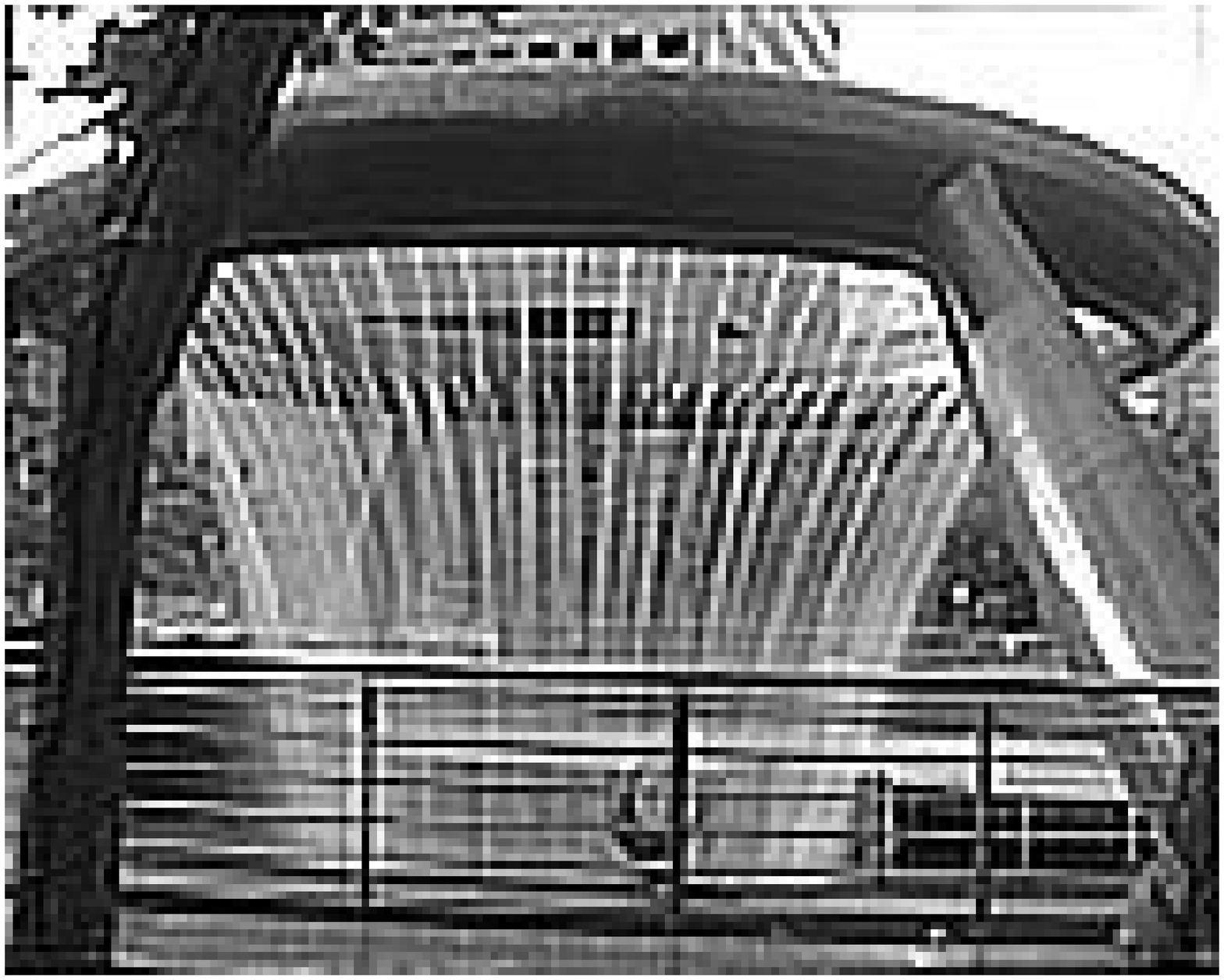}~\includegraphics[height=1.8cm,width=1.8cm]{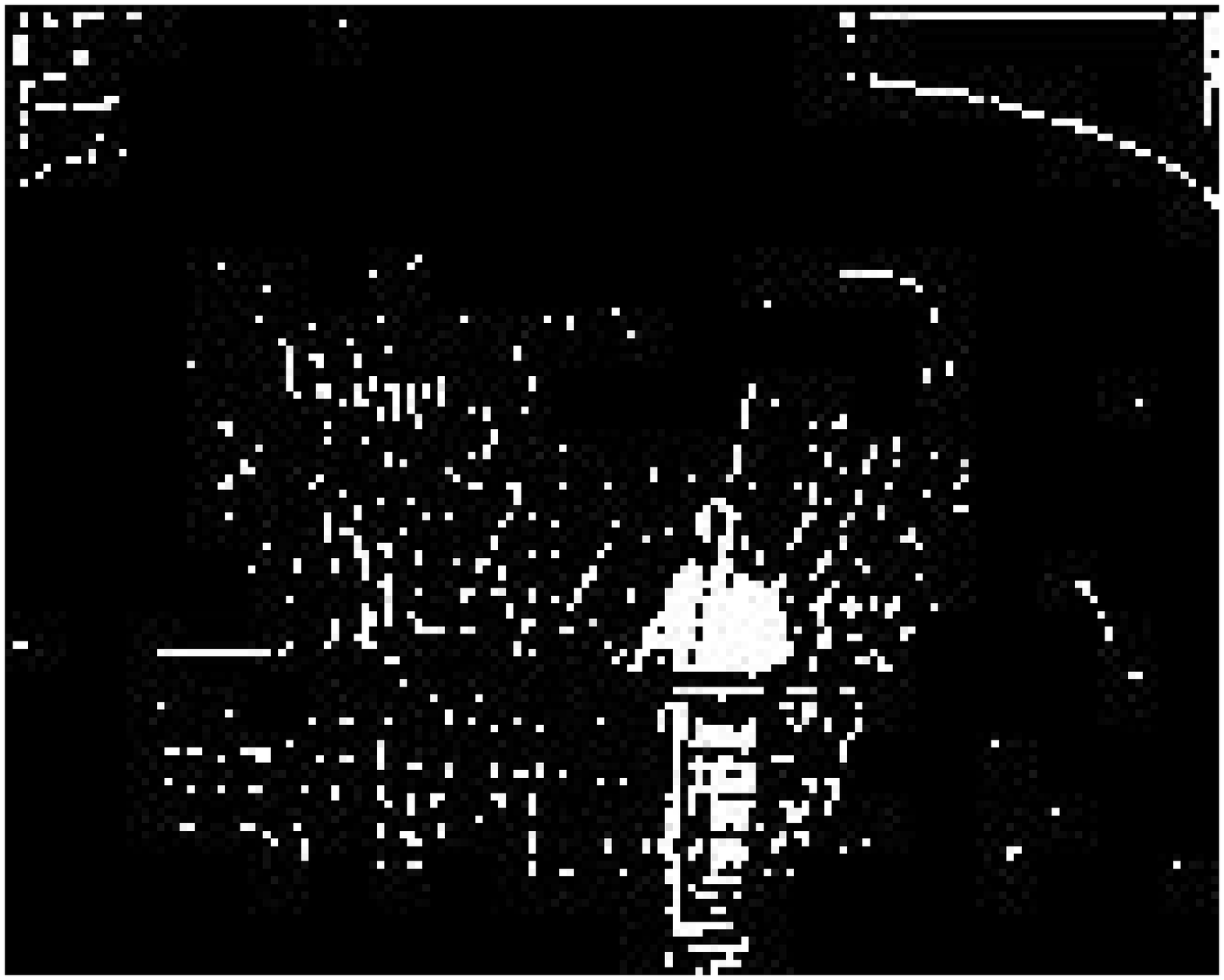}\vspace{0.5mm}\\
\includegraphics[height=1.8cm,width=1.8cm]{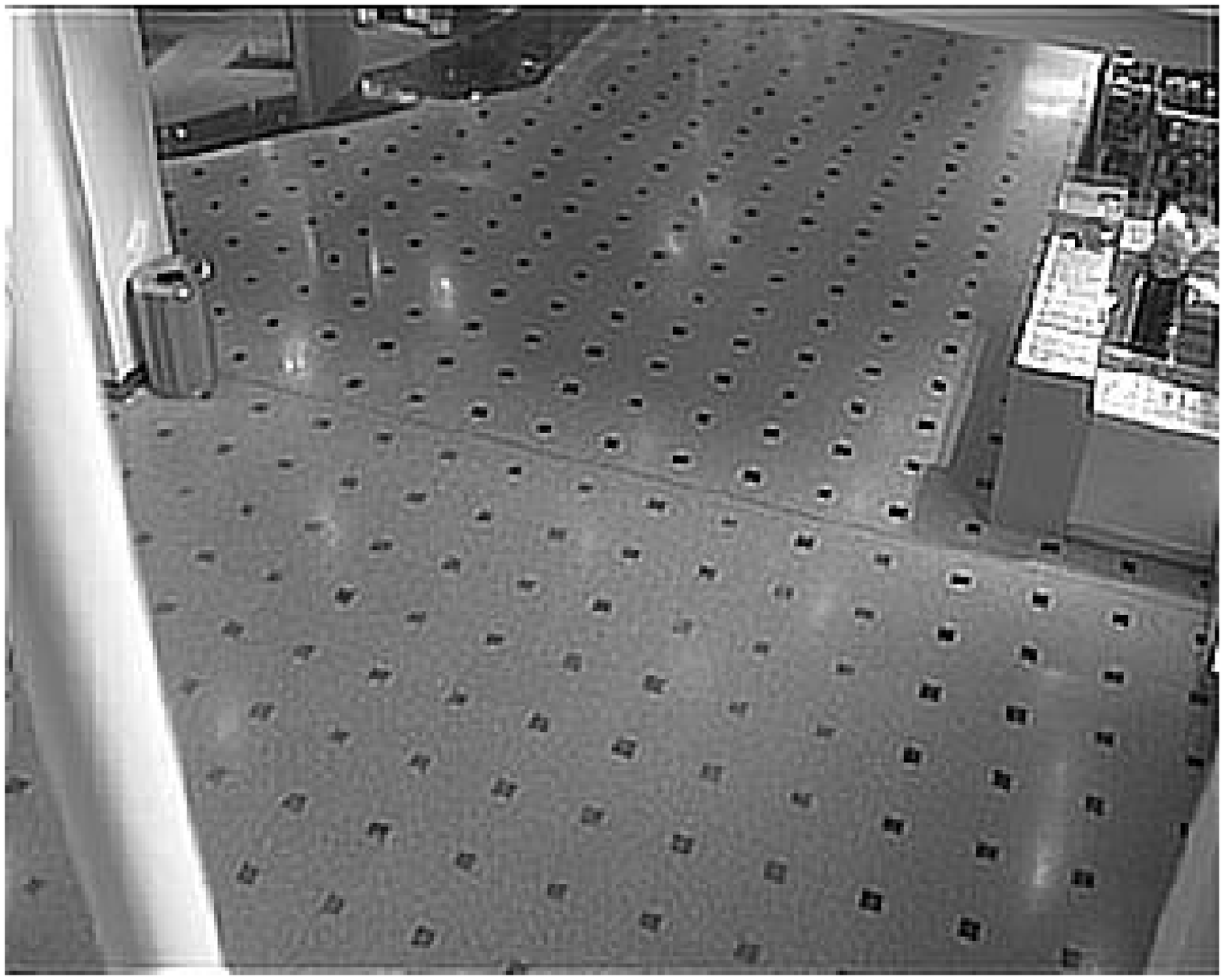}~\includegraphics[height=1.8cm,width=1.8cm]{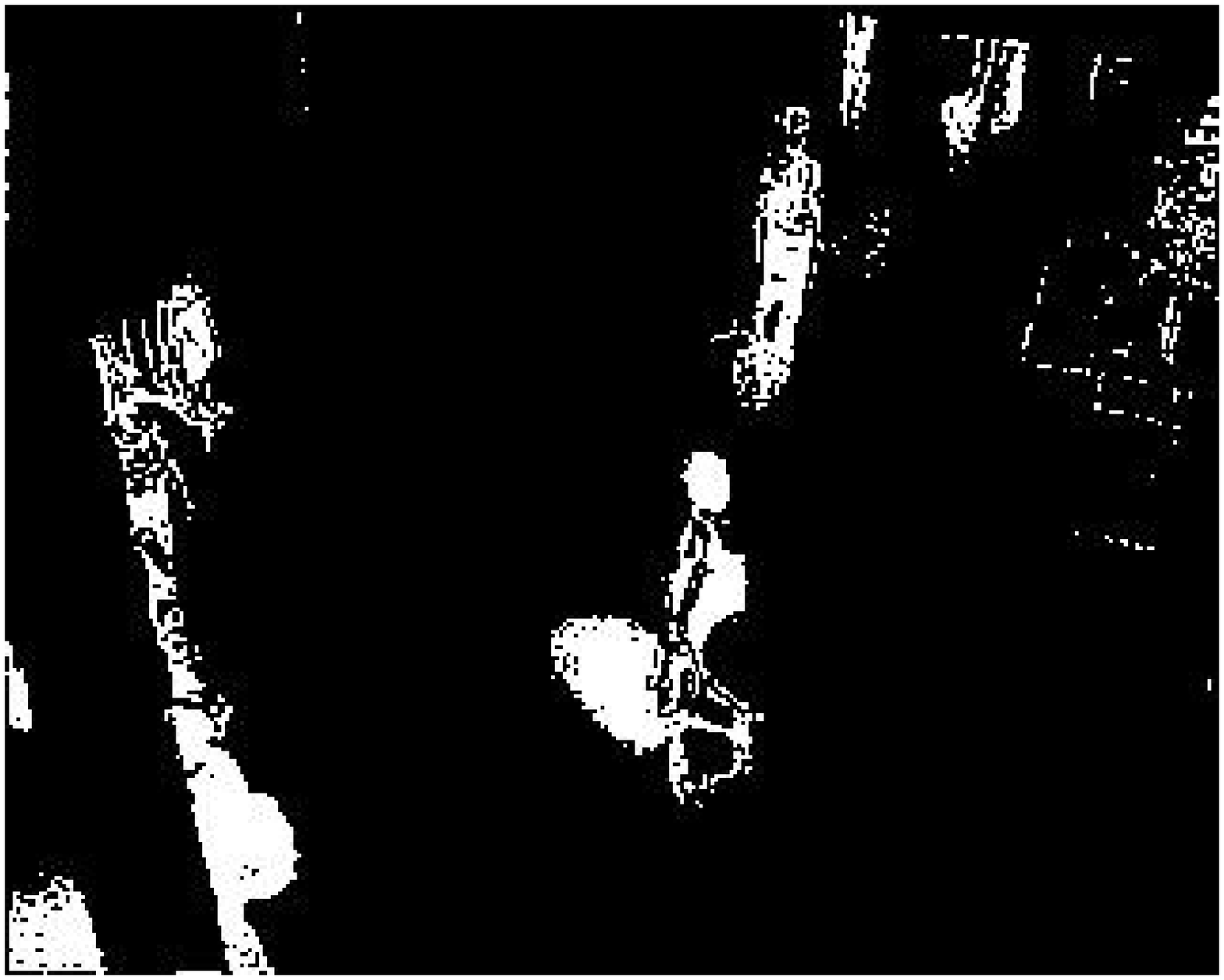}
\end{minipage}}
\caption{Extracted backgrounds and foregrounds given by the ADMM for noisy and blurred surveillance videos.}\label{fig_blur}
\end{figure}

\paragraph{Extraction from noisy and blurred surveillance videos} In this case, ${\cal A}(L+S)=H(L+S)$, $\lambda_{\max}=\lambda_{\max}(H^*H)$, $\lambda_{\min}=\lambda_{\min}(H^*H)$ and we set $\mathrm{Tol}_{A,1} = 5\times10^{-3}$, $\mathrm{Tol}_{A,2} = 10^{-2}$ and $\mathrm{Tol}_{P} = 3\times10^{-3}$. The blurring matrix $H$ is generated by the same method introduced in Subsection \ref{comparisons}. One frame of each corrupted video is shown in the second row in
Fig.\,\ref{fig_org}. We report the computational results in Table \ref{blurresult} and show the extracted backgrounds and foregrounds by the ADMM in Fig.\,\ref{fig_blur}.

\begin{table}[ht]
\caption{Numerical results for extraction from noisy and blurred surveillance videos}\label{blurresult}
\centering \tabcolsep 2pt
\begin{tabular}{|clc|cccc|cccc|}
\hline
&&& \multicolumn{4}{c|}{ADMM} & \multicolumn{4}{c|}{PALM}               \\
Data & \multicolumn{2}{c|}{regularizer} & $\mu$ & iter & time & \footnotesize{F-measure} & $\mu$ & iter & time & \footnotesize{F-measure} \\
\hline
\multirow{6}{*}{\small{Hall}}
& bri.\,$p$       &   1.0  &  5e-02  & 24  & 15.37  & 0.6801  & 5e-02  & 36  & 26.72 & 0.6626 \\
&                 &   0.5  &  1e-02  & 44  & 31.96  & 0.6358  & 1e-02  & 45  & 35.11 & 0.6357 \\
& fra.\,$\alpha$  &   1.0  &  5e-02  & 57  & 40.40  & 0.5265  & 5e-02  & 49  & 37.02 & 0.5616 \\
&                 &   2.0  &  1e-02  & 66  & 52.25  & 0.5381  & 1e-02  & 61  & 50.51 & 0.5445 \\
& log.\,$\alpha$  &   1.0  &  5e-02  & 42  & 93.40  & 0.5970  & 5e-02  & 44  & 100.84 & 0.6033 \\
&                 &   2.0  &  1e-02  & 54  & 118.58  & 0.5188  & 1e-02  & 52  & 120.68 & 0.5211 \\
\hline

\multirow{6}{*}{\small{Bootstrap}}
& bri.\,$p$      &   1.0  &  1e-01  & 22  & 10.56  & 0.7651  & 1e-01  & 50  & 26.20 & 0.7364 \\
&                &   0.5  &  1e-02  & 56  & 36.50  & 0.6692  & 1e-02  & 86  & 55.18 & 0.6589 \\
& fra.\,$\alpha$ &   1.0  &  5e-02  & 66  & 36.75  & 0.5705  & 5e-02  & 99  & 60.04 & 0.5270 \\
&                &   2.0  &  5e-02  & 73  & 41.08  & 0.5265  & 5e-02  & 111  & 61.89 & 0.4674 \\
& log.\,$\alpha$ &   1.0  &  5e-02  & 47  & 83.24  & 0.5651  & 5e-02  & 87  & 158.16 & 0.5666 \\
&                &   2.0  &  5e-02  & 73  & 131.17  & 0.4891  & 5e-02  & 115  & 213.62 & 0.4179 \\
\hline

\multirow{6}{*}{\small{Fountain}}
& bri.\,$p$      &   1.0  &  5e-02  & 28  & 12.72  & 0.7229  & 5e-02  & 64  & 36.25 & 0.6970 \\
&                &   0.5  &  1e-02  & 51  & 28.54  & 0.6881  & 1e-02  & 78  & 45.68 & 0.6606 \\
& fra.\,$\alpha$ &   1.0  &  5e-02  & 64  & 34.55  & 0.5155  & 5e-02  & 84  & 49.98 & 0.5000 \\
&                &   2.0  &  1e-02  & 62  & 37.87  & 0.4482  & 1e-02  & 87  & 56.47 & 0.4341 \\
& log.\,$\alpha$ &   1.0  &  5e-02  & 50  & 98.81  & 0.6095  & 5e-02  & 77  & 162.12 & 0.5760 \\
&                &   2.0  &  1e-02  & 53  & 105.63  & 0.4438  & 1e-02  & 80  & 169.25 & 0.4525 \\
\hline

\multirow{6}{*}{\small{ShoppingMall}}
& bri.\,$p$       &  1.0  &  5e-02  & 22  & 51.99  & 0.6431  & 5e-02  & 25  & 62.71 & 0.6411 \\
&                 &  0.5  &  5e-03  & 54  & 157.53  & 0.6271  & 5e-03  & 40  & 119.91 & 0.6328 \\
& fra.\,$\alpha$  &  1.0  &  1e-02  & 33  & 109.64  & 0.5045  & 5e-02  & 60  & 188.48 & 0.5106 \\
&                 &  2.0  &  1e-02  & 51  & 153.19  & 0.5810  & 1e-02  & 38  & 127.47 & 0.5935 \\
& log.\,$\alpha$  &  1.0  &  5e-02  & 58  & 263.00  & 0.5453  & 5e-02  & 39  & 178.60 & 0.5967 \\
&                 &  2.0  &  1e-02  & 38  & 174.88  & 0.5856  & 1e-02  & 33  & 149.36 & 0.5913 \\
\hline
\end{tabular}
\end{table}

\paragraph{Summary} From the results above, it can be seen that the ADMM with $\tau=0.8$ performs better in the sense that it takes less CPU time for solving most test problems while returning comparable F-measures. The performances of our ADMM for extraction are also promising
from Fig.\,\ref{fig_noisy} and Fig.\,\ref{fig_blur}.

%our proposed nonconvex model ($0<p<1$) is competitive with the benchmark convex model ($p=1$) in the sense that there is always at least one $p<1$ corresponding to a higher F-measure than $p=1$, meaning that the foreground recovered by the $\ell_p$-based model with a certain $0<p<1$ is better than that recovered by the $\ell_1$-based model. Particularly, for the case when $\mathscr{A}(L+S)=\mathcal{P}_{\Omega}(L+S)$, we can see from Table \ref{missresult} that all the F-measures corresponding to $0<p<1$ are higher than that corresponding to $p=1$. One can also observe readily from Fig.\,\ref{fig_miss} that the $\ell_1$-based model cannot recover the region of the foreground completely. On the other hand, in all instances, we can also observe that a sparser extracted foreground does not necessarily correspond to a higher F-measure. This is because the true foreground in the real videos may not be very sparse.

\section{Concluding remarks}\label{sec6}

In this paper, we study a general (possibly nonconvex and nonsmooth) model and adapt the ADMM with a general dual step-size $\tau$, which can be chosen in $(0,\frac{1+\sqrt{5}}{2})$, to solve it. We establish that any cluster point of the sequence generated by our ADMM gives a stationary point under some assumptions; we also give simple sufficient conditions for these assumptions. Under an additional assumption that a potential function is a Kurdyka-{\L}ojasiewicz function, we can further establish the global convergence of the whole sequence generated by our ADMM. Our computational results demonstrate the efficiency of our algorithm.

Note that our ADMM may not be beneficial when $\mathcal{B}$ or $\mathcal{C}$ has no special structure, because the corresponding subproblems of ADMM may not have closed-form solutions. Nonetheless, as in \cite{lp2014,wcx2015,wxx2014}, it may be possible to add ``proximal terms" to simplify the subproblems of our ADMM. In addition, in view of the recent work \cite{wyz2015}, it may also be possible to study the convergence of our ADMM for some specially structured nonconvex $\Psi$. These are possible future research directions.

\section*{Acknowledgments}

\noindent The authors are grateful to the editor and the anonymous referees for their valuable suggestions and comments, which helped improve this paper.

%--------------------------------------------------------------------------------------- References


\begin{thebibliography}{99}

\bibitem{ah2014}
B.P.W. Ames and M. Hong.
\newblock Alternating direction method of multipliers for penalized zero-variance discriminant analysis.
\newblock {\em arXiv preprint arXiv:1401.5492}, 2014.

\bibitem{AtBoReSo10}
H. Attouch, J. Bolte, P. Redont and A. Soubeyran.
\newblock Proximal alternating minimization and projection methods for nonconvex problems: An approach based on the Kurdyka-{\L}ojasiewicz inequality.
\newblock {\em Mathematics of Operations Research}, 35(2): 438--457, 2010.

\bibitem{bt1989}
D. Bertsekas and J.N. Tsitsiklis.
\newblock Parallel and Distributed Computation: Numerical Methods.
\newblock Prentice Hall, 1989.

\bibitem{bc2015}
W. Bian and X. Chen.
\newblock Linearly constrained non-Lipschitz optimization for image restoration.
\newblock {\em SIAM Journal on Imaging Sciences}, 8(4): 2294--2322, 2015.

\bibitem{bst2014}
J. Bolte, S. Sabach and M. Teboublle.
\newblock Proximal alternating linearized minimization for nonconvex and nonsmooth problems.
\newblock {\em Mathematical Programming}, 146(1--2): 459--494, 2014.

\bibitem{Bou11}
T. Bouwmans.
\newblock Recent advanced statistical background modelling for foreground detection: A systematic survey.
\newblock {\em Recent Patents on Computer Science}, 4(3): 147--176, 2011.

\bibitem{Bou14}
T. Bouwmans.
\newblock Traditional and recent aproaches in background modeling for foreground detection: An overview.
\newblock {\em Computer Science Review}, 11--12: 31--66, 2014.

\bibitem{bsjjz2015}
T. Bouwmans, A. Sobral, S. Javed, S.K. Jung and E.-H. Zahzah.
\newblock Decomposition into Low--rank plus Additive Matrices for Background/Foreground Separation: A Review for a Comparative Evaluation with a Large-Scale Dataset.
\newblock {\em arXiv preprint arXiv:1511.01245}, 2015.

\bibitem{bz2014}
T. Bouwmans and E.H. Zahzah.
\newblock Robust PCA via principal component pursuit: A review for a comparative evaluation in video surveillance.
\newblock {\em Computer Vision and Image Understanding}, 122: 22--34, 2014.

\bibitem{CLMW11}
E.J. Cand\`{e}s, X. Li, Y. Ma and J. Wright.
\newblock Robust principal component analysis?
\newblock {\em Journal of the ACM}, 58(3), Article 11, 2011.

\bibitem{CHYY14}
C. Chen, B. He, Y. Ye and X. Yuan.
\newblock The direct extension of ADMM for multi-block convex minimization problems is not necessarily convergent.
\newblock {\em Mathematical Programming}, 155: 57--79, 2016.

\bibitem{cz2010}
X. Chen and W. Zhou.
\newblock Smoothing nonlinear conjugate gradient method for image restoration using nonsmooth nonconvex minimization.
\newblock {\em SIAM Journal on Imaging Sciences}, 3(4): 765--790, 2010.

\bibitem{DDLZH13}
Y. Deng, Q. Dai, R. Liu, Z. Zhang and S. Hu.
\newblock Low-rank structure learning via nonconvex heuristic recovery.
\newblock {\em IEEE Transactions on Neural Networks and Learning Systems}, 24(3): 383--396, 2013.

\bibitem{dm2012}
E.D. Dolan and J.J. Mor\'{e}.
\newblock Benchmarking optimization software with performance profiles.
\newblock {\em Mathematical Programming}, 91(2): 201--213, 2012.

\bibitem{eb1992}
J. Eckstein and D. Bertsekas.
\newblock On the Douglas–Rachford splitting method and the proximal point algorithm for maximal monotone operators.
\newblock {\em Mathematical Programming}, 55, 293--318, 1992.

\bibitem{f1997}
J. Fan.
\newblock Comments on ``wavelets in statistics: A review" by A.Antoniadis.
\newblock {\em Journal of the Italian Statistical Society}, 6(2): 131--138, 1997.

\bibitem{fl2001}
J. Fan and R. Li.
\newblock Variable selection via nonconcave penalized likelihood and its oracle properties.
\newblock {\em Journal of the American Statistical Association}, 96(456): 1348--1360, 2001.

\bibitem{FPST12}
M. Fazel, T.K. Pong, D. Sun and P. Tseng.
\newblock Hankel matrix rank minimization with applications to system identification and realization.
\newblock {\em SIAM Journal on Matrix Analysis and Applications}, 34(3): 946--977, 2013.

\bibitem{GaM76}
D. Gabay and B. Mercier.
\newblock A dual algorithm for the solution of nonlinear variational problems via finite element approximations.
\newblock {\em Computers \& Mathematics with Applications}, 2(1): 17--40, 1976.

\bibitem{gr1992}
D. Geman and G. Reynolds.
\newblock Constrained restoration and the recovery of discontinuities.
\newblock {\em IEEE Transactions on Pattern Analysis and Machine Intelligence}, 14(3): 367--383, 1992.

\bibitem{GlM75}
R. Glowinski and A. Marroco.
\newblock  Sur l'approximation, par \'{e}l\'{e}ments finis d'ordre un, et la r\'{e}solution, par p\'{e}nalisation-dualit\'{e}, d'une classe de probl\`{e}mes de Dirichlet non lin\'{e}ares.
\newblock {\em Revue Francaise d'Automatique, Informatique, Recherche Op\'{e}rationelle.}, 9(R-2): 41--76, 1975.

\bibitem{gzzf2014}
S. Gu, L. Zhang, W. Zuo and X. Feng.
\newblock Weighted nuclear norm minimization with application to image denoising.
\newblock in {\em CVPR}, 2862--2869, 2014.

\bibitem{hno2006}
P.C. Hansen, J.G. Nagy and D.P. O'leary.
\newblock Deblurring Images: Matrices, Spectra, and Filtering.
\newblock Fundamentals of Algorithms 3, SIAM, Philadelphia, 2006.

\bibitem{hty2012}
B. He, M. Tao and X. Yuan.
\newblock Alternating direction method with Gaussian back substitution for separable convex programming.
\newblock {\em SIAM Journal on Optimization}, 22(2): 313--340, 2012.

\bibitem{hy2013}
B. He and X. Yuan.
\newblock Linearized alternating direction method with Gaussian back substitution for separable convex programming.
\newblock {\em Numerical Algebra, Control and Optimization}, 3(2), 247--260, 2013.

\bibitem{hlr2014}
M. Hong, Z.-Q. Luo and M. Razaviyayn.
\newblock Convergence analysis of alternating direction method of multipliers for a family of nonconvex problems.
\newblock {\em SIAM Journal on Optimization}, 26(1): 337--364, 2016.

\bibitem{hhjm2008}
J. Huang, J. L. Horowitz and S. Ma.
\newblock Asymptotic properties of bridge estimators in sparse high-dimensional regression models.
\newblock {\em Annals of Statistics}, 36(2): 587--613, 2008.

\bibitem{kf2000}
K. Knight and W. Fu.
\newblock Asymptotics for lasso-type estimators.
\newblock {\em Annals of Statistics}, 28(5): 1356--1378, 2000.

\bibitem{lhgt2004}
L. Li, W. Huang, I. Y.-H. Gu and Q. Tian.
\newblock Statistical modeling of complex backgrounds for foreground object detection.
\newblock {\em IEEE Transactions on Image Processing}, 13(11): 1459--1472, 2004.

\bibitem{lny2013}
X. Li, M. K. Ng and X. Yuan.
\newblock Median filtering--based methods for static background extraction from surveillance video.
\newblock {\em Numerical Linear Algebra with Applications}, 22: 845--865, 2015.

\bibitem{lp2014}
G. Li and T.K. Pong.
\newblock Global convergence of splitting methods for nonconvex composite optimization.
\newblock {\em SIAM Journal on Optimization}, 25(4): 2434--2460, 2015.

\bibitem{lp2014Douglas}
G. Li and T.K. Pong.
\newblock Douglas-Rachford splitting for nonconvex optimization with application to nonconvex feasibility problems.
\newblock {\em Mathematical Programming}, 159(1): 371--401, 2016.

\bibitem{lwhc2014}
L. Li, P. Wang, Q. Hu and S. Cai.
\newblock Efficient background modeling based on sparse representation and outlier iterative removal.
\newblock {\em IEEE Transactions on Circuits and Systems for Video Technology}, 26(2): 278--289, 2014.

\bibitem{llysym2013}
G. Liu, Z. Lin, S. Yan, J. Sun, Y. Yu and Y. Ma.
\newblock Robust recovery of subspace structures by low-rank representation.
\newblock {\em IEEE Transactions on Pattern Analysis and Machine Intelligence}, 35(1): 171--184, 2013.

\bibitem{lll2015}
Z. Lin, R. Liu and H. Li.
\newblock Linearized alternating direction method with parallel splitting and adaptive penalty for separable convex programs in machine learning.
\newblock {\em Machine Learning}, 99(2): 287--325, 2015.

\bibitem{lst2015}
M. Lin, D. Sun and K.-C. Toh.
\newblock A convergent 3-block semi-proximal ADMM for convex minimization problems with one strongly convex block.
\newblock {\em Asia-Pacific Journal of Operational Research}, 32(3): 1550024(19p), 2015.

\bibitem{mlfwl2014}
K. Mohan, P. London, M. Fazel, D. Witten and S.-I. Lee.
\newblock Node-based learning of multiple gaussian graphical models.
\newblock {\em Journal of Machine Learning Research}, 15(1): 445--488, 2014.

\bibitem{n2014}
M. Nikolova.
\newblock Energy Minimization Methods.
\newblock in {\em Handbook of Mathematical Methods in Imaging}, O. Scherzer, Ed. Springer, pp. 139--185, 2011.

\bibitem{nnzc2008}
M. Nikolova, M.K. Ng, S. Zhang and W.-K. Ching.
\newblock Efficient reconstruction of piecewise constant images using nonsmooth nonconvex minimization.
\newblock {\em SIAM Journal on Imaging Sciences}, 1(1): 2--25, 2008.

\bibitem{ocs2015}
R. Otazo, E.J. Cand\`{e}s and D. K. Sodickson.
\newblock Low--rank plus sparse matrix decomposition for accelerated dynamic MRI with separation of background and dynamic components.
\newblock {\em Magnetic Resonance in Medicine}, 73(3): 1125--1136, 2015.

\bibitem{pgwxm2012}
Y. Peng, A. Ganesh, J. Wright, W. Xu and Y. Ma.
\newblock Robust alignment by sparse and low-rank decomposition for linearly correlated images.
\newblock {\em IEEE Transactions on Pattern Analysis and Machine Intelligence}, 34(11): 2233--2246, 2012.

\bibitem{rw1998}
R.T. Rockafellar and R.J-B. Wets.
\newblock {\em Variational Analysis}.
\newblock Springer, 1998.

\bibitem{sty2015}
D. Sun, K.-C. Toh and L. Yang.
\newblock A convergent 3--block semi-proximal alternating direction method of multipliers for conic programming with 4--type constraints.
\newblock {\em SIAM Journal on Optimization}, 25(2): 882--915, 2015.

\bibitem{wcx2015}
F. Wang, W. Cao and Z. Xu.
\newblock Convergence of multi-block Bregman ADMM for nonconvex composite problems.
\newblock {\em arXiv preprint arXiv:1505.03063}, 2015.

\bibitem{wxx2014}
F. Wang, Z. Xu and H.-K. Xu.
\newblock Convergence of Bregman alternating direction method with multipliers for nonconvex composite problems.
\newblock {\em arXiv preprint arXiv:1410.8625}, 2014.

\bibitem{wyz2015}
Y. Wang, W. Yin and J. Zeng.
\newblock Global convergence of ADMM in nonconvex nonsmooth optimization.
\newblock {\em arXiv preprint arXiv:1511.06324}, 2015.

\bibitem{XCS2012}
H. Xu, C. Caramanis and S. Sanghavi.
\newblock Robust PCA via outlier pursuit.
\newblock {\em IEEE Transactions on Information Theory}, 58(5): 3047--3064, 2012.

\bibitem{xw2011}
M.H. Xu and T. Wu.
\newblock A class of linearized proximal alternating direction methods.
\newblock {\em Journal of Optimization Theory and Applications}, 151(2): 321--337, 2011.

\bibitem{z2010}
C.-H. Zhang.
\newblock Nearly unbiased variable selection under minimax concave penalty.
\newblock {\em Annals of Statistics}, 38(2): 894--942, 2010.

\bibitem{zrglz2015}
W. Zuo, D. Ren, S. Gu, L. Lin and L. Zhang.
\newblock Discriminative learning of iteration-wise priors for blind deconvolution.
\newblock in {\em CVPR}, 2015.











\end{thebibliography}
\end{document}